\documentclass[11pt]{article}
\usepackage[T1]{fontenc}
\usepackage{amsmath,amssymb,amsthm,mathrsfs}
\usepackage{graphicx,needspace,booktabs,longtable,enumerate}
\usepackage{float}
\usepackage{tikz,pgfplots}
\pgfplotsset{compat=1.18}
\usetikzlibrary{decorations.pathreplacing,arrows.meta,positioning,calc}
\usepackage[hidelinks]{hyperref}
\hypersetup{pdftitle={Finite-Time Blowup for Navier--Stokes with Smooth Forcing: Part I. Construction of Self-Similar Solutions with Admissible Stress and Flat Remainder},pdfauthor={Zhen Lei and Xiao Ren}}
\graphicspath{{Part-I-arxiv/}}
\newtheorem{theorem}{Theorem}[section]
\newtheorem{lemma}[theorem]{Lemma}
\newtheorem*{unnumberedlemma}{Lemma}
\newtheorem{proposition}[theorem]{Proposition}

\theoremstyle{definition}

\newtheorem{example}[theorem]{Example}
\theoremstyle{remark}
\newtheorem{remark}[theorem]{Remark}
\numberwithin{equation}{section}

\newcommand{\bR}{\mathbb R}

\begin{document}
\title{{\Large Finite-Time Blowup for Navier--Stokes with Smooth Forcing}\\[8pt]
\resizebox{\textwidth}{!}{{\large Part I. Construction of Self-Similar Solutions with Admissible Stress and Flat Remainder}}}
\author{Zhen Lei\footnote{Center for Applied Mathematics \& School of Mathematical Sciences, Fudan University, Shanghai 200433, P. R. China. Email: \texttt{zlei@fudan.edu.cn}.}
\and Xiao Ren\footnote{Center for Applied Mathematics, Fudan University, Shanghai 200433, P. R. China. Email: \texttt{xren@fudan.edu.cn}.}}
\date{}
\maketitle

\begin{abstract}
This paper gives a readable and accessible version of the
profile-construction part of OpenAI's manuscript \cite{1}.
We regard OpenAI's work as a major advance on the
Navier--Stokes Millennium Prize Problem.
We construct an axisymmetric leading flow that is regular at the axis and
sum its lower-order corrections to obtain an approximately self-similar
background field. Its Navier--Stokes
residual consists of a divergence-form term and a remainder vanishing
to infinite order in $1-t$ on fixed similarity sectors.
The associated stress and radial shear satisfy the admissible cone
condition. The stress need not be small. We introduce a new linear
model that provides a clearer explanation of the inner core construction.
The cancellation of the divergence-form residual by oscillatory pulses
is treated in the companion paper
\emph{Finite-Time Blowup for Navier--Stokes with Smooth Forcing---Part~II.
Residual Correction via Oscillatory Pulses}~\cite{PartII}.
We claim no original mathematical contribution.
\end{abstract}

\begingroup
\small
\tableofcontents
\endgroup
\section{Introduction}\label{sec:introduction}

We study the profile construction for the incompressible
Navier--Stokes equations
\begin{equation}\label{NS}
\begin{cases}
\partial_tu+(u\cdot\nabla)u+\nabla p=\Delta u+f,\\
\nabla\cdot u=0,
\end{cases}
\end{equation}
where $u$ is the velocity field and $p$ is the scalar pressure.
The leading profiles depend on the similarity variables $R\ge0$ and
$-1\le Z\le1$. Up to lower-order terms, their angular and axial
residuals are represented by a stress divergence.
The stress must vanish near the axis and in the exterior $[R_h,\infty)$.
On the intervening annulus, the stress and radial shear must satisfy
a prescribed cone condition.

The construction has two stages: the leading profile and its
lower-order corrections. The corrections remove the remainder outside the
stress divergence at every formal order. Their summation produces a
background field, smooth away from the singular point, with the required
residual and cone properties.

\subsection{Background and scope}

We place the profile problem in the self-similar setting.
Leray's work \cite{Leray} provides both the
finite-energy framework and the original motivation for studying
self-similar profiles. Exact backward self-similar solutions of the
unforced three-dimensional equations are strongly constrained.
Ne\v cas, R\r u\v zi\v cka, and \v Sver\'ak \cite{NRS} rule out
nonzero $L^3$ profiles, and Tsai \cite{Tsai} proves a corresponding
nonexistence result under local energy estimates. In the forward
direction, Jia and \v Sver\'ak \cite{JiaSverak} construct self-similar
solutions for large homogeneous initial data. 

We construct velocity fields in anisotropic similarity variables
with a nonzero admissible stress. Its magnitude need not be small.
The lower-order construction makes the remainder outside the stress divergence
flat; this does not imply smallness of the full Navier--Stokes residual.
The leading angular and axial equations retain radial viscosity; the
axial-viscosity terms have lower formal order. We construct an
approximately self-similar background with admissible stress and flat
remainder, satisfying the input conditions of Part~II.

The construction is based on OpenAI's manuscript \cite{1}, especially
its leading-profile construction and Appendices A--C.
We regard OpenAI's work as a major advance on the
Navier--Stokes Millennium Prize Problem \cite{Fefferman}.
This paper presents a readable and accessible version of this part of the construction,
retaining the explicit calculations and explaining how the pieces
determine one another. We introduce a new linear model that provides
a clearer explanation of the inner core construction.
The core argument is written in full, including the linear
model, the analytic space, the nonlinear fixed-point argument, and the
exit estimates. 

We also mention some related works. The works of Koch, Nadirashvili, Seregin and \v Sver\'ak \cite{KNSS},  Chen, Strain, Tsai and Yau\ \cite{CSTY1,CSTY2}, and Lei and Zhang \cite{LeiZhang} develop the Type I regularity theory of axi-symmetric Navier-Stokes equations.  Lei and Zhang \cite{LeiZhangCriticality} also prove regularity under a logarithmic modulus of continuity condition on $\Gamma=r u^\theta$ at the axis; see also the works of Chen, Fang, and Zhang \cite{ChenFangZhang} and Wei \cite{Wei}. Albritton, Bru\'e, and Colombo \cite{ABC} prove the non-uniqueness of Leray-Hopf weak solutions with a force singular at the initial time.  C\'ordoba, Mart\'inez-Zoroa, and Zheng \cite{CMZ} construct finite-time
singularities for forced hypodissipative Navier--Stokes equations with small dissipation orders and
 forcing in a local well-posedness class. With computer aid,  Chen and Hou \cite{ChenHouI,ChenHouII} study finite-time blowup for the three-dimensional Euler equations from smooth data in a cylinder with boundary. Recently, Alp\"oge and Buckmaster \cite{AlpogeBuckmaster} construct finite-time
blowup for the Euler equations on $\bR^3$ from smooth data with smooth forcing. OpenAI \cite{OpenAIEuler} constructs unforced Euler blowup from
smooth, compactly supported initial data. Constantin, Ignatova, and Vicol \cite{CIV} prove local regularity for a class of asymptotically axisymmetric Navier--Stokes flows under spatially analytic forcing. 

We construct the leading profiles first.
Section~\ref{sec:higher-order-correction}
then constructs the lower orders, corrects their moments, and sums the
sequence into smooth fields.
The subsequent cancellation of the divergence-form residual by oscillatory pulses
is treated in the companion paper
\emph{Finite-Time Blowup for Navier--Stokes with Smooth Forcing---Part~II.
Residual Correction via Oscillatory Pulses}~\cite{PartII}.
Part~II uses this background to construct solutions exhibiting finite-time
blowup with smooth forcing.

\subsection{The leading profile}

We state the required profile identities and inequalities.
For $0<\delta<1/200$, use cylindrical coordinates $(r,\theta,z)$ and set
\begin{equation}\label{intro:variables}
 R=\frac{r^2}{2\lambda^2},\qquad
 Z=\frac{z}{\lambda^{1-\delta}},\qquad
 1-t=\lambda^2(1-Z^2),\qquad \lambda>0.
\end{equation}
For $t<1$, these relations determine $\lambda=\lambda(t,z)$ and
$|Z|<1$. At $t=1$ and $z\ne0$, the endpoints $Z=\pm1$ describe
the smooth extension in these coordinates. The leading physical fields are
\begin{equation}\label{intro:leading-fields}
 u^\theta_{(0)}=\lambda^{-1-\delta}U^\theta_{(0)},\quad
 u^z_{(0)}=\lambda^{-1-\delta}U^z_{(0)},\quad
 u^r_{(0)}=\lambda^{-1}U^r_{(0)},\quad
 p_{(0)}=\lambda^{-2-2\delta}P_{(0)},
\end{equation}
where the profiles depend only on $(R,Z)$.
The independent profiles are $U^\theta,U^z$. Write
\[
 F(R,Z)=\frac{U^\theta(R,Z)}{\sqrt{2R}},\qquad
 L=1-Z^2,\qquad L_\delta=1-\delta Z^2.
\]
The factor $F$ is the similarity-variable
form of the regularized swirl $u^\theta/r$, also used in
\cite{HouLeiLi}. Regularity at the axis requires $F$ and $U^z$
to be smooth in $R$ at zero. The radial pressure balance and
incompressibility determine $P$ and $U^r$ through
\eqref{P-Utheta} and \eqref{Ur-Uz}.
Write $P_0(Z):=P_{(0)}(0,Z)=P(0,Z)$ for the axis value of the leading pressure profile.

We use the five moments $M^\theta,M^z,M^{\theta z},M^{z\theta},M^p$
defined in \eqref{fiveM}. The terminal conditions are
\begin{equation}\label{eq:moment-conditions}
\begin{gathered}
 M^z(\infty,Z)=M^{\theta z}(\infty,Z)=M^{z\theta}(\infty,Z)=0,\\
 \lim_{R\to\infty}\left[
 M^\theta(R,Z)-\frac{\sqrt2c_\infty}{1-\delta/2}R^{1-\delta/2}
 \right]=0,\\
 0<M^p(\infty,Z)=-P_0(Z)<\infty.
\end{gathered}
\end{equation}
The angular moment is renormalized because it does not converge at
infinity. The other conditions fix the radial velocity, mixed stresses,
and pressure normalization in the exterior. Subsection~\ref{intro:construction}
explains how the outer construction realizes the analytic axis pressure
$P_0=P_0^{\rm pre}$.

The stress has the decomposition $\mathcal T=\mathcal I+\mathcal S$,
where $\mathcal I$ is the contribution from inviscid terms, which we call
the inertial stress, and
\[
 \mathcal S=(2R\partial_RF,\sqrt{2R}\partial_RU^z)
\]
is the shear. The stress representation and the definitions of these
profiles are given in Subsections~\ref{subsec:residual-stress-shear}
and~\ref{subsec:stress-shear-radial}.
On the stress annulus $R_a\le R\le R_h$, we require
$U^\theta>0$ and $\mathcal S^\theta<0$.
Write $\mathcal S^\perp=(-\mathcal S^z,\mathcal S^\theta)$ and set
\[
 \kappa=-\frac{|\mathcal S|^2}{F\mathcal S^\theta}.
\]
The admissible cone condition is
\begin{equation}\label{eq:cone}
\begin{gathered}
 \mathcal T\cdot\mathcal S<0,\qquad \kappa>2,\\
 (\kappa-2)(\mathcal T\cdot\mathcal S^\perp)^2
       <2(\mathcal T\cdot\mathcal S)^2.
\end{gathered}
\end{equation}
The strict inequalities are imposed where the stress is nonzero. The
stress-free core $[0,R_a]$ and exterior $[R_h,\infty)$ supply the
boundary values of this annulus.
Section~\ref{sec:cone-condition} explains the geometry and introduces
the relaxed cone conditions, which are typically satisfied by the
intermediate connections. The shear modification in
Section~\ref{sec:shear-modification} strengthens them to the admissible
cone condition. Subsection~\ref{subsec:sm-completion} completes the proof
of the following theorem.

\begin{theorem}[The leading flow]\label{thm:leading}
There exist $0<\delta<1/200$, $0<R_a<R_h<\infty$, $c_\infty>0$, and
bounded leading profiles $U^\theta_{(0)}$ and $U^z_{(0)}$ on
$[0,\infty)\times[-1,1]$, with $U^\theta_{(0)}/\sqrt{2R}$ and $U^z_{(0)}$
smooth up to $R=0$ and $Z=\pm1$. Let $P_{(0)}$ and $U^r_{(0)}$ be
determined by \eqref{P-Utheta} and \eqref{Ur-Uz}, respectively, applied to
the leading profiles. Then the physical fields $(u_{(0)},p_{(0)})$
defined by \eqref{intro:leading-fields} are smooth away from $(t,x)=(1,0)$.
The velocity is exactly divergence-free, and the following properties hold.
\begin{enumerate}[(i)]
\item For $R\ge R_h$, the profiles agree with the exterior heat-flow
profiles \eqref{exterior-profile}, with the pressure given by \eqref{exterior-ansatz}.
\item The leading stress profile $\mathcal T_{(0)}$ is supported on
$[R_a,R_h]\times[-1,1]$ and is nonzero whenever $R_a<R<R_h$.
\item We have $U^\theta_{(0)}>0$ for $R>0$ and
$\mathcal S^\theta_{(0)}<0$ on $[R_a,R_h]\times[-1,1]$.
The stress and shear satisfy the admissible cone condition
\eqref{eq:cone} wherever $\mathcal T_{(0)}\ne0$.
\item The five moments satisfy \eqref{eq:moment-conditions}
for every $Z\in[-1,1]$.
\item Write $\mathsf{NS}(u,p)=\partial_tu+(u\cdot\nabla)u-\Delta u+\nabla p$
and $\mathsf T_{(0)}=\lambda^{-2-\delta}\mathcal{T}_{(0)}$.
The momentum residual satisfies
\begin{equation}\label{eq:residual-form}
\begin{aligned}
 \mathsf{NS}^\theta(u_{(0)},p_{(0)})
   &=-\left(\partial_r+\frac2r\right)\mathsf T^\theta_{(0)}
      -\partial_z^2u^\theta_{(0)},\\
 \mathsf{NS}^z(u_{(0)},p_{(0)})
   &=-\left(\partial_r+\frac1r\right)\mathsf T^z_{(0)}
      -\partial_z^2u^z_{(0)},\\
 \mathsf{NS}^r(u_{(0)},p_{(0)})&=O(\lambda^{-3}).
\end{aligned}
\end{equation}
Thus the angular and axial residuals consist of stress-divergence terms,
supported in the fixed similarity annulus $R_a\le R\le R_h$, and
the axial-viscosity remainders $-\partial_z^2u^\theta_{(0)}$ and
$-\partial_z^2u^z_{(0)}$. On fixed compact profile ranges as $\lambda\downarrow0$,
these remainders have formal order $O(\lambda^{-3+\delta})$, smaller
by a factor $\lambda^{2\delta}$ than the leading angular and axial terms.
The radial estimate in \eqref{eq:residual-form} holds on the same ranges.
\end{enumerate}
\end{theorem}

\subsection{Strategy for the leading-order construction}\label{intro:construction}

We first prescribe the complete temporary candidate. On $0<R\le R_{\rm ref}$
we prescribe
\[
 U^\theta_{\rm ref}=\frac{P_*}{1+Z^2}
              \left(\frac R{R_{\rm ref}}\right)^{1/10},
 \qquad U^z_{\rm ref}=4Z.
\]
We then prescribe its outer continuation, including every transition,
correction interval, and the heat-flow tail on $R\ge R_h$.
For each provisional waiting length $\tau$, let $\bar U^\theta_{\rm pre}$ be
the uncorrected angular profile with the heat factor $H_\delta$
replaced by $1$ and all other parameters unchanged. Matching the pre-heat
angular moment at $Z=0$ selects $\tau$ by \eqref{eq:mc-wait-choice}.
At this value of $\tau$ the axis-pressure target is
\begin{equation}\label{intro:preheat-axis-pressure}
 P_0^{\rm pre}(Z)=-\int_0^\infty
       \frac{(\bar U^\theta_{\rm pre}(R,Z))^2}{2R}\,dR.
\end{equation}
The pre-heat profile has analytic $Z$ dependence, and the integral
extends holomorphically to a complex neighborhood of $[-1,1]$ by
Lemma~\ref{lem:outer-axis-pressure}. The two angular bumps on Interval O.6
match the angular moment and restore the pressure moment.
The complete angular profile $U_o^\theta$, including its reference
continuation and the exact heat tail, therefore determines
\[
 P_0(Z)=-\int_0^\infty\frac{(U_o^\theta(R,Z))^2}{2R}\,dR
       =P_0^{\rm pre}(Z),
\]
as proved in \eqref{eq:mc-pressure-preserved}.
For fixed $P_*,\delta$ and dimensionless outer schedule, the profiles
in Lemma~\ref{lem:imc-prepared-family} share the same $\tau$ and $P_0$,
with the heat tail and moment corrections constructed at each $R_{\rm ref}$.

The inner replacement uses this exact datum.
Its regular core occupies $[0,R_a]$, where $R_a=4/\Lambda$.
We choose $U^z(0,Z)=4Z+j$ and a small, nonconstant angular amplitude $F_0(Z)$.
The linear model explains these choices and gives $\kappa_{\rm m}(R_a,Z)>3$.
The nonlinear comparison gives $\kappa(R_a,Z)\ge5/2$.
The new core is connected to the reference velocity before
$R_{\rm useless}=e^{-5}R_{\rm ref}$. A correction in $[R_m,2R_m]$, with
$R_m=e^{-6}R_{\rm ref}$, restores all five moments. The joined profile
agrees with the original outer fields on $[R_{\rm useless},\infty)$.
It satisfies the relaxed cone on the open stress annulus and the
admissible cone on a collar of each radial edge.

Finally, we modify the radial shear on $[R_{-c},R_{+c}]$ and restore
the five moments in the separate outer interval $[R_c,2R_c]$.
Together these modifications give the admissible cone throughout the open
stress annulus and preserve the core and the exact heat exterior. The
resulting leading profile has $P_{(0)}(0,Z)=-M^p(\infty,Z)=P_0^{\rm pre}(Z)$.

\subsection{Lower-order corrections and the complete background}\label{intro:lower-order}\label{intro:background-input}

Starting with the remainder in Theorem~\ref{thm:leading},
Section~\ref{sec:higher-order-correction} constructs successive corrections
and sums them to obtain the background required by
Part~II \cite[Section~1.1 and Theorem~1.1]{PartII}.

The physical stress pair is $\mathsf T_B=\lambda^{-2-\delta}\mathcal{T}_B$,
where $\mathcal{T}_B(\lambda,R,Z)$ is the normalized stress.
For a physical pair $\mathsf T=(\mathsf T^\theta,\mathsf T^z)$, put
\begin{equation}\label{intro:stress-operator}
 \mathcal D\mathsf T
 =\left(\partial_r+\frac2r\right)\mathsf T^\theta e_\theta
  +\left(\partial_r+\frac1r\right)\mathsf T^z e_z.
\end{equation}
In normalized derivatives, $(\lambda,R,Z)$ are independent variables.
Let $\rho=\sqrt{2R}$ and let $\mathscr D^{\le m}$ denote compositions
of at most $m$ operators from
$\{\lambda\partial_\lambda,\partial_\rho,\partial_Z\}$.
Thus $\partial_\rho=\sqrt{2R}\partial_R$.

In addition to the conclusions of Theorem~\ref{thm:leading}, the proof
uses three properties of the constructed profiles. Analyticity in $Z$
on a fixed inner interval extending beyond $R_a$ allows us to solve the
coefficient equations there at every order
(Lemmas~\ref{lo-core:analytic-input} and~\ref{lo-core:fixed-radius}).
The pure-power swirl and zero axial velocity on $I_2$ allow us to solve
the five moment equations there and preserve the same velocity profile
on $I_3$ (Subsection~\ref{lo-ext:moments}).
At both edges of the stress annulus, the estimates on the decay and
direction of the leading stress allow us to control the stress correction
relative to $|\mathcal T_{(0)}|$ and preserve a uniform cone margin
(Subsection~\ref{lo-cone:subsection}).

\begin{theorem}[Background with admissible stress and flat remainder]
\label{thm:background}
The leading profiles constructed in the proof of Theorem~\ref{thm:leading}
admit lower-order corrections yielding real axisymmetric fields
$(u_B,p_B,\mathsf T_B,\mathsf E_B)$ with the following properties, for some
$\lambda_*>0$ and throughout
\[
 0<\lambda\le\lambda_*,\qquad R\ge0,\qquad |Z|\le1.
\]
All fields are smooth away from $(t,x)=(1,0)$; derivatives at
$Z=\pm1$ are understood by smooth one-sided extension.
The equation is
\begin{equation}\label{intro:background-equation}
 \nabla\cdot u_B=0,\qquad
 \mathsf{NS}(u_B,p_B)=-\mathcal D\mathsf T_B+\mathsf E_B.
\end{equation}
\begin{enumerate}[(i)]
\item \emph{Flat remainder.}
For every Cartesian space--time multi-index $\nu$ and every $N>0$,
\begin{equation}\label{intro:flatness}
 |\partial_{t,x}^\nu \mathsf E_B|\le C_{\nu,N}\lambda^N.
\end{equation}
Moreover, $\mathsf E_B=0$ for $R\ge R_h$.

\item \emph{Rescaled velocity bounds and leading behavior.}
Every $\mathscr D^{\le m}$ derivative of
\[
 \lambda^{1+\delta}u_B^\theta,\qquad
 \lambda^{1+\delta}u_B^z,\qquad \lambda u_B^r
\]
is uniformly bounded. On each range $R\le K$, their differences from
the corresponding leading profiles are $O_{m,K}(\lambda^{2\delta})$
in these normalized derivatives.

\item \emph{Annular support and weighted bounds.}
The physical stress is supported in $[R_a,R_h]$. On the open annulus set
\begin{equation}\label{intro:edge-weight}
\begin{aligned}
 d_{\rm edge}(R)&=\min\{1,\log(R/R_a),\log(R_h/R)\},\\
 \omega(R)&=\exp\left[-\frac{c_a}{\log^2(R/R_a)}
                       -\frac{c_h}{\log^2(R_h/R)}\right],
 \qquad c_a,c_h>0.
\end{aligned}
\end{equation}
The constants $c_a,c_h$ are fixed by the endpoint collars, and $\omega$
is extended by zero off the annulus. There are $c>0$ and, for each
$m$, finite constants $C_m,q_m>0$ such that
\begin{equation}\label{intro:weighted-stress}
 |\mathcal{T}_B|\ge c\omega,\qquad
 |\mathscr D^{\le m}\mathcal{T}_B|\le C_m\omega d_{\rm edge}^{-q_m}.
\end{equation}
The stress extends smoothly by zero across both radial edges.

\item \emph{Uniform admissible cone.}
Define the shear of the complete background by
\begin{equation}\label{intro:background-shear}
\begin{gathered}
 F_B=\frac{\lambda^{1+\delta}u_B^\theta}{\sqrt{2R}},\qquad
 \mathcal{S}_B=\left(2R\partial_RF_B,
       \sqrt{2R}\partial_R(\lambda^{1+\delta}u_B^z)\right),\\
 \mathcal{S}_B^\perp=(-\mathcal{S}_B^z,\mathcal{S}_B^\theta),\qquad
 \kappa_B=-\frac{|\mathcal{S}_B|^2}{F_B\mathcal{S}_B^\theta}.
\end{gathered}
\end{equation}
On the closed annulus, $F_B$ and $-\mathcal{S}_B^\theta$ belong to
$[C^{-1},C]$ for a fixed $C\ge1$. For a fixed $c_{\rm cone}\in(0,1)$,
throughout the open annulus one has
\begin{equation}\label{intro:uniform-cone}
\begin{gathered}
 \kappa_B\ge2+c_{\rm cone},\qquad
 -\mathcal{T}_B\cdot\frac{\mathcal{S}_B}{|\mathcal{S}_B|}\ge c_{\rm cone}|\mathcal{T}_B|,\\
 \sqrt{\frac{\kappa_B-2}{2}}
 \left|\frac{\mathcal{T}_B\cdot \mathcal{S}_B^\perp}{\mathcal{T}_B\cdot \mathcal{S}_B}\right|
 \le1-c_{\rm cone}.
\end{gathered}
\end{equation}

\item \emph{A reserved moment-correction interval.}
There is a fixed nonempty open interval $I_3\Subset(R_a,R_h)$ on
which the complete background satisfies the exact identities
\begin{equation}\label{intro:reserved-patch}
 u_B^z=0,\qquad
 u_B^\theta=\lambda^{-1-\delta}A_*(Z)(2R)^{-1/2-\mu},
 \qquad A_*(Z)=\frac{a_*}{1+Z^2},\quad a_*,\mu>0.
\end{equation}
In particular, $A_*$ is smooth and has a positive lower bound.

\item \emph{Exact heat exterior.}
For $R\ge R_h$,
\begin{equation}\label{intro:heat-exterior}
\begin{aligned}
 u_B^r=u_B^z&=0,\qquad u_B^\theta=K(t,r),\\
 p_B(t,r)&=-\int_r^\infty K(t,s)^2\,\frac{ds}{s},\\
 K(t,r)&=c_\infty\left(\frac{r^2}{2}\right)^{-(1+\delta)/2}
 H_\delta\left(\frac{4(1-t)}{r^2}\right),\\
 H_\delta(\xi)&=\frac1{\Gamma(1+\delta/2)}
   \int_0^\infty e^{-v}v^{\delta/2}(1+\xi v)^{-\delta/2}\,dv.
\end{aligned}
\end{equation}
This exterior solves the Navier--Stokes equations exactly.
\end{enumerate}
All constants above may depend on the completed construction, the
indicated derivative orders, and, for local estimates, the chosen
compact profile range, but not on $(\lambda,R,Z)$ in their stated ranges.
\end{theorem}

At each order we gain a factor $\lambda^{2\delta}$ and solve the linear
velocity and pressure equations on the same inner interval extending beyond
$R_a$. We cut off the two velocity profiles before $R_{-c}$ and use three
angular bumps and two axial bumps in $I_2$ to solve the five affine moment
equations. The positive-order velocities then vanish for $R\ge\sup I_2$,
preserving the leading velocity on $I_3$. Since $P_{(n)}(0,Z)=0$ for $n\ge1$,
the summed pressure satisfies $\lambda^{2+2\delta}p_B(t,0,z)=P_0^{\rm pre}(Z)$.

Shrinking cutoffs in $\lambda$ sum the coefficients into smooth fields.
For the meridional velocity, we apply the cutoffs to the streamfunctions
before taking curls, which preserves incompressibility exactly.
Comparison with finite truncations proves that $\mathsf E_B$ and all its Cartesian
space--time derivatives vanish to infinite order in $\lambda$, and hence
in $1-t$ on each fixed sector $|Z|\le1-\epsilon$. The bounds for the shear
perturbation and the stress estimates relative to $\mathcal{T}_{(0)}$
up to both radial edges preserve the uniform cone margin
(Proposition~\ref{lo-cone:preservation}).

Equation~\eqref{intro:background-equation} also has the full tensor form
$\mathsf{NS}(u_B,p_B)=-\nabla\cdot\mathbb T_B+\mathsf E_B$, where
\begin{equation}\label{intro:full-tensor}
 \mathbb T_B=\mathsf T_B^\theta(e_r\otimes e_\theta+e_\theta\otimes e_r)
 +\mathsf T_B^z(e_r\otimes e_z+e_z\otimes e_r)
 +r\partial_z\mathsf T_B^z\,e_\theta\otimes e_\theta.
\end{equation}
The diagonal entry cancels the radial component of the divergence, giving
$\nabla\cdot\mathbb T_B=\mathcal D\mathsf T_B$. At each fixed $(t,z)$,
integration of the tangential equations in \eqref{intro:background-equation}
and the annular support of $\mathsf T_B$ give the force identities used in Part~II,
\begin{equation}\label{intro:force-moments}
\begin{aligned}
 \int_0^\infty r^2(\mathsf{NS}(u_B,p_B)-\mathsf E_B)^\theta\,dr&=0,\\
 \int_0^\infty r(\mathsf{NS}(u_B,p_B)-\mathsf E_B)^z\,dr&=0.
\end{aligned}
\end{equation}

\subsection{Organization of the paper}

We follow the dependencies in the construction.
Sections~\ref{sec:axi} and \ref{sec:leading-order-system} fix the
variables and profile identities. Section~\ref{sec:roadmap-construction}
gives the outer prescriptions and their requirements, followed by the
inner replacement and the shear modification.

\begin{center}
\small
\begin{tabular}{@{}p{.22\linewidth}p{.71\linewidth}@{}}
\toprule
Detailed construction & Input and output\\
\midrule
Section~\ref{sec:outer-construction} &
The exact heat exterior and the reference-plus-outer candidate.
The waiting length and the correction coefficients determine the
matched outer profile and its analytic axis pressure.\\[4pt]
Section~\ref{sec:inner-construction} &
The axis pressure gives the datum for the analytic core. A smooth
velocity connection and five moment corrections join the core to the
outer profile with the relaxed cone.\\[4pt]
Section~\ref{sec:shear-modification} &
The radial shear modification and its moment repair give an
admissible-cone profile with the same terminal moments, core, and exterior.\\[4pt]
Section~\ref{sec:higher-order-correction} &
The inner coefficient equations on a common interval, the five moment
corrections at each order, and smooth summation give the complete
background with admissible stress and a flat remainder.\\
\bottomrule
\end{tabular}
\end{center}

Appendix~\ref{app:analytic-core} proves the analytic construction of the
stress-free core and the comparison estimates at its exit.

\section{Preliminaries}\label{sec:axi}
We fix the axisymmetric equations and similarity variables.
We then identify the leading terms and the profile requirements.
We write all derivatives with explicit differentiation operators, such as
$\partial_R G$, $\partial_Z G$, $\partial_R^2G$, and
$\partial_R\partial_ZG$.
The variable of differentiation is always indicated.

\subsection{Axi-symmetric Navier-Stokes equations}
We write the equations for axisymmetric fields of the form
\[
u = u^r(r,z,t)\,e_r + u^\theta(r,z,t)\,e_\theta + u^z(r,z,t)\,e_z,\qquad p = p(r,z,t),
\]
With the force set to zero, \eqref{NS} becomes
\begin{equation}\label{eq:axiNS}
\begin{cases}
\partial_t u^r + u^r \partial_r u^r + u^z \partial_z u^r
- \dfrac{(u^\theta)^2}{r} + \partial_r p
= \Big(\Delta - \dfrac{1}{r^2}\Big) u^r, \\[1mm]
\partial_t u^\theta + u^r \partial_r u^\theta + u^z \partial_z u^\theta
+ \dfrac{u^r u^\theta}{r}
= \Big(\Delta - \dfrac{1}{r^2}\Big) u^\theta, \\[1mm]
\partial_t u^z + u^r \partial_r u^z + u^z \partial_z u^z + \partial_z p
= \Delta\, u^z, \\[1mm]
\partial_r(ru^r) + \partial_z(ru^z) = 0,
\end{cases}
\end{equation}
where all unknowns $(u^r, u^\theta, u^z, p, f)$ are independent of $\theta$ and the
axi-symmetric Laplacian is $\Delta = \partial_r^2 + \frac{1}{r}\partial_r + \partial_z^2$.

As long as the solution is smooth, it automatically satisfies
$$u^r\big|_{r = 0} = u^\theta\big|_{r = 0} = 0.$$
Indeed, $\frac{u^r}{r}$ and $\frac{u^\theta}{r}$ should be smooth functions of $r^2$ and $(t, z)$ if $u$ is smooth in $(t, x)$.

\subsection{Similarity variables}
\label{subsec:variables}

We introduce the coordinates used throughout the construction.
Let $0<\delta\ll1$ and set
\begin{equation}\label{R-Z-def}
R=\frac{r^2}{2\lambda^2},
\qquad
Z=\frac{z}{\lambda^{1-\delta}},
\end{equation}
where, for $0\le t<1$, the positive scaling parameter
$\lambda=\lambda(t,z)$ is determined implicitly by
\begin{equation}\label{lambda-def}
\lambda=\sqrt{\frac{1-t}{1-Z^2}}.
\end{equation}
We use the abbreviations
\begin{equation}\label{LZ}
L_\delta=L_\delta(Z)=1-\delta Z^2,
\qquad
L=L(Z)=1-Z^2.
\end{equation}
On $[-1,1]$, we have $1-\delta\le L_\delta\le1$ and $0\le L\le1$.
In particular, $L_\delta$ stays uniformly positive, while $L$ vanishes at $Z=\pm1$. The endpoints correspond to the limiting similarity boundary at $t=1$.

Substituting $Z=z/\lambda^{1-\delta}$ into
$1-t=\lambda^2L$ gives
\begin{equation}\label{lambda-implicit}
\lambda^2-\lambda^{2\delta}z^2=1-t.
\end{equation}
Differentiating at fixed $z$ gives
\[
\partial_\lambda
\bigl(\lambda^2-\lambda^{2\delta}z^2 - 1 + t\bigr)
=2\lambda(1-\delta  z^2\lambda^{2\delta-2})
=2\lambda L_\delta>0.
\]
Here $L_\delta\ge1-\delta>0$ on $|Z|\le1$.
Thus the left-hand side of \eqref{lambda-implicit} increases in $\lambda$
on $\lambda\ge |z|^{1/(1-\delta)}$, where $|Z|\le1$.

For each fixed $0\le t<1$, these relations define a one-to-one
change of variables
\[
(r,z)\in[0,\infty)\times\mathbb R
\quad\longleftrightarrow\quad
(R,Z)\in[0,\infty)\times(-1,1).
\]
Since the axial similarity variable $Z$ has bounded range, we need only
achieve suitable asymptotic behavior in the radial variable $R$ and glue
the profiles across radial interfaces, with the matching conditions imposed
for every $Z\in[-1,1]$.
The inverse relations are
\[
t=1-\lambda^2(1-Z^2),
\qquad
r=\lambda\sqrt{2R},
\qquad
z=\lambda^{1-\delta}Z.
\]
At $t=1$, equation \eqref{lambda-implicit} has a positive solution
precisely when $z\ne0$, and then
\[
\lambda(1,z)=|z|^{1/(1-\delta)},
\qquad
Z=\operatorname{sgn}(z).
\]
Hence the coordinates extend to $t=1$ away from $z=0$, with
$Z=\pm1$ corresponding to the terminal time. Also,
\[
\lambda(t,0)=\sqrt{1-t}.
\]
The derivative identities below show that $\lambda$ decreases
with $t$ and increases with $|z|$.
As $(t,r,z)\to(1,r_0,0)$ from $t<1$ with $r_0>0$, we have
$\lambda\to0$ and hence $R=r^2/(2\lambda^2)\to\infty$.
Near each such point the velocity and pressure therefore agree with the
heat exterior \eqref{intro:heat-exterior}, whose formulas are smooth up to
$t=1$ for $r>0$.

\begin{figure}[htbp]
\centering
\includegraphics[width=\textwidth]
{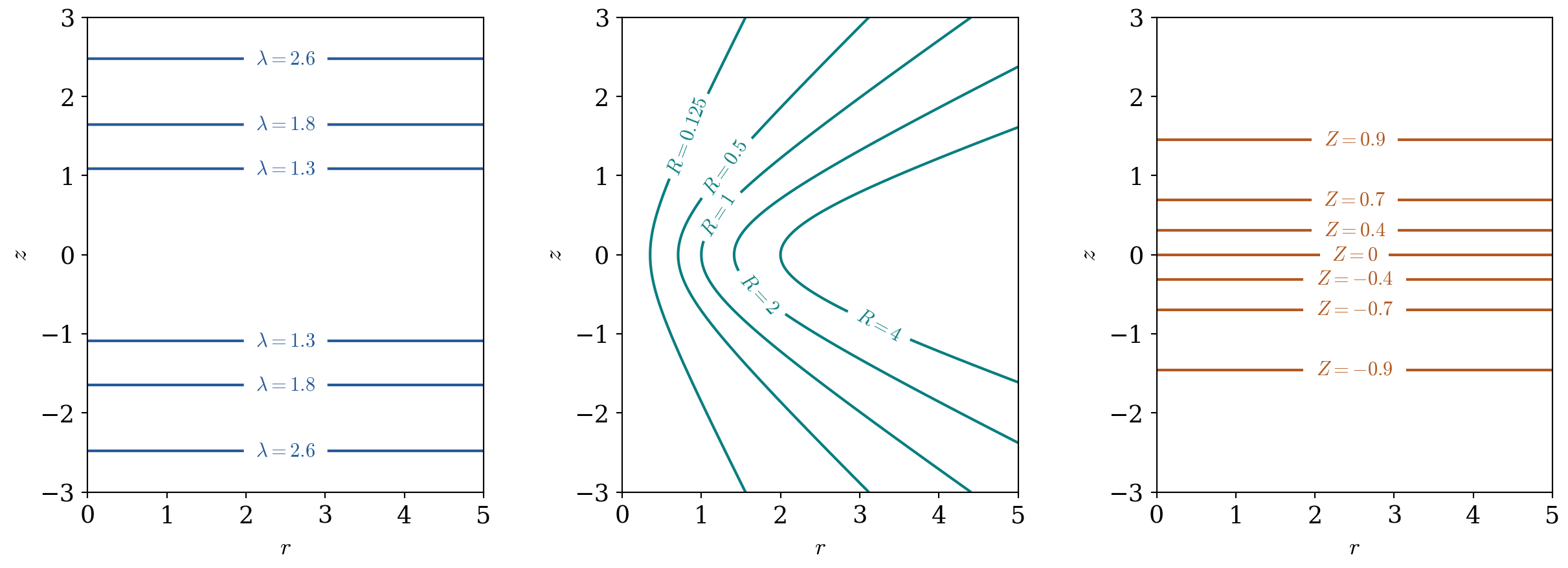}
\caption{Level lines of $\lambda$ (left), $R$ (middle),
and $Z$ (right) in the $(r,z)$-plane at $t=\tfrac12$,
with $\delta=0.01$.}
\label{fig:self-similar-level-lines}
\end{figure}

\begin{figure}[htbp]
\centering
\includegraphics[width=0.65\textwidth]
{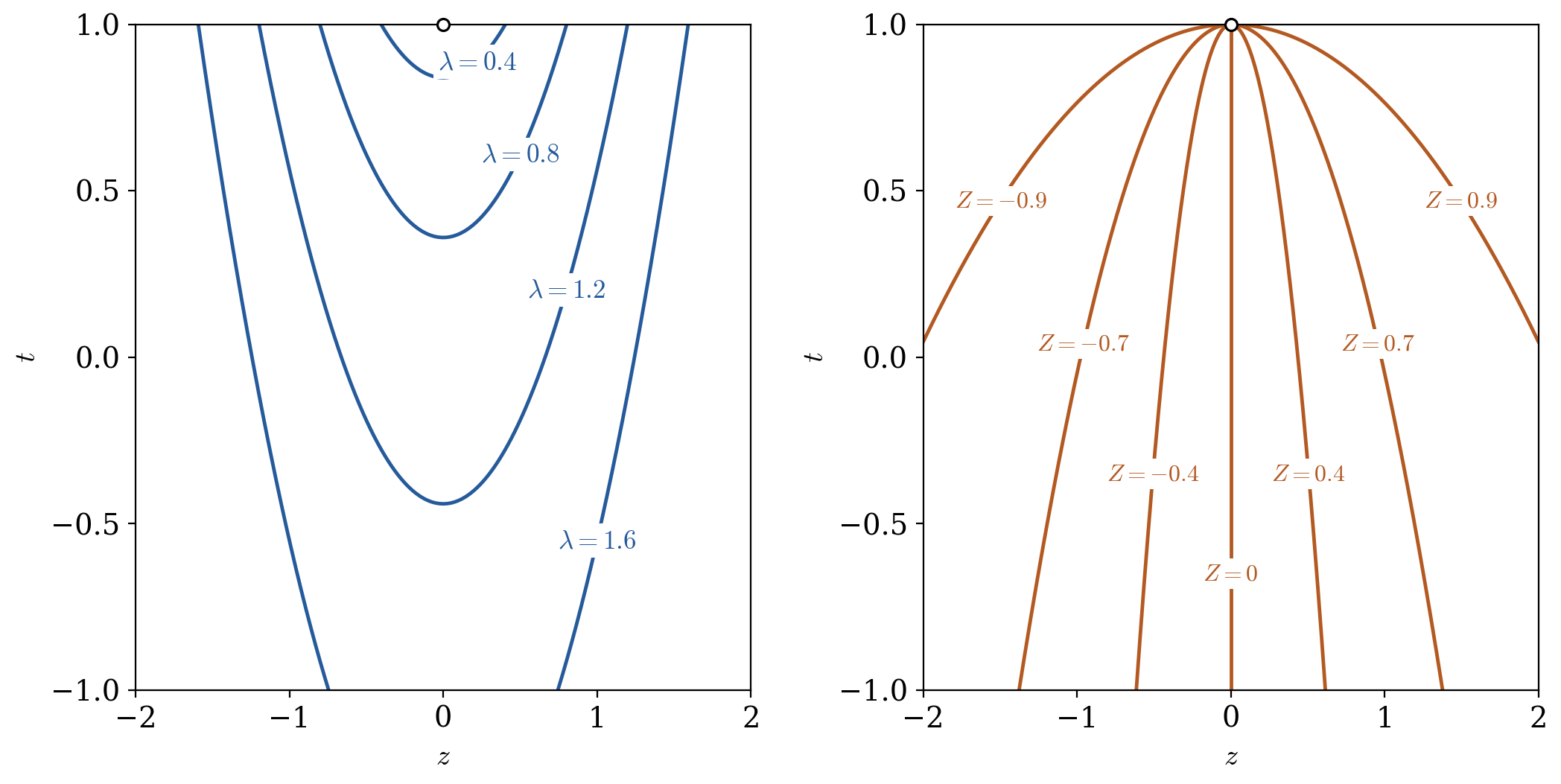}
\caption{Level lines of $\lambda$ (left) and $Z$ (right)
in the $(z,t)$-plane, with $\delta=0.01$.
The open circle marks the excluded point $(0,1)$. The $Z=\pm 1$ level lines are flat, satisfying $t=1$.}
\label{fig:lambda-Z-level-lines}
\end{figure}

\begin{figure}[htbp]
\centering
\includegraphics[width=0.65\textwidth]
{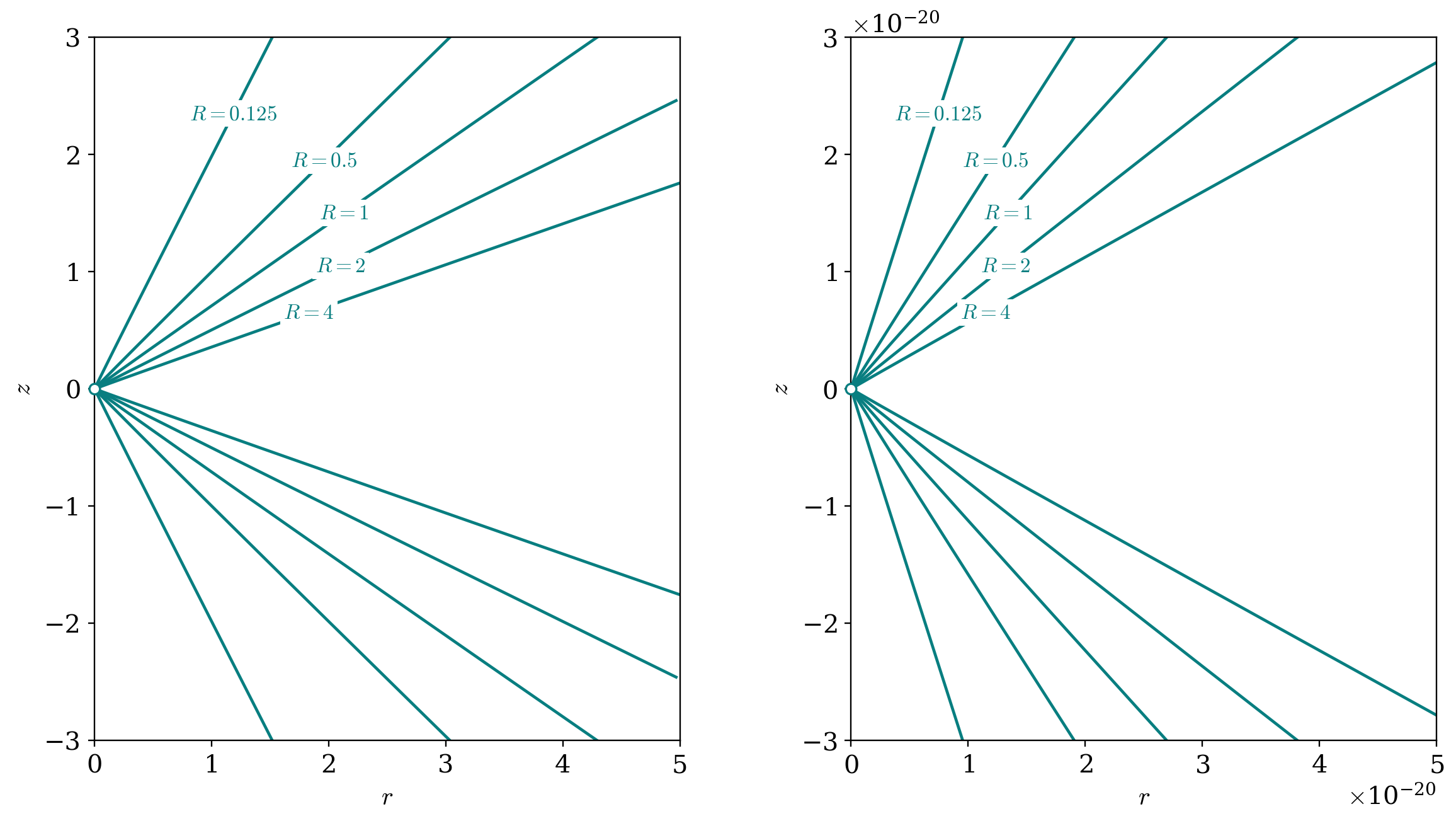}
\caption{Level lines of $R$ at $t=1$, with $\delta=0.01$:
the order-one scale (left) and the $10^{-20}$ scale (right).
The curves satisfy $r=\sqrt{2R}\,|z|^{1/(1-\delta)}$.
The open circles mark the excluded origin.}
\label{fig:R-level-lines-two-scales}
\end{figure}

For a fixed $t<1$, both $\lambda$ and $Z$ depend only on $z$,
whereas the level lines of $R$ are curved in the $(r,z)$-plane.
These level lines are illustrated in
Figure~\ref{fig:self-similar-level-lines}.
Figure~\ref{fig:lambda-Z-level-lines} shows the geometry of
$\lambda$ and $Z$ in the $(z,t)$-plane.
At $t=1$, the level lines of $R$ satisfy
$r=\sqrt{2R}\,|z|^{1/(1-\delta)}$, as illustrated in
Figure~\ref{fig:R-level-lines-two-scales}.

The leading flow is built from profiles
$(U^\theta,U^z,U^r,P)$ of $(R,Z)$ by
\begin{equation}\label{solution-profile}
\begin{cases}
u^\theta=\dfrac{U^\theta(R,Z)}{\lambda^{1+\delta}},
\qquad
u^z=\dfrac{U^z(R,Z)}{\lambda^{1+\delta}},
\\[2mm]
u^r=\dfrac{U^r(R,Z)}{\lambda},
\qquad
p=\dfrac{P(R,Z)}{\lambda^{2+2\delta}}.
\end{cases}
\end{equation}
Here we suppress the leading-order subscript $(0)$.
The two profiles to be constructed are $U^\theta$ and $U^z$;
the radial pressure balance and incompressibility determine
$P$ and $U^r$, respectively.

This scaling is supercritical relative to the Type~I rate.
Indeed, at $z=0$ we have $\lambda=\sqrt{1-t}$ and $Z=0$.
The leading angular profile satisfies $U^\theta(R,Z)>0$ for every $R>0$.
Hence, for any fixed $R_*>0$, along $r=\sqrt{2R_*}\sqrt{1-t}$,
\[
\sqrt{1-t}\,|(u^\theta,u^z)(t,r,0)|
=(1-t)^{-\delta/2}
|(U^\theta,U^z)(R_*,0)|
\longrightarrow\infty.
\]

The regular swirl variable is
\[
F(R,Z)=\frac{U^\theta(R,Z)}{\sqrt{2R}}.
\]
Indeed, \eqref{solution-profile} and $r=\lambda\sqrt{2R}$ give
$F=\lambda^{2+\delta}u^\theta/r$, identifying the regularized swirl
in the present similarity variables.
Smoothness of the physical fields at the axis requires ($\mathcal{I}^\theta$, $\mathcal{I}^z$ will be defined in \eqref{inertial-stress})
\begin{equation} \label{smooth-axis}
F,\quad U^z,\quad \frac{U^r}{\sqrt{R}},\quad P, \quad \frac{\mathcal{I}^\theta}{R}, \quad \frac{\mathcal{I}^z}{\sqrt{R}}
\end{equation}
to extend smoothly to $R=0$ as functions of $(R,Z)$.
In particular,
\[
 U^\theta(0,Z)=U^r(0,Z)=\mathcal{I}^\theta(0,Z)=\mathcal{I}^z(0,Z)=0.
\]

The chain rule for the similarity variables will be used
throughout the construction.

\begin{lemma}[Rescaling of derivatives]\label{lem:chain}
For any smooth profile $G=G(R,Z)$ and any real constant $\beta$,
\begin{equation}\label{profile-chain-rule}
\begin{aligned}
\partial_t(\lambda^\beta G)
&=\frac{\lambda^{\beta-2}}{L_\delta}
  \left(
  -\frac{\beta}{2}G
  +\frac{1-\delta}{2}Z\partial_Z G
  +R\partial_R G
  \right),\\
\partial_z(\lambda^\beta G)
&=\frac{\lambda^{\beta-1+\delta}}{L_\delta}
  \left(
  \beta ZG+L\partial_Z G-2ZR\partial_R G
  \right),\\
\partial_r(\lambda^\beta G)
&=r\lambda^{\beta-2}\partial_R G
 = \sqrt{2R}\lambda^{\beta-1}\partial_R G.
\end{aligned}
\end{equation}
\end{lemma}

\begin{proof}
Differentiating \eqref{lambda-implicit} and
\eqref{R-Z-def} gives
\begin{equation}\label{deriv}
\begin{aligned}
\partial_t\lambda&=-\frac{1}{2\lambda L_\delta},
&
\partial_z\lambda&=\frac{\lambda^\delta Z}{L_\delta},
&
\partial_r\lambda&=0,
\\
\partial_t Z&=\frac{(1-\delta)Z}{2\lambda^2L_\delta},
&
\partial_z Z&=\frac{L}{\lambda^{1-\delta}L_\delta},
&
\partial_r Z&=0,
\\
\partial_t R&=\frac{R}{\lambda^2L_\delta},
&
\partial_z R&=-\frac{2RZ}{\lambda^{1-\delta}L_\delta},
&
\partial_r R&=\frac{\sqrt{2R}}{\lambda}.
\end{aligned}
\end{equation}
For example, differentiating
$\lambda^2-\lambda^{2\delta}z^2=1-t$
with respect to $z$ at fixed $t$ and using $Z=z\lambda^{\delta-1}$ give
\[
\left(2\lambda-2\delta\lambda^{2\delta-1}z^2\right)
\partial_z\lambda
=2\lambda^{2\delta}z.
\]
Differentiating $Z=z\lambda^{\delta-1}$ yields
\[
\begin{aligned}
\partial_zZ
&=\lambda^{\delta-1}
+(\delta-1)z\lambda^{\delta-2}\partial_z\lambda\\
&=\lambda^{\delta-1}
\left(1-\frac{(1-\delta)Z^2}{L_\delta}\right) 
=\frac{L}{\lambda^{1-\delta}L_\delta}.
\end{aligned}
\]
The remaining identities follow similarly from
$R=r^2/(2\lambda^2)$ and the defining relations for $\lambda$ and $Z$.
For $q=t,z,r$, the product and chain rules yield
\[
\partial_q(\lambda^\beta G)
=
\beta\lambda^{\beta-1}(\partial_q\lambda)G
+\lambda^\beta
\bigl((\partial_R G)\partial_qR+(\partial_Z G)\partial_qZ\bigr).
\]
Substituting \eqref{deriv} gives
\eqref{profile-chain-rule}.
\end{proof}

\subsection{Formal order analysis}\label{sec:orders}

We compare powers of $\lambda$ at fixed $(R,Z)$.
Here and below,
\[
\Delta_r=\partial_r^2+\frac1r\partial_r
\]
denotes the radial part of the Laplacian, so that
$\Delta=\Delta_r+\partial_z^2$.
For fields of the form \eqref{solution-profile}, the powers of
$\lambda$ are determined by \eqref{profile-chain-rule}.

In the angular and axial equations, we have
\begin{equation}\label{orderU}
\partial_tu^\theta
+u^r\partial_ru^\theta
+u^z\partial_zu^\theta
+\frac{u^ru^\theta}{r}
-\left(\Delta_r-\frac1{r^2}\right)u^\theta
=O(\lambda^{-3-\delta}),
\end{equation}
\begin{equation}\label{orderUz}
\partial_tu^z
+u^r\partial_ru^z
+u^z\partial_zu^z
-\Delta_r u^z+\partial_zp
=O(\lambda^{-3-\delta}).
\end{equation}
In the radial equation,
\begin{equation}\label{orderP}
\partial_rp-\frac{(u^\theta)^2}{r}
=O(\lambda^{-3-2\delta}),
\end{equation}
whereas
\begin{equation}\label{orderUr}
\begin{aligned}
\partial_tu^r
+u^r\partial_ru^r
+u^z\partial_zu^r
-\left(\Delta_r-\frac1{r^2}\right)u^r
=O(\lambda^{-3}),\qquad 
\partial_z^2u^r=O(\lambda^{-3+2\delta}).
\end{aligned}
\end{equation}
The axial viscosity $\partial_z^2$ in the angular and axial equations is also
of lower order:
\begin{equation}\label{orderZZ}
\partial_z^2u^\theta=O(\lambda^{-3+\delta}),
\qquad
\partial_z^2u^z=O(\lambda^{-3+\delta}).
\end{equation}
Finally, the incompressibility constraint
\begin{equation}\label{orderIC}
\partial_r(ru^r)+\partial_z(ru^z)=0
\end{equation}
is imposed exactly, both terms having order $O(\lambda^{-1})$.

These orders dictate the construction and motivate the leading-order
Navier--Stokes system. In the angular and axial
equations, the time derivative, transport, and radial viscosity
all have the same weight $O(\lambda^{-3-\delta})$.
The axial viscosity $\partial_z^2$ gives terms smaller by a factor
$\lambda^{2\delta}$.
It is therefore left out of the leading-order profile equations.
In the radial equation, the pressure gradient $\partial_rp$ and the
centrifugal term $(u^\theta)^2/r$ each have formal order
$O(\lambda^{-3-2\delta})$ and dominate the remaining terms.

We introduce the following formal expansions at fixed similarity
coordinates $(R,Z)$:
\begin{equation}
\begin{aligned}
u^\alpha
&\sim \lambda^{-1-\delta}
\sum_{j=0}^\infty \lambda^{2j\delta}
U^\alpha_{(j)}(R,Z),
\qquad \alpha\in\{\theta,z\},\\
u^r
&\sim \lambda^{-1}
\sum_{j=0}^\infty \lambda^{2j\delta}
U^r_{(j)}(R,Z),\\
p
&\sim \lambda^{-2-2\delta}
\sum_{j=0}^\infty \lambda^{2j\delta}
P_{(j)}(R,Z).
\end{aligned}
\end{equation}
Here the coefficient profiles $U^\alpha_{(j)}$ and $P_{(j)}$
are independent of $\lambda$. The corresponding profile
expansions are
\begin{equation}
\begin{aligned}
U^\alpha
&\sim \sum_{j=0}^\infty \lambda^{2j\delta}
U^\alpha_{(j)}(R,Z),
\qquad \alpha\in\{r,\theta,z\},\\
P
&\sim \sum_{j=0}^\infty \lambda^{2j\delta}
P_{(j)}(R,Z).
\end{aligned}
\end{equation}
According to \eqref{orderP}, we impose the exact balance \emph{at leading order}:
\begin{equation}\label{radial-balance}
\partial_rp_{(0)}=\frac{(u^\theta_{(0)})^2}{r}.
\end{equation}
Together with the normalization $P_{(0)}(\infty,Z)=0$, this determines
$P_{(0)}$ from $U^\theta_{(0)}$ and leaves the radial residual
$\mathsf{NS}^r=O(\lambda^{-3})$, smaller than the leading angular
and axial terms. Incompressibility, together with regularity
at the axis, determines $U^r_{(0)}$ from $U^z_{(0)}$.
Thus it remains to construct $U^\theta_{(0)}$ and $U^z_{(0)}$.
Section~\ref{sec:leading-order-system} derives their leading equations.

\subsection{Residual, stress and shear}
\label{subsec:residual-stress-shear}
We define the momentum residual, stress, and radial shear.
For a pair $(u,p)$ of the form \eqref{solution-profile}, denote
the momentum residual and its cylindrical components by
$\mathsf{NS}=(\mathsf{NS}^r,\mathsf{NS}^\theta,\mathsf{NS}^z)$:
\begin{equation}\label{axiNSr}
\begin{cases}
\partial_tu^r+u^r\partial_ru^r+u^z\partial_zu^r
-\dfrac{(u^\theta)^2}{r}+\partial_rp
-\left(\Delta-\dfrac1{r^2}\right)u^r
=\mathsf{NS}^r,
\\[1mm]
\partial_tu^\theta+u^r\partial_ru^\theta+u^z\partial_zu^\theta
+\dfrac{u^ru^\theta}{r}
-\left(\Delta-\dfrac1{r^2}\right)u^\theta
=\mathsf{NS}^\theta,
\\[1mm]
\partial_tu^z+u^r\partial_ru^z+u^z\partial_zu^z
+\partial_zp-\Delta u^z
=\mathsf{NS}^z,
\\[1mm]
\partial_r(ru^r)+\partial_z(ru^z)=0.
\end{cases}
\end{equation}
In short, we write $\mathsf{NS} = \mathsf{NS}(u,p)$.

The aim is to express the angular and axial residuals in
divergence form, up to remainders:
\begin{equation}\label{divergenceS}
\begin{aligned}
\mathsf{NS}^\theta
&=-\left(\partial_r+\frac2r\right)\mathsf T^\theta
  +E^\theta,\\
\mathsf{NS}^z
&=-\left(\partial_r+\frac1r\right)\mathsf T^z
  +E^z.
\end{aligned}
\end{equation}
Here $\mathsf T^\theta$ and $\mathsf T^z$ denote the
$r\theta$ and $rz$ stress components, respectively.
Their profiles will be supported on a fixed annulus
$R_a\le R\le R_h$, with $0<R_a<R_h<\infty$.
This support requirement concerns the stress; the remainders
$E^\theta$ and $E^z$ need not vanish in the core $[0,R_a]$.

To explain the divergence structure, define the off-diagonal symmetric tensor
\[
\mathbb T_{\rm off}
=
\mathsf T^\theta
(e_r\otimes e_\theta+e_\theta\otimes e_r)
+
\mathsf T^z
(e_r\otimes e_z+e_z\otimes e_r).
\]
For axisymmetric coefficients,
\[
\nabla\cdot\mathbb T_{\rm off}
=
(\partial_z\mathsf T^z)e_r
+\left(\partial_r+\frac2r\right)\mathsf T^\theta e_\theta
+\left(\partial_r+\frac1r\right)\mathsf T^z e_z.
\]
Thus \eqref{divergenceS} is equivalent to
\[
\mathsf{NS}^\theta e_\theta+\mathsf{NS}^z e_z
=
-(\mathrm{Id}-e_r\otimes e_r)\nabla\cdot\mathbb T_{\rm off}
+E^\theta e_\theta+E^z e_z.
\]
Only the angular and axial components are prescribed by this
identity. The radial residual $\mathsf{NS}^r$ is treated separately since they are relatively lower order terms.

The radial residual remains part of the lower-order error; the divergence representation above concerns only the angular and axial components.

\subsection{Five integral moments}
Five cumulative moments express the residuals in divergence form.
Write
\[
\mathbf M(R,Z)=
\bigl(
M^\theta,M^z,M^{\theta z},M^{z\theta},M^p
\bigr)(R,Z),
\]
defined by
\begin{equation}\label{fiveM}
\begin{aligned}
M^\theta(R,Z)
&=\int_0^R\sqrt{2\rho}\,U^\theta(\rho,Z)\,d\rho,\\
M^z(R,Z)
&=\int_0^R U^z(\rho,Z)\,d\rho,\\
M^{\theta z}(R,Z)
&=\int_0^R\sqrt{2\rho}\,
  U^\theta(\rho,Z)U^z(\rho,Z)\,d\rho,\\
M^{z\theta}(R,Z)
&=\int_0^R
  \left((U^z)^2-\frac12(U^\theta)^2\right)(\rho,Z)\,d\rho,\\
M^p(R,Z)
&=\int_0^R\frac{(U^\theta)^2(\rho,Z)}{2\rho}\,d\rho.
\end{aligned}
\end{equation}
All five integrals start from the axis, so $\mathbf M(0,Z)=0$.
In terms of the regular variable $F$, these formulas become
\begin{equation}\label{moment-F}
\begin{aligned}
M^\theta(R,Z)
&=2\int_0^R\rho F(\rho,Z)\,d\rho,\\
M^{\theta z}(R,Z)
&=2\int_0^R\rho F(\rho,Z)U^z(\rho,Z)\,d\rho,\\
M^{z\theta}(R,Z)
&=\int_0^R
  \left((U^z(\rho,Z))^2-\rho F(\rho,Z)^2\right)\,d\rho,\\
M^p(R,Z)
&=\int_0^R F(\rho,Z)^2\,d\rho.
\end{aligned}
\end{equation}

\subsection{The profile problem and the order of construction}
We seek the five terminal moment identities and the admissible cone.
The leading profiles must also satisfy the exact centrifugal balance
\eqref{radial-balance} and
incompressibility. Their stress is supported in a fixed similarity
annulus. It depends on their cumulative moments, so it cannot be
prescribed independently.

We follow the roadmap in Section~\ref{sec:roadmap-construction}.
First we prepare the reference, its outer continuation, and the heat tail.
The pre-heat swirl fixes the axis value $P_0^{\rm pre}$, and the
outer moment correction makes the pressure vanish at infinity.
We then construct the core, connect the velocities, and restore the
five moments. The final shear modification gives the admissible cone
and preserves both the core and heat exterior.

The lower-order profiles enter through the expansion in
Section~\ref{sec:orders}. Their construction, summation, and cone
estimates are carried out in Section~\ref{sec:lower-order}.
The present detailed arguments
concern the leading profiles and their exact matching data. The radial
shear modification is itself a modification of these profiles.

\section{The leading order Navier-Stokes system}\label{sec:leading-order-system}
We derive the leading profile equations and the cone conditions.
For the leading field $(u_{(0)},p_{(0)})$, write
\[
\mathsf{NS}_{(0)}
=\mathsf{NS}(u_{(0)},p_{(0)}).
\]
The order analysis in Section~\ref{sec:orders} gives exact centrifugal
balance \eqref{radial-balance}. The angular and axial equations retain
radial viscosity and omit axial viscosity. Their residual identities are
\begin{equation}\label{eq:leading-tangential-residuals}
\begin{aligned}
\mathsf{NS}^\theta_{(0)}+\partial_z^2u^\theta_{(0)}
&=-\left(\partial_r+\frac2r\right)\mathsf T^\theta_{(0)},\\
\mathsf{NS}^z_{(0)}+\partial_z^2u^z_{(0)}
&=-\left(\partial_r+\frac1r\right)\mathsf T^z_{(0)}.
\end{aligned}
\end{equation}
Here, the leading physical stress components have the form
\begin{equation}\label{eq:leading-stress-scaling}
\mathsf T^\theta_{(0)}
=\lambda^{-2-\delta}\mathcal{T}^\theta_{(0)}(R,Z),
\qquad
\mathsf T^z_{(0)}
=\lambda^{-2-\delta}\mathcal{T}^z_{(0)}(R,Z),
\end{equation}
where
$\mathcal{T}_{(0)}
=(\mathcal{T}^\theta_{(0)},\mathcal{T}^z_{(0)})$
is independent of $\lambda$.  
At the leading-profile stage the remainders in
\eqref{divergenceS} are exactly
\begin{equation}\label{eq:error-z2}
E^\theta_{(0)}
=-\partial_z^2u^\theta_{(0)},
\qquad
E^z_{(0)}
=-\partial_z^2u^z_{(0)}.
\end{equation}
By \eqref{orderZZ}, these terms have formal order
$O(\lambda^{-3+\delta})$, smaller by a factor
$\lambda^{2\delta}$ than the leading angular and axial terms. Reducing the remaining errors to arbitrarily high order requires
higher-order corrections to the background profiles.
Section~\ref{sec:lower-order} constructs and sums these corrections.

In physical coordinates, the leading order unforced Navier-Stokes system reads 
\begin{equation}\label{leading-unforced-NS}
  \begin{cases}
  \mathsf{NS}^{\theta}(u_{(0)},p_{(0)}) +\partial_z^2u^\theta_{(0)}=   \mathsf{NS}^{z}(u_{(0)},p_{(0)})+\partial_z^2u^z_{(0)} = 0, \\
\partial_rp_{(0)}=\frac{(u^\theta_{(0)})^2}{r}, \quad \partial_r (r u^r_{(0)}) + \partial_z (r u^z_{(0)}) = 0.
  \end{cases}
\end{equation} 
The leading constructions in Subsection~\ref{sec:analytic-core} and in
Section~\ref{sec:heat-exterior}
solve this system in the core ($R\le R_a$) and exterior ($R\ge R_h$), respectively.
The exterior moment matching is verified in Sections~\ref{subsec:outer-waiting} and~\ref{subsec:outer-axial-matching}--\ref{subsec:outer-cone-closure}.
Both regions are stress-free: $\mathcal{T}_{(0)}=0$.

For brevity, we \textbf{suppress the subscript $(0)$ in the calculations below
throughout the leading-profile construction}; all profiles and derived quantities are understood
to be leading-order objects, \emph{e.g.}, $\mathbf M(R,Z) = \mathbf M_{(0)}(R,Z)$ consists of the five moments of $u_{(0)}$.

\subsection{\texorpdfstring{Determination of $P$ and $U^r$}{Determination of P and radial velocity}} \label{sec:Ur-P}

We recover pressure from the swirl and radial velocity from the axial flow.
The radial balance and chain rule give $2R\partial_R P=(U^\theta)^2$.
Using $P(\infty,Z)=0$, integration gives
\begin{equation}\label{P-Utheta}
P(R,Z)=-\int_R^\infty F^2(\rho, Z)\,d\rho = P(0,Z)
+\int_0^R F^2(\rho, Z) \,d\rho.
\end{equation}
Thus $U^\theta$ determines $P$ uniquely. On the axis,
\begin{equation}
  P(0, Z) = - M^p(\infty,Z) < 0.
\end{equation}
One can also write \eqref{P-Utheta} as
\begin{equation} \label{P-moment}
  P(R,Z)=M^p(R,Z)-M^p(\infty,Z).
\end{equation}

Next, the incompressibility constraint \eqref{orderIC} in self-similar coordinates becomes
\begin{equation} \label{Ur-radial}
L_\delta \partial_R\bigl(\sqrt{2R}\,U^r\bigr)
=
(1+\delta)Z U^z
-L\partial_ZU^z
+2ZR\partial_RU^z,
\end{equation}
where we use the abbreviations $L_\delta=1-\delta Z^2 \asymp 1$ and $L=1-Z^2 \in [0,1]$. Smoothness of the physical fields at the axis requires $U^r(0, Z) = 0$, hence integrating \eqref{Ur-radial} from the axis gives
\begin{equation}\label{Ur-Uz}
U^r(R,Z)
=\frac{
2ZR U^z(R,Z)
-(1-\delta)Z M^z(R,Z)
-L\partial_ZM^z(R,Z)
}{L_\delta\sqrt{2R}}.
\end{equation}

In conclusion, given $U^\theta$, $U^z$, the formulas \eqref{P-Utheta} and \eqref{Ur-Uz} uniquely determine $P$ and $U^r$,
with radial pressure balance and incompressibility satisfied
exactly.

\subsection{The stress, shear and radial equations}
\label{subsec:stress-shear-radial}

We separate inertial stress from radial shear.
Write $\mathcal{T}=\mathcal{I}+\mathcal{S}$.
The term $\mathcal{I}$ contains the inviscid terms and pressure gradient;
$\mathcal{S}$ contains radial viscosity.
Using the pressure balance and incompressibility  \eqref{leading-unforced-NS}$_2$, we obtain
\begin{equation}\label{inertial-stress}
\begin{aligned}
&\partial_tu^\theta
 +u^r\partial_ru^\theta+u^z\partial_zu^\theta
 +\frac{u^ru^\theta}{r} - \frac{(u^\theta)^2}{r} + \partial_r p\\
&\qquad
 =\partial_tu^\theta+\partial_z(u^zu^\theta)
  +\left(\partial_r+\frac2r\right)(u^ru^\theta)\\
&\qquad
 =-\left(\partial_r+\frac2r\right)
   \left[\lambda^{-2-\delta}\mathcal{I}^\theta(R,Z)\right],
\\[1.5ex]
&\partial_tu^z
 +u^r\partial_ru^z+u^z\partial_zu^z+\partial_zp\\
&\qquad
 =\partial_tu^z+\partial_z\bigl((u^z)^2+p\bigr)
  +\left(\partial_r+\frac1r\right)(u^ru^z)\\
&\qquad
 =-\left(\partial_r+\frac1r\right)
   \left[\lambda^{-2-\delta}\mathcal{I}^z(R,Z)\right].
\end{aligned}
\end{equation}
The radial viscous terms satisfy
\begin{equation}\label{viscous-stress}
\begin{aligned}
-\left(\Delta_r-\frac1{r^2}\right)u^\theta
&=-\left(\partial_r+\frac2r\right)
  \left[\lambda^{-2-\delta}\mathcal{S}^\theta(R,Z)\right],\\
-\Delta_ru^z
&=-\left(\partial_r+\frac1r\right)
  \left[\lambda^{-2-\delta}\mathcal{S}^z(R,Z)\right].
\end{aligned}
\end{equation}

Let $Q^\theta,Q^z$ denote the two physical inertial source terms
on the left of \eqref{inertial-stress}. Define their normalized sources by
\[
 N^j=-\lambda^{3+\delta}\sqrt{R/2}\,Q^j,
 \qquad j=\theta,z.
\]
The radial chain rule then gives \eqref{inviscid-stress-radial} below.
These sources can be expressed in terms
of the axial moments (the pressure $P$ is also determined by the $M^p$ moment through \eqref{P-moment}) as
\begin{equation}\label{inviscid-stress-sources-Ur}
\begin{aligned}
N^\theta
={}&-\frac{\sqrt{R/2}}{L_\delta}
\left[
\frac{1+\delta}{2}U^\theta
+\frac{1-\delta}{2}Z\partial_ZU^\theta
+R\partial_RU^\theta
\right]
\\
&-\frac{\sqrt{R/2}}{L_\delta}
\left[
-2(1+\delta)ZU^zU^\theta
+L\,\partial_Z(U^zU^\theta)
-2ZR\,\partial_R(U^zU^\theta)
\right]
\\
&-\left[
R\partial_R(U^rU^\theta)+U^rU^\theta
\right],
\\[2ex]
N^z
={}&-\frac{\sqrt{R/2}}{L_\delta}
\left[
\frac{1+\delta}{2}U^z
+\frac{1-\delta}{2}Z\partial_ZU^z
+R\partial_RU^z
\right]
\\
&-\frac{\sqrt{R/2}}{L_\delta}
\left[
\begin{aligned}
&-2(1+\delta)Z\bigl((U^z)^2+P\bigr)
+L\,\partial_Z\bigl((U^z)^2+P\bigr)
\\
&-2ZR\,\partial_R\bigl((U^z)^2+P\bigr)
\end{aligned}
\right]
\\
&-\left[
R\partial_R(U^rU^z)+\frac12U^rU^z
\right].
\end{aligned}
\end{equation}
The three lines of $N^z$ come respectively from
$\partial_tu^z$, $\partial_z(u^z)^2+\partial_zp$, and
$\bigl(\partial_r+\frac1r\bigr)(u^ru^z)$.
The calculation for $N^\theta$ is similar.
Both use Lemma~\ref{lem:chain}.

By \eqref{inertial-stress} and the chain rules \eqref{profile-chain-rule}, the inertial stress satisfies the radial equations
\begin{equation}\label{inviscid-stress-radial}
R\partial_R\mathcal{I}^\theta+\mathcal{I}^\theta=N^\theta,
\qquad
R\partial_R\mathcal{I}^z+\frac12\mathcal{I}^z=N^z.
\end{equation}
Regularity of the physical stress $\mathsf T$ at the axis requires (see \eqref{smooth-axis})
\[
\mathcal{I}^\theta(R,Z)=O(R),\qquad
\mathcal{I}^z(R,Z)=O(\sqrt R),
\qquad \mbox{as} \ \ R\downarrow0.
\]
This condition together with the radial equations \eqref{inviscid-stress-radial} give
\[
\mathcal{I}^\theta(R,Z)
=\frac1R\int_0^R N^\theta(\rho,Z)\,d\rho,
\qquad
\mathcal{I}^z(R,Z)
=\frac1{\sqrt R}\int_0^R\frac{N^z(\rho,Z)}{\sqrt\rho}\,d\rho,
\]
which yields
\begin{align}
\mathcal{I}^\theta
={}&
\frac{
\left(1-\frac{\delta}{2}\right)M^\theta
-\frac{1-\delta}{2}Z\partial_ZM^\theta
-R\sqrt{2R}\,U^\theta
}{2L_\delta R}
\nonumber\\
&+
\frac{
(2\delta-1)ZM^{\theta z}
-L\,\partial_ZM^{\theta z}
+2R\sqrt{2R}\,ZU^zU^\theta
}{2L_\delta R}
 - U^rU^\theta,
\label{Itheta-moments-direct}\\[2ex]
\mathcal{I}^z
={}&
\frac{
\frac{1-\delta}{2}
\left(M^z-Z\partial_ZM^z\right)-RU^z
}{L_\delta\sqrt{2R}}
\nonumber\\
&+
\frac{
2\delta ZM^{z\theta}
-L\,\partial_ZM^{z\theta}
+R\left(
2(1+\delta)ZP-L\,\partial_ZP+2Z(U^z)^2
\right)
}{L_\delta\sqrt{2R}}
 - U^rU^z.
\label{Iz-moments-direct}
\end{align}
Substituting \eqref{Ur-Uz} and using
$P(R,Z)=M^p(R,Z)-M^p(\infty,Z)$, we obtain the equivalent
moment representation
\begin{align}
\mathcal{I}^\theta
={}&\frac{U^\theta}{L_\delta\sqrt{2R}}
\left[-R+(1-\delta)ZM^z+L\,\partial_ZM^z\right]
\nonumber\\
&+\frac1{2L_\delta R}\left[
\left(1-\frac{\delta}{2}\right)M^\theta
-\frac{1-\delta}{2}Z\partial_ZM^\theta
-L\,\partial_ZM^{\theta z}
+(2\delta-1)ZM^{\theta z}
\right],
\label{Itheta-moments}\\[1ex]
\mathcal{I}^z
={}&\frac1{L_\delta\sqrt{2R}}\Bigg\{
\left[-R+(1-\delta)ZM^z+L\,\partial_ZM^z\right]U^z
\nonumber\\
&\qquad\qquad+\frac{1-\delta}{2}
\left(M^z-Z\partial_ZM^z\right)
+2\delta ZM^{z\theta}-L\,\partial_ZM^{z\theta}
\nonumber\\
&\qquad\qquad+R\left[2(1+\delta)Z-L\,\partial_Z\right]
\left[M^p(R,Z)-M^p(\infty,Z)\right]
\Bigg\}.
\label{Iz-moments}
\end{align}
The moment formulas \eqref{Itheta-moments}--\eqref{Iz-moments} are useful for matching profiles.
The radial equations \eqref{inviscid-stress-radial} propagate bounds between radii.
For instance, an upper bound for $\mathcal{I}^\theta(R_0,Z)$ and an upper
bound for $N^\theta$ on $[R_0,R_1]$ give an upper bound for
$\mathcal{I}^\theta$ throughout that interval.

On the other hand, the radial  viscosity part $\mathcal{S}$ has a relatively simple form:
\begin{equation}\label{leading-shear-profile}
  \begin{aligned}
\mathcal{S} = (\mathcal{S}^\theta, \mathcal{S}^z) 
&=\bigg(
\sqrt{2R}\,\partial_RU^\theta-\frac{U^\theta}{\sqrt{2R}},
\quad \sqrt{2R}\,\partial_RU^z
\bigg) \\
& = \left( 2R \partial_R F, \quad  \sqrt{2R} \partial_R U^z \right).
  \end{aligned}
\end{equation}
We call $\mathcal{S}$ the \emph{shear} for simplicity. 

For a positive function $G$, we call $R\partial_R\log G$ its
\emph{logarithmic radial slope}, or radial slope for short.
For example, $\sqrt R$ has radial slope $1/2$.
The angular shear involves this slope of $U^\theta$.
The axial shear is $\sqrt{2R}\,\partial_RU^z$.
It uses the radial derivative directly and does not require division by $U^z$.

\subsection{The stress-free system}\label{subsec:stress-free}

We derive the equations for the stress-free core and exterior.
Since $\mathcal{T}=\mathcal{I}+\mathcal{S}$, vanishing stress means
\[
\mathcal{S} = - \mathcal{I},
\]
or, equivalently, for $R>0$,
\begin{equation}\label{eq:stress-free-profile}
\begin{aligned}
2R\,\partial_R F 
&=-\mathcal{I}^\theta,\\
\sqrt{2R}\,\partial_RU^z
&=-\mathcal{I}^z.
\end{aligned}
\end{equation}
Together with the moment formulas for $\mathcal{I}$, these are
nonlinear nonlocal equations for $U^\theta$ and $U^z$.
 Applying $R\partial_R+1$ to the first equation and
$R\partial_R+\frac12$ to the second, and using
\eqref{inviscid-stress-radial}, we obtain
\begin{equation}\label{eq:stress-free-sources}
\begin{aligned}
2R\left(R\partial_R^2F+2\partial_RF\right)
&=-N^\theta,\\
\sqrt{2R}\left(R\partial_R^2U^z+\partial_RU^z\right)
&=-N^z.
\end{aligned}
\end{equation}
As explained in Section \ref{sec:Ur-P}, the profiles $P$ and $U^r$ are recovered from radial centrifugal
balance and incompressibility.

At the axis, we require $F = U^\theta/\sqrt{2R}$ and $U^z$ to extend
smoothly to $R=0$. Their axis values are the data for the core
construction on $[0,R_a]$, while the pressure datum $P(0,Z)$ is supplied by
the outer construction. The radial derivatives must be compatible
with \eqref{eq:stress-free-profile}; they are not independent data.

\subsection{The cone conditions} \label{sec:cone-condition}

We state the admissible and relaxed cone conditions.
On the stress annulus $R_a\le R\le R_h$, we require
$U^\theta>0$ and $\mathcal{S}^\theta<0$. Define
\[
\kappa=-\frac{\sqrt{2R}\,|\mathcal{S}|^2}
{U^\theta \mathcal{S}^\theta} = - \frac{|\mathcal{S}|^2}{F \mathcal{S}^\theta},
\]
and write
\[
\mathcal{S}^\perp=(-\mathcal{S}^z,\mathcal{S}^\theta).
\]
The \emph{admissible cone condition} that we impose on $(R_a, R_h) \times [-1, 1]$ is
\begin{equation}\label{cone-admissible}
  \boxed{
\mathcal{T}\cdot \mathcal{S}<0,\qquad
(\kappa-2)(\mathcal{T}\cdot \mathcal{S}^\perp)^2
<2(\mathcal{T}\cdot \mathcal{S})^2,\qquad
\kappa>2.
  }
\end{equation}
The last inequality is the shear-strength condition, motivated by the
pulse correction mechanism developed in Part~II and by studies of
centrifugal instability in columnar vortices and swirling jets
\cite{LeibovichStewartson,BillantGallaire}. The other inequalities restrict the
stress direction. Since \(U^\theta>0\) and \(\mathcal{S}^\theta<0\), we have
\[
\kappa>2
\quad\Longleftrightarrow\quad
\left(\mathcal{S}^\theta+F\right)^2
+(\mathcal{S}^z)^2
>F^2.
\]
Thus the shear $\mathcal{S}$ lies outside the circle centered at
\((-F,0)\), with radius
$F$ (we always have $F\ge 0$ during the construction), and remains in the left half-plane, see Figure \ref{fig:shear-circle}. 
More generally, for fixed $F>0$ and $k>0$, the level $\kappa=k$ is
\begin{equation}\label{eq:shear-level-circle}
 \left(\mathcal{S}^\theta+\frac{kF}{2}\right)^2
 +(\mathcal{S}^z)^2=\left(\frac{kF}{2}\right)^2,
 \qquad \mathcal{S}^\theta<0.
\end{equation}
Its center is $(-kF/2,0)$ and its radius is $kF/2$; the origin is
excluded. Thus the circles for different $k$ are tangent at the origin.
In the following figures, the dashed circle is the
reference level $\kappa=2$, and the solid circle is the level of the
displayed shear. All circle diagrams are in the shear plane.
For fixed shear, the middle inequality in \eqref{cone-admissible} requires the stress to lie in an open cone around
$-\mathcal{S}$, with half opening angle
$\arctan\sqrt{2/(\kappa-2)}$, see Figure \ref{fig:stress-cone}.

During the construction we first achieve the \emph{relaxed cone condition}:
\begin{equation}\label{cone-relaxed}
  \boxed{
\mathcal{T} \cdot \mathcal{S} < 0, \qquad \begin{cases}
(\kappa-2)(\mathcal{T}\cdot \mathcal{S}^\perp)^2
<2(\mathcal{T}\cdot \mathcal{S})^2, \quad {\rm if} \ \ \kappa >2, \\
-\mathcal{T}\cdot \mathcal{S}
>-F \mathcal{S}^\theta\,(2-\kappa), \quad {\rm if} \ \ \kappa \le 2.
\end{cases}
  }
\end{equation}
For $\kappa>2$, this is exactly the
admissible cone condition. For $\kappa<2$, the second inequality is
stronger than the direction condition $\mathcal{T}\cdot \mathcal{S}<0$:
it prescribes the positive lower bound $-F\mathcal{S}^\theta(2-\kappa)$
for $-\mathcal{T}\cdot \mathcal{S}$. At $\kappa=2$, it reduces to the
direction condition. See Figure \ref{fig:stress-half-plane}.
This relaxation allows us to connect profiles with
different radial growth or decay rates. A later shear modification restores
the condition $\kappa>2$.

The admissible cone controls stress direction and shear strength.
For $\kappa\le2$, the relaxed cone uses the stated lower bound on
$-\mathcal{T}\cdot \mathcal{S}$.

\begin{figure}[htbp]
\centering
\begin{tikzpicture}[x=1.15cm,y=1.15cm,>=stealth,font=\small]
  \fill[blue!8] (-4.2,-2.3) rectangle (0,2.3);
  \fill[white] (-1,0) circle (1);
  \draw[->,gray!80] (-4.4,0) -- (.8,0)
    node[right,black] {$\mathcal{S}^\theta$};
  \draw[->,gray!80] (0,-2.45) -- (0,2.5)
    node[above,black] {$\mathcal{S}^z$};
  \draw[dashed,thick,blue!65!black] (-1,0) circle (1);
  \draw[very thick,purple!75!black] (-2,0) circle (2);
  \fill (-1,0) circle (1.3pt);
  \fill[purple!75!black] (-2,0) circle (1.3pt);
  \node[below=2pt] at (-1,0) {$-F$};
  \node[below=2pt] at (-2,0) {$-2F$};
  \draw[<->,blue!65!black] (-1,.04) -- (-1,.96)
    node[midway,right] {$F$};
  \draw[<->,purple!75!black] (-2,.04) -- (-2,1.96)
    node[midway,left] {$2F$};
  \node[purple!75!black,above] at (-2,2.05) {$\kappa=4$};
  \node[blue!65!black,anchor=west] at (.5,1.35) {$\kappa=2$};
  \draw[->,blue!65!black] (.45,1.2) -- (-.35,.76);
  \node[blue!65!black,align=center] at (-3.15,-.7)
    {allowed\\shears\\$\kappa>2$};
  \draw[->,thick] (0,0) -- (-2,-2);
  \fill[purple!75!black] (-2,-2) circle (2pt);
  \node[below left=2pt] at (-2,-2)
    {$\mathcal{S}=(-2,-2)$};
  \draw[fill=white] (0,0) circle (1.6pt);
  \node[below right=2pt] at (0,0) {$0$};
\end{tikzpicture}
\caption{The shear plane, drawn with $F=1$. The dashed circle has
center $(-F,0)$ and radius $F$ ($\kappa=2$); the solid circle has
center $(-2F,0)$ and radius $2F$ ($\kappa=4$).
The displayed shear lies on the solid circle.
Admissible shears lie in the shaded region outside the dashed circle
and in $\mathcal{S}^\theta<0$; both boundaries are excluded.}
\label{fig:shear-circle}
\end{figure}

\begin{figure}[htbp]
\centering
\begin{minipage}[t]{.47\linewidth}
\centering
\textbf{Shear plane: $F=1$, $\kappa=8$}\par\smallskip
\begin{tikzpicture}[x=.44cm,y=.44cm,>=stealth,font=\footnotesize]
  \draw[->,gray!80] (-8.6,0) -- (1.2,0)
    node[right,black] {$\mathcal{S}^\theta$};
  \draw[->,gray!80] (0,-4.8) -- (0,4.8)
    node[above,black] {$\mathcal{S}^z$};
  \draw[dashed,thick,blue!65!black] (-1,0) circle (1);
  \draw[very thick,purple!75!black] (-4,0) circle (4);
  \fill (-1,0) circle (1.3pt);
  \fill[purple!75!black] (-4,0) circle (1.3pt);
  \node[below=2pt] at (-4,0) {$-4$};
  \node[blue!65!black,anchor=east] at (-2.3,1.2) {$\kappa=2$};
  \draw[->,blue!65!black] (-2.2,1.05) -- (-1.7,.7);
  \node[purple!75!black] at (-5.8,2.3) {$\kappa=8$};
  \draw[->,thick] (0,0) -- (-4,-4);
  \fill[purple!75!black] (-4,-4) circle (2pt);
  \node[below left=2pt] at (-4,-4) {$\mathcal{S}=(-4,-4)$};
  \draw[fill=white] (0,0) circle (1.5pt);
  \node[below right=2pt] at (0,0) {$0$};
\end{tikzpicture}
\end{minipage}\hfill
\begin{minipage}[t]{.49\linewidth}
\centering
\textbf{Stress plane: $2\alpha=\pi/3$}\par\smallskip
\begin{tikzpicture}[x=1.03cm,y=1.03cm,>=stealth,font=\footnotesize]
  \fill[blue!8] (0,0) -- (3.25,.8708) -- (3.25,3.25)
    -- (.8708,3.25) -- cycle;
  \draw[->,gray!80] (-.4,0) -- (3.6,0)
    node[right,black] {$\mathcal{T}^\theta$};
  \draw[->,gray!80] (0,-.4) -- (0,3.6)
    node[above,black] {$\mathcal{T}^z$};
  \draw[dashed,thick,blue!65!black] (0,0) -- (3.3,.8842);
  \draw[dashed,thick,blue!65!black] (0,0) -- (.8842,3.3);
  \draw[->,thick] (0,0) -- (2,2)
    node[above right] {direction $-\mathcal{S}$};
  \draw[->,very thick,blue!65!black] (0,0) -- (2.9,1.15)
    node[above right] {$\mathcal{T}$};
  \draw[<->] (45:1.25) arc (45:75:1.25);
  \node at (60:1.55) {$\alpha$};
  \node[blue!65!black] at (2.15,2.95) {allowed stresses};
  \draw[fill=white] (0,0) circle (1.5pt);
  \node[below left=2pt] at (0,0) {$0$};
  \path (-.4,-.85) -- (3.6,-.85);
\end{tikzpicture}
\end{minipage}
\caption{A fixed shear with $\kappa=8>2$.
The left panel distinguishes the dashed $\kappa=2$ circle from
the solid circle through $\mathcal{S}=(-4,-4)$.
The right panel shows the admissible stress cone around $-\mathcal{S}$,
with half opening angle
$\alpha=\arctan\sqrt{2/(\kappa-2)}=\pi/6$.
Its full opening is $\pi/3$.
The dashed boundary rays and the origin are excluded.}
\label{fig:stress-cone}
\end{figure}

\begin{figure}[htbp]
\centering
\begin{minipage}[t]{.47\linewidth}
\centering
\textbf{Shear plane: $F=1$, $\kappa=1$}\par\smallskip
\begin{tikzpicture}[x=1.6cm,y=1.6cm,>=stealth,font=\footnotesize]
  \draw[->,gray!80] (-2.2,0) -- (.5,0)
    node[right,black] {$\mathcal{S}^\theta$};
  \draw[->,gray!80] (0,-1.3) -- (0,1.3)
    node[above,black] {$\mathcal{S}^z$};
  \draw[dashed,thick,blue!65!black] (-1,0) circle (1);
  \draw[very thick,purple!75!black] (-.5,0) circle (.5);
  \fill (-1,0) circle (1.3pt);
  \fill[purple!75!black] (-.5,0) circle (1.3pt);
  \node[below=2pt] at (-1,0) {$-1$};
  \node[above=2pt] at (-.5,0) {$-\tfrac12$};
  \node[blue!65!black] at (-1.3,.6) {$\kappa=2$};
  \node[purple!75!black,anchor=east] at (-.85,-.75) {$\kappa=1$};
  \draw[->,thick] (0,0) -- (-.8,-.4);
  \fill[purple!75!black] (-.8,-.4) circle (2pt);
  \node[below=2pt] at (-1,-1.08)
    {$\mathcal{S}=(-\tfrac45,-\tfrac25)$};
  \draw[fill=white] (0,0) circle (1.5pt);
  \node[above right=2pt] at (0,0) {$0$};
\end{tikzpicture}
\end{minipage}\hfill
\begin{minipage}[t]{.49\linewidth}
\centering
\textbf{Stress plane: $2\mathcal{T}^\theta+\mathcal{T}^z>2$}\par\smallskip
\begin{tikzpicture}[x=1.03cm,y=1.03cm,>=stealth,font=\footnotesize]
  \fill[blue!8] (-.45,2.9) -- (-.45,3.2) -- (3.2,3.2)
    -- (3.2,-.45) -- (1.225,-.45) -- cycle;
  \draw[->,gray!80] (-.45,0) -- (3.55,0)
    node[right,black] {$\mathcal{T}^\theta$};
  \draw[->,gray!80] (0,-.45) -- (0,3.55)
    node[above,black] {$\mathcal{T}^z$};
  \draw[dashed,thick,blue!65!black] (-.45,2.9) -- (1.225,-.45);
  \node[below=2pt] at (1,0) {$1$};
  \node[left=2pt] at (0,2) {$2$};
  \draw[->,thick] (0,0) -- (2.7,1.35);
  \node[align=center] at (2.35,.65) {direction\\$-\mathcal{S}$};
  \draw[->,very thick,blue!65!black] (0,0) -- (1.1,2.45)
    node[above left] {$\mathcal{T}$};
  \node[blue!65!black] at (2.1,2.95) {allowed stresses};
  \fill (0,0) circle (1.3pt);
  \node[below left=2pt] at (0,0) {$0$};
\end{tikzpicture}
\end{minipage}
\caption{A fixed shear with $\kappa=1<2$.
The left panel shows the dashed $\kappa=2$ circle and the solid
$\kappa=1$ circle through $\mathcal{S}=(-4/5,-2/5)$.
The right panel shows the relaxed stress condition
$-\mathcal{T}\cdot \mathcal{S}>
-F\mathcal{S}^\theta(2-\kappa)=4/5$,
equivalently $2\mathcal{T}^\theta+\mathcal{T}^z>2$.
The dashed boundary line is excluded; its inward normal is $-\mathcal{S}$.}
\label{fig:stress-half-plane}
\end{figure}

\clearpage
\subsection{Five-moment matching and exterior fields}
\label{subsec:moment-gluing-principle}
Matching five moments preserves the exterior fields.

\begin{lemma}[Matching moments and matching fields]
\label{lem:moment-gluing}
Let $0<R_*<R_{**}$. Suppose two smooth angular and axial profiles
are equal on $[R_*,R_{**}]$ for every $Z\in[-1,1]$.
Suppose they use the same axis pressure $P_0$.
If their five cumulative moments are equal at $R_*$ as functions of $Z$,
then their pressure, radial velocity, inertial stress, and shear are equal
on $[R_*,R_{**}]$. If the angular and axial profiles are equal for all
$R\ge R_*$, all these fields are equal there as well.
\end{lemma}
\begin{proof}
Use the subscripts $j=1,2$ for the two profiles and their derived fields.
Write $\mathbf M_j$ for their five-moment vectors.
For $R\in[R_*,R_{**}]$, the definitions \eqref{fiveM} give
\[
\mathbf M_j(R,Z)=\mathbf M_j(R_*,Z)+\int_{R_*}^R
\begin{pmatrix}
\sqrt{2\rho}\,U_j^\theta(\rho,Z)\\
U_j^z(\rho,Z)\\
\sqrt{2\rho}\,U_j^\theta(\rho,Z)U_j^z(\rho,Z)\\
(U_j^z(\rho,Z))^2-\tfrac12(U_j^\theta(\rho,Z))^2\\
(U_j^\theta(\rho,Z))^2/(2\rho)
\end{pmatrix}\,d\rho.
\]
The initial moment vectors and the integrands are equal. Hence
\[
\mathbf M_1(R,Z)=\mathbf M_2(R,Z),\qquad
\partial_Z\mathbf M_1(R,Z)=\partial_Z\mathbf M_2(R,Z).
\]
In particular,
\[
P_1(R,Z)=P_0(Z)+M_1^p(R,Z)
        =P_0(Z)+M_2^p(R,Z)=P_2(R,Z).
\]
Substitution into \eqref{Ur-Uz}, \eqref{Itheta-moments-direct},
and \eqref{Iz-moments-direct} gives
\[
U_1^r=U_2^r,\qquad
\mathcal{I}_1^\theta=\mathcal{I}_2^\theta,\qquad
\mathcal{I}_1^z=\mathcal{I}_2^z.
\]
The angular and axial profiles have equal radial derivatives on this
interval, including the endpoints by smoothness. Thus \eqref{leading-shear-profile} gives
\[
\mathcal{S}_1^\theta=\mathcal{S}_2^\theta,\qquad
\mathcal{S}_1^z=\mathcal{S}_2^z.
\]
If the angular and axial profiles are equal for all $R\ge R_*$,
the same integral identities hold for every such $R$.
\end{proof}

\clearpage
\paragraph{Moment correction by bumps.}
In the actual constructions, we will enforce moment matching by adding finitely many smooth,
compactly supported bumps to the angular and axial velocity profiles.
Their coefficients control the changes in the moments, while these
profiles outside the bump supports remain unchanged. Since the moments
are at most quadratic in the velocities, this leads to a linear system with a quadratic
remainder. The following two elementary facts explain the mechanism.

\begin{lemma}[A generalized Vandermonde moment matrix]
\label{lem:bump-moment-matrix}
Let $\alpha_1<\cdots<\alpha_m$ be real numbers, and let
$b_1,\ldots,b_m\in C_c^\infty((0,\infty))$ be nonnegative, nonzero
bumps with strictly ordered supports:
$\sup\operatorname{supp}b_j<\inf\operatorname{supp}b_{j+1}$.
Then the matrix
\[
 A_{ij}=\int_0^\infty x^{\alpha_i}b_j(x)\,dx,
 \qquad 1\le i,j\le m,
\]
is invertible.
\end{lemma}
\begin{proof}
For $0<x_1<\cdots<x_m$, the evaluation determinant
$\det(x_j^{\alpha_i})_{i,j}$ is nonzero. Indeed, setting $x=e^t$
reduces this to the fact that a nonzero combination of $m$ exponentials
with distinct real exponents has at most $m-1$ distinct zeros.
This follows by induction: divide by the exponential with smallest exponent,
differentiate, and apply Rolle's theorem. By continuity, the determinant
has constant sign on the set of ordered positive points. Multilinearity
and Fubini's theorem give
\[
 \det A=\int_{(0,\infty)^m}
 \det(x_j^{\alpha_i})_{i,j}\prod_{j=1}^m b_j(x_j)\,
 dx_1\cdots dx_m\ne0,
\]
since the integrand has one nonzero sign wherever all bumps are positive.
\end{proof}

\begin{lemma}[A small quadratic correction]
\label{lem:quadratic-moment-correction}
Fix a norm on $\mathbb R^m$ and its induced matrix norm.
Let $A$ be invertible, put $B=\|A^{-1}\|$, and let
$Q:\mathbb R^m\times\mathbb R^m\to\mathbb R^m$ be bilinear with
$\|Q(h,k)\|\le C\|h\|\|k\|$.
If $e=\|d\|$ satisfies $8B^2Ce\le1$, then
\[
 Ah+Q(h,h)=-d
\]
has a unique solution in the closed ball $\|h\|\le2Be$.
\end{lemma}
\begin{proof}
Put $r=2Be$ and $\Phi(h)=-A^{-1}[d+Q(h,h)]$.
For $\|h\|,\|k\|\le r$,
\[
 \|\Phi(h)\|\le Be+BCr^2\le\tfrac34r,
 \qquad
 \|\Phi(h)-\Phi(k)\|\le2BCr\|h-k\|\le\tfrac12\|h-k\|.
\]
Thus $\Phi$ is a contraction of the closed ball into itself.
\end{proof}

The relevant smallness is that of the moment defect relative to the
inverse matrix bound, after any scale normalization. The first lemma
gives no uniform inverse bound as the exponents or bump supports vary.
For a linear moment system, $h=-A^{-1}d$ requires no smallness assumption.

\section{Roadmap of the leading-order construction}\label{sec:roadmap-construction}

We first prescribe the whole candidate profile by heuristic reasoning.
We use a temporary reference profile on $[0,R_{\rm ref}]$.
We then prescribe it piecewise on $[R_{\rm ref},R_h]$ and continue it to an
exact heat flow on $[R_h,\infty)$. The pre-heat swirl determines the
axis datum $P_0^{\rm pre}$ by \eqref{intro:preheat-axis-pressure}, and the
outer moment correction matches the pressure from this datum to the heat exterior.

For input scales $P_*\ge1$ and $R_{\rm ref}>1$, prescribe the reference
profile on $0<R\le R_{\rm ref}$ by
\[
 U^\theta_{\rm ref}(R,Z)=\frac{P_*}{1+Z^2}
 \left(\frac{R}{R_{\rm ref}}\right)^{1/10},
 \qquad U^z_{\rm ref}(R,Z)=4Z.
\]
Then, at $R=R_{\rm ref}$, the five moments of the candidate profiles are
\[
\begin{gathered}
 M^z=4ZR_{\rm ref},\qquad M^{\theta z}=4ZM^\theta,\qquad
 M^{z\theta}=16Z^2R_{\rm ref}-\tfrac5{12}R_{\rm ref}(U^\theta_{\rm ref})^2,\\
 M^\theta=\tfrac58R_{\rm ref}\sqrt{2R_{\rm ref}}U^\theta_{\rm ref},\qquad
 M^p=\tfrac52(U^\theta_{\rm ref})^2.
\end{gathered}
\]
All quantities in this display are evaluated at $(R_{\rm ref},Z)$.
At infinity, the moments must satisfy \eqref{eq:moment-conditions}.
We then replace the inner reference profile on $[0,R_{\rm useless}]$,
where $R_{\rm useless}=e^{-5}R_{\rm ref}$, preserving its five moments
at $R_{\rm useless}$. We retain the reference branch on
$[R_{\rm useless},R_{\rm ref}]$ and the prepared outer profile on
$[R_{\rm ref},\infty)$.
At this stage, we have the relaxed cone condition \eqref{cone-relaxed}.
To obtain the admissible cone condition \eqref{cone-admissible}, we still
need a shear modification followed by another moment repair.

\definecolor{roadblue}{RGB}{25,85,150}
\definecolor{roadorange}{RGB}{195,105,30}
\definecolor{roadpurple}{RGB}{120,65,150}
\tikzset{roadcurve/.style={very thick,roadblue},
 roadold/.style={thick,gray,dashed},
 roadguide/.style={gray!45,densely dotted},
 roadcorr/.style={very thick,roadorange}}

\providecommand{\RoadPairInfo}[4]{%
\par\smallskip
\begin{tikzpicture}[roadinfo/.style={anchor=north west,align=left,
 text width=6.65cm,inner sep=5pt,font=\footnotesize,rounded corners=2pt}]
 \node[roadinfo,draw=gray!50,fill=gray!5] at (0,0)
 {\textbf{#1}\par\smallskip #2};
 \node[roadinfo,draw=roadblue!55,fill=roadblue!4] at (7.3,0)
 {\textbf{#3}\par\smallskip #4};
\end{tikzpicture}\par\smallskip}

\begin{center}
\begin{tikzpicture}[>=stealth,
 overview/.style={draw=blue,rounded corners=2pt,fill=blue!2,
   text width=13.5cm,align=center,inner sep=6pt,font=\small},
 flow/.style={->,blue}]
 \node[overview] (candidate) {
   \textbf{Prescribe the outer candidate family}\\
   Temporary reference on $[0,R_{\rm ref}]$; piecewise continuation to the heat exterior.};
 \node[overview,anchor=north] (pressure) at ([yshift=-4mm]candidate.south) {
   \textbf{Close the outer moments with the prescribed axis pressure}\\[2pt]
   $\displaystyle P_0(Z)=-\int_0^\infty\frac{(U^\theta(R,Z))^2}{2R}\,dR=P_0^{\rm pre}(Z)$;
   impose the five conditions \eqref{eq:moment-conditions}.};
 \node[overview,anchor=north] (replace) at ([yshift=-4mm]pressure.south) {
    \textbf{Replace the inner reference on $[0,R_{\rm useless}]$}\\
    Preserve its five moments at $R_{\rm useless}$ and retain the profile for $R\ge R_{\rm useless}$.\\
   The constructed profile satisfies the relaxed cone condition \eqref{cone-relaxed}.};
 \node[overview,anchor=north] (shear) at ([yshift=-4mm]replace.south) {
   \textbf{Modify the shear}\\
   Obtain the admissible cone condition \eqref{cone-admissible}.};
 \node[overview,anchor=north] (repair) at ([yshift=-4mm]shear.south) {
   \textbf{Repair the five moments again}\\
   Restore \eqref{eq:moment-conditions} while preserving the admissible cone condition.};
 \draw[flow] (candidate.south)--(pressure.north);
 \draw[flow] (pressure.south)--(replace.north);
 \draw[flow] (replace.south)--(shear.north);
 \draw[flow] (shear.south)--(repair.north);
\end{tikzpicture}
\end{center}

Subsection~\ref{subsec:roadmap-outer} constructs the outer profile,
closes its terminal moments, and supplies the axis pressure $P_0$.
Subsection~\ref{roadin:subsection} constructs the core and connection,
then matches the five reference moments in one correction interval.
Subsection~\ref{roadin:shear-modification} modifies the shear and restores
the moments of the completed profile.

We are now ready to prescribe the candidate profiles. The radial prescriptions are
listed below. We will verify the cone conditions and carry out the moment
corrections later.

\Needspace{18cm}
\subsection{From the reference profile to the heat exterior}
\label{subsec:roadmap-outer}

The outer prescriptions are ordered by radius. Their coefficients may have
a different determination order. For example, an early correction may use
an integral of the future tail. We mark such coefficients as
\emph{pending derived data}.

We collect the fixed functions and explicit function families here.

\noindent\textbf{Flat functions and the smooth step.}
Put
\[
 \mathfrak f(t)=\begin{cases}e^{-1/t^2},&t>0,\\0,&t\le0.\end{cases}
\]
This function and all its derivatives vanish at zero. Fix the smooth step
and its primitive by
\begin{equation}\label{eq:road-sigma}
 \sigma(t)=\frac{\mathfrak f(t)}{\mathfrak f(t)+\mathfrak f(1-t)},
 \qquad \Sigma(t)=\int_0^t\sigma(s)\,ds\quad(t\ge0).
\end{equation}
Thus $\sigma=0$ on $(-\infty,0]$ and $\sigma=1$ on $[1,\infty)$.
It is increasing and flat at both endpoints.

\noindent\textbf{Normalized bumps.}
For $\ell>0$, define
\[
 \mathfrak b(t)=\begin{cases}e^{-1/(1-t^2)},&|t|<1,\\0,&|t|\ge1,\end{cases}
 \qquad \beta_\ell(t)=\frac{\mathfrak b(t/\ell)}
 {\ell\int_{-1}^1\mathfrak b(s)\,ds}.
\]
Each bump is nonnegative, has integral one, and vanishes to every order at
$t=\pm\ell$.

\noindent\textbf{The heat factor.}
For $\delta>0$ and $\xi\ge0$, use
\begin{equation}\label{eq:road-heat-functions}
 H_\delta(\xi)=\frac1{\Gamma(1+\delta/2)}
 \int_0^\infty e^{-v}v^{\delta/2}(1+\xi v)^{-\delta/2}\,dv.
\end{equation}
Here $\Gamma$ is the Euler gamma function and $H_\delta(0)=1$.
The formula fixes the entire family $H_\delta$ before $\delta$ is chosen.

The core seeds $H_0,G,F_0$ depend on the axis data. Their explicit formulas
are given in \eqref{roadin:axis-data}, after those data are available.

\paragraph{Absolute constants.}
The functions $\sigma$ and $\mathfrak b$ are fixed independently of the
input parameters. The estimates for these cutoffs and bumps enter the
sufficient lower bound for $P_*$ and upper bound for $\delta$ in the
outer construction.

We organize the construction data into input parameters, auxiliary constants,
and derived quantities. Input parameters, such as $P_*$ and $\delta$, remain
available for choice at the current stage, subject to the stated restrictions.
Auxiliary absolute constants, such as $M_1$ and $M_d$, are fixed using only
prescribed background objects and previously fixed auxiliary constants.
They are independent of the input parameters and input profiles.
Derived quantities, such as $e^{M_d}$ and $\mu=c_\mu P_*^{-4}$, are determined
by the chosen data and the stated construction rules. They are not additional
tuning parameters. This classification is relative to the construction stage:
a parameter initially treated as an input becomes derived once a matching
condition determines it.

We usually denote absolute constants by $M_1,M_2,\ldots$.
A descriptive subscript may indicate their role, as in $M_d$ for the axial
turnoff. We allow a harmless enlargement whenever the permitted dependencies
and all earlier conclusions remain valid: ``We still denote it by $M_1$;
its value may be larger than before, but this does not change any earlier
conclusions.''

A constant used in an expression defining a profile, a cutoff, or an interval
must instead be fixed at that point. If another constant must satisfy a
separate size relation with it, we introduce $M_2$, and then $M_3,\ldots$ as
needed. The subscripts distinguish constants; they do not by themselves impose
an ordering. The constants with specific names in the displayed formulas, such as
$M_d$, retain those names and their stated dependencies. Once they define
construction data, their values are not changed by the harmless-enlargement
convention.

\paragraph{The amplitude and similarity inputs.}
After fixing the required functions and auxiliary absolute constants, choose
\[
 P_*\gg1,\qquad 0<\delta\ll1.
\]
The admissible bounds for the outer construction depend only on the fixed
outer data. The upper bound for $\delta$ also depends on the selected $P_*$.
Following \eqref{eq:mc-free-parameters}, we impose
\[
 P_*>e^{T_d},\quad \mu:=c_\mu P_*^{-4}\le\mu_{\rm corr},\quad
 0<\delta<\min\{1/200,d_{\rm corr}c_\mu P_*^{-4}\},\quad T_d=e^{M_d}+10.
\]
Here $M_d>1$ is a sufficiently large fixed auxiliary constant.
The subscript $\mathrm{corr}$ marks the thresholds used for the outer
moment corrections and the cone estimates throughout the outer construction.
The constants $c_\mu,\mu_{\rm corr},d_{\rm corr}$ are fixed, sufficiently
small, and positive. They do not depend on either input parameter.
The radial placement and the core parameters will be specified later.
The outer candidate also takes the radial scale $R_{\rm ref}$ and a
provisional waiting length $\tau$ as inputs; moment matching determines
$\tau$, making it a derived quantity.

The implicit constants in the moment-matching and axis-pressure estimates
of Proposition~\ref{prop:road-outer-summary} depend only on the fixed outer
data and the indicated derivative order. In the endpoint factorization
of Proposition~\ref{prop:heat-collar}, the width $\ell$ and the constants
may depend on the selected profile parameters. The same proposition
gives the cone on the full outer collar. The data
\eqref{eq:mc-backward-normalized-data} at
$R_{\rm col}=e^{-1}R_h=e^2R_{2h}$ are uniform
for the global construction.

We may take $P_*$ arbitrarily large and $\delta$ arbitrarily small within this
admissible range. If $P_*$ is increased, $\delta$ must be decreased when needed
to retain the coupled bound. Later stages may further restrict this range
without affecting the outer estimates. Changing either input
requires recomputing its dependent quantities, including $\mu$, the waiting
length, the axis pressure, and the corresponding profiles. In the final
construction, $P_*,\delta$ are fixed before the pressure and core data are
finalized. Lemma~\ref{lem:imc-prepared-family} verifies compatibility
with the inner construction for these fixed outer data.

Proposition~\ref{prop:road-outer-summary} summarizes the outer output.
The descriptions of Intervals O.1--O.8 give its radial prescriptions and identify the few places
where moment equations are solved.

\paragraph{Order of parameter choices.}
The following diagram records the order of parameter choices for the outer
construction. The admissible input ranges are given in
\eqref{eq:mc-free-parameters}.

\begin{center}
\begin{tikzpicture}[>=stealth,
 stage/.style={draw=roadblue!65,fill=roadblue!4,rounded corners=3pt,
   text width=13.5cm,align=left,inner sep=7pt,font=\small},
 flow/.style={->,thick,roadblue!75}]
 \node[stage] (fixed) {
   \textbf{1. Fix the background functions and auxiliary absolute constants.}\\
   Choose the cutoffs, bump shapes, and auxiliary constants using only
   prescribed background data and earlier auxiliary choices,
   independently of the input parameters and profiles.};
 \node[stage,anchor=north] (input) at ([yshift=-5mm]fixed.south) {
   \textbf{2. Supply the input parameters and the reference branch.}\\
   For each admissible $P_*,\delta$ and radial scale $R_{\rm ref}$,
   define $(U^\theta_{\rm ref},U^z_{\rm ref})$ by
   \eqref{eq:road-reference} and compute its incoming moments.};
 \node[stage,anchor=north] (schedule) at ([yshift=-5mm]input.south) {
   \textbf{3. Determine the dimensionless schedule and waiting length.}\\
   Set $\mu=c_\mu P_*^{-4}$,
   $T_w=60\log(1/\mu)$, $T_s=4\log(2/\delta)$, and
   $\varepsilon=c_\varepsilon\delta$; $T_d,T_f$ are already fixed.
   Form the candidate family and determine $\tau$ by
   \eqref{eq:road-wait-choice}.
   All ratios $R_\alpha/R_{\rm ref}$ are then fixed independently
   of $R_{\rm ref}$.};
 \node[stage,anchor=north] (outer) at ([yshift=-5mm]schedule.south) {
   \textbf{4. Complete the outer profile and its moment conditions.}\\
   Use the actual heat tail in the moment targets. Determine
   $(d_1,d_2)$ first, then $(a_p,c_1,c_2)$.
   The pressure with axis value $P_0^{\rm pre}$ then vanishes at infinity.
   The angular and pressure conditions and the three axial-related
   moment conditions are all satisfied; no separate heat-compensation
   stage is needed.};
 \draw[flow] (fixed.south)--(input.north);
 \draw[flow] (input.south)--(schedule.north);
 \draw[flow] (schedule.south)--(outer.north);
\end{tikzpicture}
\end{center}

\noindent\textbf{The reference pressure moment.}
The reference contribution to the pressure moment is exactly
\[
 \int_0^{R_{\rm ref}}\frac{(U^\theta_{\rm ref}(R,Z))^2}{2R}\,dR
 =\frac{P_*^2}{2(1+Z^2)^2}\int_0^1x^{-4/5}\,dx
 =\frac{5P_*^2}{2(1+Z^2)^2}.
\]
It is independent of $R_{\rm ref}$. Changing $R_{\rm ref}$ leaves the
dimensionless schedule and restored axis pressure unchanged, but requires
recomputing the heat profile and its moment-correction coefficients.

\paragraph{Moment notation.}
Write $\mathbf M=(M^\theta,M^z,M^{\theta z},M^{z\theta},M^p)$ for the moments
in \eqref{fiveM}, including the reference branch down to $R=0$.

\begin{proposition}[The outer construction]\label{prop:road-outer-summary}
There exist absolute constants $M_d,R_{\rm corr}>1$ and
$c_\mu,c_\varepsilon,d_{\rm corr},\mu_{\rm corr}>0$, fixed as in
Section~\ref{sec:outer-construction},
such that the following hold whenever $P_*,\delta,R_{\rm ref}$ satisfy
\[
\begin{gathered}
 T_d:=e^{M_d}+10,\qquad \mu:=c_\mu P_*^{-4},\\
 R_{\rm ref}\ge R_{\rm corr},\quad P_*>e^{T_d},\quad
 \mu\le\mu_{\rm corr},\quad 0<\delta\le d_{\rm corr}\mu.
\end{gathered}
\]
The profiles are considered for $Z\in[-1,1]$.

\begin{enumerate}[(i)]
\item \textbf{Candidate family.}
For each provisional waiting length $\tau\ge0$, let
$\bar U=(\bar U^\theta,\bar U^z)$ be the uncorrected candidate described
in Intervals O.1--O.8 below, with checkpoint radii and full ansatz
given in Subsection~\ref{subsec:outer-ansatz}.
It agrees with the reference profile on $(0,R_{\rm ref}]$.

\item \textbf{Moment closure.}
For each allowed $(P_*,R_{\rm ref},\delta)$, the candidate family
in (i) contains a unique waiting length $\tau>0$ selected by
\eqref{eq:road-wait-choice}, independent of $R_{\rm ref}$.
There exist coefficients $d_1,d_2,a_p,c_1,c_2$, smooth in $Z$, such that the
angular correction \eqref{eq:road-angular-bumps} on $(R_{m\theta},R_{3/2})$
and the axial correction \eqref{eq:road-axial-pulse} on $(R_{mz},R_{0z})$
give the five moment conditions \eqref{eq:moment-conditions}
and enforce $P(\infty,Z)=0$ with the axis value in (iii).
The corrected reference-plus-outer profile satisfies the
relaxed cone on $[R_{\rm useless},R_h)$ and the admissible cone on
$[R_{+c},R_h)$.
For $R\ge R_h$ it is the exact heat exterior and $\mathcal{T}=0$.

\item \textbf{Uniform axis pressure.}
Let $\bar U^\theta_{\rm pre}$ be the candidate in (i) with $H_\delta$
replaced by $1$, and set
$P_0^{\rm pre}(Z)=-\int_0^\infty (\bar U^\theta_{\rm pre})^2/(2R)\,dR$.
The corrected axis pressure is $P_0=P_0^{\rm pre}$ at the selected
waiting length. It is even, independent of $R_{\rm ref}$, and extends
holomorphically to the fixed neighborhood $\Omega_0$ of $[-1,1]$ in
Lemma~\ref{lem:outer-axis-pressure}, with the uniform real and analytic
bounds stated there. The same bounds hold for $P_0^{\rm pre}$ for every
provisional $\tau\ge0$, with constants independent of
$P_*,\delta,R_{\rm ref},\tau$.
\end{enumerate}
\end{proposition}

\begin{proof}
The cone estimates up to $R_{mz}$ are proved in
Sections~\ref{subsec:outer-left-data}--\ref{subsec:outer-buffer}.
The moment matching and the cone estimates for the completed profile
are proved in Sections~\ref{subsec:outer-waiting}
and~\ref{subsec:outer-axial-matching}--\ref{subsec:outer-cone-closure}.
The axis-pressure assertions follow from
Lemma~\ref{lem:outer-axis-pressure} and \eqref{eq:mc-pressure-preserved}.
\end{proof}

\Needspace{7.5cm}
\paragraph{Interval O.1. The temporary reference on $[0,R_{\rm ref}]$.}
This branch supplies the incoming five moments for the outer construction.
The retained interval $[R_{\rm useless},R_{\rm ref}]$ satisfies the relaxed cone;
Section~\ref{subsec:outer-left-data} verifies it with the prescribed axis pressure.

\begin{center}
\begin{minipage}{.96\linewidth}\centering
\begin{tikzpicture}[x=1cm,y=.85cm,>=stealth,font=\small]
 \draw[->] (0,0)--(11.7,0) node[right] {$R$};
 \draw[->] (0,0)--(0,3.1) node[above] {$U^\theta_{\rm ref}$};
 \draw[roadcurve] (0,0)..controls (0,1.2) and (3,2.25)..(10.8,2.75);
 \draw[roadguide] (10.8,0)--(10.8,2.75);
 \node[below] at (0,0) {$0$};
 \node[below] at (10.8,0) {$R_{\rm ref}$};
 \node[above] at (6.7,2.75) {$U^\theta_{\rm ref}\propto R^{1/10}$};
 \node[align=center] at (6.2,1.55)
 {$U^z_{\rm ref}=4Z$\\[3pt]
 $U^r_{\rm ref}=2\sqrt{2R}\,((2+\delta)Z^2-1)/L_\delta$};
 \fill[roadpurple!18] (1.25,-.70) rectangle (2.75,-.48);
 \draw[roadpurple!75] (1.25,-.70) rectangle (2.75,-.48);
 \draw[roadpurple!55,densely dotted] (1.25,0)--(1.25,-.48);
 \draw[roadpurple!55,densely dotted] (2.75,0)--(2.75,-.48);
 \node[below,roadpurple,font=\scriptsize,align=center] at (2,-.70)
 {$[R_m,2R_m]$\\(enlarged)};
 \draw (4.2,-.06)--(4.2,.06) node[above] {$R_{\rm useless}$};
 \node[anchor=west,align=left,font=\scriptsize,text width=7.1cm]
 at (3.35,-.72)
 {later inner five-moment correction\\
  $R_m=e^{-6}R_{\rm ref}$,\quad $2R_m<R_{\rm useless}<R_{\rm ref}$};
\end{tikzpicture}
\par\small Only $[0,R_{\rm ref}]$ is shown. The curve is schematic at a fixed
$Z$; the singular slope at the axis is intentional for this temporary branch.
The interval $[R_m,2R_m]$ is enlarged for visibility.
No correction has been made there yet.
\end{minipage}
\end{center}

The subscript $\rm ref$ is reserved for this temporary branch on
$[0,R_{\rm ref}]$. For $R>R_{\rm ref}$ we write the continuation as
$U^\theta,U^r,U^z$, without that subscript. Every outer cumulative moment
starts with the reference moment at $R_{\rm ref}$ and is then accumulated
using the displayed outer components.

Recall the reference profile with input scales $P_*\ge1$ and $R_{\rm ref}>1$:
\begin{equation}\label{eq:road-reference}
 U^\theta_{\rm ref}(R,Z)=\frac{P_*}{1+Z^2}
 \left(\frac{R}{R_{\rm ref}}\right)^{1/10},
 \qquad U^z_{\rm ref}(R,Z)=4Z.
\end{equation}
This branch is temporary: its swirl has the wrong behavior at the physical
axis. Its radial moments are nevertheless finite there.

To record the radial velocity, introduce the input similarity exponent
$0<\delta<1/200$, and the derived functions $L=1-Z^2$, $L_\delta=1-\delta Z^2$.
For the later outer intervals $R\ge R_{\rm ref}$, radial velocity is
determined by
\begin{equation}\label{eq:road-radial-velocity}
\begin{aligned}
 U^r&=\frac{2ZR U^z-(1-\delta)Z M^z-L\partial_ZM^z}
 {L_\delta\sqrt{2R}},\\
 M^z(R,Z)&=4ZR_{\rm ref}+\int_{R_{\rm ref}}^R U^z(\rho,Z)\,d\rho.
\end{aligned}
\end{equation}
For the present reference branch this gives exactly
\begin{equation}\label{eq:road-reference-radial}
 M^z=4ZR,\qquad
 U^r_{\rm ref}
 =\frac{2\sqrt{2R}}{L_\delta}\bigl((2+\delta)Z^2-1\bigr)
 \qquad(0<R<R_{\rm ref}).
\end{equation}
Thus the reference branch has a nonzero radial velocity.

The five cumulative moments can already be computed for $0\le R\le R_{\rm ref}$:
\[
\begin{aligned}
M^\theta_{\rm ref}&=\tfrac58R\sqrt{2R}\,U^\theta_{\rm ref},&
M^z_{\rm ref}&=4ZR,\\
M^{\theta z}_{\rm ref}&=4Z M^\theta_{\rm ref},&
M^p_{\rm ref}&=\tfrac52(U^\theta_{\rm ref})^2,\\
M^{z\theta}_{\rm ref}&=16Z^2R-\tfrac5{12}R(U^\theta_{\rm ref})^2.
\end{aligned}
\]
All five moments vanish at $R=0$ and depend only on the reference branch.
The pressure on this branch is $P_{\rm ref}=P_0^{\rm pre}+M^p_{\rm ref}$,
where the full pre-heat swirl determines $P_0^{\rm pre}$.
The temporary profile on $[0,R_{\rm useless}]$, where
$R_{\rm useless}=e^{-5}R_{\rm ref}$, will be replaced by the inner
construction in Subsection~\ref{roadin:subsection}.

\noindent\textbf{Parameter status.}
The input data are $(P_*,R_{\rm ref},\delta)$; later stages restrict their
admissible range. The reference velocities, $L,L_\delta$, and all reference moments
are derived from these data.

\noindent\textbf{To be Checked:}
The reference profile satisfies the relaxed cone condition in a left
neighborhood of $R_{\rm ref}$. This will be checked in
Section~\ref{subsec:outer-left-data}.

\medskip
\Needspace{10cm}
\paragraph{Interval O.2. Turn the swirl downward and turn off the original axial flow.}
Turn off the original axial velocity while retaining its accumulated moments.

\begin{center}\begin{minipage}{.97\linewidth}\centering
\begin{tikzpicture}[x=1cm,y=.75cm,>=stealth,font=\small]
 \draw[->] (0,0)--(12.2,0) node[right] {$R$};
 \draw[->] (0,0)--(0,3.25) node[above] {$\overline U^\theta$};
 \draw[roadold] (0,0)..controls (0,1.5) and (1.1,2.55)..(2.8,2.75);
 \draw[roadcurve] (2.8,2.75)..controls (3.4,2.85) and (4,2.85)..(4.9,2.25)
 ..controls (6.9,1.2) and (9,.65)..(11.3,.35);
 \foreach \x/\lab in {0/0,2.8/R_{\rm ref},4.9/eR_{\rm ref},11.3/R_{-\mu}}
 {\draw[roadguide] (\x,0)--(\x,2.9);\node[below] at (\x,0) {$\lab$};}
 \node[above] at (3.9,3.1) {slope turn};
 \node[above] at (8,1.35) {$R^{-1/2}$};
 \node[align=center,font=\scriptsize,text width=2.7cm] at (1.35,1.15)
 {$U^\theta_{\rm ref}\propto R^{1/10}$\\
  $U^z_{\rm ref}=4Z$\\$U^r_{\rm ref}\not\equiv0$};
 \node[align=center,font=\scriptsize] at (3.85,1.25)
 {$B=1$\\$\overline U^z=4Z$};
 \node[align=center] at (7.7,2.7)
 {$\overline U^z:4Z\ \longrightarrow\ 0$\\$\overline U^r$ generally nonzero};
 \node[align=center,font=\scriptsize] at (10.4,1.6)
 {$B=0$\\$\overline U^z=0$};
 \draw[roadguide] (9.5,0)--(9.5,1.05);
 \draw (9.5,.06)--(9.5,-.18)
 node[below,font=\scriptsize] {$e^{-11}R_{-\mu}$};
\end{tikzpicture}
\RoadPairInfo{Left: $[0,R_{\rm ref}]$}
 {The reference amplitude is $P_* /(1+Z^2)$ at $R_{\rm ref}$.\par
  $M^z_{\rm ref}=4ZR$ and
   $U^r_{\rm ref}=\dfrac{2\sqrt{2R}}{L_\delta}[(2+\delta)Z^2-1]$.
  This radial component is determined by the accumulated axial moment.}
 {Right: $[R_{\rm ref},R_{-\mu}]$}
 {$\overline U^z=4ZB(y-1)$, $y=\log(R/R_{\rm ref})$.\par
  The accumulated quantity is $\overline M^z=4Z\Sigma_B(R)$;
  \eqref{eq:road-turnoff-radial} gives the radial component.
  It depends on the full axial history, including after $B=0$.}

\par\small Only $[0,R_{-\mu}]$ is shown. The older reference branch is dashed;
radial distances and amplitudes are schematic.
\end{minipage}\end{center}

Use $\sigma,\Sigma$ from \eqref{eq:road-sigma} and the absolute constant
$M_d$ of Proposition~\ref{prop:road-outer-summary}, and denote
\[
 T_d=e^{M_d}+10,\qquad R_{-\mu}=e^{1+T_d}R_{\rm ref}.
\]
Define the axial turnoff by
\[
 B(t)=\begin{cases}
 1,&t\le0,\\
 1-\sigma\!\left(\dfrac{\log(1+t)}{M_d}\right),&t>0.
 \end{cases}
\]
It equals zero for $t\ge e^{M_d}-1$ and is flat at both joining points.
For $y=\log(R/R_{\rm ref})$ on $0\le y\le1+T_d$, set
\begin{equation}\label{eq:road-first-turn}
\begin{aligned}
 \overline{U}^\theta(R,Z)&=\frac{P_*}{1+Z^2}
 \exp\!\left(\frac{y}{10}-\frac35\Sigma(y)\right),\\
 \overline{U}^z(R,Z)&=4Z B(y-1),
\end{aligned}
\end{equation}
Here the bar denotes the background profile before the later moment
corrections. It is used for all three velocity components and their
cumulative moments. The axial turnoff uses the function $B$ defined above.
On $[R_{\rm ref},eR_{\rm ref}]$, the logarithmic swirl slope moves
smoothly from $1/10$ to $-1/2$, while $\overline U^z=4Z$.
On $[eR_{\rm ref},R_{-\mu}]$, the swirl is an exact multiple of $R^{-1/2}$
and the axial velocity decreases smoothly to zero.
For $e^{-11}R_{-\mu}\le R\le R_{-\mu}$, $y\ge T_d-10=e^{M_d}$, so
$B(y-1)=1-\sigma(\log y/M_d)=0$ and $\overline U^z=0$.

All jets (the values of the three velocity components and their mixed
derivatives $\partial_R^j\partial_Z^k$, $j,k\ge0$) match at
$R_{\rm ref}$ because $\sigma$ is flat at zero.
The transition at $eR_{\rm ref}$ also matches all jets; there is no corner
in the swirl and the axial turnoff starts flat.
For clarity, put
\[
 \Sigma_B(R)=\int_0^R B\!\left(\log(\rho/R_{\rm ref})-1\right)\,d\rho.
\]
Then $\overline M^z=4Z\Sigma_B$ and
\begin{equation}\label{eq:road-turnoff-radial}
 \overline U^r=\frac{4\{2Z^2R B(y-1)-L_\delta\Sigma_B(R)\}}{L_\delta\sqrt{2R}}.
\end{equation}
After $\overline U^z$ vanishes, the accumulated axial moment remains: $\overline U^r$
is generally nonzero.

\noindent\textbf{Remark:}
On $[eR_{\rm ref},R_{-\mu}]$, the $R^{-1/2}$ branch gives $a=2$.
The pure axial-stress square cancels from \eqref{eq:outer-cone-test}:
\[
 a=1-2R\partial_R\log\overline U^\theta=2,\qquad
 (a-2)\left(\frac{\mathcal{I}^z}{\mathcal{I}^\theta}\right)^2=0.
\]
The slow turnoff controls the remaining axial-shear terms.
This is the reason for the first intermediate power.
Section~\ref{subsec:outer-left-data} carries out the stress estimates.

By $R_{-\mu}$, $U^z=0$, so one could already consider correcting the
accumulated axial-related moments. We postpone this correction to
Interval O.4 on $[R_{mz},R_{0z}]$. This leaves the reserved subintervals in the
pure-power buffer $[R_{+c},R_{mz}]$ of Interval O.3 on $[R_{-\mu},R_{mz}]$ available for
later surgeries. The correction on Interval O.4 sets $M^z(R_{0z},Z)=M^{\theta z}(R_{0z},Z)=0$
and enforces $M^{z\theta}(\infty,Z)=0$.

\noindent\textbf{Parameter status.}
The input data remain $(P_*,R_{\rm ref},\delta)$, now with $P_*>e^{T_d}$.
The cutoff $\sigma$ and $M_d$ are fixed; $\Sigma,B,T_d,R_{-\mu}$ are derived.
Choosing a new $P_*$ or $R_{\rm ref}$ changes the corresponding previously
displayed profiles; it does not introduce a new choice of connection function.

\noindent\textbf{To be Checked:}
The profile satisfies the relaxed cone condition on $[R_{\rm ref},R_{-\mu}]$.
This will be checked in Section~\ref{subsec:outer-left-data}.

\medskip
\Needspace{15.0cm}
\paragraph{Interval O.3. A slightly steeper power law and a long buffer.}
Here $U^z=0$. No additional moment equation is imposed.

\begin{center}\begin{minipage}{.97\linewidth}\centering
\begin{tikzpicture}[x=1cm,y=.9cm,>=stealth,font=\small]
 \draw[->] (0,0)--(12,0) node[right] {$R$};
 \draw[->] (0,0)--(0,3.2) node[above] {$U^\theta$};
 \draw[roadold] (.3,2.85)..controls (1.7,2.4) and (2.6,2.05)..(3.3,1.95);
 \draw[roadcurve] (3.3,1.95)..controls (4.1,1.8) and (4.5,1.58)..(5.5,1.35)
 ..controls (7.6,.9) and (9.6,.42)..(11.4,.25);
 \foreach \x/\lab in {.3/eR_{\rm ref},3.3/R_{-\mu},5.5/R_{+c},11.4/R_{mz}}
 {\draw[roadguide] (\x,0)--(\x,2.9);\node[below] at (\x,0) {$\lab$};}
 \node at (1.7,3.1) {$R^{-1/2}$};
 \node[align=center,font=\scriptsize,text width=2.9cm] at (1.7,1.05)
 {$U^z:4ZB\to0$\\$U^r$ generally nonzero};
 \node at (8.3,1.45) {$R^{-1/2-\mu}$};
 \node[align=center] at (7.6,2.65)
 {$U^z=0$\\[3pt]$U^r=-4\Sigma_B(R_{-\mu})/\sqrt{2R}$};
 \draw[<->] (5.5,-.65)--(11.4,-.65)
 node[midway,below] {logarithmic length $T_w$};
\end{tikzpicture}
\RoadPairInfo{Left: $[eR_{\rm ref},R_{-\mu}]$}
 {$U^r=\dfrac{4[2Z^2RB(y-1)-L_\delta\Sigma_B(R)]}{L_\delta\sqrt{2R}}$,
  where $y=\log(R/R_{\rm ref})$.\par
  The axial turnoff finishes at $e^{-11}R_{-\mu}$; its accumulated
  moment remains in the radial formula.}
 {Right: $[R_{-\mu},R_{mz}]$}
 {The logarithmic radial slope changes on $[R_{-\mu},R_{+c}]$.
  On $[R_{+c},R_{mz}]$, $\overline U^\theta\propto R^{-1/2-\mu}$.\par
  $R_{+c}=eR_{-\mu}$ and $\Sigma_B(R_{-\mu})$ is constant in $R,Z$; three correction slots are
  reserved inside the latter interval.}

\par\small Only the previous turnoff interval and the new slope/buffer
intervals are shown: $[eR_{\rm ref},R_{-\mu}]\cup[R_{-\mu},R_{mz}]$.
\end{minipage}\end{center}

\begin{center}\begin{minipage}{.97\linewidth}\centering
\begin{tikzpicture}[x=.43cm,y=.75cm,>=stealth,font=\small]
 \foreach \a/\b/\lab in {0/5/I_1,11/16/I_2,17/22/I_3}
 {\fill[roadorange!15] (\a,.05) rectangle (\b,1.45);
  \draw[roadorange!60] (\a,.05) rectangle (\b,1.45);
 }
 \fill[roadpurple!28] (1,.08) rectangle ({1+ln(2)},1.42);
 \draw[roadpurple,thick] (1,.08) rectangle ({1+ln(2)},1.42);
 \node[align=center,font=\scriptsize] at (3.5,.78) {$I_1$\\five\\moments};
 \node[align=center,font=\scriptsize] at (13.5,.78) {$I_2$\\higher\\orders};
 \node[align=center,font=\scriptsize] at (19.5,.78) {$I_3$\\phase\\mean};
 \draw[->] (-.5,0)--(26,0) node[right] {$\log(R/R_{mz})$};
 \foreach \x/\lab in {0/-25,5/-20,11/-14,16/-9,17/-8,22/-3,25/0}
 {\draw (\x,0)--(\x,-.12);}
 \foreach \x/\lab in {0/-25,5/-20,11/-14,16/-9,22/-3,25/0}
 {\node[below] at (\x,-.12) {$\lab$};}
 \node[below] at (17,-.7) {$-8$};
 \node[above] at (12.5,1.65)
 {reserved inside the exact $R^{-1/2-\mu}$ buffer};
 \draw[roadpurple,thin,->] ({1+ln(2)/2},.08)
 --({1+ln(2)/2},-1.12);
 \node[anchor=north west,align=left,font=\scriptsize,roadpurple,
 text width=7.0cm] at (0,-1.30)
 {$[R_c,2R_c]\subset I_1$,\quad $R_c=e^{-24}R_{mz}$\\
  later outer five-moment restoration};
\end{tikzpicture}
\par\footnotesize
Common to all three slots before their later corrections:
$U^\theta\propto R^{-1/2-\mu}$,
$U^z=0$, $U^r=-4\Sigma_B(R_{-\mu})/\sqrt{2R}$.
\par\small A zoom of the final $25$ logarithmic units before $R_{mz}$.
The three shaded slots are reserved, not filled by new corrections here.
The purple subinterval has the exact logarithmic endpoints
$-24$ and $-24+\log2$; it marks a later correction, with no new input endpoint.
\end{minipage}\end{center}

With the outer constants fixed as in \eqref{eq:mc-free-parameters}, set
\begin{equation}\label{eq:road-mu-scales}
 \mu=c_\mu P_*^{-4},\qquad T_w=60\log(1/\mu),\qquad
 R_{+c}=eR_{-\mu},\qquad R_{mz}=e^{T_w}R_{+c}.
\end{equation}
Here $0<\mu\le\mu_{\rm corr}$ and $0<\delta\le d_{\rm corr}\mu$.
For $t=\log(R/R_{-\mu})$, define on $[R_{-\mu},R_{mz}]$
\begin{equation}\label{eq:road-second-turn}
 \overline{U}^\theta(R,Z)=\overline{U}^\theta(R_{-\mu},Z)
 \exp\!\left(-\frac t2-\mu \Sigma(t)\right),\qquad \overline{U}^z=0.
\end{equation}
As in \eqref{eq:road-first-turn}, the bar denotes the background profile
before the moment corrections.
The interval $[R_{-\mu},R_{+c}]$ changes the slope smoothly from $-1/2$ to
$-1/2-\mu$. On $[R_{+c},R_{mz}]$ the slope is exactly $-1/2-\mu$.
The long buffer reduces the swirl amplitude before the moment corrections.
Flatness of $\sigma$ gives smooth matching at both ends of the slope change.
At $R_{-\mu}$, all jets match because $e^{-\mu \Sigma(t)}-1$ is flat at $t=0$.

Once $R_{\rm ref}$ and $M_d$ are fixed, $\Sigma_B(R_{-\mu})$ is a constant
independent of $R,Z$. The radial velocity throughout this stage is
\begin{equation}\label{eq:road-buffer-radial}
 U^r=-\frac{4\Sigma_B(R_{-\mu})}{\sqrt{2R}},\qquad U^z=0.
\end{equation}
Here $M^z$ and $M^{\theta z}$ remain equal to their incoming values.
On the pure-power buffer, their natural correction scales
$R\overline U^\theta$ and $R\sqrt{2R}(\overline U^\theta)^2$
grow like $R^{1/2-\mu}$ and $R^{1/2-2\mu}$, respectively.
The long buffer therefore makes these inherited moments small
relative to the scales used in the linear moment equations of
Interval O.4, even as the swirl amplitude decreases.
Meanwhile, $\partial_RM^{z\theta}=-(U^\theta)^2/2$ continues to
accumulate the swirl contribution.

\noindent\textbf{Three reserved correction intervals.}
In the pure-power buffer retain the three open intervals
\begin{equation}\label{eq:road-reserved-patches}
\begin{aligned}
 I_1&=(e^{-25}R_{mz},e^{-20}R_{mz}),\\
 I_2&=(e^{-14}R_{mz},e^{-9}R_{mz}),\\
 I_3&=(e^{-8}R_{mz},e^{-3}R_{mz}).
\end{aligned}
\end{equation}
They fit inside $(R_{+c},R_{mz})$ since $T_w>25$.
Here the power is fixed and the amplitude is small, with controlled cone
estimates and no overlap with the later axial, angular, or heat connections.
They correspond to the first, third, and fourth reserved intervals
after (A.9) in~\cite{1}, with new labels. Their roles are:
\begin{enumerate}[(i)]
\item $I_1$: restore five moments after the final shear modification;
\item $I_2=I_{\rm pos}$: moment corrections for higher-order backgrounds;
\item $I_3=I_{\rm mean}$: moment corrections for the phase-averaged terms.
\end{enumerate}
The original construction uses its second reserved interval to compensate
for the heat tail \cite[Proposition~A.7]{1}. Here the angular corrections
on Interval O.6 match the angular and pressure moments, and the axial
correction on Interval O.4 then uses the completed swirl to enforce
$M^{z\theta}(\infty,Z)=0$.
Their coefficients are determined as described under Interval O.8.
All three are untouched by the axial and angular corrections specified
below. They introduce no new input endpoint parameters. Their
uncorrected velocities satisfy $U^z=0$ and
$U^r=-4\Sigma_B(R_{-\mu})/\sqrt{2R}$, not $U^r=0$.

\noindent\textbf{Remark:}
On $[R_{+c},R_{mz}]$, the power $R^{-1/2-\mu}$ gives $a=\kappa=2+2\mu>2$.
With $\delta\ll\mu$, the angular inertial-stress source in
\eqref{eq:QJ-nu} has the positive lower bound
\[
 \mu\left(1-\frac{4L_\delta \Sigma_B(R_{-\mu})}R\right)-\frac\delta2
 \ge\mu(1-4e^{-11})-\frac\delta2\gtrsim\mu>0.
\]
This supplies the angular stress margin.
Section~\ref{subsec:outer-buffer} gives these estimates.

As at the end of Interval O.2, $U^z=0$ does not cancel the accumulated axial-related
moments. We leave their correction to Interval O.4 on $[R_{mz},R_{0z}]$.
It sets $M^z(R_{0z},Z)=M^{\theta z}(R_{0z},Z)=0$ and balances $M^{z\theta}$
against the future tail so that $M^{z\theta}(\infty,Z)=0$.
This keeps the reserved intervals in $[R_{+c},R_{mz}]$ available for later surgeries.

\noindent\textbf{Parameter status.}
For the given $(P_*,R_{\rm ref},\delta)$, formula \eqref{eq:road-mu-scales}
determines $\mu,T_w,R_{+c},R_{mz}$. Increasing $P_*$ changes $\mu$ and
the subsequent radii that depend on it.

\noindent\textbf{To be Checked:}
The profile satisfies the admissible cone condition on $(R_{-\mu},R_{mz}]$ and
the relaxed cone condition at $R_{-\mu}$. These will be checked in
Section~\ref{subsec:outer-buffer}.

\medskip
\Needspace{13.5cm}
\paragraph{Interval O.4. Reserve an axial pulse and two end corrections.}
Prescribe an axial pulse and two end bumps to correct the three moments
related to $U^z$.

\begin{center}\begin{minipage}{.97\linewidth}\centering
\begin{tikzpicture}[x=1cm,y=.75cm,>=stealth,font=\small]
 \draw[->] (0,0)--(12,0) node[right] {$R$};
 \draw[->] (0,0)--(0,2.8) node[above] {$U^\theta$};
 \draw[roadold] (.3,2.45)..controls (1.1,2.1) and (1.8,1.85)..(2.7,1.65);
 \draw[roadcurve] (2.7,1.65)..controls (6,1.1) and (9,.45)..(11.3,.3);
 \foreach \x/\lab in {.3/R_{+c},2.7/R_{mz},11.3/R_{0z}}
 {\draw[roadguide] (\x,0)--(\x,2.5);\node[below] at (\x,0) {$\lab$};}
 \node at (7.1,2.1) {$U^\theta=\overline{U}^\theta(R_{mz},Z)(R/R_{mz})^{-1/2-\mu}$};
 \node[align=center,font=\scriptsize,text width=2.6cm] at (1.35,.85)
 {$U^\theta\propto R^{-1/2-\mu}$\\$U^z=0$\\$U^r\not\equiv0$};
 \begin{scope}[yshift=-3.1cm]
 \draw[->] (0,0)--(12,0) node[right] {$R$};
 \node[left] at (0,1.2) {$U^z$};
 \draw[roadold] (.3,0)--(2.7,0);
 \draw[roadcorr] (2.7,0)..controls (3.5,0) and (4.2,1.45)..(5.2,1.55)
 ..controls (6.4,1.45) and (7.1,0)..(7.7,0)--(8.4,0)
 ..controls (8.5,-.45) and (8.9,-.45)..(9,0)--(9.6,0)
 ..controls (9.7,.4) and (10.1,.4)..(10.2,0)--(11.3,0);
 \node[above,align=center] at (5.2,1.55)
 {main pulse\\$M^{z\theta}(\infty,Z)=0$};
 \node[align=center] at (9.4,1.1)
 {two end bumps\\cancel $M^z,M^{\theta z}$};
 \node[align=center] at (6,-1.25)
 {$U^r$ is given by \eqref{eq:road-radial-velocity}; after exact matching,
 $U^r(R_{0z})=U^z(R_{0z})=0$.};
 \end{scope}
\end{tikzpicture}
\RoadPairInfo{Left: $[R_{+c},R_{mz}]$}
 {The inherited radial component is $-4\Sigma_B(R_{-\mu})/\sqrt{2R}$.
  Its axial moment comes from the earlier reference and turnoff.\par
  The three reserved slots remain unused.}
 {Right: $[R_{mz},R_{0z}]$}
 {The coefficients $a_p,c_1,c_2$ remain pending until the tail
  used in \eqref{eq:mc-pulse-amplitude} is complete.\par
  Inside the pulse, the radial component is recovered from the
  actual partial moment by \eqref{eq:road-radial-velocity}.}

\par\small Only $[R_{+c},R_{mz}]\cup[R_{mz},R_{0z}]$ is shown. The small end bumps
are magnified and their signs are illustrative.
\end{minipage}\end{center}

Set $R_{0z}=e^{13/\mu}R_{mz}$. On $[R_{mz},R_{0z}]$, retain the swirl power
\[
 \overline{U}^\theta(R,Z)=\overline{U}^\theta(R_{mz},Z)e^{-(1/2+\mu)t},\qquad
 t=\log(R/R_{mz}).
\]
\Needspace{3cm}
\noindent\textbf{Purpose.}
On this interval we modify $U^z$ to correct the three moments
$M^z$, $M^{\theta z}$, and $M^{z\theta}$.
The main pulse adjusts the positive integral
$\int_{R_{mz}}^{R_{0z}}(U^z(R,Z))^2\,dR$ so that
\[
 M^{z\theta}(\infty,Z)=0,\qquad
 M^{z\theta}(R_{0z},Z)=\frac12\int_{R_{0z}}^\infty(U^\theta)^2\,dR>0.
\]
The two end bumps jointly impose $M^z(R_{0z},Z)=M^{\theta z}(R_{0z},Z)=0$;
since $U^z=0$ for $R\ge R_{0z}$, this also gives $U^r=0$ there.

The outer corrections use the fixed bump $\beta=\beta_{3/20}$.
Define the main pulse and its two end bumps by
\begin{equation}\label{eq:road-pulse-shapes}
\begin{aligned}
 g_p(\xi)&=\begin{cases}
 [1-\sigma(\xi-10)]\displaystyle\int_0^\xi\sigma(50v)\,dv,&\xi\ge0,\\
 0,&\xi<0,
 \end{cases}\\
 \gamma_1(t)&=\beta(t-13/\mu+3),\qquad
 \gamma_2(t)=\beta(t-13/\mu+1).
\end{aligned}
\end{equation}
The fixed shape $g_p$ vanishes outside $[0,11]$ and is flat at both endpoints.
Its rescaling in the axial profile is $g_p(\mu t)$.
These outer $\gamma_j$ are distinct from the inner bumps defined under Interval I.4.
In the $t$-coordinate, the two support intervals have centers
$13/\mu-3$ and $13/\mu-1$, respectively, and common length $3/10$.
Since $R=R_{mz}e^t$ and $R_{0z}=R_{mz}e^{13/\mu}$, their radial intervals are
\[
 \left(e^{-3-3/20}R_{0z},e^{-3+3/20}R_{0z}\right),\qquad
 \left(e^{-1-3/20}R_{0z},e^{-1+3/20}R_{0z}\right).
\]
The corrected axial candidate on this interval is
\begin{equation}\label{eq:road-axial-pulse}
 U^z(R,Z)=\overline{U}^\theta(R,Z)
 \{a_p(Z)g_p(\mu t)+c_1(Z)\gamma_1(t)+c_2(Z)\gamma_2(t)\}.
\end{equation}
Outside this interval retain the previous axial profile. The main pulse
ends at $e^{11/\mu}R_{mz}$. Its support and those of the two end bumps are
disjoint. Every added function is flat
at its support endpoints, so all joins are smooth.

The energy condition uses the entire corrected swirl, including the later
angular corrections and heat tail. The coefficients therefore remain
pending until the matching described under Interval O.8.

For each trial pulse amplitude, the two linear moment equations determine
$c_1,c_2$ together. The energy equation then selects the positive branch
$a_p\in(0.9,1.2)$; see
\eqref{eq:mc-axial-end-coefficients} and \eqref{eq:mc-pulse-amplitude}.
Only $U^z$ changes here, so $M^\theta$ and $M^p$ are unaffected.
The third target at $R_{0z}$ is the positive tail integral displayed above,
not zero; it is $M^{z\theta}(\infty)$ that vanishes.

The swirl-energy balance requires axial energy of order
\[
 \frac{R_{mz}[\overline U^\theta(R_{mz},Z)]^2}{\mu}.
\]
A pulse with $U^z$ comparable to $U^\theta$ supplies this energy
over a logarithmic interval of length $O(1/\mu)$.
In the variable $t$, the two linear moments have weights
$e^{(1/2-\mu)t}$ and $e^{(1/2-2\mu)t}$, whereas the energy weight is
$e^{-2\mu t}$.
Together with the small incoming normalized moments from Interval O.3,
this explains why the later end bumps can have exponentially small
coefficients while cancelling the linear moments of the main pulse.

On the main pulse, the ratio $U^z/U^\theta$ varies through the slow
variable $\mu t$. Its prescribed shape and slope bound control the
stress direction; the end corrections are small together with their
first derivatives in $t$ and $Z$.
Section~\ref{subsec:outer-axial-matching} verifies the cone using these
properties and the pressure estimates for the completed swirl.

\noindent\textbf{Remark:}
The background already has $U^z=0$ here. We now insert compactly supported
axial corrections to cancel its accumulated moments. We keep $U^\theta$
fixed, including its factor $(1+Z^2)^{-1}$. Its exact radial power gives
the exponential weights in \eqref{eq:mc-axial-linear-equations} and keeps
$a=2+2\mu$ throughout the pulse. Removing the $Z$ dependence at the same
time would change these weights and the shear estimates. We therefore
leave that change to Interval O.5 on $[R_{0z},R_{m\theta}]$. After the axial matching,
$U^r=U^z=0$ there. Flattening $U^\theta$ then preserves $M^z=M^{\theta z}=0$.
Thus the correction here addresses the axial-related moments, while the construction on Interval O.5 prepares
the $Z$-independent background for the later angular corrections.

\noindent\textbf{Parameter status.}
No new scalar input parameter is introduced. The input parameter set remains
$(P_*,R_{\rm ref},\delta)$. The bump $\beta$ is fixed and $R_{0z}$ is derived.
Section~\ref{subsec:outer-axial-matching} gives the equations for
$a_p,c_1,c_2$ and their estimates for a prescribed swirl.
Section~\ref{subsec:outer-angular-matching} applies them to the final swirl.

\noindent\textbf{To be Checked:}
After the pending coefficients have been fixed, verify the following.
\begin{enumerate}[(i)]
\item \emph{Moments related to $U^z$.}
The two linear conditions and the energy condition are
\[
 M^z(R_{0z},Z)=M^{\theta z}(R_{0z},Z)=0,\qquad
 M^{z\theta}(R_{0z},Z)=\frac12\int_{R_{0z}}^\infty(U^\theta)^2\,dR.
\]
The last equality gives $M^{z\theta}(\infty,Z)=0$.
\item \emph{Cone.}
The completed profile satisfies the admissible cone condition on $[R_{mz},R_{0z}]$.
\end{enumerate}
Section~\ref{subsec:outer-axial-matching} proves both assertions
under the input estimates verified in Section~\ref{subsec:outer-angular-matching}.

\medskip
\Needspace{10.5cm}
\paragraph{Interval O.5. Remove the axial dependence of the background swirl.}
Flatten the background swirl without introducing another moment equation.
After the matching on Interval O.4, $M^z=M^{\theta z}=0$ and $U^r=U^z=0$.

\begin{center}\begin{minipage}{.97\linewidth}\centering
\begin{tikzpicture}[x=1cm,y=.86cm,>=stealth,font=\footnotesize]
 \fill[roadpurple!5] (3.2,0) rectangle (10.1,3.45);
 \draw[->] (0,0)--(12.45,0) node[right] {$R$};
 \draw[->] (0,0)--(0,3.65) node[above] {$\overline U^\theta$};
 \foreach \x/\lab in {.4/R_{mz},3.2/R_{0z},10.1/R_{m\theta}}
 {\draw[roadguide] (\x,0)--(\x,3.35);\node[below] at (\x,0) {$\lab$};}
 \draw[thick,roadpurple,dashed] (.4,3.15)..controls (1.4,3.04) and (2.4,2.91)..(3.2,2.8);
 \draw[thick,roadpurple!75,dashed] (.4,2.52)..controls (1.4,2.44) and (2.4,2.33)..(3.2,2.24);
 \draw[thick,roadpurple!50,dashed] (.4,1.58)..controls (1.4,1.53) and (2.4,1.46)..(3.2,1.4);
 \draw[very thick,roadpurple] (3.2,2.8)..controls (5.2,2.53) and (7.3,1.02)..(10.1,.68);
 \draw[very thick,roadpurple!75] (3.2,2.24)..controls (5.2,2.02) and (7.3,1.02)..(10.1,.68);
 \draw[very thick,roadpurple!50] (3.2,1.4)..controls (5.2,1.26) and (7.3,1.02)..(10.1,.68);
 \draw[roadcurve] (10.1,.68)..controls (10.8,.595) and (11.5,.52)..(12,.47);
 \node[above,roadpurple] at (1.3,3.05) {$Z=0$};
 \node[above,roadpurple!75] at (1.35,2.43) {$Z=1/2$};
 \node[above,roadpurple!60!black] at (1.35,1.51) {$Z=1$};
 \node[align=center] at (6.7,3.15)
 {$(1+Z^2)^{-1}\longrightarrow 1/2$\\flat at both endpoints};
 \node[align=center] at (10.55,3.12)
 {$U^r=U^z=0$\\after matching};
 \node[anchor=west,align=left,roadblue] at (10.25,1.65)
 {one common\\outgoing profile};
 \draw[->,roadblue] (10.75,1.2)--(11,.58);
 \draw[<->,roadpurple] (3.2,-.65)--(10.1,-.65)
 node[midway,below] {$\log(R_{m\theta}/R_{0z})=100$};
\end{tikzpicture}
\RoadPairInfo{Left: $[R_{mz},R_{0z}]$}
 {The pulse is prescribed by \eqref{eq:road-axial-pulse}; its
  three coefficients are determined only after the completed tail
  is available. The radial component retains the accumulated
  moment through the pulse.}
 {Right: $[R_{0z},R_{m\theta}]$}
 {The factor in \eqref{eq:road-flatten} is
  $[(1+Z^2)/2]^{\sigma(t/T_f)}$. Its derivative changes the
  radial slope inside the interval; both joins are flat.
  The endpoint amplitude is derived from the incoming profile.}

\par\small The left interval $[R_{mz},R_{0z}]$ contains the previously
reserved axial pulse. A short continuation beyond $R_{m\theta}$ shows
the common outgoing profile.
\end{minipage}\end{center}

The three curves show how different incoming values of $Z$ give the
same outgoing profile at $R_{m\theta}$. The shaded interval changes
the axial factor as well as the radial slope; it is not an exact
power-law interval. The flat step makes all radial derivatives join
the neighboring profiles. The horizontal spacing and amplitudes are
schematic, while the stated logarithmic length is exact.

Set $R_{m\theta}=e^{100}R_{0z}$. The logarithmic length is $T_f=100$.
For $t=\log(R/R_{0z})$, set
\begin{equation}\label{eq:road-flatten}
 \overline{U}^\theta(R,Z)=\overline{U}^\theta(R_{0z},Z)
 \left(\frac{1+Z^2}{2}\right)^{\sigma(t/T_f)}e^{-(1/2+\mu)t},
 \qquad R_{0z}\le R\le R_{m\theta}.
\end{equation}
At the left endpoint $\overline{U}^\theta(R_{0z},Z)=A_v/(1+Z^2)$ for a derived constant
$A_v>0$. At $R_{m\theta}$ the factor $1/(1+Z^2)$ has therefore become $1/2$:
\[
 \overline{U}^\theta(R_{m\theta},Z)=\frac{A_v}{2}e^{-(1/2+\mu)T_f},
\]
which is independent of $Z$. The transition is flat at both endpoints.
Its logarithmic radial slope equals
\[
 -\frac12-\mu+
 \partial_t\bigl[\sigma(t/T_f)\bigr]
 \log\!\left(\frac{1+Z^2}{2}\right).
\]
Thus it is generally not an exact power of $R$ inside this transition.

Set $U^z=0$ on this interval. Once the two linear moment equations in
Interval O.4 have been imposed, $M^z(R_{0z},Z)=0$ as an identity in $Z$;
therefore $U^r=0$ throughout this and all later outer intervals.

\noindent\textbf{Parameter status.}
The input data remain $(P_*,R_{\rm ref},\delta)$. The length $T_f$ is fixed;
$R_{m\theta},A_v$, and the amplitude at $R_{m\theta}$ are derived. There is no independent
choice of the axial flattening factor.

\noindent\textbf{To be Checked:}
The profile satisfies the admissible cone condition on $[R_{0z},R_{m\theta}]$.
Section~\ref{subsec:outer-angular-matching} verifies this condition
after both the angular and axial corrections.

\medskip
\Needspace{11.0cm}
\paragraph{Interval O.6. Relax the angular moment and reserve two angular bumps.}
The long power-law segment reduces the normalized incoming angular
mismatch; two local bumps then match the angular and pressure moments.

\begin{center}\begin{minipage}{.97\linewidth}\centering
\begin{tikzpicture}[x=1cm,y=.9cm,>=stealth,font=\small]
 \draw[->] (0,0)--(12,0) node[right] {$R$};
 \draw[->] (0,0)--(0,3.1) node[above] {$U^\theta$};
 \draw[roadold] (.4,2.6)..controls (1.8,2.25) and (2.8,1.9)..(3.8,1.7);
 \draw[roadcurve] (3.8,1.7)..controls (5.1,1.35) and (5.7,1.15)..(6.1,1.05);
 \draw[roadcorr] (6.1,1.05)..controls (6.4,.85) and (6.65,.8)..(6.9,.94);
 \draw[roadcurve] (6.9,.94)..controls (7.6,.78) and (8.1,.65)..(8.5,.58);
 \draw[roadcorr] (8.5,.58)..controls (8.8,.78) and (9.05,.75)..(9.3,.46);
 \draw[roadcurve] (9.3,.46)--(11.3,.25);
 \foreach \x/\lab in {.4/R_{0z},3.8/R_{m\theta},11.3/R_{3/2}}
 {\draw[roadguide] (\x,0)--(\x,2.7);\node[below] at (\x,0) {$\lab$};}
 \node[above] at (6.5,1.3) {$d_1\beta$};
 \node[above] at (8.95,1.05) {$d_2\beta$};
 \node[align=center] at (8,2.5)
 {background $R^{-1/2-\mu}$\\[3pt]$U^r=U^z=0$ after matching};
 \node[align=center,font=\scriptsize,text width=3.2cm] at (2,.95)
 {$U^\theta$: $Z$-flattening\\$U^r=U^z=0$ after matching};
\end{tikzpicture}
\RoadPairInfo{Left: $[R_{0z},R_{m\theta}]$}
 {The previous flattening uses \eqref{eq:road-flatten}.
  Its $Z$ dependence disappears at $R_{m\theta}$, while its radial
  slope varies inside the interval.}
 {Right: $[R_{m\theta},R_{3/2}]$}
 {$U^\theta=\overline U^\theta(1+h)$, where
  $h=d_1\beta(t+3)+d_2\beta(t+1)$. The supports lie around
  $e^{-3}R_{3/2}$ and $e^{-1}R_{3/2}$.\par
  The two coefficients depend on the completed tail.}

\par\small Only $[R_{0z},R_{m\theta}]\cup[R_{m\theta},R_{3/2}]$ is shown.
The two relative corrections are magnified; their signs are illustrative.
\end{minipage}\end{center}

Derive $R_{3/2}=\mu^{-30}R_{m\theta}$ and continue the background by
\[
 \overline{U}^\theta(R)=\overline{U}^\theta(R_{m\theta})(R/R_{m\theta})^{-1/2-\mu},\qquad U^z=0,
 \qquad R_{m\theta}\le R\le R_{3/2}.
\]
This background is independent of $Z$. Let $t=\log(R/R_{3/2})$
and add the two relative angular bumps, using the same $\beta=\beta_{3/20}$:
\begin{equation}\label{eq:road-angular-bumps}
 h(t,Z)=d_1(Z)\beta(t+3)+d_2(Z)\beta(t+1),\qquad
 U^\theta=\overline{U}^\theta(1+h).
\end{equation}

Their supports are centered at $e^{-3}R_{3/2}$ and
$e^{-1}R_{3/2}$, with logarithmic half-width $3/20$.
They lie strictly inside $(R_{m\theta},R_{3/2})$ and vanish to all orders
at their support endpoints. The corrected swirl may depend on $Z$
inside these bumps, but agrees with the $Z$-independent background near
both ends of this stage.

After the waiting length and actual heat tail are fixed, the two
coefficients impose, for every $Z$,
\begin{equation}\label{eq:road-angular-matching}
\begin{aligned}
 M^\theta(R_{3/2},Z)
   &=M_{\rm tar}^\theta(R_{3/2},Z;\tau),\\
 M^p(R_{3/2},Z)
   &=\bar M^p(R_{3/2},Z)+\Delta_H(Z).
\end{aligned}
\end{equation}
The tail target $M_{\rm tar}^\theta$ and heat deficit $\Delta_H$ are
defined under Interval O.8. The small solution is constructed in
\eqref{eq:mc-angular-coefficients}.

\noindent\textbf{Remark:}
On this interval we modify $U^\theta$ to match the incoming angular moment to
the target determined by the right tail. The waiting length on Interval O.7
is fixed first. The coefficients $d_1(Z)$ and $d_2(Z)$ then cancel the
remaining mismatch in the renormalized $M^\theta$ condition.
They also restore the total pressure moment to its pre-heat value,
$M^p(\infty,Z)=M_{\rm pre}^p(\infty,Z)$.
The angular and pressure moments give different weights to the two bump
locations. This permits the two coefficients to adjust both targets;
see \eqref{eq:mc-angular-system-unscaled}.
They leave $M^z,M^{\theta z}$ unchanged because $U^z=0$ here,
but change $M^{z\theta}$ through the swirl energy.
The pulse amplitude on Interval O.4 is therefore chosen afterwards.

\noindent\textbf{Parameter status.}
The input parameter set is unchanged. The radius $R_{3/2}$ is derived;
$d_1,d_2$ are pending derived functions until the tail is fixed.
After axial matching, $U^r=U^z=0$ throughout this stage.

\noindent\textbf{To be Checked:}
\begin{enumerate}[(i)]
\item \emph{Moments related to $U^\theta$.}
Start from the incoming moments at $R_{m\theta}$. Verify both matching conditions
in \eqref{eq:road-angular-matching} at $R_{3/2}$ for every $Z$. They give
\[
 M^\theta(R_{3/2},Z)=M_{\rm tar}^\theta(R_{3/2},Z;\tau),\qquad
 M^p(\infty,Z)=M_{\rm pre}^p(\infty,Z).
\]
\item \emph{Cone.}
After the moment matching, the profile satisfies the admissible cone
condition on $[R_{m\theta},R_{3/2}]$, including both bumps.
\end{enumerate}
Both checks are carried out in Section~\ref{subsec:outer-angular-matching}.

\medskip
\Needspace{12.0cm}
\paragraph{Interval O.7. A steep descent, a terminal power, and a waiting interval.}
The waiting length sets the baseline angular match to the future tail.
It is determined before the local corrections of Interval O.6.

\begin{center}\begin{minipage}{.97\linewidth}\centering
\begin{tikzpicture}[x=1cm,y=.9cm,>=stealth,font=\footnotesize]
 \fill[roadpurple!9] (1.65,0) rectangle (2.85,3.35);
 \fill[roadpurple!9] (6.25,0) rectangle (7.65,3.35);
 \draw[->] (0,0)--(12.45,0) node[right] {$R$};
 \draw[->] (0,0)--(0,3.7) node[above] {$\overline U^\theta$};
 \draw[roadold] (.3,3.15)..controls (.8,3.04) and (1.2,2.95)..(1.65,2.84);
 \draw[roadcurve] (1.65,2.84)..controls (2.1,2.73) and (2.35,2.43)..(2.85,2.17)
 ..controls (3.85,1.65) and (4.9,1.03)..(6.25,.78)
 ..controls (6.8,.68) and (7.1,.61)..(7.65,.56);
 \draw[roadcurve] (7.65,.56)..controls (9.2,.42) and (10.55,.30)..(11.9,.24);
 \foreach \x/\lab in {.3/R_{m\theta},1.65/R_{3/2},6.25/e^{-1}R_t,11.9/R_{2h}}
 {\draw[roadguide] (\x,0)--(\x,3.15);\node[below] at (\x,0) {$\lab$};}
 \foreach \x/\lab in {2.85/eR_{3/2},7.65/R_t}
 {\draw[roadguide] (\x,-.52)--(\x,2.85);\node[below] at (\x,-.52) {$\lab$};}
 \node[align=center,roadpurple] at (2.25,3.8) {flat slope\\transition};
 \node[align=center,roadpurple] at (6.95,3.8) {flat slope\\transition};
 \node[align=center,roadblue] at (4.5,2.95) {$R^{-3/2}$\\reduced swirl-square energy};
 \node[align=center,roadblue] at (9.9,1.8)
 {$(1-\varepsilon)c_\infty R^{-(1+\delta)/2}$\\terminal power};
 \node[align=center] at (10.05,3.15) {$U^r=U^z=0$\\after matching};
 \draw[<->] (2.85,-1.18)--(6.25,-1.18)
 node[midway,below] {logarithmic length $T_s$};
 \draw[<->] (7.65,-1.18)--(11.9,-1.18)
 node[midway,below] {waiting length $\tau$};
\end{tikzpicture}
\RoadPairInfo{Left: $[R_{m\theta},R_{3/2}]$}
 {The previous interval contains both angular bumps. They end
  before $R_{3/2}$, so the endpoint value used for the
  next prescription is independent of $Z$.}
 {Right: $[R_{3/2},R_{2h}]$}
 {The two slope transitions use \eqref{eq:road-steep-slope}.
  The steep length $T_s$ fixes $R_t$. The waiting length $\tau$
  fixes $R_{2h}=e^\tau R_t$ and is selected after the heat tail
  is specified.}

\par\small Only the preceding relaxation interval and the new descent/waiting
intervals are shown. The pure steep power holds between the two slope
transitions, on $[eR_{3/2},e^{-1}R_t]$.
\end{minipage}\end{center}

The two shaded bands are slope transitions, so the exact power
$R^{-3/2}$ holds only between $eR_{3/2}$ and $e^{-1}R_t$.
Its role is to reduce the swirl-square energy before the terminal
tail. The separate waiting interval adjusts the angular moment:
its amplitude lies below the asymptotic heat power, and varying its
length changes the backward moment target without changing
$c_\infty$. These two intervals therefore address different moment
constraints. The radial spacing in the drawing is schematic.

Introduce only the derived length $T_s=4\log(2/\delta)$ and the radii
\begin{equation}\label{eq:road-steep-radii}
 R_s=eR_{3/2},\qquad R_t=e^{T_s+1}R_s=e^{T_s+2}R_{3/2}.
\end{equation}
For $t=\log(R/R_{3/2})$ define
\begin{equation}\label{eq:road-steep-slope}
 s_7(t)=-\frac12-\mu-(1-\mu)\sigma(t)
 +(1-\delta/2)\sigma(t-1-T_s),
\end{equation}
and set
\[
 \overline{U}^\theta(R)=U^\theta(R_{3/2})\exp\!\left(\int_0^t s_7(v)\,dv\right),
 \qquad U^z=0.
\]
The two angular bumps have ended before this interval, so $U^\theta(R_{3/2})$
is independent of $Z$. On $[R_{3/2},R_s]$ the slope moves to $-3/2$;
on $[R_s,e^{-1}R_t]$ the swirl is an exact multiple of $R^{-3/2}$;
on $[e^{-1}R_t,R_t]$ its slope moves to $-(1+\delta)/2$.
All joins use the same flat step $\sigma$.

Now introduce a provisional input waiting length $\tau\ge0$ and derive
$R_{2h}=e^\tau R_t$. Continue on $[R_t,R_{2h}]$ by
\begin{equation}\label{eq:road-waiting}
 \overline{U}^\theta(R)=\overline{U}^\theta(R_t)(R/R_t)^{-(1+\delta)/2},\qquad
 U^z=0.
\end{equation}
\noindent\textbf{Remarks:}
The terminal power contributes a factor $1/\delta$ to the estimate for
$\int_{R_t}^\infty(U^\theta)^2\,dR$.
The preceding $R^{-3/2}$ segment reduces $R(U^\theta)^2$ by the factor
$(\delta/2)^8$ and absorbs this loss.
This is the swirl-square integral entering $M^{z\theta}$.
The pressure moment integrates $(U^\theta)^2/(2R)$ and has no such
$1/\delta$ loss.
On the subsequent terminal-power interval, the contribution of
$M^\theta(R_t,Z)$ relative to $R\sqrt{2R}U^\theta(R)$ decays by the factor
$(R_t/R)^{1-\delta/2}$.

The incoming angular moment is accumulated from the axis. The target is
computed backwards from the future heat-tail normalization in
\eqref{eq:road-angular-target}. They use different data and need not agree
before the matching choices. The radial powers also differ: $R^{-1/2-\mu}$
on the incoming buffer, and $R^{-(1+\delta)/2}$ in the heat asymptotics.
Accordingly, \eqref{eq:mc-incoming-angular-moment} and
\eqref{eq:mc-required-exterior-bounds} give
\[
 \frac{\bar M^\theta(R_{3/2},Z)}
 {R_{3/2}\sqrt{2R_{3/2}}\,\overline U^\theta(R_{3/2})}
 =\frac1{1-\mu}+O_{C_Z^k}(\mu^{29}),
\]
\[
 \frac{M_{\rm tar}^\theta(R,Z;\tau)}{R\sqrt{2R}\,\overline U^\theta(R,Z)}
 \longrightarrow\frac1{1-\delta/2}\qquad(R\to\infty).
\]
Moment matching compares the forward moment and the backward target at
the same radius. Between different radii, the moment continues to
accumulate by \eqref{fiveM}.

Here $\tau$ is a provisional input. Increasing $\tau$ raises the
backward pre-heat angular-moment target at $R_{3/2}$.
The matching described under Interval O.8 selects it so that this target equals the incoming background
moment at $Z=0$. This scalar matching changes a stage length.
The two local bumps on Interval O.6 subsequently solve the two
function-valued conditions \eqref{eq:road-angular-matching} for all $Z$.
These are distinct mechanisms; no additional pair of moment equations
is imposed at the junctions of the present interval.

At $R_{3/2}$, the background profile carries the angular moment
\[
\bar M^\theta(R_{3/2},Z)
=
R_{3/2}\sqrt{2R_{3/2}}\,
\bar U^\theta(R_{3/2})
\left[\frac1{1-\mu}+O(\mu^{29})\right].
\]
This quantity is determined by the preceding construction. At $R_h$, the prescribed heat field determines the velocity and shear.
After the axial moment matching, the condition
$\mathcal{T}^\theta(R_h,0)=0$ requires
\[
\left(1-\frac{\delta}{2}\right)M^\theta(R_h,0)
=
R_h\sqrt{2R_h}\,U^\theta_{\rm heat}(R_h,0)
-2R_h\mathcal{S}^\theta_{\rm heat}(R_h,0).
\]
Thus the heat field prescribes the cumulative moment that must arrive
from the left. Its leading value is
\[
M^\theta_{\rm tar}(R_h,0)
=
\frac{\sqrt2c_\infty}{1-\delta/2}R_h^{1-\delta/2}
+O(c_\infty R_h^{-\delta/2}).
\]
The required match is therefore
\[
M^\theta(R_{3/2},0)
+\int_{R_{3/2}}^{R_h}
\sqrt{2R}\,U^\theta(R,0)\,dR
=
M^\theta_{\rm tar}(R_h,0).
\]

For the present connection, the moment accumulated at $R_t$ is
\[
\begin{aligned}
\bar M^\theta(R_t,0)
&=
R_s\sqrt{2R_s}\,\bar U^\theta(R_s)\,[T_s+O(1)]\\
&\asymp T_sc_\infty R_t^{1-\delta/2}.
\end{aligned}
\]
It exceeds the heat-power primitive
\[
\frac{\sqrt2c_\infty}{1-\delta/2}R_t^{1-\delta/2}
\]
by the same order, without a small factor $\varepsilon$. On the waiting interval, we prescribe
\[
U^\theta=(1-\varepsilon)c_\infty R^{-(1+\delta)/2},
\qquad R_t\le R\le e^\tau R_t.
\]
The actual moment continues to increase. However, it increases more
slowly than the heat-power primitive. Their accumulated difference is
\[
\begin{aligned}
&\int_{R_t}^{e^\tau R_t}\sqrt{2R}
\left[c_\infty R^{-(1+\delta)/2}-U^\theta(R)\right]\,dR\\
&\qquad=
\frac{\sqrt2\,\varepsilon c_\infty R_t^{1-\delta/2}}
{1-\delta/2}
\left(e^{(1-\delta/2)\tau}-1\right).
\end{aligned}
\]
Increasing $\tau$ moves the heat connection outward while keeping
$c_\infty$ fixed. The target can then catch up with the excess
accumulated from the left. The fixed terminal connection is included
in the matching equation. We first choose $\tau$ for the pre-heat
match at $Z=0$; Interval O.6 corrects the remaining heat and
$Z$-dependent discrepancies.

The steep segment serves a different purpose. A pure terminal power
gives
\[
\int_{R_t}^{\infty}
\left[
U^\theta(R_t)(R/R_t)^{-(1+\delta)/2}
\right]^2\,dR
=
\frac{R_t[U^\theta(R_t)]^2}{\delta}.
\]
This $1/\delta$ loss enters the $M^{z\theta}$ balance in Interval O.4.
On the $R^{-3/2}$ segment,
\[
\frac{e^{-1}R_t[U^\theta(e^{-1}R_t)]^2}
{R_s[U^\theta(R_s)]^2}
=e^{-2T_s}
=(2/\delta)^{-8}
=(\delta/2)^8.
\]
Hence $T_s=4\log(2/\delta)$ absorbs the tail loss and leaves
a tail contribution of order $\delta^7$ relative to the incoming
square-integral scale. The steep segment controls this integral;
its additional angular moment is subsequently compensated
by the waiting interval.

\noindent\textbf{Parameter status.}
For the candidate family, the input data are now
$(P_*,R_{\rm ref},\delta,\tau)$. The quantities
$T_s,R_s,R_t,R_{2h},U^\theta(R_{3/2}),\overline{U}^\theta(R_t)$ are derived.
The value of $\tau$ remains pending until the matching described under Interval O.8.
After the axial matching of Interval O.4, $U^r=U^z=0$ on this entire stage.

\noindent\textbf{To be Checked:}
\begin{enumerate}[(i)]
\item \emph{Moments related to $U^\theta$.}
Verify that the waiting length has a unique choice satisfying
\[
 M_{\rm tar,pre}^\theta(R_{3/2},0;\tau)=\bar M^\theta(R_{3/2},0),
\]
as specified in \eqref{eq:road-wait-choice}. Also verify the tail-moment
bounds in \eqref{eq:mc-required-tail-bounds}, including control of the
$1/\delta$ loss in the swirl-square integral. The remaining angular and
pressure matching conditions for all $Z$ are imposed on Interval O.6.
Continue the matched moments from $R_{3/2}$ to $R_{2h}$.
They then agree with the tail targets in \eqref{eq:mc-required-moments}
at $R_{2h}$ by \eqref{fiveM}; no new matching conditions are imposed.
\item \emph{Cone.}
After moment closure, the profile satisfies the admissible cone condition
on $[R_{3/2},R_{2h}]$.
\end{enumerate}
Both checks will be carried out in Sections~\ref{subsec:outer-waiting} and~\ref{subsec:outer-cone-closure}.

\medskip
\Needspace{11.5cm}
\paragraph{Interval O.8. Join to the heat collar and the exact heat exterior.}
Complete the heat tail and determine the pending matching data.
After moment closure, $[R_{2h},R_h)$ is admissible and
$R\ge R_h$ is stress-free. These conclusions are proved in
Section~\ref{sec:heat-exterior} and Sections~\ref{subsec:outer-waiting} and~\ref{subsec:outer-axial-matching}--\ref{subsec:outer-cone-closure}.

\begin{center}\begin{minipage}{.97\linewidth}\centering
\begin{tikzpicture}[x=1cm,y=.85cm,>=stealth,font=\footnotesize]
 \fill[roadpurple!9] (3,0) rectangle (5.15,3.05);
 \draw[->] (0,0)--(12.45,0) node[right] {$R$};
 \draw[->] (0,0)--(0,3.7) node[above] {$\overline U^\theta$};
 \draw[roadold] (.4,2.75)..controls (1.2,2.57) and (2.2,2.30)..(3,2.1);
 \draw[roadcurve] (3,2.1)..controls (3.8,1.90) and (4.4,1.63)..(5.15,1.47)
 ..controls (6.15,1.25) and (7.2,1.02)..(8.25,.90);
 \draw[thick,black!55,dashed] (5.15,1.72)..controls (6.15,1.47) and (7.2,1.02)..(8.25,.90);
 \draw[very thick,black] (8.25,.90)..controls (9.4,.77) and (10.6,.58)..(11.9,.48);
 \foreach \x/\lab in {.4/R_t,3/R_{2h},8.25/R_h}
 {\draw[roadguide] (\x,0)--(\x,3.0);\node[below] at (\x,0) {$\lab$};}
 \draw[roadguide] (5.15,-.55)--(5.15,3.0);
 \node[below] at (5.15,-.55) {$eR_{2h}$};
 \node at (1.55,3.2) {waiting power};
 \node[roadpurple,align=center] at (4.1,3.55) {flat blend\\at both endpoints};
 \node[roadblue] at (6.8,2.9) {heat collar};
 \node[align=center,black!65] at (6.8,2.2) {heat continuation};
 \draw[->,black!55] (6.55,1.98)--(6.1,1.48);
 \node[align=center] at (10.25,3.25) {exact heat flow\\$\mathcal{T}=0$};
 \node[align=center] at (10.25,2.2) {$U^r=U^z=0$\\after matching};
 \draw[->,roadpurple] (8.7,1.68)--(8.25,.92);
 \node[anchor=west,roadpurple] at (8.7,1.68) {flat heat join};
 \draw[<->,roadblue] (3,-1.17)--(8.25,-1.17)
 node[midway,below] {admissible collar after matching};
 \draw[->,black] (8.25,-1.17)--(12,-1.17)
 node[midway,below] {stress-free exterior};
 \node[below] at (11.9,0) {$\longrightarrow\infty$};
\end{tikzpicture}
\RoadPairInfo{Left: $[R_t,R_{2h}]$}
 {The amplitude obeys
  $(1-\varepsilon)c_\infty=\overline U^\theta(R_t)R_t^{(1+\delta)/2}$.
  There is no heat factor on
  this interval; its logarithmic length is $\tau$.}
 {Right: $[R_{2h},\infty)$}
 {The exact exterior factor is $H_\delta(2L/R)$ in
  \eqref{eq:road-heat-exact}. The deficit
  $\varepsilon \mathfrak f((3-s)/2)$ is flat at $R_h=e^3R_{2h}$.
  The heat amplitude is independent of $\tau$.}

\par\small Only the waiting interval, terminal connection, and heat exterior
are shown. The point $R_h$ is a smooth join, with a deficit flat to every order.
\end{minipage}\end{center}

The shaded blend inserts the heat factor smoothly. In the following
collar, the solid curve differs from the dashed heat continuation by
the prescribed flat deficit. That deficit and every radial derivative
vanish at $R_h=e^3R_{2h}$, so the exact heat exterior is attached without
a corner or a residual supported at the join. The stress is nonzero
inside the active collar and vanishes in the exterior only after the
global moment conditions have been imposed. The size of the deficit
is enlarged in the drawing.

Fix a small absolute constant $c_\varepsilon>0$, derive
$\varepsilon=c_\varepsilon\delta$, and define
\begin{equation}\label{eq:road-heat-amplitude}
 c_\infty=\frac{\overline{U}^\theta(R_t)R_t^{(1+\delta)/2}}{1-\varepsilon},\qquad
 R_h=e^3R_{2h}.
\end{equation}
The heat amplitude $c_\infty$ is derived, not an input parameter. It does not depend on
$\tau$, because the product $R^{(1+\delta)/2}\overline{U}^\theta(R)$ is constant on the
waiting interval. Use the flat function $\mathfrak f$ and heat factor $H_\delta$
from \eqref{eq:road-heat-functions}.
For $s=\log(R/R_{2h})\ge0$, prescribe
\begin{equation}\label{eq:road-heat-join}
\begin{aligned}
 \overline{U}^\theta(R,Z)=c_\infty R^{-(1+\delta)/2}\bigl\{
 &(1-\sigma(s))(1-\varepsilon)\\
 &+\sigma(s)H_\delta(2L/R)[1-\varepsilon \mathfrak f((3-s)/2)]\bigr\},
 \qquad U^z=0.
\end{aligned}
\end{equation}
At $s=0$, this agrees to all orders with the waiting power law.
On $1\le s\le3$ it is exactly the heat collar
\[
 c_\infty R^{-(1+\delta)/2}H_\delta(2L/R)[1-\varepsilon \mathfrak f((3-s)/2)].
\]
At $R_h$ the deficit is flat. Thus for $R\ge R_h$ the formula becomes
\begin{equation}\label{eq:road-heat-exact}
 U^\theta=c_\infty R^{-(1+\delta)/2}H_\delta(2L/R),\qquad
 U^z=U^r=0,
\end{equation}
where the radial-velocity equality uses the completed axial moment matching.
The associated physical swirl solves the radial heat equation exactly.

\noindent\textbf{Parameter status.}
Before moment closure, the input data are $(P_*,R_{\rm ref},\delta,\tau)$.
The constant $c_\varepsilon$ and function $\mathfrak f$ are fixed;
$\varepsilon,c_\infty,R_h,H_\delta$ are derived. In particular,
$R_h$ and $c_\infty$ cannot be changed independently of the earlier data.

\noindent\textbf{Determine the pending coefficients.}
The uncorrected outer swirl $\overline{U}^\theta$ is now fully prescribed. Its angular
moment target, determined by the future tail, is
\begin{equation}\label{eq:road-angular-target}
\begin{aligned}
 M_{\rm tar}^\theta(R,Z;\tau)
 ={}&\frac{\sqrt2c_\infty}{1-\delta/2}R^{1-\delta/2}\\
 &-\int_R^\infty\sqrt{2\rho}
 [\overline{U}^\theta(\rho,Z;\tau)-c_\infty\rho^{-(1+\delta)/2}]\,d\rho.
\end{aligned}
\end{equation}
Define $\bar U^\theta_{\rm pre}$ by replacing $H_\delta$ with $1$
in \eqref{eq:road-heat-join} and its tail.
It equals the barred candidate on the earlier intervals.
Let $M_{\rm tar,pre}^\theta$ be the target in
\eqref{eq:road-angular-target} with $\bar U^\theta_{\rm pre}$ in place of
$\overline U^\theta$.
First determine the waiting length by
\begin{equation}\label{eq:road-wait-choice}
 M_{\rm tar,pre}^\theta(R_{3/2},0;\tau)
 =\bar M^\theta(R_{3/2},0).
\end{equation}
The left side increases continuously and strictly from below the fixed
right side at $\tau=0$ to $+\infty$, giving a unique $\tau>0$.
With $P_*,\delta$ and the auxiliary data fixed, scaling
$R=R_{\rm ref}x$ cancels the common factor $R_{\rm ref}^{3/2}$ from
both sides, so $\tau$ is independent of $R_{\rm ref}$.
At this selected value, the actual heat tail supplies the pressure deficit
\[
 \Delta_H(Z)=\int_{R_{2h}}^\infty
 \frac{(\bar U^\theta_{\rm pre})^2-(\overline U^\theta)^2}{2R}\,dR.
\]
First solve the two conditions \eqref{eq:road-angular-matching}
using the small branch in \eqref{eq:mc-angular-coefficients}.
With the resulting full swirl fixed, determine the axial coefficients
by the linear moment equations and the pulse-energy equation in Section~\ref{subsec:outer-axial-matching}.
Thus the coefficient order, distinct from the radial order, is
\[
 \boxed{\tau\ \longrightarrow\ (d_1,d_2)
 \ \longrightarrow\ (a_p,c_1,c_2).}
\]
Sections~\ref{subsec:outer-waiting} and~\ref{subsec:outer-axial-matching}--\ref{subsec:outer-cone-closure} prove that these choices
give all five terminal conditions \eqref{eq:moment-conditions},
with no further moment equations at the intermediate junctions.
The remaining inputs are $P_*,R_{\rm ref},\delta$;
all radii, the heat amplitude, and the correction coefficients are derived.

\noindent\textbf{Output for the inner replacement.}
The restored axis pressure is
\begin{equation}\label{roadin:axis-pressure}
 P_0(Z)=-\int_0^{R_{\rm ref}}
 \frac{(U_{\rm ref}^\theta)^2}{2R}\,dR
 -\int_{R_{\rm ref}}^\infty\frac{(U^\theta)^2}{2R}\,dR
 =P_0^{\rm pre}(Z).
\end{equation}
The angular correction on Interval O.6 enforces this identity, where
$P_0^{\rm pre}$ is computed from the uncorrected candidate at the selected
waiting length, with $H_\delta$ replaced by $1$.
By Proposition~\ref{prop:road-outer-summary}(iii), $P_0^{\rm pre}$ extends
holomorphically to a neighborhood of $[-1,1]$.
Its uniform real and analytic bounds are recorded in
Lemma~\ref{lem:outer-axis-pressure}.
It is derived from the completed outer construction and is then kept fixed
for the core.
The left replacement will match the five reference moments at $R_{\rm useless}$
and retain the velocities for $R\ge R_{\rm useless}$:
\[
 \mathbf M_{\rm new}(R_{\rm useless},Z)=\mathbf M_{\rm ref}(R_{\rm useless},Z).
\]
With the same $P_0$, this preserves all five moments and all fields
to the right, hence the prepared outer moment conditions and cones.
Only this endpoint equality is required; the individual moments
inside the replacement need not equal those of the reference.
Proposition~\ref{prop:global-relaxed-profile} proves this preservation.

\subsection{From the axis to the reference profile}
\label{roadin:subsection}

We replace the temporary reference profile on $[0,R_{\rm useless}]$, where
$R_{\rm useless}=e^{-5}R_{\rm ref}$, by a regular core and a connection.
Retain the completed outer profile of Proposition~\ref{prop:road-outer-summary}
and its axis pressure $P_0=P_0^{\rm pre}$ from
\eqref{roadin:axis-pressure}. Write $M_{\rm ref}^j$ for the reference
moments on $[0,R_{\rm ref}]$ and $P_{\rm ref}=P_0+M_{\rm ref}^p$.
The reference velocity is retained on $[R_{\rm useless},R_{\rm ref}]$.

The outer construction is used here as a family with scale parameter
$R_{\rm ref}$. First fix $P_*,\delta$ and the auxiliary data, including
the gluing tolerance and axis shift $j$ specified in Interval I.1.
The resulting $\tau$ and $P_0$ are independent of $R_{\rm ref}$.
Next fix $\Lambda$ large enough for the core and inner connection.
Finally take $C_*$ sufficiently large and determine the scale parameter
by $R_{\rm ref}=110(C_*P_*)^{10}$, selecting the corresponding outer
profile with the same $\tau$ and $P_0$. With the preceding data fixed,
all remaining compatibility and inner moment-correction conditions
hold for every sufficiently large $C_*$ by
Lemma~\ref{lem:imc-prepared-family}.

The construction is organized by Intervals I.1--I.4: the core, its smooth
continuation, the return to the reference velocities, and the inner
moment-correction patch. The constructions on Intervals I.1--I.3 verify
regularity, smooth joining and the stated cone bounds. Exact reference-moment
matching is imposed only on Interval I.4, a patch inside the terminal part
of Interval I.3. It gives the same five moments at $R_{\rm useless}$; with velocity
agreement near $R_{\rm useless}$ and the same $P_0$, the prepared outer fields are
then recovered exactly.
Subsection~\ref{roadin:shear-modification} performs the final shear
modification on this completed profile.

Write $F=U^\theta/\sqrt{2R}$, $L=1-Z^2$, and $L_\delta=1-\delta Z^2$.
All actual moments are integrated continuously from the axis by
\eqref{fiveM}: they are never reset at a join or replaced by
comparison moments. Their estimates enter the pressure and cone
bounds, but no new moment identities are required at intermediate
joins. Throughout the construction,
\begin{equation}\label{roadin:recovery}
\begin{aligned}
P&=P_0+M^p,\\
U^r&=\frac{2ZR\,U^z-\mathcal A M^z}{L_\delta\sqrt{2R}},
\qquad \mathcal A q=(1-\delta)Zq+L\partial_Zq.
\end{aligned}
\end{equation}
Thus $U^r$ is recovered from $U^z$ and its accumulated moment,
rather than interpolated independently. Before the correction on Interval I.4, the
forward pressure need not yet have its final normalization at infinity.

\providecommand{\RoadInStrip}[8]{%
\begin{center}
\begin{tikzpicture}[x=1cm,y=.72cm,>=stealth,
 every node/.style={font=\footnotesize}]
 \fill[blue!4] (0,0) rectangle (12.4,2.1);
 \draw[blue!25] (0,0) rectangle (12.4,2.1);
 \node[blue!65!black,align=center,text width=11.8cm] at (6.2,1.66) {#3};
 \node[red!65!black,align=center,text width=11.8cm] at (6.2,1.02) {#4};
 \node[black!75,align=center,text width=11.8cm] at (6.2,.39) {#5};
 \draw[->] (0,-.20)--(12.6,-.20);
 \draw (0,-.14)--(0,-.26) node[below right] {$#1$};
 \draw (12.4,-.14)--(12.4,-.26) node[below left] {$#2$};
 \begin{scope}[shift={(.45,2.45)}]
  \draw[->,black!50] (0,0)--(11.55,0);
  \draw[->,black!50] (0,0)--(0,2.1);
 \node[anchor=west,blue!65!black] at (.2,2.0)
    {$U^\theta$ at a fixed $Z$ (schematic)};
  \node[anchor=south west,align=left,text width=5.3cm]
    at (0,2.30) {#7};
  \node[anchor=south east,align=left,text width=5.3cm]
    at (11.25,2.30) {#8};
  #6
 \end{scope}
\end{tikzpicture}
\end{center}}
The local plots show schematic angular profiles, or their regularized
form $F$, at a fixed $Z$;
orange marks identify the part being changed. Short annotations in
each figure give the velocity components on its two pieces. The
boxes below retain longer moment formulas and stress information.
The scales are not quantitative.

\Needspace{8cm}
\begin{center}\begin{minipage}{.97\linewidth}\centering
\begin{tikzpicture}[x=1cm,y=1cm,>=stealth,font=\small]
 \node[anchor=west,font=\bfseries\small] at (0,1.55)
  {Inner replacement: location of the five-moment correction};
 \fill[roadpurple!15] (6.1,.05) rectangle (7.55,.42);
 \draw[roadpurple,thick] (6.1,.05) rectangle (7.55,.42);
 \draw[->] (0,0)--(12.25,0) node[right] {$R$};
 \foreach \x/\lab in {.15/0,1.35/R_a,6.1/R_m,7.55/2R_m,9.25/R_{\rm useless},11.6/R_{\rm ref}}
  {\draw (\x,-.06)--(\x,.07);\node[below] at (\x,-.06) {$\lab$};}
 \node[align=center,font=\footnotesize,text width=5.8cm,roadpurple] at (6.5,.90)
  {Inner moment correction (I.4)\\$[R_m,2R_m]$};
 \node[align=center,font=\footnotesize] at (10.4,.73)
  {reference profile\\after matching};
 \begin{scope}[yshift=-3.1cm]
  \node[anchor=west,font=\bfseries\small] at (0,1.55)
   {Final modification: shear first, moment restoration afterwards};
  \fill[orange!15] (1.7,.05) rectangle (5.8,.42);
  \draw[orange!80!black,thick] (1.7,.05) rectangle (5.8,.42);
  \fill[roadpurple!15] (8,.05) rectangle (9.25,.42);
  \draw[roadpurple,thick] (8,.05) rectangle (9.25,.42);
  \draw[->] (0,0)--(12.25,0) node[right] {$R$};
  \foreach \x/\lab in {.35/R_a,1.7/R_{-c},5.8/R_{+c},8/R_c,9.25/2R_c,11.6/R_h}
   {\draw (\x,-.06)--(\x,.07);\node[below] at (\x,-.06) {$\lab$};}
  \node[align=center,font=\footnotesize,text width=4.5cm,orange!80!black] at (3.75,.90)
   {Shear modification\\$[R_{-c},R_{+c}]$};
  \node[align=center,font=\footnotesize,text width=5.8cm,roadpurple] at (8.8,.90)
   {Outer moment restoration\\$[R_c,2R_c]\subset I_1$};
 \end{scope}
\end{tikzpicture}
\par\small The two rows use separate schematic radial scales;
the lower row anticipates Subsection~\ref{roadin:shear-modification}.
The labels anticipate the later definitions and are subject to the
scale and endpoint conditions described below. The shear-modification
interval $[R_{-c},R_{+c}]$ is not one of the three reserved outer slots.
\end{minipage}\end{center}

\Needspace{11.5cm}
\paragraph{Interval I.1. Choose the axis data and construct the core: construction on $[0,R_a]$.}
\leavevmode\par
\begin{figure}[H]
\centering
\begin{tikzpicture}[x=1cm,y=1cm,>=stealth,font=\small]
 \begin{scope}
  \node[font=\bfseries\small] at (2.8,3.45) {Regular angular profile};
  \draw[->] (0,0)--(5.8,0) node[right] {$R$};
  \draw[->] (0,0)--(0,3.05) node[above] {$F/F_0$};
  \draw[gray!55,densely dotted] (5,0)--(5,2.65);
  \draw[very thick,blue!70!black] (0,2.5)
    ..controls (1.6,1.78) and (3.2,1.1)..(5,.65);
  \draw[thick,blue!60!black] (0,2.5)
    ..controls (1.7,2.17) and (3.5,1.7)..(5,1.45);
  \draw[thick,purple!80!black,dashed] (0,2.5)
    ..controls (1.8,2.47) and (3.5,2.36)..(5,2.24);
  \node[left] at (0,2.5) {$1$};
  \node[below left] at (0,0) {$0$};
  \node[below] at (5,0) {$R_a$};
  \node[font=\footnotesize,align=center] at (2.8,2.9)
    {different fixed $Z$; schematic};
  \node[blue!70!black,align=center] at (2.6,.7)
    {$F>0$\\$\partial_RF<0$};
 \end{scope}
 \begin{scope}[shift={(7.5,1.5)}]
  \node[font=\bfseries\small] at (2.3,1.95) {Analytic domain and exit};
  \fill[blue!5] (0,-1.25) rectangle (4.55,1.25);
  \draw[blue!55,densely dashed] (-.15,-1.4) rectangle (4.7,1.4);
  \draw[blue!65!black] (0,-1.25) rectangle (4.55,1.25);
  \draw[->,gray!70] (-.35,0)--(5.25,0) node[right,black] {$R$};
  \draw[->,gray!70] (0,-1.55)--(0,1.65) node[left,black] {$Z$};
  \draw[very thick,purple!80!black] (4.3,-1.25)--(4.3,1.25);
  \node[left] at (0,1.25) {$1$};
  \node[left] at (0,-1.25) {$-1$};
  \node[below left] at (0,0) {$0$};
  \node[blue!70!black,align=center,font=\footnotesize] at (2.3,.7)
    {$F,U^z$: jointly analytic\\$\mathcal{T}=0$};
  \node[align=center,font=\footnotesize] at (2.05,-.68)
    {exit: $\kappa(R_a,Z)\ge5/2$\\for every $|Z|\le1$};
  \draw[->,purple!80!black] (3.45,-.68)--(4.27,-.68);
  \node[below,purple!80!black] at (4.3,-1.45) {$R_a$};
  \node[below,font=\footnotesize,align=center] at (2.3,-1.87)
    {analytic continuation to $4.1/\Lambda>R_a$};
 \end{scope}
\end{tikzpicture}
\caption{Two views of the stress-free core.
The left panel plots the regular swirl $F/F_0$, rather than $U^\theta$;
all curves start at one but their radial decay depends on $Z$.
The right panel marks the exit $R_a=4/\Lambda$ inside the analytic
continuation domain. The exit bound is a shear margin: $\mathcal{T}=0$ there.}
\label{fig:road-core-exit}
\end{figure}

Figure~\ref{fig:road-core-exit} separates two requirements on the axis
data. The prescribed $Z$ dependence of $F_0$ creates angular decay;
where this decay is weak, the axial shear supplies the exit margin.
The continuation past $R_a$ leaves room for the smooth comparison used
in Interval I.2. It does not extend the nonzero-stress annulus to the
left of $R_a$.
\RoadPairInfo{Axis regularity and pressure}
 {$U^\theta/\sqrt{2R}\to F_0$; the regular variables extend to the axis.\par
  $M^j(0,Z)=0$ and $P(0,Z)=P_0(Z)$.}
 {Exit stress and radial velocity}
 {Recover $U_c^r$ from the actual core moment by \eqref{roadin:recovery}.\par
  $\mathcal{T}=0$, $\kappa\ge5/2$: an exit shear margin, not a strict stress cone.}

\emph{Inherited:} the prepared outer family with fixed $P_*,\delta$
and its selected waiting length, together with $P_0=P_0^{\rm pre}$.
The restored pressure has the required analytic extension.
Its complex neighborhood is chosen before the core parameters.

\emph{Fix the gluing accuracy before the axis shift.}
Use the absolute pressure constant $C_{\rm pres}$ in \eqref{eq:inner-axis-pressure}
and set $K_N=1000(1+C_{\rm pres})$, $\epsilon_0=[10^6(1+K_N)]^{-1}$.
Fix the bump and its translates in Interval I.4, and hence the absolute
constants $C_A,C_Q,C_S$ and thresholds $t_*,\mathfrak e_*$ in
\eqref{eq:imc-thresholds}, before choosing the core.
\emph{Fixed axis choices:} take
\begin{equation}\label{roadin:j-budget}
\eta_{\rm tol}=\min\{\epsilon_0/4,\mathfrak e_*/100\},\qquad
j=\eta_{\rm tol}/8,\qquad \sigma_0=j/500.
\end{equation}
Thus $j$ and $\sigma_0$ are fixed before the analytic neighborhood
and the large core parameters. Define
\begin{equation}\label{roadin:axis-data}
\begin{aligned}
U_0^z(Z)&=4Z+j, & H_0(Z)&=\tfrac{1-\delta}{2}Z+L\,U_0^z(Z),\\
G(Z)&=\int_{Z_0}^Z\frac{L_\delta(w)H_0(w)}{H_0(w)^2+\sigma_0^2}\,dw,
&F_0(Z)&=C_*^{-1}e^{-\Lambda G(Z)},
\end{aligned}
\end{equation}
where $Z_0$ is the unique zero of $H_0$.
The shift $j>0$ moves $Z_0$ away from $Z=0$, where
$\partial_ZP_0=0$, allowing the pressure to supply axial shear
near $Z_0$. The choice of $G$ gives
\[
 H_0\partial_Z\log F_0
 =-\Lambda L_\delta\frac{H_0^2}{H_0^2+\sigma_0^2},
\]
which strengthens angular decay away from $Z_0$.
In the linear model of Section~\ref{sec:analytic-core}, fixing
$\Lambda$ and increasing $C_*$ leaves the angular contribution
to the model shear intensity unchanged, while the axial contribution
is proportional to $C_*^2$. This supplies the exit margin where
angular decay is weak. The nonlinear comparison in that section
transfers the margin to the actual core.

Choose a common complex neighborhood $\Omega$ for these coefficients
and $P_0$, and an analytic width $h_{\rm an}>0$ inside it. The derived
quantity $A_\Omega=\max_{\overline\Omega}(-\operatorname{Re}G)$ is then fixed.
The neighborhood, its width and the core threshold $\Lambda_0$ may
depend on this fixed $j$ and $\sigma_0$; no uniformity as $j\to0$ is claimed.
Proposition~\ref{prop:continuation-core} applies for every fixed $0<j\le1/20$,
$\Lambda\ge\Lambda_0$ and $C_*\ge\Lambda^2e^{\Lambda A_\Omega}$.
With $\widetilde R=\Lambda R$, the analytic
core has the form
\begin{equation}\label{roadin:core}
F_c=F_0\Phi(\widetilde R,Z),\qquad
U_c^z=U_0^z+\Lambda^{-1}\Psi(\widetilde R,Z),\qquad
\Phi(0,Z)=1,\quad\Psi(0,Z)=0.
\end{equation}
These formulas give the jointly analytic solution of the stress-free core equations. They satisfy
$F_c>0$, $\partial_R F_c<0$, and $\mathcal{T}=0$ on
$0\le R\le4.1/\Lambda$. The derived exit radius is
$R_a=4/\Lambda$, where $\kappa\ge5/2$.
With $g$ defined in \eqref{eq:g}, the axial comparison estimates required at this exit are
\begin{equation}\label{roadin:axial-budget}
\begin{aligned}
U_c^z(R_a,Z)&=4Z+j-\frac{2g(Z)}{\Lambda L_\delta(Z)}
                         +O_{C^2_Z}(\Lambda^{-2}),\\
M_c^z(R_a,Z)/R_a&=4Z+j-\frac{g(Z)}{\Lambda L_\delta(Z)}
                         +O_{C^2_Z}(\Lambda^{-2}).
\end{aligned}
\end{equation}
\noindent\textbf{Verification.}
Proposition~\ref{prop:continuation-core} gives axis regularity, the stress-free
core and the stated exit shear margin. For each fixed $j$,
\eqref{roadin:axial-budget} bounds each deviation from $4Z+j$ by $B/\Lambda$
in $C_Z^2$; their maximum is $\rho_{\rm core}$ in \eqref{eq:inner-core-axial-budget}.
The derivative estimate at the exit is \eqref{eq:inner-core-input}, proved
in Appendix~\ref{app:core-exit}. The connection uses these bounds with
the choice of $j$ in \eqref{roadin:j-budget}.

\Needspace{4cm}
\paragraph{Interval I.2. Leave the core smoothly: construction on $[R_a,100]$.}
\leavevmode\par
Use the exit data at $R_a$ to prepare a smooth comparison, then
integrate the prescribed shear direction to obtain the actual
connection. The comparison is defined up to $110$, since it also
supplies the first short switch on Interval I.3. The actual connection
starts directly at $R_a$ and matches the core smoothly.

\medskip\noindent\textbf{(a) Exit data and a smooth comparison.}

The analytic core contains $[0,R_a e^{1/100}]$
since $4e^{1/100}<4.1$.
At the exit write $f=F_c(R_a,\cdot)$, $v=U_c^z(R_a,\cdot)$ and
$m_j=M_c^j(R_a,\cdot)$ for $j\in\{\theta,z,\theta z,z\theta,p\}$.
\emph{Inherited accuracy:} use $C_{\rm pres},K_N,\epsilon_0,\eta_{\rm tol}$
and the fixed $j$ from Interval I.1. For the connection, take $\Lambda$
sufficiently large that the core estimate \eqref{eq:inner-core-axial-budget} gives
\begin{equation}\label{roadin:lambda-budget}
 \rho_{\rm core}\le B/\Lambda\le\eta_{\rm tol}/4.
\end{equation}
Then $\|v-4Z\|_{C_Z^2}+\|m_z/R_a-4Z\|_{C_Z^2}
\le2(j+\rho_{\rm core})<\epsilon_0$.
Set $\chi_c=H_0^2/(H_0^2+\sigma_0^2)$ with $H_0$ from
\eqref{roadin:axis-data}. The exit estimate \eqref{eq:inner-core-input} reads
\begin{equation}\label{roadin:exit-tests}
 |\partial_Z\log f+\Lambda G'|\le B(\sqrt{\chi_c}+\Lambda^{-1})
 \qquad (|Z|\le1).
\end{equation}
Where $\chi_c\le99/100$, the proposition gives the axial exit bound
$-2L_\delta\partial_R U_c^z(R_a,Z)\ge3j/10$.

\emph{Derived data sizes:} with core norms in $(R/R_a,Z)$, put
\begin{equation}\label{roadin:sizes}
\begin{aligned}
 A={}&10+\|\log(C_*F_c)\|_{C^3}+\|U_c^z\|_{C^3},\\
 K={}&10^6+100+R_a^{-1}+C_*+P_*+A
        +\|F_c\|_{C^3}+\|1/F_c\|_{C^3}
        +\|P_0\|_{C^3_Z}+\sum_j\|M_c^j\|_{C^3}.
\end{aligned}
\end{equation}
The exit constant $B$ depends only on $j,\delta,P_0,\Omega,h_{\rm an}$.
With these data and $\Lambda$ fixed, $A=O_\Lambda(1)$ and $K=O_\Lambda(C_*)$
as $C_*$ increases.
Use the absolute constants $K_1,c_*$ fixed in
Subsection~\ref{sec:inner-leaving-core}, with
\[
0<c_*\le\min\left\{
\frac{\epsilon_0}{10^{10}K_1},
\frac{\eta_{\rm tol}}{12K_1}\right\}.
\]
These are independent of the supplied profile.
Set $h_b=c_*K^{-100}$. Over a logarithmic interval of length $h_b$,
the shear multiplier decreases from one to $h_b$.
The second restriction on $c_*$ reserves the axial accuracy needed
after the bridge and short switches.
Use the fixed step $\sigma$ from \eqref{eq:road-sigma}.
For $y=\log(R/R_a)$ define the smooth comparison by
\begin{equation}\label{roadin:comparison}
\begin{aligned}
\alpha(y)&=1-\sigma((y-h_b)/h_b),\\
\partial_y\log\bar F&=\alpha\partial_y\log F_c,
&\partial_y\bar U^z&=\alpha\partial_y U_c^z,
&(\bar F,\bar U^z)|_{y=0}&=(f,v).
\end{aligned}
\end{equation}
For $0\le y\le h_b$ the comparison agrees with the exact core.
The cutoff decreases to zero on $[h_b,2h_b]$, with $2h_b<1/100$.
At $r_*=R_a e^{2h_b}$ the regular variables attain values
$f_*=\bar F(r_*,\cdot)$ and $v_*=\bar U^z(r_*,\cdot)$ which they retain
up to $R=110$. Continue the comparison's own moments from the core,
using $\bar P=P_0+\bar M^p$ and $\bar U^r$ from \eqref{Ur-Uz} in the
stress formulas \eqref{Itheta-moments-direct}--\eqref{Iz-moments-direct}. Write
\[
 \mathbf q=(\bar D,\bar E)=\bar{\mathcal I}/\bar F,\qquad
 \bar H=|\mathbf q|^2/\bar D .
\]
Use $\bar D=R\bar Q/L_\delta$ and $\bar N=(\bar U^\theta)^2\bar J$ in
the radial equations \eqref{eq:normalized-Q-J}. The stress-free core
supplies $\bar D(R_a)=-2R_a\,\partial_R\log F_c(R_a)$ and
$\bar N(R_a)=-2L_\delta\,\partial_RU_c^z(R_a)$. Integration across the
short transition $[R_a,r_*]$ gives the values at $r_*$.
On $[r_*,110]$, $\bar a=\bar b=0$ and $\bar U^\theta=\sqrt{2R}\,f_*$.
Writing $H_{v_*}=(1-\delta)Z/2+Lv_*$, the two equations reduce to
\[
\begin{aligned}
 \partial_R(R\bar D)&= R S_Q/L_\delta,\qquad
 S_Q=-\bar W-\tfrac{\delta}{2}(1-2Zv_*)-H_{v_*}\partial_Z\log f_*,\\
 \partial_y\bar N+\bar N
 &=(\bar U^\theta)^2(Z+\bar\Pi)
   -\tfrac{1+\delta}{2}(1-2Zv_*)v_*-H_{v_*}\partial_Zv_* .
\end{aligned}
\]
Since $H_0G'=L_\delta\chi_c$ and $v=4Z+j+O_{C_Z^2}(\Lambda^{-1})$,
\eqref{roadin:exit-tests} gives the
\hyperlink{inner-exit-transport}{weighted angular slope bound} in
Lemma~\ref{prop:inner-comparison}. On the constant segment both
$v_*$ and the axial average $\bar M^z/R$ differ from $4Z+j$ in $C_Z^1$
by $O(\Lambda^{-1})+O_\Lambda(h_b)$. Substitution in the displayed
angular equation gives
\[
 S_Q\ge 3-\tfrac{\delta}{2}+jZ+\tfrac34\Lambda L_\delta\chi_c
            -O(\Lambda^{-1})-O_\Lambda(h_b).
\]
Choose $\Lambda$ sufficiently large and then, with $\Lambda$ fixed,
$C_*$ sufficiently large to make $h_b$ small. This gives $S_Q\ge2$
on the constant segment and $R S_Q/L_\delta\ge3\chi_c>9/4$ where
$\chi_c\ge99/100$, since $R\ge R_a=4/\Lambda$.
Lemma~\ref{prop:inner-comparison} propagates these bounds from $r_*$.
Where $\chi_c\ge99/100$ the angular equation gives $\bar D\ge9/4$.
Where $\chi_c\le99/100$ the axial equation gives $\bar E\ge9/8$,
with $\bar E=R\bar N/(L_\delta\bar U^\theta)$. Together they give
$\bar D>0$ and $\bar H\ge9/4$ on $[R_a,110]$, with $\bar D\ge4$ on
$[100,110]$.

\Needspace{11.5cm}
\noindent\textbf{(b) Smooth departure along the prescribed shear direction.}
\par
\begin{figure}[H]
\centering
\begin{tikzpicture}[x=1cm,y=1cm,>=stealth,font=\small]
 \fill[blue!6] (0,0.2) rectangle (3.0,1.15);
 \fill[purple!9] (3.0,.2) rectangle (3.65,1.15);
 \fill[orange!8] (3.65,.2) rectangle (5.25,1.15);
 \fill[gray!6] (5.25,.2) rectangle (13,.1+1.05);
 \draw[gray!35] (0,.2) rectangle (13,1.15);
 \foreach \x in {3,3.65,5.25}
   {\draw[gray!45,densely dotted] (\x,-.05)--(\x,1.15);}
 \draw[->] (0,0)--(13.4,0) node[right] {$R$};
 \foreach \x/\lab in {0/0,3/R_a,3.65/R_{ac},5.25/R_ae^{h_b},13/100}
   {\draw (\x,-.06)--(\x,.06);\node[below] at (\x,-.09) {$\lab$};}
 \node[align=center] at (1.5,.67) {analytic core\\$\mathcal{T}=0$};
 \node[align=center,font=\footnotesize] at (9.1,.67)
   {small $|\mathcal{S}|/F$; $\mathcal{S}$ nearly opposite to $\mathcal{I}$\\$\chi_b=h_b$, $\kappa<1$};
 \node[align=center,purple!80!black,font=\footnotesize] at (4.12,1.8)
   {reduce the shear on a short logarithmic interval\\$0<h_b\ll1$};
 \draw[->,purple!80!black] (4.12,1.4)--(4.12,1.05);
 \draw[<->,purple!80!black] (3,-.9)--(3.65,-.9);
 \node[purple!80!black,align=center,font=\footnotesize,anchor=east]
   at (2.8,-.9) {admissible\\inner collar};
 \draw[<->,blue!65!black] (3.65,-.9)--(13,-.9);
 \node[above,blue!65!black,font=\footnotesize] at (8.4,-.9)
   {relaxed cone on the continuation};
\end{tikzpicture}
\caption{Departure from the stress-free core. The radial order is
schematic, and the short intervals near $R_a$ are enlarged. The
admissible collar ends at $R_{ac}<R_ae^{h_b}$; the formula
$\kappa<1$ applies only after the shear reduction is complete.}
\label{fig:road-core-departure}
\end{figure}

The inner collar and the long continuation in
Figure~\ref{fig:road-core-departure} have different roles.
The small collar retains the core's strong
shear while the stress turns on from zero. Farther out, reducing the
shear makes the relaxed cone available for the long connection.
The comparison prescribes the direction, but the stress and cone
must still be computed from the actual moments of the connection.
\RoadPairInfo{Collar: actual stress and core matching}
 {Both integration equations are \eqref{roadin:bridge}; match every
  core derivative at $R_a$.\par
  Admissible cone for $R_a<R\le R_{ac}$; at $R_a$ the stress is zero.}
 {Continuation: change in cone strength}
 {Use actual moments in \eqref{roadin:recovery}, not the comparison moments.\par
  Relaxed cone throughout; $\chi_b=h_b$ and $\kappa<1$
  only once $R\ge R_ae^{h_b}$.}

\emph{Inherited and fixed:} the comparison $\mathbf q$, $h_b$,
and all core moments. Define the radial multiplier
\begin{equation}\label{roadin:bridge}
\begin{aligned}
\chi_b(y)&=1-(1-h_b)\sigma(y/h_b),\\
\partial_y\log F&=-\tfrac12\chi_b\bar D,
&\partial_y U^z&=-\chi_b\sqrt{R/2}\,F\bar E,
\qquad R_a<R\le100.
\end{aligned}
\end{equation}
The cutoff $\chi_b$ changes from one to $h_b$ on $[0,h_b]$.
The other connecting cutoffs are translates or rescalings of $\sigma$.
Integrate the angular equation with $F(R_a)=f$, then the axial
equation with $U^z(R_a,Z)=U_c^z(R_a,Z)$. Thus this is an explicit integration of known
coefficients, not a forward problem for an unknown exterior stress.
The shear is exactly $\mathcal{S}=-\chi_bF\mathbf q$. Flatness at
$y=0$ matches every derivative of the core. The comparison error
satisfies
\[
|\mathcal{I}/F-\mathbf q|
 \le 2K_1K^{20}h_b(1-\chi_b).
\]
Writing $\mathbf e=\mathcal I/F-\mathbf q$, the bound gives
$|\mathbf e\cdot\mathbf q|/\bar D\le1/16$ and
$|\mathbf e|\le(1-\chi_b)|\mathbf q|/(10\sqrt{1+\bar H})$.
Thus $\mathcal T/F=(1-\chi_b)\mathbf q+\mathbf e$ has a controlled
direction even as it vanishes at the core. The cone is admissible on the collar
\[
R_a<R\le R_{ac},\qquad
R_{ac}=R_a\exp(h_b t_{ac}),\qquad
t_{ac}=\sigma^{-1}\!\left(\frac1{16K^{10}}\right).
\]
Indeed, $\bar H\le K^{10}$ gives $(1-\chi_b)\bar H\le1/16$ there,
so $\kappa=\chi_b\bar H\ge35/16$. The relaxed cone holds on the
remaining interval. For $y\ge h_b$ the shear is small:
$\kappa=h_b\bar H<1$.

\noindent\textbf{Verification.}
The comparison estimates above supply the prescribed shear direction in
Subsection~\ref{sec:inner-leaving-core}. That section proves smooth
matching of every core jet, positivity of $F$, the displayed bridge
error bound, the relaxed cone on $(R_a,100]$, and the admissible cone
on $(R_a,R_{ac}]$. At $R=100$ put $(D,E)=\mathcal I/F$.
The bound $\bar D\ge4$ and the comparison error give $D>3$ and
\[
 D+\bar E E/\bar D
 =\bar H+\mathbf e\cdot\mathbf q/\bar D\ge35/16.
\]
These are the two endpoint margins for the test
$D+t\bar E E/\bar D>2$, $0\le t\le1$, used in the first switch of
Interval I.3.

\Needspace{6cm}
\begin{remark}[The two shear adjustments]
We first prepare the comparison by \eqref{roadin:comparison}.
It agrees with the core near $R_a$, and its
regular variables become constant in $R$ after a short transition.
The actual profile is then constructed by integrating
$\partial_y\log F=-\chi_b\bar D/2$ and the axial equation in
\eqref{roadin:bridge}. This prescription changes the background
shear itself. Together with the error bound above, it produces
nonzero stress immediately to the right of $R_a$ and the admissible
inner collar $(R_a,R_{ac}]$.
The later oscillatory shear modification in
Subsection~\ref{roadin:shear-modification} assumes that the profile
already satisfies the relaxed cone and has admissible collars at
the ends of the stress annulus. It uses the inner collar constructed
here and therefore cannot be applied directly in place of this
initial departure from the stress-free core.
\end{remark}

\Needspace{8cm}
\paragraph{Interval I.3. Restore the reference angular and axial velocities: construction on $[100,R_{\rm useless}]$.}
\leavevmode\par
This piece first restores the radial power of the swirl, then its reference
$Z$ dependence and normalization, and finally the axial reference velocity.
The three stages use the actual endpoint values of the preceding stage.

\begin{center}
\begin{tikzpicture}[x=1cm,y=1cm,>=stealth,
 every node/.style={font=\footnotesize}]
 \fill[roadorange!7] (0,.18) rectangle (3,2.55);
 \fill[roadblue!5] (3,.18) rectangle (7.65,2.55);
 \fill[roadblue!9] (7.65,.18) rectangle (13.3,2.55);
 \draw[gray!35] (0,.18) rectangle (13.3,2.55);
 \draw[gray!35] (3,.18)--(3,2.55);
 \draw[gray!35] (7.65,.18)--(7.65,2.55);
 \node[align=center,text width=2.7cm] at (1.5,1.38)
 {\textbf{(a) Power}\\[5pt]
  $b\to0$; $a\to4/5$\\[4pt]
  $U^\theta\propto R^{1/10}$\\
  $U^z\to U^z(110,Z)$};
 \node[align=center,text width=4.25cm] at (5.325,1.38)
 {\textbf{(b) Angular reference}\\[5pt]
  reshape and normalize\\[4pt]
  $U^\theta=U_{\rm ref}^\theta$ from $R_{\theta\mathrm{ok}}$\\
  $U^z=U^z(110,Z)$};
 \node[align=center,text width=5.2cm] at (10.475,1.38)
 {\textbf{(c) Axial reference}\\[5pt]
  $U^z\to4Z$ by $eR_{4Z}$\\[4pt]
  $U^\theta=U_{\rm ref}^\theta$ throughout\\
  continue both components to $R_{\rm useless}$};
 \draw[->] (0,0)--(13.65,0) node[right] {$R$};
 \foreach \x/\lab in {0/100,3/110,5.2/R_{\theta\mathrm{ok}},7.65/R_{4Z},9.35/eR_{4Z},13.3/R_{\rm useless}}
  {\draw (\x,-.06)--(\x,.08);
   \node[below] at (\x,-.08) {$\lab$};}
 \fill[roadpurple!18] (10.2,-1.1) rectangle (11.7,-.9);
 \draw[roadpurple] (10.2,-1.1) rectangle (11.7,-.9);
 \draw[roadpurple!55,densely dotted] (10.2,0)--(10.2,-.9);
 \draw[roadpurple!55,densely dotted] (11.7,0)--(11.7,-.9);
 \node[anchor=east,align=right,text width=6.4cm,roadpurple]
  at (9.8,-1) {Interval I.4: moment correction\\inside this same piece};
 \node[below,roadpurple] at (10.95,-1.1) {$[R_m,2R_m]$};
\end{tikzpicture}
\par\small Radial order is schematic; the correction interval is enlarged.
The restored components are $U^\theta$ and $U^z$.
\end{center}

\emph{(a) Remove the axial shear, then restore the angular power on $[100,110]$.}
Use the shear coordinates
\[
\mathcal{S}=F(-a,b),\qquad
a=-2\partial_y\log F=1-2\partial_y\log U^\theta,
\qquad b=2(\partial_y U^z)/U^\theta.
\]
The same inherited width $h_b$ controls both short switches; no new
input parameter is introduced. On $100\le R\le100e^{h_b}$ prescribe
\[
a=h_b\bar D,\qquad
b=-h_b\bar E
 [1-\sigma(\log(R/100)/h_b)].
\]
On $100e^{h_b}\le R\le100e^{2h_b}$, put
$t=\log(R/(100e^{h_b}))/h_b$ and prescribe
\[
b=0,\qquad
a=(1-\sigma(t))h_b\bar D+\tfrac45\sigma(t).
\]
Integrate $\partial_y\log F=-a/2$ and $\partial_y U^z=b\,U^\theta/2$,
and keep $(a,b)=(4/5,0)$ from the second endpoint to $R=110$.
Thus $U^\theta$ has radial power $1/10$, while $U^z$ is constant in $R$.
At $110$, the inherited angular inertial stress satisfies
$\mathcal{I}^\theta/F>3$. For the bridge and these switches, the axial
change from the core exit value obeys
\[
\rho_{\rm br}\le3K_1c_*K^{-80}\le\eta_{\rm tol}/4.
\]
The resulting axial accuracy, also used on Interval I.4, is
\begin{equation}\label{roadin:eta-budget}
\eta:=\|U^z(110,\cdot)-4Z\|_{C^1_Z}
\le j+\rho_{\rm core}+\rho_{\rm br}
\le5\eta_{\rm tol}/8<\mathfrak e_*/100.
\end{equation}

\emph{(b) Restore the angular $Z$ dependence and reference normalization,
then continue to $R_{4Z}$.}
Set the derived length $T=400A$ and radius $R_{\theta\mathrm{ok}}=110e^T$.
For $y=\log(R/110)$, $0\le y\le T$, prescribe
\begin{equation}\label{roadin:reshape}
U^\theta=e^{y/10}\bigl[U^\theta(110,Z)\bigr]^{1-\sigma(y/T)}
 \left[\frac1{C_*(1+Z^2)}\right]^{\sigma(y/T)},
\qquad U^z=U^z(110,Z).
\end{equation}
This interpolates the logarithm of the positive swirl. With
$B=\log\!\bigl(C_*\,U^\theta(110,Z)(1+Z^2)\bigr)$, its angular shear is
$a=4/5+2B\partial_y[\sigma(y/T)]$, so $7/10\le a\le9/10$.
The long interval controls the derivative cost of changing the
$Z$ dependence. Impose the compatibility conditions
\begin{equation}\label{roadin:compatibility}
\begin{gathered}
R_{\rm ref}=110(C_*P_*)^{10},\qquad C_*\ge e^{4A},\\
R_{\rm ref}\ge110e^{T+10}(1+A)^{10}.
\end{gathered}
\end{equation}
At $R_{\theta\mathrm{ok}}$ these give
\[
U^\theta=\frac{(R/110)^{1/10}}{C_*(1+Z^2)}
 =\frac{P_*}{1+Z^2}\left(\frac R{R_{\rm ref}}\right)^{1/10}
 =U_{\rm ref}^\theta.
\]
Keep $U^\theta=U_{\rm ref}^\theta$ and $U^z=U^z(110,Z)$ on
$[R_{\theta\mathrm{ok}},R_{4Z}]$.

\emph{(c) Restore $U^z=4Z$ on $[R_{4Z},eR_{4Z}]$ and continue to $R_{\rm useless}$.}
The derived radii are $R_{4Z}=e^{-8}R_{\rm ref}$ and
$R_{\rm useless}=e^{-5}R_{\rm ref}$. Impose the additional scale separation
\begin{equation}\label{roadin:stress-separation}
R_{4Z}\ge110(1+K)^2/P_*^2,
\end{equation}
which controls the inherited axial-stress contribution. On
$R_{4Z}\le R\le eR_{4Z}$ prescribe
\begin{equation}\label{roadin:axial}
U^z=U^z(110,Z)+(4Z-U^z(110,Z))\sigma(\log(R/R_{4Z})),
\qquad U^\theta=U_{\rm ref}^\theta,
\end{equation}
and keep $U^z=4Z$ afterwards. The actual mismatch
$U^z(110,Z)-4Z$, bounded by \eqref{roadin:eta-budget}, controls
the axial shear. Both prescribed components agree with the reference
on $[eR_{4Z},R_{\rm useless}]=[e^{-7}R_{\rm ref},R_{\rm useless}]$.

\noindent\textbf{Verification and output of the whole piece.}
Verify smooth joining at every switch, positivity of $U^\theta$,
and the relaxed cone on $[100,R_{\rm useless}]$, together with the pressure bound.
The flat cutoff gives the endpoint jets; the remaining estimates use
the inherited stress and the actual moments in \eqref{roadin:recovery}.
Check the axial budget \eqref{roadin:eta-budget}, the compatibility
conditions \eqref{roadin:compatibility}, and the stress separation
\eqref{roadin:stress-separation}, retaining the inherited smallness
and width bounds. In particular, the stages and the later correction
have the radial order
\[
100e^{2h_b}\le110<R_{\theta\mathrm{ok}}<R_{4Z}<eR_{4Z}<R_m<2R_m<R_{\rm useless},
\qquad R_m=e^{-6}R_{\rm ref}.
\]
Subsection~\ref{subsec:inner-reference-connection} verifies the switches and their
axial budget. The same subsection, in
its angular-stress and pressure estimates, proves
the relaxed cone and pressure estimates through the angular reshaping
and continuation; its axial-stress estimate controls the axial change
and its continuation to $R_{\rm useless}$.

Only the angular and axial velocities have been restored. Once both
prescribed velocities agree with the reference, the five differences
$\Delta_j=M^j-M_{\rm ref}^j$ are constant in $R$, but may be nonzero;
in particular,
\[
P-P_{\rm ref}=\Delta_p,\qquad
U^r-U_{\rm ref}^r=-\frac{\mathcal A\Delta_z}{L_\delta\sqrt{2R}}.
\]
The correction on Interval I.4 removes these accumulated defects on $[R_m,2R_m]$, inside
the terminal reference-velocity interval of this piece. This is a
correction to the profile already constructed here, not an additional
interval to the right of $R_{\rm useless}$.

\Needspace{11.5cm}
\paragraph{Interval I.4. Correct all five moments on an inner reference interval: construction on $[R_m,2R_m]$.}
\leavevmode\par
\begin{figure}[H]
\centering
\begin{tikzpicture}[x=1cm,y=1cm,>=stealth,font=\small]
 \draw[->] (0,0)--(13.4,0) node[right] {$R$};
 \foreach \x/\lab in {0/eR_{4Z},2.5/R_m,10.1/2R_m,12.7/R_{\rm useless}}
   {\draw[gray!60,densely dotted] (\x,0)--(\x,3.2);
    \node[below] at (\x,-.05) {$\lab$};}
 \draw[blue!65!black,thick] (0,2.1)--(12.9,2.7);
 \draw[purple!80!black,very thick] (0,2.1)--(3.0,2.2395)
   ..controls (3.5,2.263) and (3.55,2.75)..(4.0,2.75)
   ..controls (4.45,2.75) and (4.5,2.31)..(5.0,2.3326)
   --(5.3,2.3465)
   ..controls (5.75,2.367) and (5.8,1.92)..(6.3,1.92)
   ..controls (6.8,1.92) and (6.85,2.419)..(7.3,2.4395)
   --(7.6,2.4535)
   ..controls (8.05,2.474) and (8.1,2.98)..(8.6,2.98)
   ..controls (9.1,2.98) and (9.2,2.53)..(9.65,2.5488)
   --(12.9,2.7);
 \node[purple!80!black,above,font=\footnotesize] at (6.5,3.2)
   {three angular bumps};
 \node[blue!65!black,anchor=west,font=\footnotesize] at (.15,1.72)
   {baseline: $U^\theta_{\rm ref}$};
 \draw[red!65!black,thick] (0,.6)--(12.9,.6);
 \draw[red!75!black,very thick] (0,.6)--(3.0,.6)
   ..controls (3.5,.6) and (3.55,1.16)..(4,1.16)
   ..controls (4.45,1.16) and (4.5,.6)..(5,.6)
   --(7.6,.6)
   ..controls (8.05,.6) and (8.1,.16)..(8.6,.16)
   ..controls (9.1,.16) and (9.2,.6)..(9.65,.6)
   --(12.9,.6);
 \node[red!75!black,above,font=\footnotesize] at (10.8,.77)
   {two axial bumps};
 \node[red!65!black,anchor=west,font=\footnotesize] at (.15,.9)
   {baseline: $4Z$};
 \draw[<->,purple!80!black] (2.5,-.65)--(10.1,-.65);
 \node[below,align=center,font=\footnotesize,purple!80!black] at (6.3,-.65)
   {five moment corrections; smooth, flat joins\\
    both velocity components equal the reference outside the supports};
\end{tikzpicture}
\caption{The inner five-moment correction, with the angular and
axial changes shown separately. The rows are schematic and use
independent vertical offsets. The two axial bumps share the first
and third radial supports of the angular bumps; the displayed signs
illustrate a possible correction, not a prescribed sign pattern.}
\label{fig:road-inner-five-bumps}
\end{figure}

In Figure~\ref{fig:road-inner-five-bumps}, all five coefficients are
chosen jointly to restore the reference moments at $2R_m$ while
leaving $P_0$ unchanged. Equality of the velocity components outside
the bump supports alone would not recover the radial velocity or
pressure: this also uses equality of the accumulated moments.
The long preceding wait reduces some defects, but the axial-energy,
mixed, and linear axial moments also require the small inherited
$C_Z^1$ discrepancy between $U^z$ and $4Z$.
\RoadPairInfo{Before the first bump: inherited defects}
 {$\widehat U^r=U_{\rm ref}^r-\mathcal A\Delta_z/(L_\delta\sqrt{2R})$.\par
  On $R_m\le R\le(49/40)R_m$, the correction is zero;
  the original five moment defects remain.}
 {After all five corrections: exact field match}
 {For $2R_m\le R\le R_{\rm useless}$, all five moments and the pressure
  match exactly.\par
  The cone remains relaxed; no strict admissible cone is claimed on this patch.}

\emph{Inherited:} the completed velocity connection and its actual
moments. \emph{Derived:} let $R_m=e^{-6}R_{\rm ref}$, $x=R/R_m$,
and $A_m^\theta=e^{-3/5}P_* /(1+Z^2)$, so that $U^\theta=A_m^\theta x^{1/10}$ and $U^z=4Z$
on $[R_m,R_{\rm useless}]$. Write $\Delta_j=M^j(R_{\rm useless})-M_{\rm ref}^j(R_{\rm useless})$ and define
\begin{equation}\label{roadin:defect}
\begin{aligned}
\mathbf d=\bigg(&\frac{\Delta_z}{R_m},\
\frac{\Delta_{\theta z}-4Z\Delta_\theta}{\sqrt2R_m^{3/2}A_m^\theta},\
\frac{\Delta_\theta}{\sqrt2R_m^{3/2}A_m^\theta},\\[-1mm]
&\frac{\Delta_{z\theta}-8Z\Delta_z}{R_m(A_m^\theta)^2},\
\frac{\Delta_p}{(A_m^\theta)^2}\bigg),\qquad
\mathfrak e=\|\mathbf d\|_{C^1_Z,\ell^1}.
\end{aligned}
\end{equation}
\emph{Fixed correction data:} the inner moment correction uses
$\beta=\beta_{1/40}$. This differs from the outer bump $\beta_{3/20}$;
the symbol $\beta$ is local to each construction. For the inner correction, fix
\[
\gamma_j(x)=\beta_{1/40}(x-s_j),\quad(s_1,s_2,s_3)=(5/4,3/2,7/4),
\qquad (b_1,b_2)=(\gamma_1,\gamma_3).
\]
The three supports are disjoint and lie inside $(1,2)$. These shapes
are fixed before any input parameter, and before $j$ in Interval I.1.
Their coefficients are determined by the moment equations.
For five coefficient functions $\mathbf c=(c_1,c_2,\xi_1,\xi_2,\xi_3)$,
the exact profile correction is
\begin{equation}\label{roadin:bumps}
\widehat U^\theta=A_m^\theta\left(x^{1/10}+f\right),\quad
\widehat U^z=4Z+g,\qquad
f=\sum_{j=1}^3\xi_j\gamma_j,\quad g=\sum_{i=1}^2c_ib_i.
\end{equation}
It is supported in $R_m<R<2R_m$. Proposition~\ref{prop:imc}
solves the exact five-dimensional moment map,
including its quadratic terms, to impose
\begin{equation}\label{roadin:five-equations}
\widehat M^j(R,Z)=M_{\rm ref}^j(R,Z),\qquad
2R_m\le R\le R_{\rm useless},\quad
j\in\{\theta,z,\theta z,z\theta,p\}.
\end{equation}
Use the constants $C_A,C_Q,C_S$ and thresholds $t_*,\mathfrak e_*$
fixed before the core in Interval I.1. The inputs are $R_m\ge16$, the
inherited average, pressure, and stress bounds, and
$\mathfrak e\le\mathfrak e_*$. A sufficient defect test is
\begin{equation}\label{roadin:defect-test}
\mathfrak e_0+16(R_{\theta\mathrm{ok}}/R_m)^{1/5}
 +3\eta+40\eta^2/P_*^2\le\mathfrak e_*.
\end{equation}
Here $\eta=\|U^z(110,\cdot)-4Z\|_{C^1_Z}$ and $\mathfrak e_0$
is the norm in \eqref{roadin:defect} evaluated with the moment
differences at $R=110$, using the same denominators $R_m,A_m^\theta$.
Lemma~\ref{lem:imc-prepared-family} establishes this full defect
budget for the family $R_{\rm ref}=110(C_*P_*)^{10}$ with sufficiently
large $C_*$, after fixing $P_*,\delta,j,\Lambda,\tau,P_0$.
It uses $A(C_*)=O(1)$ and $K(C_*)=O(C_*)$ from Interval I.2 and the
velocity bound up to $R=110$ from Interval I.3 to obtain
$\mathfrak e_0+(R_{\theta\mathrm{ok}}/R_m)^{1/5}=O(C_*^{-2})$.

\noindent\textbf{Verification and output.}
Proposition~\ref{prop:imc} verifies solvability of the exact five-moment
system, positivity of the corrected swirl and preservation of the
relaxed cone on the patch. After the patch, the reference velocity,
all five moments and $P_0$ agree, hence so do $P,U^r,\mathcal{I}$ and
all radial jets near $R_{\rm useless}$. We may therefore append the prepared
outer fields without changing their moment conditions or cone.
These are finite-radius reference targets, not five zero moments.

\Needspace{4cm}
\subsection{Shear modification and restoration of the five moments}
\label{roadin:shear-modification}

After the correction on Interval I.4, the joined profile already satisfies all five global
moment conditions \eqref{eq:moment-conditions}. The remaining
operation acts on the completed annulus, not just the inner connection.
The aim is to make the cone admissible throughout the stress annulus
while retaining the core, the axis pressure and the heat exterior.
We first realize an admissible shear by a small velocity modulation,
then remove its five moment defects in the reserved outer interval.

\Needspace{11.5cm}
\noindent The shear is modified on $[R_{-c},R_{+c}]$, with moment restoration
in the reserved outer interval $[R_c,2R_c]\subset I_1$.
\par
\begin{figure}[H]
\centering
\begin{tikzpicture}[x=1cm,y=1cm,>=stealth,font=\footnotesize,
  input/.style={blue!65!black,line width=.75pt},
  modulated/.style={orange!85!black,line width=.55pt}]
 \node[anchor=west,inner sep=0pt] at (0,4.15) {Angular velocity at fixed $Z$};
 \node[anchor=east,inner sep=0pt] at (13.1,4.15)
   {$O(N^{-1})$ change in velocity, $O(1)$ change in shear};
 \draw[input] (.45,3.65)--(1.05,3.65) node[right,black] {$U^\theta$};
 \draw[modulated] (2.05,3.65)--(2.65,3.65) node[right,black] {$U_N^\theta$};
 \draw[->,gray!75] (.45,2.45)--(.45,3.4);
 \draw[->,gray!75] (.45,2.45)--(13.1,2.45);
 \draw[input] (.45,2.76)--(12.55,3.28);
 \draw[modulated] plot[domain=0:1,samples=321,variable=\t]
   ({.45+12.1*\t},
    {2.76+.52*\t + .115*sin(360*12*(\t-.05)/.90)
       *sin(90*max(0,min(1,(\t-.05)/.10)))^2
       *sin(90*max(0,min(1,(.95-\t)/.10)))^2});
 \foreach \x/\lab in {.45/0,12.55/{\log(R_{+c}/R_{-c})}}
   {\draw[gray!75] (\x,2.51)--(\x,2.39) node[below,black] {$\lab$};}
 \node[below] at (6.5,2.40) {$y=\log(R/R_{-c})$};

 \fill[orange!13] (1.05,0) rectangle (5.45,.57);
 \fill[purple!13] (6.80,0) rectangle (9.75,.57);
 \node[orange!75!black] at (3.25,.285) {shear modulation};
 \node[purple!80!black] at (8.275,.285) {five-moment repair};
 \node[align=center,anchor=south] at (5.45,.85)
   {$U_N^\theta,U_N^z$ match the input\\moments may still differ};
 \draw[gray!75,->] (5.45,.86)--(5.45,.57);
 \node[align=center,anchor=south] at (11.25,.85)
   {all five moments and\\outer fields match the input};
 \draw[gray!75,->] (9.75,.62)--(12.95,.62);
 \draw[gray!75] (9.75,.52)--(9.75,.72);
 \draw[->] (0,0)--(13.4,0) node[right] {$R$};
 \foreach \x/\lab in {0/R_a,1.05/R_{-c},5.45/R_{+c},6.80/R_c,9.75/2R_c,11.25/R_{mz},12.95/R_h}
   {\draw (\x,.05)--(\x,-.07) node[below] {$\lab$};}
\end{tikzpicture}
\caption{Shear modification and five-moment repair. The upper panel enlarges
$[R_{-c},R_{+c}]$ and shows the modulated angular velocity rejoining the
input near both endpoints. The input curve is drawn schematically, without
tracking its radial monotonicity. The lower panel locates the separate repair
interval $[R_c,2R_c]\subset I_1$ and the exact recovery for $R\ge 2R_c$.}
\label{fig:road-shear-repair}
\end{figure}

Figure~\ref{fig:road-shear-repair} distinguishes the two effects of
the modulation. Its radial derivatives follow the admissible shear
loop, whereas its values and moment integrands change by
$O(N^{-1})$. Although the modulated angular and axial velocities
already agree with the input near $R_{+c}$, their moments can still
carry errors. The separate repair on $[R_c,2R_c]$ removes those
errors and restores the pressure, radial velocity, and stress as well.
\RoadPairInfo{Entry: exact input agreement}
 {The modulation vanishes near $R_{-c}$; all incoming moments and
  the existing admissible cone are preserved.}
 {Exit, followed by the reserved correction patch}
 {Near $R_{+c}$, $M_N^z$ and hence $U_N^r$ may still carry a moment error.\par
  After restoration at $2R_c$, all three velocity components,
  the five moments and $P$ equal the input again.}

\emph{Inherited:} strict relaxed cone margins and admissible collars
at both annular ends. The outer construction reserves
$I_1=(e^{-25}R_{mz},e^{-20}R_{mz})$ for the final moment restoration.
\emph{Derived correction radius:} choose
\[
R_c=e^{-24}R_{mz},\qquad [R_c,2R_c]\subset I_1.
\]
This outer correction interval is distinct from the inner five-moment
patch $[R_m,2R_m]$ of Interval I.4 and from the shear-modification
interval $[R_{-c},R_{+c}]$.
The upper endpoint is the radius $R_{+c}=eR_{-\mu}$
already fixed in Interval O.3. Choose the lower endpoint $R_{-c}$
inside the inner admissible collar so that
\[
R_a<R_{-c}<R_{+c}<R_c<2R_c<R_h.
\]
The input is admissible on $(R_a,R_{-c}]\cup[R_{+c},R_h)$ and relaxed
on the whole stress annulus $(R_a,R_h)$. On the unused correction
patch it has the exact form
\[
U^\theta(R,Z)=U^\theta(R_c,Z)(R/R_c)^{-1/2-\mu},\qquad
U^z=0,\qquad 0<\mu\le\tfrac12.
\]
These are the geometric and profile hypotheses of
Proposition~\ref{prop:sm}; they follow from the preceding inner
collar and the reserved outer interval.

\paragraph{The Poisson kernel and the admissible shear loop.}
Recall the shear coordinates $\mathcal{S}=F(-a,b)$, with
$\kappa=(a^2+b^2)/a$. First hold the input $F$ and inertial stress
$\mathcal{I}$ fixed. The geometric point is that shears satisfying the
admissible cone at every phase can have the original relaxed shear as
their mean; see Example~\ref{ex:sm-shear-mean} and
Figure~\ref{fig:sm-loop}. This mean constraint permits a small velocity
modulation. The full relaxed cone inequalities give a
strict margin for varying the shear direction and, where an upgrade
is needed, taking its length parameter slightly above $\kappa=2$.
The two scalar inequalities in \eqref{eq:sm-scalar-cone} determine
which directions give the admissible cone.

To construct the loop explicitly, hold $(R,Z)$ fixed and put
$p=\mathcal{I}/F=(p_1,p_2)$ and $t_0=-b/a$.
The positive scales $d_*,\eta$ are determined from the cone margins
of this completed input by \eqref{eq:sm-loop-scales}--\eqref{eq:sm-eta}.
Here $\eta$ is local to the loop construction and is unrelated to the
axial error in \eqref{roadin:eta-budget}.
For $|r|<1$, the period-$2\pi$ Poisson kernel is
\[
 \mathsf P_r(\psi)=\frac{1-r^2}{1-2r\cos\psi+r^2}
 =1+2\sum_{n=1}^\infty r^n\cos(n\psi),\qquad
 \frac1{2\pi}\int_0^{2\pi}\mathsf P_r(\psi)\,d\psi=1.
\]
Set
\[
\begin{aligned}
 q&=\begin{cases}
 \displaystyle\sigma\!\left(\frac{2+\eta-\kappa}{\eta}\right)
 \sqrt{\frac{2+2\eta-\kappa}{2a}},&\kappa<2+\eta,\\
 0,&\kappa\ge2+\eta,
 \end{cases}
 &v&=\kappa+2aq^2,\\
 u&=\frac{p_2q}{d_*},&h&=\sqrt{1+u^2},\qquad r=\frac uh.
\end{aligned}
\]
Here $r$ is the kernel parameter, not a radial coordinate. The direction
and the change of phase are
\[
 t(\psi)=t_0+\frac{2q}{h}\frac{\cos\psi-r}{1-2r\cos\psi+r^2},
 \qquad \phi(\psi)=\frac{a}{2\pi v}\int_0^\psi[1+t(s)^2]\,ds.
\]
For $p_2\ne0$, the direction is equivalently
$t-t_0=d_*(\mathsf P_r-1)/p_2$; the displayed formula also covers $p_2=0$.
The phase is strictly increasing and satisfies
$\phi(\psi+2\pi)=\phi(\psi)+1$. With its inverse $\psi=\psi(\phi)$, define
\[
 (a_L,-b_L)(\phi)=\frac{v}{1+t(\psi(\phi))^2}(1,t(\psi(\phi))).
\]
The phase change weights directions by $1+t^2$ to preserve the shear mean.
The resulting smooth period-one loop satisfies
\[
 \int_0^1(a_L,b_L)(R,Z,\phi)\,d\phi=(a,b)(R,Z),\qquad
 \mathcal{S}_L=F(-a_L,b_L),\quad
 \mathcal{T}_L=\mathcal{I}+\mathcal{S}_L.
\]
The pair $(\mathcal{S}_L,\mathcal{T}_L)$ satisfies the admissible cone
at every phase, with a positive margin on the compact modification
region. Near its radial boundaries $q=0$, so the loop is constant and
equals the input. Thus the change of shear has zero mean and can be
integrated into a small velocity change. Subsection~\ref{subsec:sm-loop}
proves this loop construction and the displayed mean identity.

\paragraph{Zero-mean periodic primitives.}
For a period-one function $f$ with $\int_0^1 f=0$, put
\[
(Qf)(\phi)=\int_0^\phi f(s)\,ds
 -\int_0^1\int_0^t f(s)\,ds\,dt.
\]
Then $Qf$ has period one, has mean zero, and $\partial_\phi Qf=f$.
Use this periodic loop to define, with $(R,Z)$ held fixed,
\[
\mathcal A_L=-\tfrac12 Q(a_L-a),\qquad
\mathcal B_L=\tfrac12 U^\theta Q(b_L-b).
\]
These are the zero-mean periodic primitives used in the shear modification. They satisfy
\[
\partial_\phi\mathcal A_L=-\tfrac12(a_L-a),\qquad
\partial_\phi\mathcal B_L=\tfrac12U^\theta(b_L-b).
\]
They vanish near the radial boundaries of the modification interval
$(R_{-c},R_{+c})$. \emph{Last input choice:} after computing the actual
input norms and cone margins, choose the integer frequency $N$
sufficiently large and set, with $y=\log(R/R_{-c})$,
\begin{equation}\label{roadin:upgrade}
U_N^\theta=U^\theta\exp\!\left[N^{-1}\mathcal A_L(R,Z,Ny)\right],
\qquad
U_N^z=U^z+N^{-1}\mathcal B_L(R,Z,Ny).
\end{equation}
\paragraph{Small velocity changes and the actual stress.}
The values change little while their radial derivatives realize
the admissible shear loop. Indeed, at $\phi=Ny$ the radial
operator is $R\partial_R=\partial_y+N\partial_\phi$.
If $(a_N,b_N)$ are the shear coordinates of the modulated velocity,
differentiation gives
\[
 a_N=a_L-\frac{2\partial_y\mathcal A_L}{N},\qquad
 b_N=e^{-\mathcal A_L/N}
 \left(b_L+\frac{2\partial_y\mathcal B_L}{NU^\theta}\right).
\]
Here all terms on the right are evaluated at $\phi=Ny$, and
$\partial_y$ holds the phase fixed. Thus the factor $N$ from
radial differentiation cancels the modulation amplitude $N^{-1}$,
while the remaining errors are small for large $N$.
The phase is independent of $Z$, so the velocity changes and their
first $Z$ derivatives are also small. The five moment integrands
contain no radial derivatives of the velocity; their changes, and hence
the induced inertial-stress change, can therefore be controlled.
We use the actual moments integrated from the axis and the same
pressure convention $P_N=P_0+M_N^p$. The perturbation estimates \eqref{eq:sm-first-estimates},
together with \eqref{Ur-Uz} and the stress formulas
\eqref{Itheta-moments-direct}--\eqref{Iz-moments-direct}, quantify this control.
The dependence on the fixed input and loop is recorded explicitly
by the derived norm \eqref{eq:sm-K}; the frequency is chosen
after these data.

\paragraph{Restoring the moments and retaining the cone.}
For the final moment repair, use the three angular bump shapes from
Interval I.4 and the two axial bumps $(\beta_1,\beta_2)=(\gamma_1,\gamma_3)$
in the coordinate $x=R/R_c$. These shapes are fixed before any input
parameter. Their coefficients are determined by the moment equations.
Recompute all moments with the same $P_0$, then use the exact moment
map \eqref{eq:sm-moment-map} to restore them by five bumps on
$[R_c,2R_c]$. On this reserved power-law interval, the linearized
five-moment map is invertible for $\mu>0$. Lemma~\ref{lem:sm-moments}
therefore solves the small nonlinear defect equations with small
coefficients; the loss as $\mu$ decreases is included in the final
frequency requirement. The target is now equality with
the completed profile \emph{before the shear modification} for
every $R\ge2R_c$, rather than the reference targets in Interval I.4.
This restores all five global moment conditions and the exact heat
fields. It remains to pass from the frozen inertial stress used for
the loop to the actual stress after both modifications. The complete
error estimate \eqref{eq:sm-final-errors} and the stability bound \eqref{eq:sm-stability} give this passage, including
on the correction patch, which already has an admissible margin.
Choose $N$ according to \eqref{eq:sm-frequency-choice}, after the
input norms, cone margins and moment-correction losses are fixed.
This is one finite choice; no limit in $N$ is required.
Proposition~\ref{prop:sm} then proves that the final cone is
admissible throughout $(R_a,R_h)$, including the correction patch.
The stress-free core, the axis pressure and the heat exterior are
preserved; no strict cone is imposed on $[0,R_a]$ or $[R_h,\infty)$.

Lemma~\ref{lem:imc-prepared-family} and Subsection~\ref{subsec:sm-completion} make all the choices above
simultaneously: after the outer dimensionless data and $\Lambda$
are fixed, increase $C_*$ with $R_{\rm ref}=110(C_*P_*)^{10}$, and choose
the modulation frequency last. The resulting profile has the five
global moment conditions, a regular stress-free core, an exact heat
exterior and the admissible cone on the intervening stress annulus,
as asserted in Theorem~\ref{thm:leading}.

\label{roadmap:last-page}

\newcommand{\MergedHeatFigure}{%
\begin{figure}[htbp]
\centering
\begin{tikzpicture}[x=1cm,y=.85cm,>=stealth,font=\footnotesize]
 \fill[roadpurple!9] (3,0) rectangle (5.15,3.05);
 \draw[->] (0,0)--(12.45,0) node[right] {$R$};
 \draw[->] (0,0)--(0,3.7) node[above] {$\overline U^\theta$};
 \draw[roadold] (.4,2.75)..controls (1.2,2.57) and (2.2,2.30)..(3,2.1);
 \draw[roadcurve] (3,2.1)..controls (3.8,1.90) and (4.4,1.63)..(5.15,1.47)
 ..controls (6.15,1.25) and (7.2,1.02)..(8.25,.90);
 \draw[thick,black!55,dashed] (5.15,1.72)..controls (6.15,1.47) and (7.2,1.02)..(8.25,.90);
 \draw[very thick,black] (8.25,.90)..controls (9.4,.77) and (10.6,.58)..(11.9,.48);
 \foreach \x/\lab in {.4/R_t,3/R_{2h},8.25/R_h}
 {\draw[roadguide] (\x,0)--(\x,3.0);\node[below] at (\x,0) {$\lab$};}
 \draw[roadguide] (5.15,-.55)--(5.15,3.0);
 \node[below] at (5.15,-.55) {$eR_{2h}$};
 \node at (1.55,3.2) {waiting power};
 \node[roadpurple,align=center] at (4.1,3.55) {flat blend\\at both endpoints};
 \node[roadblue] at (6.8,2.9) {heat collar};
 \node[align=center,black!65] at (6.8,2.2) {heat continuation};
 \draw[->,black!55] (6.55,1.98)--(6.1,1.48);
 \node[align=center] at (10.25,3.25) {exact heat flow\\$\mathcal{T}=0$};
 \node[align=center] at (10.25,2.2) {$U^r=U^z=0$\\after matching};
 \draw[->,roadpurple] (8.7,1.68)--(8.25,.92);
 \node[anchor=west,roadpurple] at (8.7,1.68) {flat heat join};
 \draw[<->,roadblue] (3,-1.17)--(8.25,-1.17)
 node[midway,below] {admissible collar after matching};
 \draw[->,black] (8.25,-1.17)--(12,-1.17)
 node[midway,below] {stress-free exterior};
 \node[below] at (11.9,0) {$\longrightarrow\infty$};
 \draw[roadguide] (6.7,0)--(6.7,1.2);
 \node[below] at (6.7,0) {$R_{\rm col}=e^2R_{2h}$};
\end{tikzpicture}
\caption{The exact heat exterior and its flat inward connection (Interval O.8).
The shaded blend inserts $H_\delta$. The dashed curve is the direct heat
continuation and the deficit is enlarged. Proposition~\ref{prop:heat-collar}
proves admissibility for the right-side data on $[R_{\rm col},R_h)$.
Section~\ref{subsec:outer-waiting}
extends the cone to the full marked interval $[R_{2h},R_h)$.
Moment matching identifies these fields with the completed profile. At $R_h$ the deficit
vanishes to every order and the stress-free exterior begins.}
\label{fig:outer-collar}
\end{figure}
}

\newcommand{\MergedTailFigure}{%
\begin{figure}[htbp]
\centering
\begin{tikzpicture}[x=1cm,y=.9cm,>=stealth,font=\footnotesize]
 \fill[roadpurple!9] (1.65,0) rectangle (2.85,3.35);
 \fill[roadpurple!9] (6.25,0) rectangle (7.65,3.35);
 \draw[->] (0,0)--(12.45,0) node[right] {$R$};
 \draw[->] (0,0)--(0,3.7) node[above] {$\overline U^\theta$};
 \draw[roadold] (.3,3.15)..controls (.8,3.04) and (1.2,2.95)..(1.65,2.84);
 \draw[roadcurve] (1.65,2.84)..controls (2.1,2.73) and (2.35,2.43)..(2.85,2.17)
 ..controls (3.85,1.65) and (4.9,1.03)..(6.25,.78)
 ..controls (6.8,.68) and (7.1,.61)..(7.65,.56);
 \draw[roadcurve] (7.65,.56)..controls (9.2,.42) and (10.55,.30)..(11.9,.24);
 \foreach \x/\lab in {.3/R_{m\theta},1.65/R_{3/2},6.25/e^{-1}R_t,11.9/R_{2h}}
 {\draw[roadguide] (\x,0)--(\x,3.15);\node[below] at (\x,0) {$\lab$};}
 \foreach \x/\lab in {2.85/eR_{3/2},7.65/R_t}
 {\draw[roadguide] (\x,-.52)--(\x,2.85);\node[below] at (\x,-.52) {$\lab$};}
 \node[align=center,roadpurple] at (2.25,3.8) {flat slope\\transition};
 \node[align=center,roadpurple] at (6.95,3.8) {flat slope\\transition};
 \node[align=center,roadblue] at (4.5,2.95) {$R^{-3/2}$\\reduced swirl-square energy};
 \node[align=center,roadblue] at (9.9,1.8)
 {$(1-\varepsilon)c_\infty R^{-(1+\delta)/2}$\\terminal power};
 \node[align=center] at (10.05,3.15) {$U^r=U^z=0$\\after matching};
 \draw[<->] (2.85,-1.18)--(6.25,-1.18)
 node[midway,below] {logarithmic length $T_s$};
 \draw[<->] (7.65,-1.18)--(11.9,-1.18)
 node[midway,below] {waiting length $\tau$};
\end{tikzpicture}
\caption{Steepening and waiting serve different moment conditions (Interval O.7).
The exact $R^{-3/2}$ segment has logarithmic length $T_s=4\log(2/\delta)$;
it absorbs the $1/\delta$ cost of the tail energy. The lower terminal
power waits for length $\tau$ to adjust the angular moment.
The two shaded slope transitions are flat at their endpoints.}
\label{fig:outer-terminal-stage}
\end{figure}
}

\newcommand{\MergedTurnoffFigure}{%
\begin{figure}[htbp]
\centering
\begin{tikzpicture}[x=1cm,y=.75cm,>=stealth,font=\small]
 \draw[->] (0,0)--(12.2,0) node[right] {$R$};
 \draw[->] (0,0)--(0,3.25) node[above] {$\overline U^\theta$};
 \draw[roadold] (0,0)..controls (0,1.5) and (1.1,2.55)..(2.8,2.75);
 \draw[roadcurve] (2.8,2.75)..controls (3.4,2.85) and (4,2.85)..(4.9,2.25)
 ..controls (6.9,1.2) and (9,.65)..(11.3,.35);
 \foreach \x/\lab in {0/0,2.8/R_{\rm ref},4.9/eR_{\rm ref},11.3/R_{-\mu}}
 {\draw[roadguide] (\x,0)--(\x,2.9);\node[below] at (\x,0) {$\lab$};}
 \node[above] at (3.9,3.1) {slope turn};
 \node[above] at (8,1.35) {$R^{-1/2}$};
 \node[align=center,font=\scriptsize,text width=2.7cm] at (1.35,1.15)
 {$U^\theta_{\rm ref}\propto R^{1/10}$\\
  $U^z_{\rm ref}=4Z$\\$U^r_{\rm ref}\not\equiv0$};
 \node[align=center,font=\scriptsize] at (3.85,1.25)
 {$B=1$\\$\overline U^z=4Z$};
 \node[align=center] at (7.7,2.7)
 {$\overline U^z:4Z\ \longrightarrow\ 0$\\$\overline U^r$ generally nonzero};
 \node[align=center,font=\scriptsize] at (10.4,1.6)
 {$B=0$\\$\overline U^z=0$};
 \draw[roadguide] (9.5,0)--(9.5,1.05);
 \draw (9.5,.06)--(9.5,-.18)
 node[below,font=\scriptsize] {$e^{-11}R_{-\mu}$};
\end{tikzpicture}
\caption{Turning the swirl and switching off the axial velocity (Interval O.2).
First $\mathcal{S}^z=0$ while the swirl slope changes; then $a=2$ throughout
the slow axial turnoff. With the original choice $T_d=e^{M_d}+10$,
the last eleven logarithmic units have $U^z=0$.
The cumulative axial moment, and hence $U^r$, need not vanish.}
\label{fig:outer-axial-turnoff}
\end{figure}
}

\newcommand{\MergedBufferFigure}{%
\begin{figure}[H]
\centering
\begin{tikzpicture}[x=1cm,y=.9cm,>=stealth,font=\small]
 \draw[->] (0,0)--(12,0) node[right] {$R$};
 \draw[->] (0,0)--(0,3.2) node[above] {$U^\theta$};
 \draw[roadold] (.3,2.85)..controls (1.7,2.4) and (2.6,2.05)..(3.3,1.95);
 \draw[roadcurve] (3.3,1.95)..controls (4.1,1.8) and (4.5,1.58)..(5.5,1.35)
 ..controls (7.6,.9) and (9.6,.42)..(11.4,.25);
 \foreach \x/\lab in {.3/eR_{\rm ref},3.3/R_{-\mu},5.5/R_{+c},11.4/R_{mz}}
 {\draw[roadguide] (\x,0)--(\x,2.9);\node[below] at (\x,0) {$\lab$};}
 \node at (1.7,3.1) {$R^{-1/2}$};
 \node[align=center,font=\scriptsize,text width=2.9cm] at (1.7,1.05)
 {$U^z:4ZB\to0$\\$U^r$ generally nonzero};
 \node at (8.3,1.45) {$R^{-1/2-\mu}$};
 \node[align=center] at (7.6,2.65)
 {$U^z=0$\\[3pt]$U^r\not\equiv0$};
 \draw[<->] (5.5,-.65)--(11.4,-.65)
 node[midway,below] {logarithmic length $T_w$};
\end{tikzpicture}
\par\medskip
\begin{tikzpicture}[x=.43cm,y=.75cm,>=stealth,font=\small]
 \foreach \a/\b/\lab in {0/5/I_1,11/16/I_2,17/22/I_3}
 {\fill[roadorange!15] (\a,.05) rectangle (\b,1.45);
  \draw[roadorange!60] (\a,.05) rectangle (\b,1.45);
 }
 \fill[roadpurple!28] (1,.08) rectangle ({1+ln(2)},1.42);
 \draw[roadpurple,thick] (1,.08) rectangle ({1+ln(2)},1.42);
 \node[align=center,font=\scriptsize] at (3.5,.78) {$I_1$\\five\\moments};
 \node[align=center,font=\scriptsize] at (13.5,.78) {$I_2$\\higher\\orders};
 \node[align=center,font=\scriptsize] at (19.5,.78) {$I_3$\\phase\\mean};
 \draw[->] (-.5,0)--(26,0) node[right] {$\log(R/R_{mz})$};
 \foreach \x/\lab in {0/-25,5/-20,11/-14,16/-9,17/-8,22/-3,25/0}
 {\draw (\x,0)--(\x,-.12);}
 \foreach \x/\lab in {0/-25,5/-20,11/-14,16/-9,22/-3,25/0}
 {\node[below] at (\x,-.12) {$\lab$};}
 \node[below] at (17,-.7) {$-8$};
 \node[above] at (12.5,1.65)
 {reserved inside the exact $R^{-1/2-\mu}$ buffer};
 \draw[roadpurple,thin,->] ({1+ln(2)/2},.08)
 --({1+ln(2)/2},-1.12);
 \node[anchor=north west,align=left,font=\scriptsize,roadpurple,
 text width=7.0cm] at (0,-1.30)
 {$[R_c,2R_c]\subset I_1$,\quad $R_c=e^{-24}R_{mz}$\\
  later outer five-moment restoration};
\end{tikzpicture}
\caption{The long buffer (Interval O.3) and its reserved subintervals.
The extra slope $-\mu$ produces a positive angular source even at $Z=0$.
The uncorrected $M^z,M^{\theta z}$ stay constant while their repair
scales grow. The lower panel resolves the final $25$ logarithmic units;
$[R_c,2R_c]\subset I_1$ is reserved for the later five-moment repair,
with $R_c=e^{-24}R_{mz}$. None of these slots is used by the present corrections.}
\label{fig:outer-power-stage}
\end{figure}
}

\newcommand{\MergedRemoveZFigure}{%
\begin{figure}[H]
\centering
\begin{tikzpicture}[x=1cm,y=.86cm,>=stealth,font=\footnotesize]
 \fill[roadpurple!5] (3.2,0) rectangle (10.1,3.45);
 \draw[->] (0,0)--(12.45,0) node[right] {$R$};
 \draw[->] (0,0)--(0,3.65) node[above] {$\overline U^\theta$};
 \foreach \x/\lab in {.4/R_{mz},3.2/R_{0z},10.1/R_{m\theta}}
 {\draw[roadguide] (\x,0)--(\x,3.35);\node[below] at (\x,0) {$\lab$};}
 \draw[thick,roadpurple,dashed] (.4,3.15)..controls (1.4,3.04) and (2.4,2.91)..(3.2,2.8);
 \draw[thick,roadpurple!75,dashed] (.4,2.52)..controls (1.4,2.44) and (2.4,2.33)..(3.2,2.24);
 \draw[thick,roadpurple!50,dashed] (.4,1.58)..controls (1.4,1.53) and (2.4,1.46)..(3.2,1.4);
 \draw[very thick,roadpurple] (3.2,2.8)..controls (5.2,2.53) and (7.3,1.02)..(10.1,.68);
 \draw[very thick,roadpurple!75] (3.2,2.24)..controls (5.2,2.02) and (7.3,1.02)..(10.1,.68);
 \draw[very thick,roadpurple!50] (3.2,1.4)..controls (5.2,1.26) and (7.3,1.02)..(10.1,.68);
 \draw[roadcurve] (10.1,.68)..controls (10.8,.595) and (11.5,.52)..(12,.47);
 \node[above,roadpurple] at (1.3,3.05) {$Z=0$};
 \node[above,roadpurple!75] at (1.35,2.43) {$Z=1/2$};
 \node[above,roadpurple!60!black] at (1.35,1.51) {$Z=1$};
 \node[align=center] at (6.7,3.15)
 {$(1+Z^2)^{-1}\longrightarrow 1/2$\\flat at both endpoints};
 \node[align=center] at (10.55,3.12)
 {$U^r=U^z=0$\\after matching};
 \node[anchor=west,align=left,roadblue] at (10.25,1.65)
 {one common\\outgoing profile};
 \draw[->,roadblue] (10.75,1.2)--(11,.58);
 \draw[<->,roadpurple] (3.2,-.65)--(10.1,-.65)
 node[midway,below] {$\log(R_{m\theta}/R_{0z})=100$};
\end{tikzpicture}
\caption{Removing the background $Z$-dependence (Interval O.5).
The flat transition turns $(1+Z^2)^{-1}$ into $1/2$ over logarithmic
length $100$. The outgoing background swirl is independent of $Z$;
the accumulated left angular moment is only nearly independent of $Z$.
Its remaining variation will be corrected on Interval O.6.}
\label{fig:outer-remove-Z-stage}
\end{figure}
}

\newcommand{\MergedAngularFigure}{%
\begin{figure}[htbp]
\centering
\begin{tikzpicture}[x=1cm,y=.9cm,>=stealth,font=\small]
 \draw[->] (0,0)--(12,0) node[right] {$R$};
 \draw[->] (0,0)--(0,3.1) node[above] {$U^\theta$};
 \draw[roadold] (.4,2.6)..controls (1.8,2.25) and (2.8,1.9)..(3.8,1.7);
 \draw[roadcurve] (3.8,1.7)..controls (5.1,1.35) and (5.7,1.15)..(6.1,1.05);
 \draw[roadcorr] (6.1,1.05)..controls (6.4,.85) and (6.65,.8)..(6.9,.94);
 \draw[roadcurve] (6.9,.94)..controls (7.6,.78) and (8.1,.65)..(8.5,.58);
 \draw[roadcorr] (8.5,.58)..controls (8.8,.78) and (9.05,.75)..(9.3,.46);
 \draw[roadcurve] (9.3,.46)--(11.3,.25);
 \foreach \x/\lab in {.4/R_{0z},3.8/R_{m\theta},11.3/R_{3/2}}
 {\draw[roadguide] (\x,0)--(\x,2.7);\node[below] at (\x,0) {$\lab$};}
 \node[above] at (6.5,1.3) {$d_1\beta$};
 \node[above] at (8.95,1.05) {$d_2\beta$};
 \node[align=center] at (8,2.5)
 {background $R^{-1/2-\mu}$\\[3pt]$U^r=U^z=0$ after matching};
 \node[align=center,font=\scriptsize,text width=3.2cm] at (2,.95)
 {$U^\theta$: $Z$-flattening\\$U^r=U^z=0$ after matching};

 \draw[roadguide] (6.5,0)--(6.5,.9);
 \draw[roadguide] (8.95,0)--(8.95,.6);
 \node[below,font=\scriptsize] at (6.5,-.55) {$e^{-3}R_{3/2}$};
 \node[below,font=\scriptsize] at (8.95,-1.05) {$e^{-1}R_{3/2}$};
\end{tikzpicture}
\caption{Two angular bumps match two moments (Interval O.6).
The centers are $e^{-3}R_{3/2}$ and $e^{-1}R_{3/2}$, as in the
coefficient equations below. Their weights give independent control
of $M^\theta$ and $M^p$. They match the right angular target and enforce
$P(\infty,Z)=0$ with $P(0,Z)=P_0^{\rm pre}(Z)$. The bumps are shown with illustrative sizes and
signs, without preserving the radial monotonicity of the full profile $U^\theta$.
The coefficients satisfy $|d_1|+|d_2|\lesssim\mu^{29}$.}
\label{fig:merged-angular}
\end{figure}
}

\newcommand{\MergedAxialFigure}{%
\begin{figure}[htbp]
\centering
\begin{tikzpicture}[x=1cm,y=.75cm,>=stealth,font=\small]
 \draw[->] (0,0)--(12,0) node[right] {$R$};
 \draw[->] (0,0)--(0,2.8) node[above] {$U^\theta$};
 \draw[roadold] (.3,2.45)..controls (1.1,2.1) and (1.8,1.85)..(2.7,1.65);
 \draw[roadcurve] (2.7,1.65)..controls (6,1.1) and (9,.45)..(11.3,.3);
 \foreach \x/\lab in {.3/R_{+c},2.7/R_{mz},11.3/R_{0z}}
 {\draw[roadguide] (\x,0)--(\x,2.5);\node[below] at (\x,0) {$\lab$};}
 \node at (7.1,2.1) {$U^\theta=\overline{U}^\theta(R_{mz},Z)(R/R_{mz})^{-1/2-\mu}$};
 \node[align=center,font=\scriptsize,text width=2.6cm] at (1.35,.85)
 {$U^\theta\propto R^{-1/2-\mu}$\\$U^z=0$\\$U^r\not\equiv0$};
 \begin{scope}[yshift=-3.1cm]
 \draw[->] (0,0)--(12,0) node[right] {$R$};
 \node[left] at (0,1.2) {$U^z$};
 \draw[roadold] (.3,0)--(2.7,0);
 \draw[roadcorr] (2.7,0)..controls (3.5,0) and (4.2,1.45)..(5.2,1.55)
 ..controls (6.4,1.45) and (7.1,0)..(7.7,0)--(8.4,0)
 ..controls (8.5,-.45) and (8.9,-.45)..(9,0)--(9.6,0)
 ..controls (9.7,.4) and (10.1,.4)..(10.2,0)--(11.3,0);
 \node[above,align=center] at (5.2,1.55)
 {main pulse\\$M^{z\theta}(\infty,Z)=0$};
 \node[align=center] at (9.4,1.1)
 {two end bumps\\cancel $M^z,M^{\theta z}$};
 \node[align=center] at (6,-1.25)
 {After exact matching,
 $U^r(R_{0z})=U^z(R_{0z})=0$.};
 \end{scope}
\end{tikzpicture}
\caption{One main axial pulse and two late end corrections (Interval O.4).
The pulse supplies the positive axial energy in the $M^{z\theta}$ balance;
the late bumps cancel $M^z$ and $M^{\theta z}$.
Their amplitudes are chosen only after the complete angular profile,
including Interval O.6 and the heat tail, has been fixed.
The small end bumps are enlarged, and their signs are illustrative.}
\label{fig:merged-axial}
\end{figure}
}

\section{The outer construction: pressure, moment matching, and the cone}
\label{sec:outer-construction}\label{sec:leading-order-construction}\label{sec:outer-profile}
We construct an \emph{outer connecting profile}
$(U^\theta,U^z)$ on $[R_{\rm ref},\infty)\times[-1,1]$,
joining the reference velocity at $R_{\rm ref}$ to the exterior heat flow.
This section proves the outer construction summarized in
Proposition~\ref{prop:road-outer-summary}.
We first specify the candidate family and its axis-pressure target.
The exact heat exterior supplies data from the right, while the reference
profile supplies cumulative moments from the left. We treat Intervals
O.8 and O.7 first, then O.2--O.6. Matching the two sets of data at
$R_{3/2}$ determines the waiting length and the correction coefficients.
We then verify the cone on the completed outer profile.

\subsection{The outer connecting profile ansatz}\label{subsec:outer-ansatz}

For the cumulative moment calculations, define the \emph{reference velocity profile} on
$[0,R_{\rm ref}]\times[-1,1]$ by
\begin{equation}\label{reference-profile}
U^\theta_{\rm ref}(R,Z)=\frac{P_*}{1+Z^2}
\left(\frac{R}{R_{\rm ref}}\right)^{1/10},
\qquad U^z_{\rm ref}(R,Z)=4Z.
\end{equation}
The corresponding moments $\mathbf M_{\rm ref}(R,Z)$ are well-defined
by \eqref{fiveM} and vanish at $R=0$. The reference profile alone
does not determine a pressure by integration from infinity;
extending its displayed power law to infinity would give a divergent
pressure integral.

For each provisional waiting length $\tau\ge0$, let
$\bar U=(\bar U^\theta,\bar U^z)$ denote the uncorrected candidate.
On $[0,R_{\rm ref}]$ it equals the reference profile \eqref{reference-profile}.
For $R\ge R_{\rm ref}$ the ansatz below continues it smoothly to the
exact heat exterior \eqref{exterior-profile} on $R\ge R_h$, with amplitude
$c_\infty$. We write $U$ for the profile obtained after the angular and
axial moment corrections. In calculations involving only the uncorrected
candidate we sometimes omit the bar.

The candidate moments $\mathbf M=(M^\theta,M^z,M^{\theta z},M^{z\theta},M^p)$
are computed from $\bar U$ with its reference continuation down to the
axis. In the outward propagation we use $P=P_0^{\rm pre}+M^p$ to
compute the axial inertial stress, with $P_0^{\rm pre}$ defined in
\eqref{eq:outer-preheat-pressure}.
These cumulative moments are propagated across every interface, without
resetting them or imposing new terminal conditions on each interval.
Their estimates below are used to verify the cone; the required global
identities are imposed by the angular and axial corrections below.

Fix a smooth increasing step function $\sigma$, flat at both $y=0$ and $y=1$:
\[
\sigma(y)=\frac{e^{-1/y^2}}{e^{-1/y^2}+e^{-1/(1-y)^2}}
\quad(0<y<1),\qquad
\sigma=0\ (y\le0),\quad \sigma=1\ (y\ge1).
\]
Then $\|\partial_y\sigma\|_{L^\infty}=8$. This fixed function is used
throughout the outer connecting ansatz.

\paragraph{Parameters.}
With the step $\sigma$ and the bump for the moment corrections fixed,
take $M_d$ sufficiently large for the axial turnoff and set $T_d=e^{M_d}+10$.
Fix small positive constants $c_\mu,c_\varepsilon$ for the power-law
buffer and the heat transition, respectively. The proof then fixes
$d_{\rm corr}$ for the pulse direction, $\mu_{\rm corr}$ for the error
estimates, and finally $R_{\rm corr}$ as a common lower bound on
$R_{\rm ref}$ for the cone estimates on both the uncorrected and
corrected intervals.
These choices are independent of the inputs. With $P_*$ chosen first,
the outer construction applies throughout the range
\begin{equation}\label{eq:mc-free-parameters}
 P_*>e^{T_d},\qquad \mu=c_\mu P_*^{-4}\le\mu_{\rm corr},\qquad
 0<\delta\le d_{\rm corr}\mu,\qquad R_{\rm ref}\ge R_{\rm corr}.
\end{equation}
The logarithmic lengths and the relative deficit below the terminal power law are
\begin{equation}\label{eq:outer-derived-parameters}
 T_w=60\log(1/\mu),\qquad T_f=100,\qquad
 T_s=4\log(2/\delta),\qquad \varepsilon=c_\varepsilon\delta.
\end{equation}
The candidate family allows any waiting length $\tau\ge0$.
The matching equation \eqref{eq:mc-wait-choice} selects its unique value
for the completed profile. In the uniform estimates below, implicit
constants depend only on the previously fixed auxiliary data and,
for derivative estimates, on the order of differentiation.
The $C^1_Z$ norms are taken on $[-1,1]$.

\paragraph{Checkpoint radii.}
For the reference-region estimates, set $R_{\rm useless}=e^{-5}R_{\rm ref}$. Define the outer checkpoint radii by
\[
\begin{aligned}
R_{-\mu}&=e^{1+T_d}R_{\rm ref},
&
R_{+c}&=eR_{-\mu},
\\
R_{mz}&=e^{T_w}R_{+c},
&
R_{0z}&=e^{13/\mu}R_{mz},
\\
R_{m\theta}&=e^{T_f}R_{0z},
&
R_{3/2}&=\mu^{-30}R_{m\theta},
\\
R_{2h}&=e^{T_s+2+\tau}R_{3/2},
&
R_h&=e^3R_{2h}.
\end{aligned}
\]
The radii through $R_{3/2}$ are fixed multiples of $R_{\rm ref}$ once
$P_*$ is chosen. The two remaining ratios $R_{2h}/R_{\rm ref}$ and
$R_h/R_{\rm ref}$ also depend on $\delta$ and the waiting length $\tau$.

\paragraph{The ansatz for $\bar U^\theta$.}
Use $y=\log(R/R_{\rm ref})$. Write $y_d$, $y_v$, $y_{\rm rel}$ and
$y_{\rm tail}$ for the values of $y$ at $R_{-\mu}$, $R_{0z}$,
$R_{3/2}$ and $R_{2h}$, respectively. Set
\[
A(y)=P_*\exp\left(\int_0^y s(v)\,dv\right),
\]
where
\[
\begin{aligned}
s(y)
={}&\frac1{10}-\frac35\sigma(y)
-\mu\sigma(y-y_d)
\\
&-(1-\mu)\sigma(y-y_{\rm rel})
+\left(1-\frac{\delta}{2}\right)
 \sigma(y-y_{\rm rel}-1-T_s).
\end{aligned}
\]
The four steps change the logarithmic slope of $A$ successively from
$1/10$ to $-1/2$, then to $-1/2-\mu$, $-3/2$, and finally
$-(1+\delta)/2$. For $R \in [R_{\rm ref}, R_{2h}]$, define
\[
\bar U^\theta(R,Z)=
\frac{A(y)}{1+Z^2}
\left(\frac{1+Z^2}{2}\right)^{
\sigma((y-y_v)/T_f)}.
\]
For $R \in [R_{2h}, +\infty)$, put
\[
c_\infty=
\frac{A(y_{\rm tail})R_{2h}^{(1+\delta)/2}}
{2(1-\varepsilon)}.
\]
This amplitude is independent of $\tau$ since $A(y)e^{(1+\delta)y/2}$
is constant for $y\ge y_{\rm rel}+T_s+2$. Set
\begin{equation}\label{eq:connecting-swirl}
\begin{aligned}
\bar U^\theta(R,Z)=
{}&c_\infty R^{-(1+\delta)/2}
\Bigl[
(1-\sigma(y-y_{\rm tail}))(1-\varepsilon)
\\
&\qquad+\sigma(y-y_{\rm tail})
H_\delta\!\left(\frac{2L}{R}\right)
\bigl(1-\varepsilon \mathfrak f((3-y+y_{\rm tail})/2)\bigr)
\Bigr].
\end{aligned}
\end{equation}
For $y\ge y_{\rm tail}+1$, this ansatz gives the heat profile of
Section~\ref{sec:heat-exterior} multiplied by the flat collar factor
studied in Section~\ref{subsec:outer-waiting}.

\paragraph{The axial profile ansatz.}
Define
\[
B(t)=
\begin{cases}
1,&t\le0,\\[1mm]
\displaystyle
1-\sigma\!\left(\frac{\log(1+t)}{M_d}\right),
&t>0.
\end{cases}
\]
With the reference continuation included, for $R>0$ we have
\begin{equation}\label{eq:connecting-axial}
\bar U^z(R,Z)=4ZB(y-1).
\end{equation}
At the axis, $\bar U^z(0,Z)=4Z$.

\paragraph{Moment corrections.}
The corrections have support above $R_{mz}$, so the reference velocities
and cumulative velocity moments at $R_{\rm ref}$ remain unchanged.

Write $\mathbf M_{\mathrm{left}}$ for cumulative moments obtained from the axis,
and $\mathbf M_{\mathrm{right}}$ for the terminal data propagated inward from
the heat exterior. These subscripts distinguish the direction of
calculation; the left data here come from the temporary reference
profile, not from the analytic core constructed later.
When no comparison is being made, we retain the notation $\mathbf M$ of
\eqref{fiveM}. The right-side definitions are given in
\eqref{eq:mc-required-moments}.

There are two kinds of corrections. Two small angular bumps adjust the
angular moment and enforce $P(\infty,Z)=0$ for the prescribed axis
pressure $P_0^{\rm pre}$. An axial
pulse, followed by two much smaller axial bumps, corrects the other three
moments. The pulse has a bounded relative amplitude; it is not a small
relative perturbation of the zero axial velocity. Its absolute amplitude
is small because the swirl has already decayed before $R_{mz}$.
The three intervals reserved in Interval O.3 remain available for the
later construction; the corrections in this section do not use them.
The axial pulse is a smooth component of the stationary profile.

\paragraph{The desired output.}
Use the pre-heat profile and pressure defined below in
\eqref{eq:outer-preheat-pressure}. Set
\[
M_{\rm pre}^p(\infty,Z)=-P_0^{\rm pre}(Z).
\]
The angular correction restores this pressure moment while retaining
the actual heat tail.
We impose the exact identities
\begin{equation}\label{eq:mc-targets}
\begin{gathered}
M^z(\infty,Z)=M^{\theta z}(\infty,Z)
=M^{z\theta}(\infty,Z)=0,\\
\int_0^\infty\sqrt{2R}
\left[U^\theta(R,Z)-c_\infty R^{-(1+\delta)/2}\right]\,dR=0,\\
\int_0^\infty
\frac{(U^\theta)^2-(\bar U^\theta_{\rm pre})^2}{2R}\,dR=0.
\end{gathered}
\end{equation}
The second line is the renormalized angular moment condition.
Its integral converges at the axis since $\delta<2$, and at infinity
since $H_\delta(2L/R)-1=O(\delta/R)$ and $\delta>0$.
The last line realizes the analytic pre-heat axis pressure,
evaluated at the waiting length selected below.

The order of determination is
\[
\tau\longrightarrow(d_1,d_2)\longrightarrow(a_p,c_1,c_2).
\]
The waiting length on Interval O.7 matches the pre-heat angular moment at
$Z=0$. The two angular bumps on Interval O.6 then match the angular and
pressure moments for every $Z$. Although the axial pulse lies earlier in
radius, its energy balance uses the entire future swirl, including these
bumps and the heat tail. Section~\ref{subsec:outer-axial-matching}
therefore gives the axial construction for a prescribed swirl, and
Section~\ref{subsec:outer-angular-matching} applies it after the two angular corrections.

\subsection{The pressure estimates and the axis datum}\label{subsec:outer-axis-data}

The axis datum $P_0(Z)$ is the value $P(0,Z)$ of the completed leading
pressure profile, with the normalization $P(\infty,Z)=0$ from
\eqref{P-Utheta}. It is determined by the full angular profile.
The scalar $P_*$ is the amplitude parameter of the reference swirl;
it is not the value of the axis pressure.

In the following lemma, let $P$ denote the tail integral
\begin{equation}\label{eq:outer-pressure}
P(R,Z)=-\int_R^\infty
\frac{(U^\theta(\rho,Z))^2}{2\rho}\,d\rho
\end{equation}
for either $U^\theta=\bar U^\theta$ or its pre-heat profile
$U^\theta=\bar U^\theta_{\rm pre}$, obtained by setting $H_\delta=1$
in \eqref{eq:connecting-swirl}. The former gives the right-side
pressure for $R\ge R_{3/2}$. The latter gives $P_{\rm pre}$ and
its axis value $P_0^{\rm pre}$, from which the left pressure is recovered below.
The estimates hold throughout the range \eqref{eq:mc-free-parameters}.

\begin{lemma}[Pressure bounds]\label{lem:outer-pressure}
For every $R\ge0$ and $Z\in[-1,1]$,
\[
|P(R,Z)|\lesssim\frac{P_*^2}{(1+Z^2)^2},
\qquad
|\partial_Z P(R,Z)|\lesssim
\frac{P_*^2|Z|}{(1+Z^2)^3}.
\]
For $R\ge eR_{\rm ref}$, the sharper estimates hold:
\begin{equation}\label{eq:outer-pressure-sharp}
|P(R,Z)|\le\frac12(U^\theta(R,Z))^2,
\qquad
|\partial_Z P(R,Z)|\le
\frac{2|Z|}{1+Z^2}(U^\theta(R,Z))^2.
\end{equation}
\end{lemma}

\begin{proof}
We first verify the pointwise bounds
\begin{equation}\label{eq:outer-swirl-pressure-bounds}
\partial_y\log U^\theta\le-\frac12
\quad(R\ge eR_{\rm ref}),
\qquad
|\partial_Z\log U^\theta|
\le\frac{2|Z|}{1+Z^2}
\quad(R>0).
\end{equation}
Before the terminal transition, these follow directly from the
ansatz: the factor removing the $Z$-dependence is nonincreasing
in $y$, and
\[
\partial_Z\log U^\theta
=-\frac{2Z}{1+Z^2}
\left[1-\sigma\left(\frac{y-y_v}{T_f}\right)\right].
\]

On the terminal transition, write
\[
U^\theta=c_\infty R^{-(1+\delta)/2}K(x,Z),
\qquad x=\log(R/R_{2h}),
\]
where $K$ is the bracket in \eqref{eq:connecting-swirl}, with
$H_\delta$ replaced by $1$ for the pre-heat profile.
The definition of $H_\delta$ gives
\[
0<H_\delta(\xi)\le1,\qquad
0\le1-H_\delta(\xi)\lesssim\delta\xi,\qquad
|\partial_\xi H_\delta(\xi)|\lesssim\delta
\quad(0\le\xi\le1).
\]
Consequently, uniformly for $R\ge R_{2h}$,
\[
K\gtrsim1,\qquad
\partial_y\log K\lesssim\varepsilon+\frac{\delta}{R},
\qquad
|\partial_Z\log K|\lesssim\frac{\delta|Z|}{R}.
\]
Since $\varepsilon=c_\varepsilon\delta$, we take the fixed
$c_\varepsilon\le10^{-3}$ sufficiently small and $R_{\rm corr}$ sufficiently
large that $\partial_y\log K\le\delta/2$ and
$|\partial_Z\log K|\le2|Z|/(1+Z^2)$. This proves
\eqref{eq:outer-swirl-pressure-bounds} also on this region.

The radial slope bound in \eqref{eq:outer-swirl-pressure-bounds} implies
\[
U^\theta(\rho,Z)
\le U^\theta(R,Z)(\rho/R)^{-1/2}
\qquad(\rho\ge R\ge eR_{\rm ref}).
\]
Substitution into \eqref{eq:outer-pressure} yields
\[
-P(R,Z)\le\frac12(U^\theta(R,Z))^2.
\]
The $Z$-derivative bound gives, for every $R\ge0$,
\begin{equation}\label{eq:outer-pressure-derivative}
|\partial_Z P(R,Z)|
\le\frac{4|Z|}{1+Z^2}\,[-P(R,Z)].
\end{equation}
Together these prove \eqref{eq:outer-pressure-sharp}.

Finally, the common reference branch on $[0,R_{\rm ref}]$ contributes
\[
\int_0^{R_{\rm ref}}
\frac{(U^\theta_{\rm ref})^2}{2R}\,dR
=\frac{5P_*^2}{2(1+Z^2)^2}.
\]
On $0\le y\le1$, we have
$U^\theta\le P_*e^{y/10}/(1+Z^2)$, while the remaining
contribution from $[e R_{\rm ref}, +\infty)$ is controlled by \eqref{eq:outer-pressure-sharp}.
Thus
\[
-P(0,Z)\lesssim \frac{5 P_*^2}{(1+Z^2)^2}.
\]
The global estimates follow from the monotonicity of $P$ in $R$ and
\eqref{eq:outer-pressure-derivative}.
\end{proof}

\paragraph{The pre-heat pressure target.}
For each provisional waiting length $\tau$, the pre-heat profile retains
the reference branch and the flat deficit multiplier. Define
\begin{equation}\label{eq:outer-preheat-pressure}
\begin{aligned}
P_{\rm pre}(R,Z;\tau)&=-\int_R^\infty
 \frac{(\bar U^\theta_{\rm pre}(\rho,Z;\tau))^2}{2\rho}\,d\rho,\\
P_0^{\rm pre}(Z;\tau)&=P_{\rm pre}(0,Z;\tau).
\end{aligned}
\end{equation}
For the uncorrected candidate $\bar U$, we propagate from the axis with
\[
 P_{\mathrm{left}}(R,Z)=P_0^{\rm pre}(Z)+\bar M_{\mathrm{left}}^p(R,Z).
\]
Here and below we suppress $\tau$ when it is fixed. Since
$\bar U^\theta=\bar U^\theta_{\rm pre}$ on $[0,R_{2h}]$, we have
$P_{\mathrm{left}}=P_{\rm pre}$ throughout this interval.
Extending the same cumulative moment to infinity gives
$P_{\mathrm{left}}(\infty,Z)=-\Delta_H(Z)$, where
\begin{equation}\label{eq:outer-heat-pressure-defect}
\Delta_H(Z)=\int_{R_{2h}}^\infty
 \frac{(\bar U^\theta_{\rm pre}(\rho,Z))^2-(\bar U^\theta(\rho,Z))^2}
 {2\rho}\,d\rho.
\end{equation}
Thus $P_{\mathrm{right}}-P_{\mathrm{left}}=\Delta_H$ on $R\ge R_{3/2}$.
Once the angular matching in Section~\ref{subsec:outer-waiting} fixes $\tau$,
the angular correction increases the total pressure moment by $\Delta_H$.
The corrected profile then satisfies $P(\infty,Z)=0$ with the prescribed axis value
\begin{equation}\label{eq:outer-final-axis-pressure}
P_0(Z):=P(0,Z)=-M^p(\infty,Z)=P_0^{\rm pre}(Z).
\end{equation}
This identity is proved in \eqref{eq:mc-pressure-preserved}. The axial pulse leaves the
angular profile unchanged and hence preserves this pressure.

\begin{lemma}[The axis pressure]\label{lem:outer-axis-pressure}
For every input satisfying \eqref{eq:mc-free-parameters} and every
provisional waiting length $\tau\ge0$, the axis datum $P_0^{\rm pre}$
is even and satisfies, for $Z\in[-1,1]$,
\[
\begin{aligned}
\frac{5P_*^2}{2(1+Z^2)^2}
&\le-P_0^{\rm pre}(Z)\le\frac{KP_*^2}{(1+Z^2)^2},\\
\frac{10P_*^2Z^2}{(1+Z^2)^3}
&\le Z\partial_Z P_0^{\rm pre}(Z)\le\frac{KP_*^2Z^2}{(1+Z^2)^3}.
\end{aligned}
\]
It extends holomorphically to the fixed rectangle
\[
 \Omega_0=\{z\in\mathbb C:|\operatorname{Re}z|<5/4,\quad
 |\operatorname{Im}z|<1/4\},
\]
with
\[
 \sup_{\Omega_0}|P_0^{\rm pre}|\le KP_*^2,\qquad
 \sup_{[-1,1]}|\partial_Z^kP_0^{\rm pre}|\le K\,8^k k!\,P_*^2
 \quad(k\ge0).
\]
Here $K>1$ depends only on the fixed auxiliary data and is independent
of $P_*,\delta,R_{\rm ref},\tau$ and $k$.
With $P_*,\delta$ and the dimensionless schedule, including $\tau$,
fixed, $P_0^{\rm pre}$ is independent of $R_{\rm ref}$.
\end{lemma}

\begin{proof}
Evenness follows from the ansatz. The lower bound for the
pressure magnitude follows from the reference branch, and the
upper bounds follow from Lemma~\ref{lem:outer-pressure} applied to
$\bar U^\theta_{\rm pre}$.
The pre-heat swirl has the form
\[
\bar U^\theta_{\rm pre}(R,Z)
=\mathcal A(y)(1+Z^2)^{-\vartheta(y)},\qquad
0\le\vartheta(y)\le1.
\]
In the tail, the choice of amplitude in Subsection~\ref{subsec:outer-ansatz}
gives
\[
 c_\infty R_{\rm ref}^{-(1+\delta)/2}
 =\frac{A(y_{\rm tail})}{2(1-\varepsilon)}
   \left(\frac{R_{2h}}{R_{\rm ref}}\right)^{(1+\delta)/2}.
\]
For fixed $P_*,\delta$ and dimensionless schedule, the right-hand side
is independent of $R_{\rm ref}$. Thus $\mathcal A$ and $\vartheta$
depend only on $y$ throughout the pre-heat profile.
The change of variables $R=R_{\rm ref}e^y$ gives
\[
 P_0^{\rm pre}(z)=-\frac12\int_{\mathbb R}
 \mathcal A(y)^2(1+z^2)^{-2\vartheta(y)}\,dy,
\]
which shows the asserted independence of $R_{\rm ref}$.
On $\Omega_0$, $\operatorname{Re}(1+z^2)>15/16$, so the principal
logarithm gives $|(1+z^2)^{-2\vartheta}|\le(16/15)^2$.
Since $\mathcal A^2$ is integrable, uniform integration of these
holomorphic functions proves the analyticity of $P_0^{\rm pre}$ and gives
\[
\sup_{\Omega_0}|P_0^{\rm pre}|
\le(16/15)^2[-P_0^{\rm pre}(0)]\le KP_*^2.
\]
Cauchy's estimate on disks of radius $1/8$ centered on $[-1,1]$ gives
\[
\sup_{[-1,1]}|\partial_Z^kP_0^{\rm pre}|
\le K\,8^k k!\,P_*^2.
\]
The constant $K$ may be enlarged to cover all the displayed bounds;
it is independent of $P_*,\delta,R_{\rm ref},\tau$ and $k$.

Every radial interval contributes nonnegatively to
$Z\partial_Z P_0^{\rm pre}(Z)$. The contribution on
$[0,eR_{\rm ref}]$ gives the stronger lower bound
\begin{equation}\label{eq:outer-axis-pressure-margin}
Z\partial_Z P_0^{\rm pre}(Z)
\ge\frac{[10+2(1-e^{-1})]P_*^2Z^2}{(1+Z^2)^3}.
\end{equation}
\end{proof}

We also record the size of the pressure defect \eqref{eq:outer-heat-pressure-defect}.
On $[R_{2h},\infty)$ the pre-heat swirl is independent
of $Z$, $\bar U^\theta\le\bar U^\theta_{\rm pre}$, and
$Z\partial_Z\bar U^\theta\ge0$. Thus $\Delta_H\ge0$ and
$Z\partial_Z\Delta_H\le0$.
The estimates
\[
1-H_\delta\left(\frac{2L}R\right)
\lesssim\frac{\delta(1-Z^2)}R,
\qquad
\left|\partial_ZH_\delta\left(\frac{2L}R\right)\right|
\lesssim\frac{\delta|Z|}R
\]
yield, by integration,
\begin{equation}\label{eq:outer-axis-pressure-heat-error}
\begin{aligned}
0\le\Delta_H(Z)
&\lesssim\frac{\delta(1-Z^2)U_{\rm tail}^2}{R_{2h}},\\
0\le-Z\partial_Z\Delta_H(Z)
&\lesssim\frac{\delta Z^2U_{\rm tail}^2}{R_{2h}}.
\end{aligned}
\end{equation}
Here $U_{\rm tail}=\bar U^\theta(R_{2h},Z)$ is
independent of $Z$. These estimates are uniform in $\tau\ge0$.

\paragraph{Compatibility with the axis pressure.}
If a left profile is recovered from an axis value $P_0(Z)$ by
$P(R,Z)=P_0(Z)+\int_0^R F(\rho,Z)^2\,d\rho$, matching its pressure
to this exterior requires
\[
\int_0^{R_h}F(R,Z)^2\,dR
=-P_0(Z)-\int_{R_h}^\infty
\frac{(U^\theta_{\rm heat}(R,Z))^2}{2R}\,dR.
\]
The left-hand integral includes the fixed collar contribution.
This pressure identity is one of the matching requirements; it
does not replace the other moment identities needed for
$U^r=0$ and \eqref{eq:outer-stress-boundary-data} in the exterior.
In the two-sided construction, the common axis datum in this identity
is the selected $P_0^{\rm pre}$. Equality of the pressure moments at
$R_{3/2}$ will enforce the same pressure normalization throughout the
unchanged right-hand profile.

\subsection{The exact heat exterior and right-side moments}
\label{sec:heat-exterior}

We construct the exact heat exterior in \eqref{eq:road-heat-exact}
and determine its right-side moment data. The inward collar of
\eqref{eq:road-heat-join} is treated with the inward propagation in
Section~\ref{subsec:outer-waiting}.
In the global construction, the amplitude $c_\infty$ and interface
radius $R_h$ are those derived in Section~\ref{subsec:outer-ansatz}.

Fix $0<\delta\le1/2$ and $c_\infty>0$. In the exterior we use
a pure swirl that is independent of the axial variable:
\begin{equation}\label{exterior-ansatz}
u^r=u^z=0,\qquad u^\theta=u^\theta(t,r),\qquad
p(t,r)=-\int_r^\infty
\frac{(u^\theta(t,\tilde r))^2}{\tilde r}\,d\tilde r.
\end{equation}
The pressure gives exact radial balance.
The axial equation and incompressibility hold automatically.
The angular equation reduces to
\begin{equation}\label{K-heat}
\partial_tu^\theta
=\left(\partial_r^2+\frac1r\partial_r-\frac1{r^2}\right)u^\theta.
\end{equation}
The axial viscosity terms vanish as well, since the physical
velocity and pressure are independent of $z$.

\paragraph{The self-similar heat factor.}
We seek an exterior solution with the homogeneous terminal trace
\begin{equation}\label{eq:heat-terminal-trace}
u^\theta(1-0,r)
=c_\infty\left(\frac{r^2}{2}\right)^{-(1+\delta)/2},
\qquad r>0.
\end{equation}
The scaling compatible with this trace is
\[
u^\theta(1-\lambda^2(1-t),\lambda r)
=\lambda^{-1-\delta}u^\theta(t,r).
\]
Accordingly, put
\begin{equation}\label{eq:heat-physical-profile}
u^\theta(t,r)
=c_\infty\left(\frac{r^2}{2}\right)^{-(1+\delta)/2}
H_\delta\left(\frac{4(1-t)}{r^2}\right).
\end{equation}
Writing $\xi=4(1-t)/r^2$, direct substitution into
\eqref{K-heat} gives
\begin{equation}\label{eq:heat-factor-ode}
\xi^2\partial_\xi^2 H_\delta
+\bigl(1+(2+\delta)\xi\bigr)\partial_\xi H_\delta
+\frac{\delta}{2}\left(1+\frac{\delta}{2}\right)H_\delta=0.
\end{equation}
The terminal normalization is $H_\delta(0)=1$.

The solution used here is
\begin{equation}\label{eq:heat-factor-integral}
H_\delta(\xi)
=\frac1{\Gamma(1+\delta/2)}
\int_0^\infty e^{-v}v^{\delta/2}
(1+\xi v)^{-\delta/2}\,dv,\qquad \xi\ge0.
\end{equation}
For $\xi>0$, this is the classical Kummer integral representation
\cite[\S13.4, equation (13.4.4)]{DLMF}, after rescaling its
integration variable. We verify the equation directly below.
For a direct verification, temporarily put $a=\delta/2$.
Differentiation under the integral yields
\[
\partial_\xi^n H_\delta(\xi)
=\frac{(-1)^n(a)_n}{\Gamma(1+a)}
\int_0^\infty e^{-v}v^{a+n}(1+\xi v)^{-a-n}\,dv,
\qquad n\ge0,
\]
where $(a)_0=1$ and $(a)_n=a(a+1)\cdots(a+n-1)$ for $n\ge1$.
The integrable majorants are uniform for $\xi$ in every compact
subset of $[0,\infty)$, so $H_\delta$ is smooth up to $\xi=0$.
Substitution of its first two derivatives into the left-hand side
of \eqref{eq:heat-factor-ode} gives
\[
\frac a{\Gamma(1+a)}
\int_0^\infty
\partial_v\!\left[e^{-v}v^{a+1}(1+\xi v)^{-a-1}\right]\,dv=0.
\]
Both boundary terms vanish. The Gamma integral also gives
$H_\delta(0)=1$, completing the verification of the heat solution.

The same representation shows that $H_\delta$ is positive,
strictly decreasing and strictly convex. In particular,
\begin{equation}\label{eq:heat-factor-basic-bounds}
0<H_\delta(\xi)\le1,\qquad
|\partial_\xi H_\delta(\xi)|\le a(1+a),\qquad
0\le1-H_\delta(\xi)\le a(1+a)\xi.
\end{equation}
Thus $H_\delta=1+O_\delta(\xi)$ in the small-argument regime
used in the exterior construction. Figure~\ref{fig:heat-factor}
illustrates this factor.

Smoothness at $\xi=0$ must be distinguished from analyticity.
Indeed,
\begin{equation}\label{eq:heat-factor-origin-derivatives}
(\partial_\xi^n H_\delta)(0)=(-1)^n(a)_n(1+a)_n.
\end{equation}
For $a>0$, the ratio of successive absolute Taylor coefficients
is $(a+n)(1+a+n)/(n+1)$, which tends to infinity.
Consequently the Taylor series at zero has radius of convergence
zero. The heat profile and the inward collar below require only
smoothness, not analyticity at this endpoint.
\begin{figure}[htbp]
\centering
\begin{tikzpicture}
\begin{axis}[
    width=9cm, height=5.5cm,
    xlabel={$\xi$},
    ylabel={$H_\delta(\xi)$},
    xmin=0, xmax=10,
    ymin=0.989, ymax=1.0005,
    xtick={0,2,4,6,8,10},
    ytick={0.990,0.992,0.994,0.996,0.998,1.000},
    yticklabel style={
        /pgf/number format/fixed,
        /pgf/number format/precision=3,
        /pgf/number format/zerofill
    },
    scaled y ticks=false,
    axis lines=left,
    tick align=outside,
]
\addplot[blue!70!black, thick, smooth, no marks]
table {
0.0  1.000000000
0.5  0.998188187
1.0  0.997012970
1.5  0.996116053
2.0  0.995383388
2.5  0.994760990
3.0  0.994218416
3.5  0.993736607
4.0  0.993302760
4.5  0.992907823
5.0  0.992545138
5.5  0.992209653
6.0  0.991897443
6.5  0.991605389
7.0  0.991330969
7.5  0.991072116
8.0  0.990827113
8.5  0.990594514
9.0  0.990373094
9.5  0.990161802
10.0 0.989959731
};
\end{axis}
\end{tikzpicture}
\caption{The heat-flow factor $H_\delta$ for $\delta=0.01$.
It is positive, decreasing, and convex, with $H_\delta(0)=1$.}
\label{fig:heat-factor}
\end{figure}

\paragraph{The exterior profile and its shear.}
In the similarity coordinates, \eqref{eq:heat-physical-profile}
becomes
\begin{equation}\label{exterior-profile}
U^\theta_{\rm heat}(R,Z)
=c_\infty R^{-(1+\delta)/2}
H_\delta\left(\frac{2L}R\right),
\qquad U^r_{\rm heat}=U^z_{\rm heat}=0.
\end{equation}
Here $2L/R=4(1-t)/r^2$. The profile is smooth on
$R>0$, $-1\le Z\le1$, including the endpoints in $Z$, and
\[
U^\theta_{\rm heat}(R,Z)
=c_\infty R^{-(1+\delta)/2}\bigl(1+O_\delta(R^{-1})\bigr)
\qquad(R\to\infty).
\]
Its pressure is
\[
P_{\rm heat}(R,Z)
=-\int_R^\infty
\frac{(U^\theta_{\rm heat}(\rho,Z))^2}{2\rho}\,d\rho.
\]
The global profile will agree with these fields for $R\ge R_h$,
equivalently $r\ge\sqrt{2R_h}\lambda(t,z)$.
This is agreement of the fields and moments needed to recover
zero stress; velocity agreement alone does not ensure it.

For $\xi=2L/R$,
\[
F_{\rm heat}
=\frac{U^\theta_{\rm heat}}{\sqrt{2R}}
=\frac{c_\infty}{\sqrt2}R^{-1-\delta/2}H_\delta(\xi),
\]
and hence
\begin{equation}\label{eq:heat-shear}
\begin{aligned}
\mathcal{S}^\theta_{\rm heat}
&=2R\partial_RF_{\rm heat}\\
&=-\sqrt2\,c_\infty R^{-1-\delta/2}
\left[\left(1+\frac\delta2\right)H_\delta(\xi)
+\xi\partial_\xi H_\delta(\xi)\right],\\
\mathcal{S}^z_{\rm heat}&=0.
\end{aligned}
\end{equation}
The bracket is strictly positive, because
\[
\left(1+\frac\delta2\right)H_\delta(\xi)+\xi\partial_\xi H_\delta(\xi)
=\frac1{\Gamma(1+\delta/2)}
\int_0^\infty e^{-v}v^{\delta/2}
\frac{1+\delta/2+\xi v}{(1+\xi v)^{1+\delta/2}}\,dv>0.
\]
Therefore $\mathcal{S}^\theta_{\rm heat}<0$, and
\begin{equation}\label{eq:heat-kappa}
\kappa_{\rm heat}
=-\frac{\mathcal{S}^\theta_{\rm heat}}{F_{\rm heat}}
=2+\delta+2\xi\frac{\partial_\xi H_\delta(\xi)}{H_\delta(\xi)}.
\end{equation}
On $e^{-1}R_h\le R\le R_h$, the basic bounds give
$\kappa_{\rm heat}=2+\delta+O_\delta(R_h^{-1})$, uniformly in $Z$.

The exact physical heat flow has zero stress. The corresponding
inviscid stresses are
\[
\mathcal{I}^\theta_{\rm heat}=-\mathcal{S}^\theta_{\rm heat},
\qquad \mathcal{I}^z_{\rm heat}=0,
\qquad N^z_{\rm heat}=0.
\]
They satisfy \eqref{inviscid-stress-radial}; this also follows by
substituting \eqref{eq:heat-factor-ode}.

\paragraph{Compatibility with the terminal moments.}
\label{par:heat-moment-compatibility}
Suppose a global profile satisfies \eqref{eq:moment-conditions} and has
$U^\theta=U^\theta_{\rm heat}$ and $U^z=0$ for $R\ge R_h$.
Then its stress recovered from the axis vanishes on this exterior.
Indeed, the terminal conditions give $M^z=M^{\theta z}=U^r=0$ and
determine the remaining moments from their terminal normalizations.
In particular, the heat-factor bounds give
\[
M^\theta(R,Z)
=\frac{\sqrt2c_\infty}{1-\delta/2}R^{1-\delta/2}
+O_{C_Z^1}(c_\infty R^{-\delta/2}).
\]
The leading parts of $(1-\delta/2)M^\theta$ and
$R\sqrt{2R}U^\theta$ therefore cancel in \eqref{Itheta-moments}.
Similarly,
\[
M^{z\theta}(R,Z)
=\frac12\int_R^\infty (U^\theta_{\rm heat}(\rho,Z))^2\,d\rho
=\frac{c_\infty^2}{2\delta}R^{-\delta}
\left[1+O\!\left(\frac{\delta L}{R}\right)\right].
\]
The factor $1/\delta$ is the tail-energy cost addressed in Interval O.7.
Using the same heat asymptotics in
\eqref{Itheta-moments}--\eqref{Iz-moments} and adding the shear gives
\[
R\mathcal{T}^\theta=O(c_\infty R^{-\delta/2})\longrightarrow0,\qquad
\sqrt R\,\mathcal{T}^z=O(c_\infty^2R^{-\delta})\longrightarrow0
\quad\text{as }R\to\infty.
\]
On the exact heat exterior, the homogeneous radial stress equations give
$\mathcal{T}^\theta=C_\theta(Z)/R$ and
$\mathcal{T}^z=C_z(Z)/\sqrt R$.
The weighted limits force $C_\theta=C_z=0$, proving the assertion.

\paragraph{The right-side target moments.}
These targets are defined for each provisional waiting length $\tau$,
which will be selected in Section~\ref{subsec:outer-waiting}.
In view of the requirements \eqref{eq:mc-targets}, we define the right-side target moments by integration from infinity: for $R \ge R_{3/2}$,
\begin{equation}\label{eq:mc-required-moments}
\begin{aligned}
M_{\mathrm{right}}^\theta(R,Z)
&=\frac{\sqrt2c_\infty}{1-\delta/2}R^{1-\delta/2}
-\int_R^\infty\sqrt{2\rho}
\left[
\bar U^\theta(\rho,Z)-c_\infty\rho^{-(1+\delta)/2}
\right]\,d\rho,\\
M_{\mathrm{right}}^z(R,Z)&=M_{\mathrm{right}}^{\theta z}(R,Z;\tau)=0,\\
M_{\mathrm{right}}^{z\theta}(R,Z)
&=\frac12\int_R^\infty
(\bar U^\theta(\rho,Z))^2\,d\rho,\\
M_{\mathrm{right}}^p(R,Z)
&=-P_0^{\rm pre}(Z;\tau)-\int_R^\infty
\frac{(\bar U^\theta(\rho,Z))^2}{2\rho}\,d\rho.
\end{aligned}
\end{equation}
For $R\ge R_{3/2}$, we write $M_{\mathrm{right}}^\theta$ for the target
$M_{\rm tar}^\theta$ in \eqref{eq:road-angular-target} and $M_{\mathrm{right,pre}}^\theta$
for its pre-heat counterpart $M_{\rm tar,pre}^\theta$ in \eqref{eq:road-wait-choice}.

All moment corrections are supported below $R_{3/2}$. Hence for
$R\ge R_{3/2}$ the corrected velocities are $U^\theta=\bar U^\theta$
and $U^z=\bar U^z=0$. At the selected waiting length, the terminal
conditions therefore give exactly \eqref{eq:mc-required-moments} as
the right-side moments of the completed profile.

The left-side moments $\mathbf M_{\mathrm{left}}$ are the cumulative
integrals \eqref{fiveM} from the axis, including the reference branch.
We use $\bar{\mathbf M}_{\mathrm{left}}$ for the uncorrected candidate.
The right-side moments instead encode the terminal conditions.
In particular, the angular moment is renormalized; its value at infinity
is not zero. The pressure moment uses the common axis-pressure target:
\[
P_{\mathrm{right}}(R,Z)
=P_0^{\rm pre}(Z;\tau)+M_{\mathrm{right}}^p(R,Z)
=-\int_R^\infty\frac{(\bar U^\theta(\rho,Z))^2}{2\rho}\,d\rho.
\]
The full matching problem is
\begin{equation}\label{eq:outer-two-sided-matching}
\mathbf M_{\mathrm{left}}(R_{3/2},Z)
=\mathbf M_{\mathrm{right}}(R_{3/2},Z)
\qquad(-1\le Z\le1).
\end{equation}
This is equality as smooth functions of $Z$, so their $Z$ derivatives
match as well. With the common pressure target and the same velocity
formulas near the matching radius, the reconstruction identities give
matching pressure, radial velocity, and stress. No additional independent
stress-matching conditions are imposed.

Thus matching the five moments at $R_{3/2}$ enforces the global
requirements \eqref{eq:mc-targets}.

The heat-factor estimate
\[
\left\|H_\delta\left(\frac{2L}R\right)-1
\right\|_{C_Z^k}\le C_k\frac{\delta}{R}
\]
gives, for $R\ge R_h$, uniformly for $\tau\ge0$,
\begin{equation}\label{eq:mc-required-exterior-bounds}
\begin{aligned}
M_{\mathrm{right}}^\theta
&=\frac{\sqrt2c_\infty}{1-\delta/2}R^{1-\delta/2}
+O_{C_Z^k}(c_\infty R^{-\delta/2}),\\
M_{\mathrm{right}}^{z\theta}
&=\frac{c_\infty^2}{2\delta}R^{-\delta}
+O_{C_Z^k}(c_\infty^2\delta R^{-1-\delta}),\\
M_{\rm pre}^p(\infty,Z)-M_{\mathrm{right}}^p
&=\frac{c_\infty^2}{2(1+\delta)}R^{-1-\delta}
+O_{C_Z^k}(c_\infty^2\delta R^{-2-\delta}).
\end{aligned}
\end{equation}
Here $M_{\rm pre}^p(\infty,Z)=-P_0^{\rm pre}(Z)$, so the last line estimates
the pressure-moment tail of the actual heat profile.

\subsection{Intervals O.7--O.8: The inward collar, steep descent, and waiting}
\label{subsec:outer-right-data}\label{subsec:outer-waiting}

We first record a cone test and radial equations used for propagation
from either side. Starting from the heat exterior, we construct the
inward collar and propagate the right-side data. The pre-heat angular
matching selects the waiting length and gives the positive angular
stress remaining at the inner endpoint of the steep segment.

\paragraph{A test for the cone condition.}
We repeatedly changed the radial growth/decay slope of $U^\theta$ in its defining  ansatz. To keep the estimates short, we
introduce three dimensionless scalar quantities.
They are defined using the outer
profile for $R\ge R_{\rm ref}$; for $0<R<R_{\rm ref}$,
the same definitions, including those of $F$, $\mathcal{S}$, and
$\mathcal{T}$, use $U_{\rm ref}$ instead.
\begin{equation}\label{eq:outer-cone-variables}
a=1-2R\partial_R\log U^\theta,\qquad
b=\frac{2R\partial_RU^z}{U^\theta},\qquad
w=\frac{\mathcal{I}^z}{\mathcal{I}^\theta}.
\end{equation}
The stress, shear and shear intensity can be written as
\[
\mathcal{S}=F(-a,b),\qquad
\mathcal{T}=F
\left(\frac{\mathcal{I}^\theta}{F}-a,
      \frac{\mathcal{I}^\theta}{F}w+b\right),
\qquad \kappa=a+\frac{b^2}{a}.
\]
In particular, $a>0$ is exactly $\mathcal{S}^\theta<0$.

\begin{lemma}\label{lem:outer-cone-test}
Suppose $U^\theta>0$, $\mathcal{I}^\theta>0$, $a>0$, and, uniformly on a fixed compact domain,
\begin{equation}\label{eq:outer-cone-test}
a-bw>0 \ \, (\Leftrightarrow \,  \mathcal{I} \cdot \mathcal{S} < 0),\qquad
2bw+\frac{b^2}{a}+(a-2)w^2<2.
\end{equation}
\begin{itemize}
    \item If $\frac{\mathcal{I}^\theta}{F}$ is sufficiently large (with  $a, b, w$ given and fixed),
then the relaxed cone condition holds.
\item If  $b=0$ and $0<a\le2$, the relaxed cone condition is equivalent to
\begin{equation}\label{eq:outer-cone-axial-zero}
\frac{\mathcal{I}^\theta}{F}>2.
\end{equation}
\end{itemize}
\end{lemma}

\begin{proof}
Substitute the displayed formulas for $\mathcal{S},\mathcal{T}$ into
the relaxed cone inequalities \eqref{cone-relaxed}. First suppose $\kappa\le2$. Since $a^2+b^2=a\kappa$,
\[
-\mathcal{T}\cdot \mathcal{S}
=\frac{(U^\theta)^2}{2R}
\left[
\frac{\sqrt{2R}\,\mathcal{I}^\theta}{U^\theta}(a-bw)
-a\kappa
\right].
\]
Thus the relaxed cone inequality $-\mathcal{T} \cdot \mathcal{S} > \frac{-U^\theta \mathcal{S}^\theta}{\sqrt{2R}} (2-\kappa)$ is equivalent to
\begin{equation}\label{eq:inner-small-shear-test}
\frac{\mathcal{I}^\theta}{F}(a-bw)>2a.
\end{equation}
This holds when $\mathcal{I}^\theta/F$
is sufficiently large, since $a-bw$ has a uniform positive
lower bound. The  assertion with $b=0$
is also clear from the above inequality.

For $\kappa>2$, the same identity gives $\mathcal{T}\cdot \mathcal{S}<0$
when $\sqrt{2R}\,\mathcal{I}^\theta/U^\theta$ is sufficiently large.
The admissible cone inequality $(\kappa -2) (\mathcal{T} \cdot \mathcal{S}^\perp)^2 < 2 (\mathcal{T} \cdot \mathcal{S})^2$ can be written as
\[
2\left(a-bw-\frac{a\kappa}{\mathcal{I}^\theta/ F}\right)^2
-(\kappa-2)(aw+b)^2>0.
\]
Since
\[
2(a-bw)^2-(\kappa-2)(aw+b)^2
=(a^2+b^2)\left(2-2bw-\frac{b^2}{a}-(a-2)w^2\right)>0
\]
uniformly by assumption, the cone inequality follows by taking $\mathcal{I}^\theta/F$ sufficiently large.
\end{proof}

\paragraph{Normalized radial equations.}
We write radial equations for the normalized stresses and moments.
Write
\[
y=\log(R/R_{\rm ref}),\qquad \partial_y=R\partial_R.
\]
Recall that $a=1-2\partial_y\log U^\theta$ and define
\[
\zeta=\partial_Z\log U^\theta.
\]
Introduce the normalized inertial stresses
\begin{equation}\label{eq:normalized-inertial-stresses}
Q=\frac{L_\delta \mathcal{I}^\theta}{\sqrt{R/2}\,U^\theta},
\qquad
J=\frac{L_\delta \mathcal{I}^z}{\sqrt{R/2}\,(U^\theta)^2}.
\end{equation}
The quantities in the cone test then satisfy
\[
w=\frac{U^\theta J}{Q},\qquad
b=\frac{2\partial_yU^z}{U^\theta},\qquad
\frac{\mathcal{I}^\theta}{F}=\frac{RQ}{L_\delta}.
\]
Set
\begin{equation}\label{eq:normalized-transport-pressure}
W=1-\frac{(1-\delta)ZM^z+L\,\partial_ZM^z}{R},
\qquad
\Pi=\frac{2(1+\delta)ZP-L\partial_Z P}{(U^\theta)^2}.
\end{equation}
After removing the common similarity factor, the radial part of the
material derivative is $W\partial_y$, since \eqref{Ur-Uz} gives
$W=1-2ZU^z+L_\delta\sqrt{2/R}\,U^r$.
The radial pressure balance $2\partial_yP=(U^\theta)^2$
gives the pressure contribution to the axial stress equation in the form
\[
2Z\partial_yP+2(1+\delta)ZP-L\partial_ZP=(U^\theta)^2(Z+\Pi).
\]
This explains the square in the normalization of $J$ in
\eqref{eq:normalized-inertial-stresses}.

The radial equations for the inertial stresses become
\begin{equation}\label{eq:normalized-Q-J}
\begin{aligned}
\partial_yQ+\left(2-\frac a2\right)Q
={}&-W\left(1-\frac a2\right)
-\frac{\delta}{2}(1-2ZU^z)\\
&-\left(\frac{1-\delta}{2}Z+LU^z\right)\zeta,
\\[1ex]
\partial_yJ+(2-a)J
={}&Z+\Pi
-\frac{1}{(U^\theta)^2}
\Bigg[
W\partial_yU^z
+\frac{1+\delta}{2}(1-2ZU^z)U^z\\
&\hspace{35mm}
+\left(\frac{1-\delta}{2}Z+LU^z\right)\partial_ZU^z
\Bigg].
\end{aligned}
\end{equation}
For the inward propagation on $R\ge R_{3/2}$,
\eqref{eq:mc-required-moments} gives $M_{\mathrm{right}}^z=0$ and hence
$W_{\mathrm{right}}=1$. The right-side pressure is the tail integral
\eqref{eq:outer-pressure}, since all moment corrections are supported
below $R_{3/2}$. Thus
\[
\Pi_{\mathrm{right}}(R,Z)
=-\frac{1}{(\bar U^\theta(R,Z))^2}
\int_R^\infty
\left[(1+\delta)Z-L\,\partial_Z\log\bar U^\theta(\rho,Z)\right]
\frac{(\bar U^\theta(\rho,Z))^2}{\rho}\,d\rho.
\]
The bounds \eqref{eq:outer-pressure-sharp} give
$|\Pi_{\mathrm{right}}|\lesssim |Z|$. Since $U^z=0$ here,
the only source in the equation for $J$ is $Z+\Pi_{\mathrm{right}}$.
On the exact heat exterior $\Pi_{\rm heat}=-Z$, so this source vanishes.

\paragraph{The terminal scales.}

The right-side moments are defined by \eqref{eq:mc-required-moments}.
By the ansatz in Section~\ref{subsec:outer-ansatz},
$\bar U^\theta(R,Z)$ is independent of $Z$ for
$R\in[R_{3/2},R_{2h}]$, so we simply write
$\bar U^\theta(R)$ on this interval. Recall that
\begin{equation}\label{eq:mc-tail-scales}
\begin{aligned}
R_t&=e^{T_s+2}R_{3/2}\asymp\delta^{-4}R_{3/2},
&\bar U^\theta(R_t)&\asymp\bar U^\theta(R_{3/2})\delta^6,\\
R_{2h}&=e^\tau R_t,
&\bar U^\theta(R_{2h})
&=e^{-(1+\delta)\tau/2}\bar U^\theta(R_t),
\end{aligned}
\end{equation}
and
\begin{equation}\label{eq:mc-exterior-amplitude}
c_\infty
=\frac{\bar U^\theta(R_t)R_t^{(1+\delta)/2}}{1-\varepsilon}
\asymp
\bar U^\theta(R_{3/2})R_{3/2}^{(1+\delta)/2}
\delta^{4-2\delta}.
\end{equation}
We emphasize that $R_{2h}$ and $\bar U^\theta(R_{2h})$
depend on $\tau$, which is still an input in Section~\ref{subsec:outer-ansatz}
and will be fixed below. The amplitude $c_\infty$ is independent
of $\tau$, since $R^{(1+\delta)/2}\bar U^\theta(R)$ is constant
on the waiting interval $[R_t,R_{2h}]$.
At fixed $y=\log(R/R_{\rm ref})$, the heat factor still depends on
$R_{\rm ref}$ through its argument $2L/R$.

\paragraph{Steepening and waiting serve different moment conditions.}
For $R_{3/2}\le R\le R_{2h}$, the swirl is independent of $Z$ and
$\bar U^z=0$. Its slope changes to $-3/2$ over one logarithmic unit,
remains there for $T_s=4\log(2/\delta)$ units, and then changes to
$-(1+\delta)/2$ over one unit. Put
\[
R_s=eR_{3/2},\qquad R_t=e^{T_s+1}R_s.
\]
The amplitude and its radial energy weight satisfy
\[
\frac{\bar U^\theta(e^{T_s}R_s)}{\bar U^\theta(R_s)}
=(\delta/2)^6,
\qquad
\frac{e^{T_s}R_s[\bar U^\theta(e^{T_s}R_s)]^2}
{R_s[\bar U^\theta(R_s)]^2}
=(\delta/2)^8.
\]
In units of $R[U^\theta(R,Z)]^2$ at the start of each power law,
the swirl energy on Interval O.4 is of order $1/\mu$, while the
integral of the final power law costs $1/\delta$. The reduction
by $(\delta/2)^8$ thus prevents a loss of order $\mu/\delta$ in
the axial pulse balance. In particular,
\begin{equation}\label{eq:outer-tail-energy}
\begin{aligned}
\int_{R_{3/2}}^\infty(\bar U^\theta)^2\,dR
&\lesssim R_{3/2}[\bar U^\theta(R_{3/2})]^2,\\
\int_{R_{0z}}^\infty(\bar U^\theta)^2\,dR
&\lesssim[1+\log(1/\mu)]R_{0z}[\bar U^\theta(R_{0z},Z)]^2.
\end{aligned}
\end{equation}
The logarithm in the second estimate comes from Interval O.6, where
$\bar U^\theta\propto R^{-1/2-\mu}$ and $R_{3/2}=\mu^{-30}R_{m\theta}$.
With $t=\log(R/R_{m\theta})$, we have
\[
\frac{\displaystyle\int_{R_{m\theta}}^{R_{3/2}}(\bar U^\theta)^2\,dR}
     {R_{m\theta}[\bar U^\theta(R_{m\theta})]^2}
=\int_0^{30\log(1/\mu)}e^{-2\mu t}\,dt
=\frac{1-\mu^{60\mu}}{2\mu}
\sim 30\log(1/\mu)\qquad(\mu\to0).
\]
Since $R(\bar U^\theta)^2$ is nonincreasing on $[R_{0z},R_{3/2}]$,
the preceding interval of fixed logarithmic length $T_f=100$ and the
tail bounded by the first estimate each contribute at most a constant
multiple of $R_{0z}[\bar U^\theta(R_{0z},Z)]^2$.
The constants are uniform for small $\delta$ and all
$\tau\ge0$; the terminal multiplier and heat factor remain
bounded by absolute constants.

On the steep segment the angular-moment weight has a different behavior:
\begin{equation}\label{eq:outer-steep-angular-cost}
\int_{R_s}^{e^{-1}R_t}\sqrt{2R}\,\bar U^\theta(R)\,dR
=T_sR_s\sqrt{2R_s}\,\bar U^\theta(R_s).
\end{equation}
The extra factor $T_s$ is the angular-moment cost of the steep segment.
Thus steepening controls the tail energy, while waiting below the
heat power law (see the $1-\varepsilon$ factor in the waiting ansatz)
will remove the main angular-moment mismatch.
The waiting length is selected below; the two angular bumps in
Section~\ref{subsec:outer-angular-matching} then correct only the
remaining small discrepancies.
\MergedTailFigure

At $R_{3/2}$ we also have, uniformly for $\tau\ge0$,
\begin{equation}\label{eq:mc-required-tail-bounds}
\begin{aligned}
\|M_{\mathrm{right}}^{z\theta}(R_{3/2}, \cdot)\|_{C_Z^k}
&\lesssim_k R_{3/2}[\bar U^\theta(R_{3/2})]^2,\\
\|P_0^{\rm pre}(\cdot;\tau)+M_{\mathrm{right}}^p(R_{3/2}, \cdot)\|_{C_Z^k}
&\lesssim_k [\bar U^\theta(R_{3/2})]^2,\\
|\partial_ZM_{\mathrm{right}}^{z\theta}(R_{3/2},Z)|
&\le C\delta|Z|[\bar U^\theta(R_{2h})]^2,\\
\left|\partial_Z\left[P_0^{\rm pre}(Z;\tau)+M_{\mathrm{right}}^p(R_{3/2},Z)\right]\right|
&\le \frac{C\delta|Z|}{R_{2h}}
[\bar U^\theta(R_{2h})]^2.
\end{aligned}
\end{equation}
Indeed, the steep segment gives an energy integral bounded as
$\lesssim R_{3/2}[\bar U^\theta(R_{3/2})]^2$, while its entire
subsequent tail costs at most
\[
C\frac{R_t[\bar U^\theta(R_t)]^2}{\delta}
\le CR_{3/2}[\bar U^\theta(R_{3/2})]^2\delta^7.
\]
By \eqref{eq:mc-required-moments}, $P_0^{\rm pre}+M_{\mathrm{right}}^p$
is the right-side pressure, whose integral has the additional weight $1/R$.
Only the heat factor introduces $Z$-dependence beyond $R_{3/2}$;
differentiating it and integrating gives the last two estimates.

\paragraph{The stress computed from the right.}

We determine the right-side target stresses on this interval from the exterior,
using the normalized radial equations.
No stress estimate propagated from below $R_{3/2}$ is used.
The scalar matching condition at $R_{3/2}$ will identify
the selected waiting length within this backward solution.

\paragraph{The exact exterior and the last collar.}
The prescribed velocity is $U^\theta=U^\theta_{\rm heat}$, $U^z=0$
for $R\ge R_h$. We prescribe the right-side stress there by
\begin{equation}\label{eq:mc-exact-exterior}
\mathcal{T}_{\mathrm{right}}=0,\qquad
\mathcal{I}^\theta_{\mathrm{right}}=-\mathcal{S}^\theta_{\rm heat},\qquad
\mathcal{I}^z_{\mathrm{right}}=0
\qquad(R\ge R_h).
\end{equation}
These are terminal data, not yet an assertion about the stress recovered
from the reference branch. They are consistent with the right-side moments
\eqref{eq:mc-required-moments} by the
\hyperref[par:heat-moment-compatibility]{terminal-moment calculation}
in Section~\ref{sec:heat-exterior}.
In the backward calculations below, $\mathbf M,\mathcal{I},\mathcal{T},Q,J,W,\Pi$ denote the
right-side data, with the subscript suppressed; $U^\theta$ is the
prescribed terminal swirl, unchanged by the future corrections.
After \eqref{eq:outer-two-sided-matching} is established in
Section~\ref{subsec:outer-cone-closure}, these data will be identified
with the fields recovered from the axis.

\paragraph{A flat inward modification.}
We now extend the exterior slightly inward, with nonzero stress
at every point to the left of $R_h$. Fix $0<\varepsilon\le1/2$ and
write $y_h=\log(R_h/R)$. Set $R_{\rm col}=e^{-1}R_h=e^2R_{2h}$.
On $R_{\rm col}\le R<R_h$, prescribe
\begin{equation}\label{eq:outer-collar}
U^\theta(R,Z)
=c_\infty R^{-(1+\delta)/2}
H_\delta\left(\frac{2L}R\right)
\left(1-\varepsilon e^{-4/y_h^2}\right),
\qquad U^r=U^z=0.
\end{equation}
For $R\ge R_h$, retain the exact heat profile. Normalize the
pressure by the backward integral
\[
P(R,Z)=-\int_R^\infty\frac{(U^\theta(\rho,Z))^2}{2\rho}\,d\rho.
\]
The multiplier is flat at $y_h=0$, so the velocity and pressure
match the heat fields to every order at $R_h$.
Figure~\ref{fig:outer-collar} shows the modification.
In the notation of Interval O.8, $s=\log(R/R_{2h})$ and
$R_h=e^3R_{2h}$, so $y_h=3-s$ and
$e^{-4/y_h^2}=\mathfrak f((3-s)/2)$.
Thus this is the terminal part of \eqref{eq:road-heat-join},
where the global construction takes $\varepsilon=c_\varepsilon\delta$.
\MergedHeatFigure

\begin{proposition}[Admissible inward heat collar]\label{prop:heat-collar}
Fix $0<\delta\le1/2$, $c_\infty>0$ and
$0<\varepsilon\le1/2$. For $R_h$ sufficiently large, there exists
$\ell\in(0,1]$ such that the profile \eqref{eq:outer-collar},
with the pressure above and the inviscid-stress boundary values
\eqref{eq:outer-stress-boundary-data}, satisfies
\[
U^\theta>0,\qquad \mathcal{S}^\theta<0,\qquad
\mathcal{T}^\theta>0,\qquad \kappa>2+\delta/2
\quad(e^{-\ell}R_h\le R<R_h).
\]
Its stress satisfies the admissible cone condition
\eqref{cone-admissible} throughout this half-open collar and
\[
\frac{\mathcal{T}^z}{\mathcal{T}^\theta}
=O\!\left(\log(R_h/R)^6\right),\qquad
\frac{\mathcal{T}}{|\mathcal{T}|}\longrightarrow e_\theta
\quad(R\uparrow R_h),
\]
uniformly in $Z\in[-1,1]$.
Writing $y_h=\log(R_h/R)$, the leading stress on $0<y_h\le\ell$ satisfies
\[
|\mathcal T|\asymp e^{-4/y_h^2}y_h^{-3},\qquad
|\mathscr D^{\le m}\mathcal T|\le C_m e^{-4/y_h^2}y_h^{-3-3m}
\quad(m\ge0).
\]
This gives $c_h=4$ in the radial edge weight \eqref{intro:edge-weight}
used in Part~II. Near $R_h$ the inner factor of $\omega$ is bounded
above and below by positive constants, and $d_{\rm edge}=y_h$.
All fields match the heat fields smoothly at $R_h$.
The stress vanishes to every order there and is zero on the
exterior; no strict stress cone is asserted at $R_h$ or outside it.
The width $\ell$ and the constants in these local estimates may depend
on $\delta,c_\infty,\varepsilon,R_h$, and $C_m$ also on $m$. Realizing
the displayed stress boundary values by an axis-based profile
is a separate global moment-matching requirement.

For the parameters of the outer construction in \eqref{eq:mc-free-parameters},
with $\varepsilon=c_\varepsilon\delta$, these sign and cone conditions hold
on the full interval $[R_{\rm col},R_h)$.
\end{proposition}
\begin{proof}
Throughout the proof, we use the temporary notation
\[
y=\log(R_h/R),\qquad
\xi(y,Z)=\frac{2Le^y}{R_h},
\]
and write
\[
\Delta P=P-P_{\rm heat},\qquad
\Delta N^\theta
=N^\theta-N^\theta_{\rm heat},\qquad
\Delta \mathcal{I}^\theta
=\mathcal{I}^\theta-\mathcal{I}^\theta_{\rm heat}.
\]
Thus $R=R_he^{-y}$ and $R\partial_R=-\partial_y$.
We work on $0\le y\le1$ and choose $\ell\le1$ below for the
local endpoint estimates.
All bounds are uniform in $Z\in[-1,1]$. We keep
$\delta,c_\infty,\varepsilon,R_h$ fixed in the local estimates.
For the global construction the quantitative collar estimates below
are uniform in $P_*,R_{\rm ref},\delta$ and the provisional waiting
length $\tau\ge0$ in the range \eqref{eq:mc-free-parameters}.

On the outer collar, set
\[
f(y)=1-\varepsilon \mathfrak f(y/2),\qquad
U^\theta=fU^\theta_{\rm heat},\qquad
U^r=U^z=0,
\]
where
\[
U^\theta_{\rm heat}
=c_\infty R^{-(1+\delta)/2}H_\delta(\xi),
\qquad
F_{\rm heat}
=\frac{c_\infty}{\sqrt2}
R^{-1-\delta/2}H_\delta(\xi).
\]
Then $F=fF_{\rm heat}$ and
\[
\partial_y f(y)=-8\varepsilon \mathfrak f(y/2)y^{-3}.
\]
We take $0<\varepsilon\le1/2$, so that $f\ge1/2$.

The pressure is normalized by
$P(R_h,Z)=P_{\rm heat}(R_h,Z)$ and satisfies
$\partial_RP=F^2$.
For the inviscid stress, we use the matching conditions
\begin{equation}\label{eq:outer-stress-boundary-data}
\mathcal{I}^\theta(R_h,Z)
=\mathcal{I}^\theta_{\rm heat}(R_h,Z),
\qquad
\mathcal{I}^z(R_h,Z)=0.
\end{equation}
These are the boundary conditions for the backward radial
integration below; in the global construction, they must
be ensured by the moment matching conditions.

\medskip
\noindent\textbf{Shear and pressure estimates.}
Since $\mathcal{S}^\theta=2R\partial_RF$, we have
\begin{equation}\label{eq:outer-shear-exact}
\begin{aligned}
\mathcal{S}^\theta-\mathcal{S}^\theta_{\rm heat}
&=(f-1)\mathcal{S}^\theta_{\rm heat}
-2F_{\rm heat}\partial_y f(y)\\
&=\varepsilon \mathfrak f(y/2)y^{-3}
\left(16F_{\rm heat}
-y^3\mathcal{S}^\theta_{\rm heat}\right).
\end{aligned}
\end{equation}
The heat profile and its shear satisfy
\[
|F_{\rm heat}|+|\mathcal{S}^\theta_{\rm heat}|
=O\left(c_\infty R_h^{-1-\delta/2}\right)
\]
on this collar. Consequently,
\[
\mathcal{S}^\theta-\mathcal{S}^\theta_{\rm heat}
=O\left(
\varepsilon c_\infty R_h^{-1-\delta/2}
\mathfrak f(y/2)y^{-3}\right).
\]
Here the term containing $f-1$ is included in the estimate,
since $\mathfrak f(y/2)\le \mathfrak f(y/2)y^{-3}$ for $0<y\le1$.
Also, $\mathcal{S}^z=0$.

For the pressure, the common boundary value at $R_h$ gives
\[
\Delta P(R,Z)
=\int_R^{R_h}
\left[1-f\left(\log(R_h/\rho)\right)^2\right]
F_{\rm heat}(\rho,Z)^2\,d\rho.
\]
Changing variables to $s=\log(R_h/\rho)$ yields
\begin{equation}\label{eq:outer-pressure-integral}
\Delta P
=\frac{\varepsilon c_\infty^2}{2}R_h^{-1-\delta}
\int_0^y
\mathfrak f(s/2)\bigl(2-\varepsilon \mathfrak f(s/2)\bigr)e^{(1+\delta)s}
H_\delta(\xi(s,Z))^2\,ds.
\end{equation}

To justify both the size and the smooth factorization of
this integral, we record the following identity.
For $m=-3$ or $m=0$ and a smooth function $A$, by changing variable
\[t=\frac4{s^2}-\frac4{y^2},\qquad s=\frac{y}{\sqrt{1+y^2t/4}},\qquad \mathfrak f(s/2)=\mathfrak f(y/2)e^{-t},\]
one has 
\[
\mathfrak f(s/2)s^m\,ds
=-\frac18 \mathfrak f(y/2)y^{m+3}e^{-t}
\left(1+\frac{y^2t}{4}\right)^{-(m+3)/2}dt,
\]
hence
\begin{equation}\label{eq:outer-flat-integral}
\int_0^y \mathfrak f(s/2)s^mA(s,Z)\,ds
=\mathfrak f(y/2)y^{m+3}\mathcal A_m(y,Z),
\end{equation}
where the substitution $t=4/s^2-4/y^2$ gives
\[
\mathcal A_m(y,Z)
=\frac18\int_0^\infty e^{-t}
\left(1+\frac{y^2t}{4}\right)^{-(m+3)/2}
A\left(
\frac{y}{\sqrt{1+y^2t/4}},Z
\right)\,dt.
\]
This formula defines a smooth function up to $y=0$, with
\[
\mathcal A_m(0,Z)=\frac18A(0,Z).
\]
Indeed, every derivative of the integrand is bounded by
$e^{-t}$ times a polynomial in $t$, uniformly on the closed
collar, so differentiation under the integral is justified.

Applying \eqref{eq:outer-flat-integral} with $m=0$ to
\eqref{eq:outer-pressure-integral}, we obtain
\begin{equation}\label{eq:outer-pressure-factorization}
\Delta P
=\varepsilon c_\infty^2R_h^{-1-\delta}
\mathfrak f(y/2)y^3p(y,Z),
\end{equation}
where $p$ is smooth and bounded, together with its
$Z$ derivatives. In particular,
\[
|\Delta P|+|\partial_Z\Delta P|
=O\left(
\varepsilon c_\infty^2R_h^{-1-\delta}
\mathfrak f(y/2)y^3\right).
\]
The radial derivative is given exactly by
\begin{equation}\label{eq:outer-pressure-radial}
R\partial_R\Delta P
=R(f^2-1)F_{\rm heat}^2
=-\varepsilon \mathfrak f(y/2)\bigl(2-\varepsilon \mathfrak f(y/2)\bigr)
RF_{\rm heat}^2.
\end{equation}

\medskip
\noindent\textbf{Inviscid sources and backward integration.}
Since $U^r=U^z=0$, the angular equation in
\eqref{inviscid-stress-sources-Ur} reduces to
\[
N^\theta
=-\frac{\sqrt{R/2}}{L_\delta}
\left[
\frac{1+\delta}{2}U^\theta
+\frac{1-\delta}{2}Z\partial_ZU^\theta
+R\partial_RU^\theta
\right].
\]
Substituting $U^\theta=fU^\theta_{\rm heat}$ gives
\begin{equation}\label{eq:outer-angular-source-exact}
\Delta N^\theta
=(f-1)N^\theta_{\rm heat}
+\frac{\sqrt{R/2}}{L_\delta}U^\theta_{\rm heat}\partial_y f(y).
\end{equation}
Likewise, since $N^z_{\rm heat}=0$,
\begin{equation}\label{eq:outer-axial-source-exact}
N^z
=\frac{\sqrt{R/2}}{L_\delta}
\left[
2(1+\delta)Z\Delta P
-L\partial_Z\Delta P
+2ZR\partial_R\Delta P
\right].
\end{equation}
Equations \eqref{eq:outer-pressure-factorization}--%
\eqref{eq:outer-axial-source-exact} therefore imply the
exact factorizations
\begin{equation}\label{eq:outer-source-factorization}
\begin{aligned}
\Delta N^\theta
&=\varepsilon c_\infty R_h^{-\delta/2}
\mathfrak f(y/2)y^{-3}n_\theta(y,Z),\\
N^z
&=\varepsilon c_\infty^2R_h^{-1/2-\delta}
\mathfrak f(y/2)n_z(y,Z),
\end{aligned}
\end{equation}
with smooth bounded coefficients $n_\theta,n_z$.
For example,
\[
n_\theta
=-\frac{y^3N^\theta_{\rm heat}}
{c_\infty R_h^{-\delta/2}}
-\frac{4\sqrt2}{L_\delta}e^{\delta y/2}H_\delta(\xi).
\]
In the axial source, the terms involving $\Delta P$ and
$\partial_Z\Delta P$ contain a factor $y^3$; the term
$R\partial_R\Delta P$ has only the factor $\mathfrak f(y/2)$.
This explains the different powers of $y$ in
\eqref{eq:outer-source-factorization}.

In the variable $y$, the radial equations
\eqref{inviscid-stress-radial} read
\[
-\partial_y\Delta \mathcal{I}^\theta
+\Delta \mathcal{I}^\theta
=\Delta N^\theta,\qquad
-\partial_y\mathcal{I}^z+\frac12\mathcal{I}^z
=N^z.
\]
Using \eqref{eq:outer-stress-boundary-data}, we obtain
\[
\begin{aligned}
\Delta \mathcal{I}^\theta(y,Z)
&=-e^y\int_0^y e^{-s}
\Delta N^\theta(s,Z)\,ds,\\
\mathcal{I}^z(y,Z)
&=-e^{y/2}\int_0^y e^{-s/2}
N^z(s,Z)\,ds.
\end{aligned}
\]
Applying \eqref{eq:outer-flat-integral}, with $m=-3$
in the first integral and $m=0$ in the second, gives
\begin{equation}\label{eq:outer-inviscid-factorization}
\begin{aligned}
\Delta \mathcal{I}^\theta
&=\varepsilon c_\infty R_h^{-\delta/2}
\mathfrak f(y/2)a_\theta(y,Z),\\
\mathcal{I}^z
&=\varepsilon c_\infty^2R_h^{-1/2-\delta}
\mathfrak f(y/2)y^3a_z(y,Z),
\end{aligned}
\end{equation}
where $a_\theta,a_z$ are smooth and bounded up to $y=0$.

\medskip
\noindent\textbf{Stress direction and the admissible cone.}
The heat stress vanishes:
\[
\mathcal{I}^\theta_{\rm heat}
+\mathcal{S}^\theta_{\rm heat}=0,\qquad
\mathcal{I}^z_{\rm heat}=\mathcal{S}^z_{\rm heat}=0.
\]
Combining \eqref{eq:outer-shear-exact} with
\eqref{eq:outer-inviscid-factorization}, we find
\begin{equation}\label{eq:outer-stress-factorization}
\mathcal{T}^\theta=\mathfrak f(y/2)y^{-3}b_\theta(y,Z),\qquad
\mathcal{T}^z=\mathfrak f(y/2)y^3b_z(y,Z),
\end{equation}
where
\[
\begin{aligned}
b_\theta
&=\varepsilon\left[
16F_{\rm heat}
-y^3\mathcal{S}^\theta_{\rm heat}
+c_\infty R_h^{-\delta/2}y^3a_\theta
\right],\\
b_z
&=\varepsilon c_\infty^2R_h^{-1/2-\delta}a_z.
\end{aligned}
\]
Both coefficients are smooth and bounded. Moreover,
\[
b_\theta(0,Z)
=16\varepsilon F_{\rm heat}(R_h,Z)
=8\sqrt2\,\varepsilon c_\infty R_h^{-1-\delta/2}
H_\delta\left(\frac{2L}{R_h}\right)>0.
\]
For fixed parameters, compactness of $[-1,1]$ gives
\[
m_\theta:=\min_{Z\in[-1,1]}b_\theta(0,Z)>0.
\]
After reducing $\ell$, we therefore have
$b_\theta(y,Z)\ge m_\theta/2$ throughout the collar.
It follows that
\begin{equation}\label{eq:outer-stress-direction}
\mathcal{T}^\theta>0,\qquad
\frac{\mathcal{T}^z}{\mathcal{T}^\theta}
=y^6\frac{b_z}{b_\theta}=O(y^6),
\qquad 0<y\le\ell.
\end{equation}

Since $\mathcal S^z=0$, the shear intensity satisfies
\[
\kappa=\kappa_{\rm heat}+\frac{2\partial_y f}{f},
\qquad
\kappa_{\rm heat}
=2+\delta+2\xi\frac{\partial_\xi H_\delta}{H_\delta}
=2+\delta+O(R_h^{-1}).
\]
For fixed $\delta>0$, take $R_h$ sufficiently large and then $\ell$
sufficiently small. The flatness of $\partial_y f$ at zero gives
$\kappa>2+\delta/2$ and $\mathcal S^\theta=-\kappa F<0$ on this collar.
By \eqref{eq:outer-stress-direction},
\[
\mathcal T\cdot\mathcal S<0,\qquad
(\kappa-2)\left(\frac{\mathcal T^z}{\mathcal T^\theta}\right)^2
=O(y^{12})<2
\]
after a further reduction of $\ell$. These are the admissible cone
inequalities when $\mathcal S^z=0$.
Also $\mathcal T/|\mathcal T|\to e_\theta$ uniformly in $Z$.
The smooth factors in \eqref{eq:outer-stress-factorization} show
that $|\mathcal T|\asymp e^{-4/y^2}y^{-3}$ for $0<y\le\ell$.
Each $y$ derivative costs at most a factor $y^{-3}$, while $Z$
derivatives act only on the smooth coefficients. Since
$\sqrt{2R}\,\partial_R=-\sqrt{2/R}\,\partial_y$ and the leading stress
is independent of $\lambda$, this proves the stated derivative bounds
and smooth vanishing at $R_h$.

\medskip
\noindent\textbf{Quantitative bounds on the full collar.}
We now use $\varepsilon=c_\varepsilon\delta$ and the global scales
\eqref{eq:mc-tail-scales}. In particular,
$R_{\rm col}\gtrsim\delta^{-4}R_{3/2}$ and
$U^\theta_{\rm heat}(R_{\rm col})\lesssim
\delta^6 U^\theta(R_{3/2})$.
Subtracting the heat equation and integrating backward from
\eqref{eq:mc-exact-exterior} gives the exact identity
\begin{equation}\label{eq:mc-terminal-positive-stress}
\begin{aligned}
\mathcal{T}^\theta(R,Z)
={}&-\sqrt{\frac2R}\,U^\theta_{\rm heat}(R,Z)\partial_y f(y)\\
&-\frac1R\int_R^{R_h}
\left[
\frac{\sqrt{\rho/2}\,U^\theta_{\rm heat}(\rho,Z)}{L_\delta}
-\mathcal{S}^\theta_{\rm heat}(\rho,Z)
\right]
(\partial_y f)\left(\log\frac{R_h}{\rho}\right)\,d\rho.
\end{aligned}
\end{equation}
Both terms are positive since $\partial_y f<0$.
Since the remaining logarithmic length is at most one,
the radii and heat amplitudes in the integral are comparable
by fixed constants. Thus
\[
\mathcal{T}^\theta
\ge c\varepsilon\sqrt R\,U^\theta_{\rm heat}\mathfrak f(y/2).
\]
In \eqref{eq:outer-pressure-integral}, the $Z$ dependence comes from
$H_\delta(2L/\rho)$, whose $Z$ derivative is $O(\delta |Z|/\rho)$.
Hence
\[
 |Z\Delta P|+|\partial_Z\Delta P|
 \le C\varepsilon |Z|[U^\theta_{\rm heat}]^2\mathfrak f(y/2).
\]
Substituting this and \eqref{eq:outer-pressure-radial} into
\eqref{eq:outer-axial-source-exact} and integrating backward gives
\[
|\mathcal{T}^z|
\le C\varepsilon|Z|\sqrt R\,
[U^\theta_{\rm heat}]^2\mathfrak f(y/2),
\qquad
\left|\frac{\mathcal{T}^z}{\mathcal{T}^\theta}\right|
\le C|Z|U^\theta_{\rm heat}.
\]
The heat estimates and the flat multiplier give
\[
0\le2+\delta-\kappa_{\rm heat}\le C\delta/R_{\rm col},
\qquad
\left|\frac{2\partial_y f}{f}\right|\le Cc_\varepsilon\delta.
\]
Taking the fixed $c_\varepsilon$ small and $R_{\rm corr}$ large gives
\[
b=0,\qquad
2+\frac\delta2<a=\kappa\le2+\delta.
\]
Consequently $\mathcal{T}^\theta>0$, $\mathcal{S}^\theta<0$, and
\[
(\kappa-2)
\left(\frac{\mathcal{T}^z}{\mathcal{T}^\theta}\right)^2
\le C\delta[U^\theta_{\rm heat}]^2<2.
\]
This proves the cone for the right-side data on $R_{\rm col}\le R<R_h$,
with vanishing stress at $R_h$.
\end{proof}

\paragraph{Initial data for inward propagation.}
At $R_{\rm col}$ the collar coordinate $y=\log(R_h/R)$ equals one.
The integral of $-\partial_y f$ from $0$ to $1$ is $\varepsilon e^{-4}$. On this unit
interval the radii are comparable to $R_{\rm col}$ and the heat
amplitudes to $U^\theta(R_{\rm col},Z)$. Thus the integral term in
\eqref{eq:mc-terminal-positive-stress} is comparable to
$\varepsilon\sqrt{R_{\rm col}}\,U^\theta(R_{\rm col},Z)$, and its
first term is $O(\varepsilon U^\theta(R_{\rm col},Z)/\sqrt{R_{\rm col}})$.
Together with the axial bound above, this gives
\begin{equation}\label{eq:mc-backward-collar-data}
\begin{aligned}
\mathcal{T}^\theta(R_{\rm col},Z)
&\asymp\varepsilon\sqrt{R_{\rm col}}\,U^\theta(R_{\rm col},Z),\\
|\mathcal{T}^z(R_{\rm col},Z)|
&\le C\varepsilon|Z|\sqrt{R_{\rm col}}\,[U^\theta(R_{\rm col},Z)]^2.
\end{aligned}
\end{equation}
Since $\mathcal S^\theta=-\kappa U^\theta/\sqrt{2R}$ and
$\mathcal S^z=0$, \eqref{eq:normalized-inertial-stresses} gives
$Q=L_\delta\mathcal T^\theta/(\sqrt{R/2}\,U^\theta)+L_\delta\kappa/R$
and $|J|\lesssim\varepsilon |Z|$ at $R_{\rm col}$.
Here $R_{\rm col}^{-1}\lesssim\delta^4/R_{3/2}\ll\varepsilon$ by
\eqref{eq:mc-tail-scales}. For the pressure, subtracting the heat identity
$\Pi_{\rm heat}=-Z$ yields
\[
 (U^\theta)^2(\Pi+Z)
 =2(1+\delta)Z\Delta P-L\partial_Z\Delta P
   +Z(f^2-1)(U^\theta_{\rm heat})^2.
\]
The collar pressure bound and $|f(1)^2-1|\le2\varepsilon e^{-4}$
therefore give the normalized initial data
\begin{equation}\label{eq:mc-backward-normalized-data}
\begin{gathered}
c\varepsilon\le Q(R_{\rm col},Z)\le C\varepsilon,\\
|J(R_{\rm col},Z)|+|\Pi(R_{\rm col},Z)+Z|\le C\varepsilon|Z|.
\end{gathered}
\end{equation}

\paragraph{The equations used for backward propagation.}
On $[R_{3/2},R_{\rm col}]$, the right-side moment data give
$U^z=U^r=M^z=M^{\theta z}=0$, hence $W=1$.
Recall that $\partial_y=R\partial_R$, $a=1-2\partial_y\log U^\theta$
and $\zeta=\partial_Z\log U^\theta$. The normalized equations are
\begin{equation}\label{eq:mc-backward-Q-J}
\begin{aligned}
\partial_yQ+\left(2-\frac a2\right)Q
&=\frac a2-1-\frac\delta2-\frac{1-\delta}{2}Z\zeta,\\
\partial_yJ-(a-2)J&=Z+\Pi.
\end{aligned}
\end{equation}
The pressure integral also gives
\begin{equation}\label{eq:mc-backward-pressure}
\partial_y(\Pi+Z)-(a-1)(\Pi+Z)
=(2+\delta-a)Z-L\zeta.
\end{equation}
We solve these equations inward from $R_{\rm col}$.

\paragraph{From $R_{\rm col}$ to $R_{2h}$.}
This interval has logarithmic length two.
The terminal ansatz satisfies
\[
2+\frac\delta2\le a\le2+\delta,\qquad
Z\zeta\ge0,\qquad
|2+\delta-a|\le C\varepsilon,\qquad
|\zeta|\le C\frac{\delta|Z|}{R_{2h}}.
\]
The small heat terms are absorbed here using
$R_{2h}\ge R_t\asymp\delta^{-4}R_{3/2}$.
In particular, the source in the $Q$ equation is nonpositive.
Backward integration of \eqref{eq:mc-backward-Q-J} and
\eqref{eq:mc-backward-pressure} over these two units gives
\begin{equation}\label{eq:mc-heat-entry-data}
Q\asymp\varepsilon,\qquad
|J|+|\Pi+Z|\le C\varepsilon|Z|
\qquad(R_{2h}\le R\le R_{\rm col}).
\end{equation}

\paragraph{The waiting interval: exact formulas without a bound on $\tau$.}
For $R_t\le R\le R_{2h}$, put
\[
u=\log(R_{2h}/R)\in[0,\tau].
\]
Here $a=2+\delta$ and $\zeta=0$.
Solving the three equations backward gives exactly
\begin{equation}\label{eq:mc-backward-waiting}
\begin{aligned}
Q(R,Z)
&=Q(R_{2h},Z)e^{(1-\delta/2)u},\\
\Pi(R,Z)+Z
&=[\Pi(R_{2h},Z)+Z]e^{-(1+\delta)u},\\
J(R,Z)
&=e^{-\delta u}\Big[
J(R_{2h},Z)
-[\Pi(R_{2h},Z)+Z](1-e^{-u})
\Big].
\end{aligned}
\end{equation}
Thus
\[
|J(R,Z)|\le C\varepsilon|Z|e^{-\delta u}.
\]
Recall that $w=\mathcal I^z/\mathcal I^\theta=U^\theta J/Q$.
Using the bound for $J$, the formula for $Q$, and
$U^\theta(R)=U^\theta(R_{2h})e^{(1+\delta)u/2}$, we obtain
\begin{equation}\label{eq:mc-waiting-direction}
|w(R,Z)|
\le C U^\theta(R_{2h})|Z|e^{-u/2}
\le C U^\theta(R_{3/2})\delta^6|Z|.
\end{equation}
The last inequality uses $\tau\ge0$ and \eqref{eq:mc-tail-scales}.
No upper bound on $\tau$ has been used.

\paragraph{The slope transitions and the steep segment.}
Recall
\[
R_s=eR_{3/2},\qquad
R_t=e^{T_s+1}R_s.
\]
There is one logarithmic unit of slope change on either side
of the steep segment $[R_s,e^{T_s}R_s]$.
Throughout $[R_{3/2},R_t]$, $\zeta=0$ and $2+\delta\le a\le4$.
The backward pressure equation gives $|\Pi|\le C|Z|$.
Backward propagation across either unit transition therefore
preserves $|J|\le C|Z|$.
On $[R_s,e^{-1}R_t]$ the prescribed power $U^\theta\propto R^{-3/2}$
gives $a=1-2(-3/2)=4$. With $t=\log(R/R_s)\in[0,T_s]$, the axial
equation gives
\[
\begin{aligned}
J(R_se^t,Z)
={}&e^{-2(T_s-t)}J(e^{T_s}R_s,Z)\\
&-\int_t^{T_s}e^{-2(v-t)}
[Z+\Pi(R_se^v,Z)]\,dv.
\end{aligned}
\]
Consequently,
\begin{equation}\label{eq:mc-backward-J-bound}
|J(R,Z)|\le C|Z|
\qquad(R_{3/2}\le R\le R_t),
\end{equation}
with no loss depending on $T_s$ or $\tau$.

On the pure steep segment, $a=4$ and $\zeta=0$, so the angular equation
in \eqref{eq:mc-backward-Q-J} reduces to $\partial_y Q=1-\delta/2$.
For any provisional waiting length $\tau$ this gives
\[
Q(R,Z)=Q(e^{-1}R_t,Z)
-\left(1-\frac{\delta}{2}\right)\log\frac{e^{-1}R_t}{R},
\qquad R_s\le R\le e^{-1}R_t.
\]
Thus $Q$ decreases by $(1-\delta/2)T_s$ when propagated inward to $R_s$.
Choosing $\tau$ by angular-moment matching below gives the endpoint bound
$Q(R_{3/2},Z)\ge c\mu$ in \eqref{eq:mc-backward-compatibility}
and hence the lower bound throughout $[R_{3/2},R_t]$
in \eqref{eq:outer-steep-positive-propagation}.

\paragraph{Choosing the waiting length.}
We choose $\tau$ by matching the angular moments of $\bar U^\theta_{\rm pre}$
computed from the left and from the right at $(R,Z)=(R_{3/2},0)$.
The two angular bumps then correct the remaining
$Z$-dependent and heat discrepancies and restore the pressure moment.
Since $\bar U^\theta_{\rm pre}=\bar U^\theta$ for $R\le R_{3/2}$,
the incoming angular moment is $\bar M_{\mathrm{left}}^\theta$ of the
uncorrected reference-plus-outer candidate.

The right-side pre-heat angular target is
\begin{equation}\label{eq:mc-preheat-angular-target}
\begin{aligned}
M_{\mathrm{right,pre}}^\theta(R,Z;\tau)
={}&\frac{\sqrt2c_\infty}{1-\delta/2}R^{1-\delta/2}\\
&-\int_R^\infty\sqrt{2\rho}
 [\bar U^\theta_{\rm pre}(\rho,Z;\tau)
       -c_\infty\rho^{-(1+\delta)/2}]\,d\rho.
\end{aligned}
\end{equation}
For $R\ge R_{3/2}$ it is independent of $Z$.
The matching equation is
\begin{equation}\label{eq:mc-wait-choice}
\boxed{M_{\mathrm{right,pre}}^\theta(R_{3/2},0;\tau)
=\bar M_{\mathrm{left}}^\theta(R_{3/2},0)}.
\end{equation}
The left-side incoming moment is independent of $\tau$.
For the following estimates take $\mu_{\rm corr}\le1/60$ and
$d_{\rm corr}\le1/4$ in \eqref{eq:mc-free-parameters}, so that
$\mu\le1/60$ and $\delta\le\mu/4$.

\paragraph{The incoming moment.}
Set
\[
X(R,Z)=
\frac{\bar M_{\mathrm{left}}^\theta(R,Z)}
{R\sqrt{2R}\bar U^\theta(R,Z)}.
\]
With $\partial_y=R\partial_R$, the moment equation gives
\[
\partial_yX+
\left(\frac32+\partial_y\log\bar U^\theta\right)X=1.
\]
Direct integration of the reference profile gives
$X(R_{\rm ref},Z)=5/8$.
Up to $R_{m\theta}$, the coefficient in this equation is at least $4/5$,
and its $Z$-derivatives are bounded by fixed constants.
The equation and its differentiated equations therefore give
\[
\|X(R_{m\theta},\cdot)\|_{C_Z^k}\le C_k.
\]

On $[R_{m\theta},R_{3/2}]$, the coefficient is exactly $1-\mu$.
Since $R_{3/2}/R_{m\theta}=\mu^{-30}$, we obtain
\[
X(R_{3/2},Z)
=\frac1{1-\mu}
+\mu^{30(1-\mu)}
\left[X(R_{m\theta},Z)-\frac1{1-\mu}\right].
\]
Thus
\begin{equation}\label{eq:mc-incoming-angular-moment}
\left\|
\frac{\bar M_{\mathrm{left}}^\theta(R_{3/2},\cdot)}
{R_{3/2}\sqrt{2R_{3/2}}\bar U^\theta(R_{3/2})}
-\frac1{1-\mu}
\right\|_{C_Z^k}
\le C_k\mu^{29}.
\end{equation}
In particular,
\[
\bar M_{\mathrm{left}}^\theta(R_{3/2},0)
\asymp R_{3/2}\sqrt{2R_{3/2}}\bar U^\theta(R_{3/2}),
\]
and the incoming angular moment is almost independent of $Z$.

\paragraph{Monotonicity of the target.}
For each fixed $R$, the pre-heat profile
$\bar U^\theta_{\rm pre}(R,0;\tau)$ is nonincreasing in $\tau$.
Indeed, increasing $\tau$ decreases
$s=\log(R/R_t)-\tau$, and its terminal multiplier is
\[
(1-\sigma(s))(1-\varepsilon)
 +\sigma(s)(1-\varepsilon \mathfrak f((3-s)/2)).
\]
Both $\sigma(s)$ and $1-\varepsilon \mathfrak f((3-s)/2)$ are
nondecreasing, and the latter is at least $1-\varepsilon$.
Thus the multiplier is nondecreasing in $s$, and increasing
$\tau$ strictly decreases the profile on a nonempty interval.
Since $c_\infty$ is independent of $\tau$, the backward formula
\eqref{eq:mc-preheat-angular-target} shows that
$M_{\mathrm{right,pre}}^\theta(R_{3/2},0;\tau)$ is continuous
and strictly increasing in $\tau$.
\paragraph{Existence of the matching waiting length.}
At $\tau=0$, the terminal formula and
\eqref{eq:mc-tail-scales} give
\[
\begin{aligned}
M_{\mathrm{right,pre}}^\theta(R_{3/2},0;0)
={}&-\int_{R_{3/2}}^{R_t}
\sqrt{2R}\bar U^\theta(R)\,dR\\
&+O\left(
R_{3/2}\sqrt{2R_{3/2}}\bar U^\theta(R_{3/2})
\right).
\end{aligned}
\]
On the steep segment,
$R\sqrt{2R}\bar U^\theta(R)$ is constant and comparable to
$R_{3/2}\sqrt{2R_{3/2}}\bar U^\theta(R_{3/2})$.
Its logarithmic length is $T_s$.
The two unit transition intervals contribute only a fixed
multiple of the same scale.
Consequently,
\begin{equation}\label{eq:mc-zero-wait-discrepancy}
\frac{
M_{\mathrm{right,pre}}^\theta(R_{3/2},0;0)
-\bar M_{\mathrm{left}}^\theta(R_{3/2},0)}
{R_{3/2}\sqrt{2R_{3/2}}\bar U^\theta(R_{3/2})}
\le C-cT_s<0,
\end{equation}
after decreasing the fixed upper bound $\mu_{\rm corr}$ if needed.

On the other hand, the target tends to infinity as $\tau\to\infty$.
Indeed, on the waiting interval $[R_t,e^\tau R_t]$,
\[
\bar U^\theta(R)
-c_\infty R^{-(1+\delta)/2}
=-\varepsilon c_\infty R^{-(1+\delta)/2}.
\]
Its contribution to the backward moment grows without bound.
The terminal part also contributes nonnegatively, since
$\bar U^\theta_{\rm pre}(R,0;\tau)\le c_\infty R^{-(1+\delta)/2}$
for $R\ge R_t$.
More explicitly,
\[
\frac{
M_{\mathrm{right,pre}}^\theta(R_{3/2},0;\tau)
-\bar M_{\mathrm{left}}^\theta(R_{3/2},0)}
{R_{3/2}\sqrt{2R_{3/2}}\bar U^\theta(R_{3/2})}
\ge
c\varepsilon\left(e^{(1-\delta/2)\tau}-1\right)
-C(1+T_s)
\longrightarrow+\infty.
\]
Continuity and strict monotonicity therefore give a unique
$\tau>0$ satisfying \eqref{eq:mc-wait-choice}.
We fix this value from now on. After division by
$R_{3/2}\sqrt{2R_{3/2}}\bar U^\theta(R_{3/2})$,
the matching equation contains only the dimensionless schedule.
Hence this choice of $\tau$ is independent of $R_{\rm ref}$,
and so is $P_0^{\rm pre}$.
The angular correction in Section~\ref{subsec:outer-angular-matching}
matches the pressure from this axis datum to the prescribed right-side pressure.

\paragraph{Positivity of the right-side angular stress.}
We verify the lower bound for the right-side solution at the selected
waiting length. Since $M^z=M^{\theta z}=0$, the angular moment formula
\eqref{Itheta-moments} gives
\[
Q_{\mathrm{right}}(R,Z)
=-1+
\frac{(1-\delta/2)M_{\mathrm{right}}^\theta(R,Z)
-\frac{1-\delta}{2}Z\partial_ZM_{\mathrm{right}}^\theta(R,Z)}
{R\sqrt{2R}\,U^\theta(R)}
\qquad(R_{3/2}\le R\le R_t).
\]
Let $Q_{\mathrm{pre}}$ be the same expression with
$M_{\mathrm{right}}^\theta$ replaced by $M_{\mathrm{right,pre}}^\theta$.
The latter moment is independent of $Z$ on this interval.
By \eqref{eq:mc-wait-choice} and \eqref{eq:mc-incoming-angular-moment},
\[
Q_{\mathrm{pre}}(R_{3/2})
=(1-\delta/2)X(R_{3/2},0)-1
=\frac{\mu-\delta/2}{1-\mu}+O(\mu^{29})
\ge c\mu.
\]
Here we use $\delta\le\mu/4$ and the fixed smallness choice of
$\mu_{\rm corr}$.

The heat replacement improves this bound. On $R_{3/2}\le R\le R_{2h}$,
the difference $M_{\mathrm{right}}^\theta-M_{\mathrm{right,pre}}^\theta$
is independent of $R$ and equals
\[
\Delta M^\theta(Z)
=\int_{R_{2h}}^\infty\sqrt{2\rho}\,
  [\bar U^\theta_{\rm pre}(\rho,Z)-\bar U^\theta(\rho,Z)]\,d\rho.
\]
In \eqref{eq:connecting-swirl}, the difference in brackets is a
nonnegative factor independent of $Z$ times $1-H_\delta(2L/\rho)$.
Since $H_\delta\le1$ and $Z\partial_ZH_\delta(2L/\rho)\ge0$,
we have $\Delta M^\theta\ge0$ and $Z\partial_Z\Delta M^\theta\le0$.
Consequently,
\[
Q_{\mathrm{right}}-Q_{\mathrm{pre}}
=\frac{(1-\delta/2)\Delta M^\theta
-\frac{1-\delta}{2}Z\partial_Z\Delta M^\theta}
{R\sqrt{2R}\,U^\theta(R)}
\ge0.
\]
In particular,
\begin{equation}\label{eq:mc-backward-compatibility}
Q_{\mathrm{right}}(R_{3/2},Z)\ge c\mu
\qquad(-1\le Z\le1).
\end{equation}

For the inward propagation through the steep segment and its two
slope transitions, the right-side data have $W=1$ and $\zeta=0$.
The angular equation \eqref{eq:mc-backward-Q-J} gives
\[
\begin{aligned}
R\sqrt{2R}\,U^\theta(R)Q(R,Z)
={}&R_t\sqrt{2R_t}\,U^\theta(R_t)Q(R_t,Z)\\
&-\int_R^{R_t}\sqrt{2\rho}\,U^\theta(\rho)
  \left(\frac{a(\rho)}2-1-\frac\delta2\right)\,d\rho.
\end{aligned}
\]
The integrand is nonnegative, so the total loss is largest at
$R=R_{3/2}$. By \eqref{eq:mc-backward-compatibility}, the remaining
quantity there is at least
$c\mu R_{3/2}\sqrt{2R_{3/2}}U^\theta(R_{3/2})$.
The factor $R\sqrt{2R}\,U^\theta(R)$ is constant on the steep segment
and increases by at most $e^2$ over the two unit transitions.
Therefore
\begin{equation}\label{eq:outer-steep-positive-propagation}
Q_{\mathrm{right}}(R,Z)\ge c\mu
\qquad(R_{3/2}\le R\le R_t,\ -1\le Z\le1).
\end{equation}

\paragraph{The cone for the right-side data.}
By \eqref{eq:outer-steep-positive-propagation}, $Q\ge c\mu$ on
$[R_{3/2},R_t]$. There \eqref{eq:mc-backward-J-bound} gives
$|w|=U^\theta|J|/Q\le CU^\theta(R_{3/2})/\mu\ll1$. On
$[R_t,R_{\rm col}]$, the waiting and heat estimates
\eqref{eq:mc-backward-waiting}--\eqref{eq:mc-waiting-direction} and
\eqref{eq:mc-heat-entry-data} give $Q\ge c\varepsilon$ and $|w|\ll1$.
Consequently
\[
 \frac{RQ}{L_\delta}
 \ge c\min\{\mu R_{3/2},\varepsilon R_t\}\ge c R_{\rm ref}
 \qquad(R_{3/2}\le R\le R_{\rm col})
\]
by \eqref{eq:mc-tail-scales}. Taking $R_{\rm ref}\ge R_{\rm corr}$ with
$R_{\rm corr}$ sufficiently large makes $\mathcal I^\theta/F=RQ/L_\delta$
uniformly large on $[R_{3/2},R_{\rm col}]$. Since $b=0$,
$2<\kappa=a\le4$ and $|w|\ll1$, Lemma~\ref{lem:outer-cone-test}
gives the admissible cone there. Proposition~\ref{prop:heat-collar} gives the cone on the remaining
collar up to $R_h$. Thus the right-side data satisfy the admissible
cone throughout $[R_{3/2},R_h)$ for the selected waiting length.

\subsection{Interval O.2: Turning the swirl and switching off the axial flow}
\label{subsec:outer-left-data}\label{subsec:outer-o2}

We now return to the uncorrected reference-plus-outer candidate $\bar U$.
In this subsection, $\mathbf M=\bar{\mathbf M}_{\mathrm{left}}$ consists of its cumulative
moments from the axis, and we write $U=\bar U$. Since
$\bar U^\theta=\bar U^\theta_{\rm pre}$ up to $R_{2h}$, the pressure
obtained from the prescribed axis value satisfies
\[
 P(R,Z)=P_0^{\rm pre}(Z)+\bar M_{\mathrm{left}}^p(R,Z)
       =P_{\rm pre}(R,Z)\qquad(0\le R\le R_{2h}).
\]
With this pressure we propagate the stress and moments outwards
to verify the relaxed cone.
The starting data come from the temporary reference profile, not from
the analytic core that will be constructed later.

\begin{lemma}\label{lem:outer-turnoff-cone}
    The profile $(\bar U^\theta, \bar U^z)$ specified in Section~\ref{subsec:outer-ansatz} satisfies the relaxed cone condition on $[R_{\rm useless}, R_{-\mu}]$.
\end{lemma}

\begin{proof}
We prove it by propagating the inviscid stress $\mathcal{I}$ using its radial equations \eqref{inviscid-stress-radial} in increasing $R$. We also collect useful estimates on the moments along the way.

\paragraph{Initial data at $R_{\rm ref}$.}
As already mentioned, the five moments at $R_{\rm ref}$ are obtained by integrating
the reference profile from the axis.
For the uncorrected candidate, \eqref{eq:connecting-axial} gives
\[
\bar M_{\mathrm{left}}^z(R,Z)=4Z\Sigma_B(R),\qquad
W(R,Z)=1-\frac{4L_\delta\Sigma_B(R)}{R},\qquad
1-4L_\delta\le W(R,Z)\le1,
\]
where $\Sigma_B(R)=\int_0^R B(\log(\rho/R_{\rm ref})-1)\,d\rho$
lies between $0$ and $R$.

A direct computation using \eqref{Itheta-moments} gives
\begin{equation}\label{eq:outer-reference-Q}
Q(R_{\rm ref},Z)
=\frac98-\frac{5\delta}{16}+\delta Z^2
+\frac{5Z^2(9-\delta-8Z^2)}{8(1+Z^2)}
>1
\end{equation}
for sufficiently small $\delta$. Since $\Sigma_B(R_{\rm ref})=R_{\rm ref}$,
\begin{equation}\label{eq:outer-reference-W}
W(R_{\rm ref},Z)=1-4L_\delta=-3+4\delta Z^2.
\end{equation}
For the axial stress, the moment formula \eqref{Iz-moments-direct} yields
\begin{equation}\label{eq:outer-reference-J}
\begin{aligned}
J(R_{\rm ref},Z)
={}&\Pi(R_{\rm ref},Z)
-\frac{5Z}{6}\left(\delta+\frac{2L}{1+Z^2}\right)\\
&+\frac{4Z\bigl((8+4\delta)Z^2-5\bigr)(1+Z^2)^2}{P_*^2}.
\end{aligned}
\end{equation}
On $[R_{\rm ref},R_{2h}]$, Lemma~\ref{lem:outer-pressure} applied to
$P_{\rm pre}$ gives $|\Pi(R,Z)|\lesssim |Z|$. Here we use
$U^\theta\asymp P_*/(1+Z^2)$ on $[R_{\rm ref},eR_{\rm ref}]$
and \eqref{eq:outer-pressure-sharp} on $[eR_{\rm ref},R_{2h}]$. Hence
\[
|J(R_{\rm ref},Z)|
\lesssim |Z|(1+P_*^{-2})
\lesssim |Z|.
\]

On the reference branch $[R_{\rm useless},R_{\rm ref}]$,
the  normalized quantities satisfy
\[
Q(R,Z)=Q(R_{\rm ref},Z),\qquad
a=\frac45,\qquad b=0,\qquad \kappa=\frac45.
\]
In particular,
\[
\frac{\mathcal{I}^\theta}{F}
=\frac{RQ_{\rm ref}}{L_\delta}
>\frac{R}{L_\delta}.
\]
Lemma~\ref{lem:outer-cone-test} therefore gives the relaxed
cone condition throughout this interval provided
\[
R_{\rm useless}>2,
\qquad\text{equivalently}\qquad R_{\rm ref}>2e^5.
\]
We include $R_{\rm corr}>2e^5$ in the choice of the radial threshold.
No estimate for $J$ is needed for this particular cone check,
although its initial bound will be used in the subsequent
propagation. The admissible cone condition fails on the
reference interval because $\kappa=4/5<2$.

\paragraph{Changing the slope of $U^\theta$ on $[R_{\rm ref}, eR_{\rm ref}]$.}

Here, the ansatz reads
\begin{equation}\label{eq:outer-first-turn}
U^\theta=\frac{P_*}{1+Z^2}
\exp\left(\frac y{10}-\frac35\int_0^y\sigma(v)\,dv\right),
\qquad U^z=4Z,
\end{equation}
as shown in Figure~\ref{fig:outer-axial-turnoff}. Here
$a=4/5+(6/5)\sigma(y)$ increases from $4/5$ to $2$ while $b=0$,
so the relaxed cone test \eqref{eq:outer-cone-axial-zero} reduces to
$\mathcal{I}^\theta/F=RQ/L_\delta>2$.
We have $W=1-4L_\delta$ and $\zeta=-2Z/(1+Z^2)$.
The radial equation for $Q$ becomes
\begin{equation}\label{eq:outer-first-Q}
\begin{aligned}
\partial_yQ+\left(\frac85-\frac35\sigma(y)\right)Q
={}&\frac35(3-4\delta Z^2)(1-\sigma(y))
-\frac\delta2(1-8Z^2)\\
&+\frac{(9-\delta-8Z^2)Z^2}{1+Z^2}.
\end{aligned}
\end{equation}
The right-hand side is at least $-\delta/2$.
Since $Q(R_{\rm ref}, Z)>1$, integration gives
\[
Q(R,Z)\ge e^{-8/5}-\frac\delta2>0
\qquad(R_{\rm ref}\le R\le eR_{\rm ref}).
\]
On this interval $U^\theta\asymp P_*/(1+Z^2)$, so the axial
source in \eqref{eq:normalized-Q-J} is $O(|Z|)$.
Together with $|J_{\rm ref}|\lesssim|Z|$, this gives
$|J(R,Z)|\lesssim|Z|$ throughout the interval.
The condition~\eqref{eq:outer-cone-axial-zero} again gives the relaxed cone
throughout this first turn.

\paragraph{Slow axial turnoff on $[eR_{\rm ref},R_{-\mu}]$.}
The power $U^\theta\propto R^{-1/2}$ gives $a=2$, removing the
$w^2$ term from \eqref{eq:outer-cone-test}. Thus the slow turnoff
only needs to control the axial contributions $bw$ and $b^2$ there.
Put $t=y-1\in[0,T_d]$. The ansatz reads
\[
U^\theta(R,Z)=U^\theta(eR_{\rm ref},Z)e^{-t/2},
\qquad U^z(R,Z)=f(t)Z,
\]
where
\[
f(t)=4\left[1-\sigma\left(\frac{\log(1+t)}{M_d}\right)\right],
\qquad
|\partial_t f(t)|\le\frac{4\|\partial_y\sigma\|_\infty}{M_d(1+t)} = \frac{32}{M_d(1+t)}.
\]
In particular, $f=0$ on $[T_d-11,T_d]$. See Figure \ref{fig:outer-axial-turnoff} for an illustration. Here $a=2$ and $\zeta=-2Z/(1+Z^2)$, so the normalized
propagation equations become
\begin{equation}\label{eq:outer-turnoff-propagation}
\begin{aligned}
\partial_tQ+Q
&=-\frac\delta2(1-2fZ^2)
+\left(\frac{1-\delta}{2}+Lf\right)
 \frac{2Z^2}{1+Z^2},\\
\partial_tJ
&=Z+\Pi-\frac{Z}{(U^\theta)^2}
\left[W\partial_t f+f+\bigl(1-(2+\delta)Z^2\bigr)f^2\right].
\end{aligned}
\end{equation}
In the last paragraph, propagation across the preceding unit interval
$[R_{\rm ref},eR_{\rm ref}]$ has given us
\[
Q(eR_{\rm ref},Z)\ge e^{-8/5}-\delta/2\gtrsim1,
\qquad |J(eR_{\rm ref},Z)|\lesssim |Z|.
\]
The right-hand side in the first line of \eqref{eq:outer-turnoff-propagation}
is bounded below by $cZ^2-\delta/2$.
Since $P_*>e^{T_d}$ and $\delta\le d_{\rm corr}\mu
=d_{\rm corr}c_\mu P_*^{-4}$, we have $\delta\ll e^{-T_d}$.
Using $|W|\lesssim1$ and $|\Pi|\lesssim|Z|$, integration gives
\begin{equation}\label{eq:outer-turnoff-estimates}
Q\gtrsim Z^2+e^{-t},\qquad
|J|\lesssim |Z|\left(1+t+\frac{e^t}{P_*^2}\right)
\lesssim |Z|(1+t).
\end{equation}
Consequently,
\begin{equation}\label{eq:outer-slow-bw}
\begin{aligned}
    |bw| &= \left|\frac{2\partial_y U^z}{U^\theta}\right|\left|\frac{U^\theta J}{Q}\right|=\left|\frac{2Z(\partial_t f)J}{Q}\right|
\lesssim\frac1{M_d}\frac{Z^2}{Z^2+e^{-t}}
\lesssim M_d^{-1},\\
b^2 &= \left|\frac{2\partial_y U^z}{U^\theta}\right|^2 \lesssim\frac{e^{T_d}}{M_d^2P_*^2}\ll1.
\end{aligned}
\end{equation}
Taking $M_d$ sufficiently large verifies the two algebraic
inequalities in Lemma~\ref{lem:outer-cone-test}, since $a=2$.
Since $R=e^{1+t}R_{\rm ref}$ and $Q\gtrsim e^{-t}$,
$\mathcal{I}^\theta/F=RQ/L_\delta\gtrsim R_{\rm ref}$.
Thus the same lemma gives the relaxed cone for
$R_{\rm ref}\ge R_{\rm corr}$ with an absolute threshold.

At the exit $R_{-\mu}$ the axial velocity has been zero for eleven
logarithmic units. Thus $\Sigma_B(R_{-\mu})\le e^{-11}R_{-\mu}$ and
$W(R_{-\mu},Z)\ge1-4e^{-11}>3/4$.
Hence when the next interval reaches slope $-1/2-\mu$, its angular
source at $Z=0$ satisfies $\mu W-\delta/2\gtrsim\mu$.

\MergedTurnoffFigure
\end{proof}

\subsection{Interval O.3: The power-law buffer and reserved intervals}
\label{subsec:outer-buffer}

We continue the outward propagation for the uncorrected candidate,
using the cumulative moments and $P=P_{\rm pre}$ as in
Section~\ref{subsec:outer-left-data}.

\begin{lemma}\label{lem:outer-buffer-cone}
The uncorrected candidate $\bar U$ satisfies the admissible cone condition
on $(R_{-\mu},R_{mz}]$.
\end{lemma}

\begin{proof}\leavevmode\par
\paragraph{The slope transition and power law on $[R_{-\mu},R_{mz}]$.}
The axial velocity remains zero throughout this interval, while the angular velocity changes its slope, 
see Figure \ref{fig:outer-power-stage} for an illustration.
Its accumulated moment, however, need not vanish. Define
\[
m_\infty=4\int_0^{R_{-\mu}}
B\left(\log(\rho/R_{\rm ref})-1\right)\,d\rho.
\]
Thus $m_\infty/R_{\rm ref}$ depends only on $M_d$, by the
change of variables $\rho=R_{\rm ref}r$.
Since the axial cutoff is supported below $e^{-11}R_{-\mu}$,
\begin{equation}\label{eq:outer-W-after-turnoff}
M^z(R,Z)=Zm_\infty,\qquad
W(R,Z)=1-\frac{L_\delta m_\infty}{R},\qquad
0\le m_\infty\le4e^{-11}R_{-\mu}
\quad(R\ge R_{-\mu}).
\end{equation}
Thus the lower bound $W>3/4$ from the end of Interval O.2 persists
throughout the buffer.

On $[R_{-\mu},R_{+c}]$, put $t=y-y_d\in[0,1]$ and
$\nu(t)=\mu\sigma(t)$. Then
\[
\partial_y\log U^\theta=-\frac12-\nu,
\qquad a=2+2\nu, \qquad \zeta = \partial_Z \log U^\theta = - \frac{2Z}{1+Z^2}
\]
and
\begin{equation} \label{eq:QJ-nu}
\partial_tQ+(1-\nu)Q
=\nu W-\frac\delta2+\frac{(1-\delta)Z^2}{1+Z^2},
\qquad
\partial_tJ-2\nu J=Z+\Pi.    
\end{equation}
The slightly steeper power supplies the positive source
$\mu W-\delta/2\gtrsim\mu$ at $Z=0$ once $\nu=\mu$.
This is the reason for using $-1/2-\mu$, with $\delta\ll\mu$,
before the heat tail: it provides the angular stress margin used below.
More generally, the source on the ensuing power-law interval satisfies
\[
\mu W-\frac\delta2+\frac{(1-\delta)Z^2}{1+Z^2}
\gtrsim\mu+Z^2.
\]
The preceding bounds and suitable choice of constants yield
\[
Q(R_{+c},Z)\gtrsim Z^2+\mu+e^{-T_d},\qquad
|J(R_{+c},Z)|\lesssim(1+T_d)|Z|.
\]

On $[R_{+c},R_{mz}]$, put $t=\log(R/R_{+c})\in[0,T_w]$. Here
\[
U^\theta(R,Z)=U^\theta(R_{+c},Z)e^{-(1/2+\mu)t},
\qquad a=2+2\mu,
\]
and the same equations hold with $\nu=\mu$. Therefore
\begin{equation}\label{eq:outer-power-estimates}
\begin{aligned}
Q&\gtrsim Z^2+\mu+e^{-T_d}e^{-(1-\mu)t},\\
|J|&\lesssim |Z|(1+T_d+t)e^{2\mu t}.
\end{aligned}
\end{equation}
Moreover,
\begin{equation}\label{eq:outer-power-w}
\begin{aligned}
\sqrt\mu\,|w|
&=\sqrt\mu\,\frac{U^\theta|J|}{Q}\\
&\lesssim P_*\sqrt\mu(1+T_d+t)e^{\mu t/2}\ll1
\qquad(0\le t\le T_w).
\end{aligned}
\end{equation}
For this estimate, use $|Z|/(Z^2+h)\le1/(2\sqrt h)$ with
$h=e^{-T_d}e^{-(1-\mu)t}$. The factor $h^{-1/2}$ cancels
$e^{-T_d/2}e^{-(1-\mu)t/2}$ in the product of the swirl
and stress bounds, leaving $e^{\mu t/2}$.
Finally, $T_w=60\log(1/\mu)$ bounds this exponential uniformly.
With $\mu=c_\mu P_*^{-4}$ and $T_d<\log P_*$, the inequality
$(\log P_*)/P_*\le1/e$ gives
\begin{equation}\label{eq:outer-parameter-smallness}
 P_*\sqrt\mu\,(1+T_d+T_w)
 \le\sqrt{c_\mu}\left[1+60\log(1/c_\mu)+\frac{241}{e}\right].
\end{equation}
Taking $c_\mu$ sufficiently small therefore gives $(a-2)w^2=2\mu w^2<2$
uniformly in $P_*$.
On the preceding unit transition, $Q\gtrsim e^{-T_d}$,
$|w|\lesssim P_*(1+T_d)$, and $0\le a-2\le2\mu$, giving the same
strict inequality. Since $b=0$, the cone is admissible wherever
$a>2$, and relaxed at $R_{-\mu}$, for $R_{\rm ref}\ge R_{\rm corr}$.

At $R_{mz}$, integration of the explicit ansatz also gives
\begin{equation}\label{eq:outer-power-end-moments}
\begin{gathered}
\|U^\theta(R_{mz},\cdot)\|_{C^1_Z}
\lesssim P_*e^{-T_d/2}\mu^{30},
\qquad
\left\|\frac{M^z(R_{mz},\cdot)}
{R_{mz}U^\theta(R_{mz},\cdot)}\right\|_{C^1_Z}
\lesssim\mu^{29},\\
\left\|\frac{M^{z\theta}(R_{mz},\cdot)}
{R_{mz}[U^\theta(R_{mz},\cdot)]^2}\right\|_{C^1_Z}
\lesssim1+T_d+\log(1/\mu).
\end{gathered}
\end{equation}
These estimates are useful later for moment corrections.
The accumulated $M^z$ and $M^{\theta z}$ stay fixed after the axial
turnoff, whereas the natural repair scales $RU^\theta$ and
$R\sqrt{2R}(U^\theta)^2$ grow along the buffer. This makes their
relative sizes small. The three intervals $I_1,I_2,I_3$ reserved in
Interval O.3 remain unused here; in particular,
$R_c=e^{-24}R_{mz}$ and $[R_c,2R_c]\subset I_1$ are retained for the
later shear modification and its moment repair.

The same absolute threshold $R_{\rm corr}$ suffices for all cone checks
in Section~\ref{subsec:outer-left-data} and the present subsection,
through $R_{mz}$. Indeed, the preceding lower bounds give
\[
\frac{\mathcal{I}^\theta}{F}=\frac{RQ(R,Z)}{L_\delta}\gtrsim R_{\rm ref}
\qquad(R_{\rm useless}\le R\le R_{mz}).
\]
The quantity $a-bw$ has absolute positive lower and upper bounds,
and $2-2bw-b^2/a-(a-2)w^2$ has an absolute positive lower bound,
while $a$ and $b^2$ have absolute upper bounds.
Thus the thresholds in Lemma~\ref{lem:outer-cone-test} are
independent of the input parameters on this interval.
\end{proof}
\MergedBufferFigure

\subsection{Interval O.4: An axial pulse and its two end corrections}
\label{subsec:outer-axial-matching}

We construct the axial correction under the tail energy assumption
\eqref{eq:mc-pulse-future-energy} on the future swirl $U^\theta$ for $R\ge R_{0z}$.
The waiting length has been fixed in Section~\ref{subsec:outer-waiting};
the two angular corrections will be constructed in
Section~\ref{subsec:outer-angular-matching}. The result below will be
applied to that completed swirl at the end of Interval O.6.
We change only the axial velocity. Recall the radii
\[
\begin{aligned}
R_{mz}&=\mu^{-60}R_{+c},&\qquad
R_{0z}&=e^{13/\mu}R_{mz},\\
R_{m\theta}&=e^{T_f}R_{0z},&
R_{3/2}&=\mu^{-30}R_{m\theta},\qquad T_f=100.
\end{aligned}
\]
The input axial velocity satisfies $\bar U^z=0$ on $[R_{mz},\infty)$.
For this subsection, fix a swirl $U^\theta$ that is smooth and positive
for $R>0$ and agrees with $\bar U^\theta$ for $0<R\le R_{0z}$.
The angular bumps lie in $(R_{m\theta},R_{3/2})$, so the swirl to be
constructed has this property. In particular, on $[R_{mz},R_{0z}]$,
\[
U^\theta(R,Z)=\bar U^\theta(R,Z)
=U^\theta(R_{mz},Z)\left(\frac R{R_{mz}}\right)^{-1/2-\mu},
\qquad
\partial_Z\log U^\theta=-\frac{2Z}{1+Z^2}.
\]

\MergedAxialFigure

\paragraph{The main pulse and the end bumps.}
Write
\[
t=\log(R/R_{mz}),\qquad \xi=\mu t,\qquad
E(t,Z)=U^\theta(R_{mz}e^t,Z).
\]
Define the fixed pulse shape
\begin{equation}\label{eq:mc-pulse-shape}
g_p(\xi)=
\begin{cases}
\displaystyle
\bigl[1-\sigma(\xi-10)\bigr]\int_0^\xi\sigma(50v)\,dv,
&\xi\ge0,\\
0,&\xi<0.
\end{cases}
\end{equation}
It is smooth, supported in $[0,11]$, and satisfies
\[
0\le g_p(\xi)\le\xi,\qquad \partial_\xi g_p(\xi)\le1
\qquad(\xi\ge0).
\]
Its derivatives of each fixed order are bounded by fixed constants.
For the two end bumps, use the smooth nonnegative bump
$\beta=\beta_{3/20}$ of unit integral, supported in $[-3/20,3/20]$.
Set
\[
\gamma_1(t)=\beta\left(t-\frac{13}{\mu}+3\right),
\qquad
\gamma_2(t)=\beta\left(t-\frac{13}{\mu}+1\right).
\]
They are nonzero precisely on the respective radial intervals
\[
\left(e^{-3-3/20}R_{0z},e^{-3+3/20}R_{0z}\right),
\qquad
\left(e^{-1-3/20}R_{0z},e^{-1+3/20}R_{0z}\right).
\]
With $\mu_{\rm corr}\le1/60$, the main pulse ends at $t=11/\mu$,
before the first end bump starts at $t=13/\mu-3-3/20$.
The two end bumps have disjoint supports as well.

For $0\le t\le13/\mu$, set
\begin{equation}\label{eq:mc-axial-pulse}
U^z=E B_{\mathrm{axial}},\qquad
B_{\mathrm{axial}}(t,Z)=a_p(Z)g_p(\mu t)
+c_1(Z)\gamma_1(t)+c_2(Z)\gamma_2(t),
\end{equation}
and retain $\bar U^z$ outside $[R_{mz},R_{0z}]$.
In particular, $U^z=0$ after $R_{0z}$.

The logarithmic length is of order $1/\mu$, so a pulse with
$U^z$ comparable to $U^\theta$ supplies energy of order
$R_{mz}[U^\theta(R_{mz},Z)]^2/\mu$.
The amplitude $a_p$ supplies the axial energy needed to make
$M^{z\theta}(\infty,Z)=0$.
The two exponentially small coefficients $c_1,c_2$ cancel the linear
moments $M^z$ and $M^{\theta z}$.
They are placed near $R_{0z}$ because these two moments give much
larger weights to later radii, whereas the energy weight changes
only on the longer scale $1/\mu$.
On the main pulse, the slow variation of the ratio $B_{\mathrm{axial}}=U^z/U^\theta$ controls the
stress direction even where $U^z$ and $U^\theta$ have comparable size.

\paragraph{The incoming data.}
The preceding power-law estimates give
\begin{equation}\label{eq:mc-pulse-data}
\begin{aligned}
\|U^\theta(R_{mz},\cdot)\|_{C_Z^1}
&\le C P_*e^{-T_d/2}\mu^{30},\\
\left\|\frac{\bar M^z(R_{mz},\cdot)}
{R_{mz}U^\theta(R_{mz},\cdot)}\right\|_{C_Z^1}
+\left\|\frac{\bar M^{\theta z}(R_{mz},\cdot)}
{R_{mz}\sqrt{2R_{mz}}[U^\theta(R_{mz},\cdot)]^2}\right\|_{C_Z^1}
&\le C\mu^{29},\\
\left\|\frac{\bar M^{z\theta}(R_{mz},\cdot)}
{R_{mz}[U^\theta(R_{mz},\cdot)]^2}\right\|_{C_Z^1}
&\le C[1+T_d+\log(1/\mu)].
\end{aligned}
\end{equation}
For the mixed linear moment, integrate
$\sqrt{2R}\,\bar U^\theta\bar U^z$ only over the support
of the original axial velocity, below $R_{-\mu}$.
The estimate $\sqrt R\,\bar U^\theta\lesssim
P_*\sqrt{R_{\rm ref}}$ there, together with
$R_{mz}/R_{+c}=\mu^{-60}$, gives the second bound.
These incoming moments are unchanged by the angular bumps.

\begin{proposition}[Axial correction for a prescribed swirl]
\label{prop:outer-conditional-axial-correction}
Suppose that the prescribed swirl has the properties above, and that
its normalized tail energy is smooth in $Z$ and satisfies
\begin{equation}\label{eq:mc-pulse-future-energy}
\left\|
\frac{\displaystyle\int_{R_{0z}}^\infty[U^\theta(\rho,\cdot)]^2\,d\rho}
{R_{mz}[U^\theta(R_{mz},\cdot)]^2}
\right\|_{C_Z^1}
\le C[1+\log(1/\mu)].
\end{equation}
Here $C$ depends only on the fixed auxiliary data and is fixed before
$\mu_{\rm corr}$ is chosen.
With the incoming bounds \eqref{eq:mc-pulse-data} and the fixed
smallness choices, there exist smooth coefficients $a_p,c_1,c_2$
in \eqref{eq:mc-axial-pulse} such that
\[
0.9<a_p<1.2,\qquad
\|c_1\|_{C_Z^1}+\|c_2\|_{C_Z^1}\le C e^{-1/(2\mu)},
\]
and
\[
M^z(R_{0z},Z)=M^{\theta z}(R_{0z},Z)=M^{z\theta}(\infty,Z)=0.
\]
The angular velocity and pressure are unchanged. If the pressure of
the prescribed swirl satisfies \eqref{eq:outer-pressure-sharp} on
$[R_{mz},R_{0z}]$, the corrected profile satisfies the admissible
cone there, with $Q(R,Z)\ge c(\mu+Z^2)$.
\end{proposition}

The integral in \eqref{eq:mc-pulse-future-energy} uses the completed
swirl, including the two angular bumps. Interval O.6 verifies both
this bound and the pressure estimate for that swirl.

\begin{proof}
We first solve the two linear moment equations for $c_1,c_2$ in terms of
a trial amplitude, then determine $a_p$ from the energy balance and verify
the cone condition using the angular and axial stress estimates.

\paragraph{The two linear equations.}
For a scalar trial amplitude $a\in[0.9,1.2]$, first solve
\[
M^z(R_{0z},Z)=M^{\theta z}(R_{0z},Z)=0
\]
for $c_1,c_2$. Define
\[
\begin{aligned}
m_1(Z)&=\frac{\bar M^z(R_{mz},Z)}
{R_{mz}U^\theta(R_{mz},Z)},\\
m_2(Z)&=\frac{\bar M^{\theta z}(R_{mz},Z)}
{R_{mz}\sqrt{2R_{mz}}[U^\theta(R_{mz},Z)]^2},
\end{aligned}
\qquad
\lambda_1=\frac12-\mu,\quad
\lambda_2=\frac12-2\mu.
\]
Changing variables to $t$, the equations are exactly
\begin{equation}\label{eq:mc-axial-linear-equations}
\sum_{j=1}^2c_j\int_0^{13/\mu}e^{\lambda_i t}\gamma_j(t)\,dt
=-m_i(Z)-a\int_0^{13/\mu}e^{\lambda_i t}g_p(\mu t)\,dt,
\qquad i=1,2.
\end{equation}
Thus the exponents in the two linear weights differ by $\mu$.

Divide row $i$ by $e^{13\lambda_i/\mu}$.
The resulting matrix has entries
\[
\begin{aligned}
A_{i1}&=e^{-3\lambda_i}
\int_{\mathbb R}e^{\lambda_i s}\beta(s)\,ds,\\
A_{i2}&=e^{-\lambda_i}
\int_{\mathbb R}e^{\lambda_i s}\beta(s)\,ds.
\end{aligned}
\]
On the support of $\beta$, the exponential factors lie between
fixed positive constants. In particular,
\[
\frac12\le\int e^{\lambda_i s}\beta(s)\,ds\le2.
\]
The determinant is explicitly
\[
\det A
=-2e^{-2+6\mu}\sinh\mu
\prod_{i=1}^2\int e^{\lambda_i s}\beta(s)\,ds.
\]
Consequently,
\begin{equation}\label{eq:mc-axial-matrix-bounds}
|\det A|\asymp\mu,\qquad \|A^{-1}\|\le C\mu^{-1}.
\end{equation}

The main pulse ends at $11/\mu$, so
\[
e^{-13\lambda_i/\mu}
\int_0^{13/\mu}e^{\lambda_i t}g_p(\mu t)\,dt
\le C e^{-2\lambda_i/\mu}
\le C e^{-1/\mu}.
\]
The normalized incoming term satisfies the same bound in $C_Z^1$,
by \eqref{eq:mc-pulse-data}.
Thus the linear system \eqref{eq:mc-axial-linear-equations} has a unique solution $c_j=c_j(a,Z)$,
with
\begin{equation}\label{eq:mc-axial-end-coefficients}
\sup_{0.9\le a\le1.2}
\left(
\|c_j(a,\cdot)\|_{C_Z^1}
+\|\partial_a c_j(a,\cdot)\|_{C_Z^1}
\right)
\le C\mu^{-1}e^{-1/\mu}
\le C e^{-1/(2\mu)}.
\end{equation}
Here $Z$ derivatives are taken with $a$ held fixed.
Since the matrix is independent of $Z$ and the right-hand side
is affine in $a$, we can write
\[
c_j(a,Z)=u_j(Z)+v_j a,
\]
where $v_j$ depends only on $\mu$ and the fixed pulse shapes.

\paragraph{The quadratic equation and the choice of $a_p$.}
Define
\[
K_p=\int_0^{13}e^{-2\xi}g_p(\xi)^2\,d\xi,
\qquad
K_j(\mu)=\int_0^{13/\mu}e^{-2\mu t}\gamma_j(t)^2\,dt.
\]
The first number is a fixed constant.
The other two are derived coefficients. More explicitly,
\[
\begin{aligned}
K_1(\mu)&=e^{-26+6\mu}
\int e^{-2\mu s}\beta(s)^2\,ds,\\
K_2(\mu)&=e^{-26+2\mu}
\int e^{-2\mu s}\beta(s)^2\,ds.
\end{aligned}
\]
In particular,
\[
\frac12e^{-26}\|\beta\|_{L^2}^2
\le K_j(\mu)\le2e^{-26}\|\beta\|_{L^2}^2.
\]
For $K_p$, symmetry of $\sigma$ gives
$g_p(\xi)=\xi-1/100$ on $[1/50,10]$.
Integrating there bounds $K_p$ below, while using
$g_p\le\xi$ on $[0,1/50]$ and
$g_p\le\xi-1/100$ afterwards bounds it above.
These elementary integrals show that the fixed number
\[
a_0=\left(\frac{1-e^{-26}}{4K_p}\right)^{1/2}
\]
satisfies $1<a_0<1.1$.

The three supports are disjoint, so
\[
B_{\mathrm{axial}}^2=a^2g_p(\mu t)^2+
c_1(a,Z)^2\gamma_1(t)^2+c_2(a,Z)^2\gamma_2(t)^2.
\]
Also
\[
R[U^\theta(R,Z)]^2
=R_{mz}[U^\theta(R_{mz},Z)]^2e^{-2\mu t}.
\]
It follows that $M^{z\theta}(\infty,Z)=0$ is exactly the
quadratic equation
\begin{equation}\label{eq:mc-pulse-balance}
\begin{aligned}
K_pa^2+\mu\sum_{j=1}^2K_j(\mu)c_j(a,Z)^2
={}&\frac{1-e^{-26}}4
-\frac{\mu\bar M^{z\theta}(R_{mz},Z)}
{R_{mz}[U^\theta(R_{mz},Z)]^2}\\
&+\frac{\mu}{2R_{mz}[U^\theta(R_{mz},Z)]^2}
\int_{R_{0z}}^\infty[U^\theta(\rho,Z)]^2\,d\rho.
\end{aligned}
\end{equation}
All quantities on the right are fixed by the prescribed swirl and
the incoming moments.
In particular, the future energy includes the two angular bumps.

Substituting $c_j=u_j+v_ja$ makes the coefficients of $a^2$
and $a$ on the left explicit:
\[
K_p+\mu\sum_jK_jv_j^2,\qquad
2\mu\sum_jK_ju_jv_j.
\]
The first differs from $K_p$ by at most $C\mu e^{-1/\mu}$.
The second, and the constant term $\mu\sum_jK_ju_j^2$,
have $C_Z^1$ norm at most $C\mu e^{-1/\mu}$.

The future-energy bound \eqref{eq:mc-pulse-future-energy}, together
with \eqref{eq:mc-pulse-data}, shows that the right
side of \eqref{eq:mc-pulse-balance} is
\[
\frac{1-e^{-26}}4+
O_{C_Z^1}\bigl(\mu[1+T_d+\log(1/\mu)]\bigr).
\]

Since $1<a_0<1.1$, for sufficiently small $\mu_{\rm corr}$
the left side minus the right side is negative at $a=0.9$
and positive at $a=1.1$.
Its derivative on this interval is
\[
2K_pa+2\mu\sum_jK_jc_j(a,Z)v_j\ge0.3
\]
by the fixed smallness choice of $\mu_{\rm corr}$.
There is therefore a unique root $a=a_p(Z)$ in this interval.
The equation and its $Z$ derivative give
\begin{equation}\label{eq:mc-pulse-amplitude}
0.9<a_p<1.1,\qquad
\|a_p-a_0\|_{C_Z^1}
\le C\mu[1+T_d+\log(1/\mu)].
\end{equation}
Smoothness follows from the nonvanishing derivative in $a$.
After composition with $a_p(Z)$, the coefficients $c_j$ retain
the bound \eqref{eq:mc-axial-end-coefficients}.

The three axial moment conditions are now exact.
For every $R\ge R_{0z}$,
\begin{equation}\label{eq:mc-backward-moments}
\begin{gathered}
U^z=U^r=M^z=M^{\theta z}=0,\\
M^{z\theta}(R,Z)
=\frac12\int_R^\infty[U^\theta(\rho,Z)]^2\,d\rho.
\end{gathered}
\end{equation}
The identity for $U^r$ follows from the reconstruction formula
using $U^z=M^z=0$.
Changing only $U^z$ does not change the pressure. The radius $R_{0z}$
is the intermediate matching point for these three moments; their
identities propagate to the common matching point $R_{3/2}$ and beyond.

\paragraph{The angular stress on the pulse.}
Set
\[
 q_0(Z)=\mu-\frac\delta2+\frac{(1-\delta)Z^2}{1+Z^2}.
\]
At $R_{mz}$ the incoming angular stress is fixed by the velocities
through that radius, independently of the pressure. Indeed, on the
preceding power-law interval, \eqref{eq:QJ-nu} and
\eqref{eq:outer-W-after-turnoff} give
\[
 \partial_yQ+(1-\mu)Q=q_0+\mu(W-1),\qquad
 |W-1|\le C R_{+c}/R .
\]
The initial $Q(R_{+c},Z)$ is bounded by the angular equation on
Intervals O.2--O.3. Integrating through $T_w=60\log(1/\mu)$ gives
$Q(R_{mz},Z)=q_0/(1-\mu)+O(e^{-(1-\mu)T_w})
=q_0/(1-\mu)+O(\mu^{59})$.

On the pulse $a=2+2\mu$ and $\zeta=-2Z/(1+Z^2)$.
By \eqref{eq:mc-pulse-data} and the bounds on $B_{\mathrm{axial}}$ and $\partial_Z B_{\mathrm{axial}}$,
integrating $\partial_R M^z=E B_{\mathrm{axial}}$ gives $|W-1|\le C E$. Indeed,
the integral of the power law $U^\theta$ from $R_{mz}$ to $R$
is bounded by $C R E$, also after one $Z$ derivative. Thus the source in
\eqref{eq:normalized-Q-J} is $q_0+O(E)$, with damping
$1-\mu\ge1/2$. Integration gives
\begin{equation}\label{eq:mc-pulse-angular-stress}
 Q=\frac{q_0}{1-\mu}
 +O\bigl(P_*e^{-T_d/2}\mu^{30}+\mu^{59}\bigr)
 \ge c(\mu+Z^2)
 \qquad (R_{mz}\le R\le R_{0z}).
\end{equation}

\paragraph{The axial stress direction on the pulse.}
We now use the pressure hypothesis to verify the cone. The ratio
$B_{\mathrm{axial}}=U^z/U^\theta$ can be of order one, so we estimate the two cone
expressions directly.
Put
\[
C(t,Z)=\frac{M^z(R_{mz}e^t,Z)}{R_{mz}e^t E(t,Z)},
\qquad \partial_\xi B_{\mathrm{axial}}=\mu^{-1}\partial_t B_{\mathrm{axial}}.
\]
The moment equation gives
\[
\partial_tC+\lambda_1C=B_{\mathrm{axial}},\qquad \lambda_1=\frac12-\mu.
\]
Two integrations by parts in its exponential convolution yield
\begin{equation}\label{eq:mc-pulse-average}
C=\frac{B_{\mathrm{axial}}}{\lambda_1}
-\frac{\mu\partial_\xi B_{\mathrm{axial}}}{\lambda_1^2}
+C(0,Z)e^{-\lambda_1t}+O_{C_Z^1}(\mu^2).
\end{equation}
The main pulse is flat at $t=0$ and both end bumps vanish there, so
$B_{\mathrm{axial}}(0,Z)=\partial_t B_{\mathrm{axial}}(0,Z)=0$. The remainder is
$\lambda_1^{-2}\int_0^t e^{-\lambda_1(t-s)}\partial_s^2 B_{\mathrm{axial}}(s,Z)\,ds$.
In \eqref{eq:mc-axial-pulse}, two derivatives on $g_p(\mu s)$ supply
$\mu^2$, and the end bumps are translates of a fixed function. Thus
\eqref{eq:mc-axial-end-coefficients} and \eqref{eq:mc-pulse-amplitude}
give $\sup_s\|\partial_s^2 B_{\mathrm{axial}}(s,\cdot)\|_{C_Z^1}
\le C(\mu^2+e^{-1/(2\mu)})$. Since $\lambda_1\ge29/60$ and the
kernel is independent of $Z$, the stated bound holds throughout the pulse.

The linear axial terms in \eqref{Iz-moments} give
\begin{equation}\label{eq:mc-pulse-axial-stress}
\frac{L_\delta \mathcal{I}^z}{\sqrt{R/2}\,E}
=-B_{\mathrm{axial}}+\frac{1-\delta}{2}
\bigl[(1-Z\zeta)C-Z\partial_Z C\bigr]+O(E/\mu).
\end{equation}
For the remaining terms we use the energy cancellation
$M^{z\theta}(\infty,Z)=0$ just proved, which gives
\[
M^{z\theta}(R,Z)=
\int_R^\infty\left[\frac12(U^\theta)^2-(U^z)^2\right]\,d\rho.
\]
Up to $R_{0z}$ the power law gives
$\int_R^{R_{0z}}(U^\theta)^2\,d\rho\le RE^2/(2\mu)$.
Beyond $R_{0z}$, \eqref{eq:mc-pulse-future-energy} controls the
swirl energy since $RE^2\ge e^{-26}R_{mz}[U^\theta(R_{mz},Z)]^2$.
Both estimates hold after a $Z$ derivative. With $B_{\mathrm{axial}},\partial_Z B_{\mathrm{axial}}$
bounded this gives
\[
\frac{|M^{z\theta}|+|\partial_ZM^{z\theta}|}{R E^2}
\le C\mu^{-1}.
\]
By \eqref{eq:outer-pressure-sharp}, $|P|+|\partial_ZP|\le C E^2$.
These bounds and $|W-1|\le C E$ give the error in
\eqref{eq:mc-pulse-axial-stress}.

Insert \eqref{eq:mc-pulse-average} in \eqref{eq:mc-pulse-axial-stress}.
The coefficient of $B_{\mathrm{axial}}$ is
$-1+(1-\delta)(1-Z\zeta)/(2\lambda_1)=q_0/\lambda_1$.
Dividing by $Q$ from \eqref{eq:mc-pulse-angular-stress} and using
$b=2\partial_t(E B_{\mathrm{axial}})/E$ gives
\begin{equation}\label{eq:mc-pulse-directions}
a=2+2\mu,\qquad
b=-B_{\mathrm{axial}}+O(\mu),\qquad
w=2B_{\mathrm{axial}}-K(Z)\partial_\xi B_{\mathrm{axial}}+O(\mu^{1/4}),
\end{equation}
where
\[
K(Z)=
\frac{\mu(1-\delta)
\left(\frac12+\frac{Z^2}{1+Z^2}\right)(1-\mu)}
{(\frac12-\mu)^2q_0(Z)}.
\]
For clarity, the error before simplification is bounded by
\[
C\left[
\sqrt\mu\,[1+T_d+\log(1/\mu)]
+P_*e^{-T_d/2}\mu^{28}+\mu
\right].
\]
The first term comes from differentiating $a_p$ and using
$|Z|/(\mu+Z^2)\le C\mu^{-1/2}$.
The displayed bound is $O(\mu^{1/4})$, since $T_d,c_\mu$
are fixed and $P_*=(c_\mu/\mu)^{1/4}$.

The function $K$ decreases with $Z^2/(1+Z^2)$.
Taking $d_{\rm corr}\le1/100$ and $\mu_{\rm corr}\le1/100$ in the initial
parameter choices gives
\[
0\le K(Z)\le K(0)
=\frac{2(1-\delta)(1-\mu)}
{(1-2\mu)^2(1-\delta/(2\mu))}
\le2.1.
\]
By \eqref{eq:mc-pulse-amplitude} and $\partial_\xi g_p\le1$,
apart from the exponentially small end corrections,
$B_{\mathrm{axial}}\ge0$ and $\partial_\xi B_{\mathrm{axial}}\le1.1$.
All quantities in \eqref{eq:mc-pulse-directions} are bounded by
fixed constants. Therefore
\[
\begin{aligned}
bw
&\le\sup_{s\ge0}(-2s^2+2.4s)+O(\mu^{1/4})
<0.8,\\
2bw+\frac{b^2}{a}+(a-2)w^2
&\le\sup_{s\ge0}(-3.5s^2+4.8s)+O(\mu^{1/4})
<1.8.
\end{aligned}
\]
Indeed, $2.1\cdot1.1<2.4$, and completing the square gives maxima
strictly below $0.8$ and $1.7$ for the two quadratics.
Their fixed gaps allow the displayed error terms to be absorbed
by the choice of $\mu_{\rm corr}$.
A negative derivative on the falling side improves the upper bounds for both leading terms.
Hence $a-bw>1.2$, and the second cone inequality retains
a fixed positive margin.
Also
\[
\frac{RQ}{L_\delta}\ge cR_{mz}\mu
=c e^{T_d+2}\mu^{-59}R_{\rm ref}\ge cR_{\rm ref}.
\]
Lemma~\ref{lem:outer-cone-test} therefore applies with the fixed threshold
$R_{\rm corr}$. Since $a>2$, the pulse satisfies the admissible cone.
\end{proof}

\subsection{Intervals O.5--O.6: Removing \texorpdfstring{$Z$}{Z}-dependence and matching moments}
\label{subsec:outer-remove-z}
\label{subsec:outer-angular-matching}\label{sec:outer-moment-corrections}

The waiting length has been fixed in Section~\ref{subsec:outer-waiting}.
On Interval O.5 the background swirl becomes independent of $Z$.
We use two angular bumps on Interval O.6 to match the angular and pressure
moments at $R_{3/2}$ with the axis pressure fixed at $P_0^{\rm pre}$.
We then use the completed swirl to apply
Proposition~\ref{prop:outer-conditional-axial-correction} and verify
the cone on O.5--O.6 for the resulting fields.

\paragraph{Removing the background $Z$-dependence on Interval O.5.}
On $[R_{0z},R_{m\theta}]$, put $t=y-y_v\in[0,100]$.
The background swirl is
\[
\bar U^\theta(R,Z)
=\bar U^\theta(R_{0z},Z)e^{-(1/2+\mu)t}
\left(\frac{1+Z^2}{2}\right)^{\sigma(t/100)}.
\]
The angular bumps are supported beyond $R_{m\theta}$, so the completed
swirl has the same formula on O.5. Consequently,
\[
\begin{aligned}
a&=2+2\mu-\frac{2(\partial_y\sigma)(t/100)}{100}
\log\left(\frac{1+Z^2}{2}\right),\\
\zeta&=-\frac{2Z}{1+Z^2}[1-\sigma(t/100)].
\end{aligned}
\]
In particular, $2+2\mu\le a<3$ and $Z\zeta\le0$ for small $\mu$.
At $R_{m\theta}$ the factor $(1+Z^2)^{-1}$ has become $1/2$.
The outgoing background swirl is independent of $Z$ and retains slope
$-1/2-\mu$ up to $R_{3/2}$, as in
Figure~\ref{fig:outer-remove-Z-stage}. This permits one waiting length
$\tau$ to match the main pre-heat angular moment at $Z=0$ by
\eqref{eq:mc-wait-choice}. The earlier contribution to the angular
moment still depends on $Z$ and is corrected on Interval O.6.

\MergedRemoveZFigure

\paragraph{The remaining heat and angular mismatches.}
The final targets still use the actual heat tail. At the selected
$\tau$, put $E_{\rm rel}=\bar U^\theta(R_{3/2})$ and define
\begin{equation}\label{eq:mc-angular-rhs-definition}
s_\theta(Z)=\frac{M_{\mathrm{right}}^\theta(R_{3/2},Z)
                 -\bar M_{\mathrm{left}}^\theta(R_{3/2},Z)}
 {R_{3/2}\sqrt{2R_{3/2}}E_{\rm rel}}.
\end{equation}
Split $s_\theta=s_{\theta,\rm pre}+s_{\theta,H}$, where $s_{\theta,\rm pre}$ uses
$M_{\mathrm{right,pre}}^\theta$ instead of $M_{\mathrm{right}}^\theta$.
The waiting-length equation gives $s_{\theta,\rm pre}(0)=0$, so the
remaining discrepancy at $Z=0$ comes entirely from the heat tail.
The incoming estimate \eqref{eq:mc-incoming-angular-moment}
gives $\|s_{\theta,\rm pre}\|_{C_Z^k}\le C_k\mu^{29}$.
The heat contribution is exactly
\[
s_{\theta,H}(Z)=\frac{\displaystyle\int_{R_{2h}}^\infty
 \sqrt{2R}\,[\bar U^\theta_{\rm pre}(R,Z)-\bar U^\theta(R,Z)]\,dR}
 {R_{3/2}\sqrt{2R_{3/2}}E_{\rm rel}}.
\]
Its integrand is nonnegative. The heat-factor bound
$\|H_\delta(2L/R)-1\|_{C_Z^k}\le C_k\delta/R$ and
\eqref{eq:mc-exterior-amplitude} imply
\[
\|s_{\theta,H}\|_{C_Z^k}
\le\frac{C_kc_\infty\delta
 \int_{R_{2h}}^\infty R^{-1-\delta/2}\,dR}
 {R_{3/2}\sqrt{2R_{3/2}}E_{\rm rel}}
\le\frac{C_k\delta^{4-2\delta}}{R_{3/2}}.
\]
Both contributions are even in $Z$. Since $R_{3/2}^{-1}\le\mu^{90}$,
their sum satisfies
\begin{equation}\label{eq:mc-angular-rhs}
\|s_\theta\|_{C_Z^k}\le C_k\mu^{29}.
\end{equation}

\paragraph{The pressure moment discrepancy.}
By \eqref{eq:outer-heat-pressure-defect}, the required increase in the
total pressure moment is $\Delta_H$. In units of $E_{\rm rel}^2$, it is
\begin{equation}\label{eq:mc-preheat-pressure-rhs}
s_H(Z)=\frac{\Delta_H(Z)}{E_{\rm rel}^2}
=\frac{M_{\rm pre}^p(\infty,Z)-\bar M_{\mathrm{left}}^p(\infty,Z)}
       {E_{\rm rel}^2}.
\end{equation}
This function is smooth, even, and nonnegative. It vanishes
at $Z=\pm1$ but need not vanish at zero. The heat estimates
and \eqref{eq:mc-tail-scales} give
\begin{equation}\label{eq:mc-pressure-rhs-bounds}
\|s_H\|_{C_Z^k}
\le C_k\frac{\delta U_{\rm tail}^2}
                   {R_{2h}E_{\rm rel}^2}
\le \frac{C_k\delta^{17}}{R_{3/2}}.
\end{equation}
All constants are uniform for $\tau\ge0$. The power $17$
comes from $U_{\rm tail}/E_{\rm rel}\lesssim\delta^6$ and
$R_{2h}\gtrsim\delta^{-4}R_{3/2}$.

\paragraph{The two angular bumps.}
\MergedAngularFigure

Use the fixed nonnegative smooth bump $\beta$, supported in
$[-3/20,3/20]$ and normalized by $\int\beta=1$.
Put
\[
t=\log(R/R_{3/2}),\qquad
\beta_1(t)=\beta(t+3),\qquad
\beta_2(t)=\beta(t+1),
\]
and set
\begin{equation}\label{eq:mc-angular-bumps}
h(t,Z)=d_1(Z)\beta_1(t)+d_2(Z)\beta_2(t),
\qquad
U^\theta=\bar U^\theta(1+h),\qquad U^z=\bar U^z.
\end{equation}
The two supports are disjoint and lie in $(R_{m\theta},R_{3/2})$.
On this interval,
\[
\bar U^\theta(R)
=\bar U^\theta(R_{3/2})e^{-(1/2+\mu)t},
\qquad
\bar U^z=0.
\]
We choose $d_1,d_2$ to match the angular moment and restore
the pre-heat pressure moment:
\[
M_{\mathrm{left}}^\theta(R_{3/2},Z)=M_{\mathrm{right}}^\theta(R_{3/2},Z),
\qquad
M_{\mathrm{left}}^p(R_{3/2},Z)=\bar M_{\mathrm{left}}^p(R_{3/2},Z)+\Delta_H(Z).
\]
The right side of the second identity is equal to the prescribed pressure
moment $M_{\mathrm{right}}^p(R_{3/2},Z)$ in \eqref{eq:mc-required-moments}.
The future heat swirl is unchanged.

The angular moment weights the later bump more strongly,
whereas the pressure moment weights the earlier bump more strongly.
This gives independent control of the two moment discrepancies.

\paragraph{The coefficient equations.}
Changing variables from $R$ to $t$ gives
\[
\begin{aligned}
\frac{M_{\mathrm{left}}^\theta(R_{3/2},Z)-\bar M_{\mathrm{left}}^\theta(R_{3/2},Z)}
{R_{3/2}\sqrt{2R_{3/2}}\bar U^\theta(R_{3/2})}
&=\int_{-4}^0e^{(1-\mu)t}h(t,Z)\,dt,\\
\frac{M_{\mathrm{left}}^p(R_{3/2},Z)-\bar M_{\mathrm{left}}^p(R_{3/2},Z)}
{[\bar U^\theta(R_{3/2})]^2}
&=\int_{-4}^0e^{-(1+2\mu)t}
\left(h+\frac12h^2\right)\,dt.
\end{aligned}
\]
Define the derived weights
\[
\begin{aligned}
A_\mu&=\int_{-3/20}^{3/20}e^{(1-\mu)s}\beta(s)\,ds,\\
B_\mu&=\int_{-3/20}^{3/20}e^{-(1+2\mu)s}\beta(s)\,ds,\\
D_\mu&=\int_{-3/20}^{3/20}e^{-(1+2\mu)s}\beta(s)^2\,ds.
\end{aligned}
\]
Translation of the bumps gives the explicit equations
\begin{equation}\label{eq:mc-angular-system-unscaled}
\begin{gathered}
A_\mu\left(e^{-3(1-\mu)}d_1+e^{-(1-\mu)}d_2\right)=s_\theta(Z),\\
\begin{aligned}
&B_\mu\left(e^{3(1+2\mu)}d_1+e^{1+2\mu}d_2\right)\\
&\qquad+\frac{D_\mu}{2}
\left(e^{3(1+2\mu)}d_1^2+e^{1+2\mu}d_2^2\right)=s_H(Z).
\end{aligned}
\end{gathered}
\end{equation}
Since the bumps $\beta_1$ and $\beta_2$ have disjoint supports, the
coefficient of $d_1d_2$ in the second equation vanishes.

For $0\le\mu\le1/60$, both exponential weights on the support
of $\beta$ lie between $4/5$ and $6/5$.
Hence
\begin{equation}\label{eq:mc-bump-weight-bounds}
\begin{gathered}
\frac45\le A_\mu,B_\mu\le\frac65,\\
\frac45\|\beta\|_{L^2}^2
\le D_\mu\le\frac65\|\beta\|_{L^2}^2,\\
\frac13\|\beta\|_{L^2}^2
\le\frac{D_\mu}{2B_\mu}
\le\frac34\|\beta\|_{L^2}^2.
\end{gathered}
\end{equation}
These bounds are independent of all input parameters.

\paragraph{Solving the equations.}\leavevmode\par
Divide the first equation by $A_\mu e^{-(1-\mu)}$
and the second by $B_\mu e^{3(1+2\mu)}$.
The resulting system is
\begin{equation}\label{eq:mc-angular-system}
\begin{aligned}
e^{-2(1-\mu)}d_1+d_2
&=\frac{e^{1-\mu}}{A_\mu}s_\theta(Z),\\
d_1+e^{-2(1+2\mu)}d_2
&=\frac{e^{-3(1+2\mu)}}{B_\mu}s_H(Z)
-\frac{D_\mu}{2B_\mu}
\left(d_1^2+e^{-2(1+2\mu)}d_2^2\right).
\end{aligned}
\end{equation}
The determinant of its linear part is
\[
e^{-4-2\mu}-1,
\qquad
|e^{-4-2\mu}-1|\ge1-e^{-4}.
\]
Its inverse has maximum-row-sum norm less than $2$,
while the quadratic coefficients are bounded by
$\frac34\|\beta\|_{L^2}^2$.

Put $D_{\rm rhs}=\|s_\theta\|_{C_Z^0}+\|s_H\|_{C_Z^0}$.
By \eqref{eq:mc-angular-rhs} and \eqref{eq:mc-pressure-rhs-bounds},
$D_{\rm rhs}\le C\mu^{29}$ with $C$ depending only on the fixed auxiliary data.
We take $\mu_{\rm corr}$ small enough that for every $\mu\le\mu_{\rm corr}$,
\[
D_{\rm rhs}
\le\frac1{1000(1+\|\beta\|_{L^2}^2)}.
\]

Applying the linear inverse gives a fixed-point map on
\[
\max\{|d_1|,|d_2|\}\le10D_{\rm rhs}.
\]
Its linear contribution has size at most $7D_{\rm rhs}$.
On this ball, the nonlinear contribution has size at most
\[
200\|\beta\|_{L^2}^2D_{\rm rhs}^2,
\]
and its Lipschitz constant is at most
\[
40\|\beta\|_{L^2}^2D_{\rm rhs}<\frac1{20}.
\]
Thus the map preserves the ball and is a contraction.

There is consequently a unique small solution. The same argument at each
$Z$ bounds it in terms of the two pointwise discrepancies. The contraction
bound gives a uniformly invertible linearization, so differentiation yields,
for each fixed $k\ge0$,
\begin{equation}\label{eq:mc-angular-coefficients}
\begin{aligned}
|d_1(Z)|+|d_2(Z)|&\le20(|s_\theta(Z)|+|s_H(Z)|),\\
\|d_1\|_{C_Z^k}+\|d_2\|_{C_Z^k}
&\le C_k\mu^{29}.
\end{aligned}
\end{equation}
More explicitly, inversion of the linear part gives
\[
\begin{aligned}
d_1
&=\frac{e^{-3(1+2\mu)}s_H/B_\mu-e^{-1-5\mu}s_\theta/A_\mu}
{1-e^{-4-2\mu}}
+O_{C_Z^k}(\mu^{58}),\\
d_2
&=\frac{e^{1-\mu}s_\theta/A_\mu-e^{-5-4\mu}s_H/B_\mu}
{1-e^{-4-2\mu}}
+O_{C_Z^k}(\mu^{58}).
\end{aligned}
\]
All constants in these estimates depend only on the fixed
auxiliary data, $\beta$, and the indicated derivative order.

Both right sides are even, so uniqueness makes each $d_j$ even.
Thus $d_j'(0)=0$, and the $C_Z^2$ bound in \eqref{eq:mc-angular-coefficients} gives
\[
|d_j'(Z)|\le |Z|\,\|d_j''\|_\infty\le C\mu^{29}|Z|,
\qquad j=1,2.
\]

\paragraph{The corrected swirl and pressure.}
The coefficient bounds give
\[
\|\partial_t^j h\|_{C_Z^k}\le C_{j,k}\mu^{29}.
\]
Since the bumps have disjoint supports,
$\|h\|_\infty\le 10\|\beta\|_\infty D_{\rm rhs}\le 1/2$
for sufficiently small $\mu_{\rm corr}$. Thus $1+h\ge1/2$.

Since $a=1-2\partial_t\log U^\theta$, we have
\[
a-(2+2\mu)=-\frac{2\partial_t h}{1+h},
\qquad
\partial_Z\log U^\theta=\frac{\partial_Zh}{1+h}.
\]
Using $\partial_t h=\sum_{j=1}^2 d_j\beta_j'$ and the bound for $d_j'(Z)$ gives
\begin{equation}\label{eq:mc-angular-perturbation-bounds}
|a-(2+2\mu)|\le C\mu^{29},
\qquad
|\partial_Z\log U^\theta|\le C\mu^{29}|Z|.
\end{equation}
Taking $C\mu_{\rm corr}^{28}\le 1/8$ in these two bounds gives
$a\ge 2+\mu$ and $|\partial_Z\log U^\theta|\le\mu |Z|/8$
on the angular supports.
At fixed axis pressure, the change in $P$ equals the change in the
cumulative pressure moment. Integration of $(2h+h^2)(\bar U^\theta)^2/(2R)$
over the bumps up to radius $R$ gives
\[
\|M_{\mathrm{left}}^p(R,\cdot)-\bar M_{\mathrm{left}}^p(R,\cdot)\|_{C_Z^k}
\le C_k[\bar U^\theta(R_{3/2})]^2\mu^{29}.
\]
Before the first bump this change is zero, so $P=P_{\rm pre}$ there.
After the second bump it equals $\Delta_H(Z)$. The total pressure
moment is now $M_{\rm pre}^p(\infty,Z)$, so $P(\infty,Z)=0$ and
the pressure beyond the bumps agrees with the prescribed tail integral.
In particular,
\begin{equation}\label{eq:mc-pressure-preserved}
P_0(Z)=P(0,Z)=-M^p(\infty,Z)=P_0^{\rm pre}(Z).
\end{equation}
It is holomorphic near $[-1,1]$ by
Lemma~\ref{lem:outer-axis-pressure}. It is independent of
$R_{\rm ref}$ after the pre-heat waiting length has been fixed.
The global pressure bounds also survive. On the angular supports,
the bounds for $a$ and $\partial_Z\log U^\theta$ above give
$\partial_t\log U^\theta\le-1/2$ and
$|\partial_Z\log U^\theta|\le2|Z|/(1+Z^2)$.
Every other radial interval is unchanged. The proof of
Lemma~\ref{lem:outer-pressure} therefore gives the same global
and sharp pressure bounds for the corrected profile.

The angular and pressure moments now have their required values
at $R_{3/2}$: $M_{\mathrm{left}}^\theta=M_{\mathrm{right}}^\theta$
and $M_{\mathrm{left}}^p=M_{\mathrm{right}}^p$ there, as functions of $Z$.
The bumps leave $M^z$ and $M^{\theta z}$ unchanged, but change
$M^{z\theta}$ through the angular energy.
In Proposition~\ref{prop:outer-conditional-axial-correction} we constructed
the axial correction that enforces the remaining three moment conditions
under the tail energy assumption \eqref{eq:mc-pulse-future-energy}.

\paragraph{Applying the axial correction.}
\label{par:mc-final-tail-energy}
The corrected swirl satisfies \eqref{eq:mc-pulse-future-energy}.
Indeed, \eqref{eq:outer-tail-energy} gives
\[
\int_{R_{0z}}^\infty(\bar U^\theta)^2\,dR
\le C[1+\log(1/\mu)]R_{0z}[\bar U^\theta(R_{0z},Z)]^2.
\]
Extend $h$ by zero outside the angular supports. Then
$U^\theta=\bar U^\theta(1+h)$ on $[R_{0z},\infty)$, and
\eqref{eq:mc-angular-coefficients} gives
$\|h\|_{C_Z^1}\le C\mu^{29}$ uniformly in $R$.
Together with the logarithmic derivative bound in
\eqref{eq:outer-swirl-pressure-bounds}, this gives the same estimate
for the integral and its first $Z$ derivative. On the pulse interval
the swirl has not changed, and its pure power law gives exactly
\[
R_{0z}[U^\theta(R_{0z},Z)]^2
=e^{-26}R_{mz}[U^\theta(R_{mz},Z)]^2.
\]
The logarithmic derivative of the denominator is also bounded, so
division gives \eqref{eq:mc-pulse-future-energy}. The corresponding
bounds at each fixed derivative order show that the tail integral
is smooth in $Z$. The steep segment absorbs the factor $1/\delta$ in the final
power-law energy integral.

The angular bumps lie beyond $R_{0z}$. Thus the incoming moments
\eqref{eq:mc-pulse-data} and stresses at $R_{mz}$ are those of
Interval O.3, and on the pulse interval $P=P_{\rm pre}$ satisfies
\eqref{eq:outer-pressure-sharp}. Together with
\eqref{eq:mc-pulse-future-energy}, this verifies the hypotheses of
Proposition~\ref{prop:outer-conditional-axial-correction}.

Apply that proposition to choose $a_p,c_1,c_2$. The pulse satisfies
the admissible cone and gives $Q(R_{0z},Z)\ge c(\mu+Z^2)$ by
\eqref{eq:mc-pulse-angular-stress}. It imposes the three axial
moment conditions in \eqref{eq:mc-backward-moments} and preserves
the matched values of $M^\theta$ and $M^p$ at $R_{3/2}$.
Since $M^z=M^{\theta z}=U^z=0$ for $R\ge R_{0z}$,
\eqref{Ur-Uz} gives $U^r=0$ on and beyond Interval O.5.
Thus the coefficients have been chosen in the order
\[
 \tau\ \longrightarrow\ (d_1,d_2)\ \longrightarrow\ (a_p,c_1,c_2).
\]

\paragraph{The corrected cone after the pulse.}
On $[R_{0z},R_{3/2}]$ we now have $U^z=M^z=0$, hence $b=0$ and $W=1$.
The angular equation \eqref{eq:normalized-Q-J} becomes
\[
 \partial_yQ+\left(2-\frac a2\right)Q
 =\frac a2-1-\frac\delta2-\frac{1-\delta}{2}Z\zeta ,
 \qquad \zeta=\partial_Z\log U^\theta .
\]
On Interval O.5, $a\ge2+2\mu$ and $Z\zeta\le0$. On Interval O.6,
\eqref{eq:mc-angular-perturbation-bounds} gives $a\ge2+\mu$ and
$|Z\zeta|\le\mu/8$. Since $\delta\le\mu/4$, the source is
bounded below by $c\mu$ on both intervals. With $2+\mu\le a<3$,
a sufficiently small constant multiple of $\mu$ is a subsolution.
The pulse exit bound therefore gives
\begin{equation}\label{eq:mc-postpulse-Q}
 Q(R,Z)\ge c\mu
 \qquad(R_{0z}\le R\le R_{3/2},\ |Z|\le1).
\end{equation}

For the axial direction, \eqref{eq:mc-backward-moments} and
\eqref{Iz-moments-direct} give
\[
 M^{z\theta}(R,Z)=\frac12\int_R^\infty (U^\theta)^2\,d\rho,\qquad
 J=\Pi+\frac{2\delta ZM^{z\theta}-L\partial_Z M^{z\theta}}
                  {R(U^\theta)^2}.
\]
On $[R,R_{3/2}]$ the slope bound $a\ge2+\mu$ gives
$[U^\theta(\rho,Z)]^2\le [U^\theta(R,Z)]^2(\rho/R)^{-1-\mu}$.
Integrating this bound and using \eqref{eq:outer-tail-energy} beyond
$R_{3/2}$ yields $M^{z\theta}\lesssim R(U^\theta)^2/\mu$.
The corrected logarithmic derivative bound
\eqref{eq:outer-swirl-pressure-bounds} gives
$|\partial_Z M^{z\theta}|\lesssim |Z|M^{z\theta}$.
Since $|\Pi|\lesssim|Z|$, it follows that $|J|\lesssim |Z|/\mu$.
Consequently
\[
 |w|=\frac{U^\theta |J|}{Q}\le CP_*e^{-T_d/2}\mu^{28}|Z|\ll1,
 \qquad
 \frac{\mathcal I^\theta}{F}=\frac{RQ}{L_\delta}\ge cR_{\rm ref}.
\]
Here we used \eqref{eq:outer-power-end-moments} and the decay of the
swirl after $R_{mz}$. As $b=0$ and $2+\mu\le a<3$, the second
expression in \eqref{eq:outer-cone-test} is $(a-2)w^2<1/2$.
Lemma~\ref{lem:outer-cone-test} gives the admissible cone through
$R_{3/2}$, including both angular bumps. At $R_{3/2}$ all five
moments of these left-side fields have their prescribed right-side values.

\subsection{The completed outer construction}
\label{subsec:outer-cone-closure}

At the waiting length selected by \eqref{eq:mc-wait-choice}, the two
angular corrections match $M^\theta$ and $M^p$ at $R_{3/2}$. For this
completed swirl the axial correction gives the three remaining moments
by \eqref{eq:mc-backward-moments}. Thus
\begin{equation}\label{eq:outer-left-right-matching}
\mathbf M_{\mathrm{left}}(R_{3/2},Z)
=\mathbf M_{\mathrm{right}}(R_{3/2},Z),
\qquad -1\le Z\le1.
\end{equation}
Beyond $R_{3/2}$ the two sides have the same velocities, so their
moment derivatives in $R$ agree. The five moments therefore coincide
as functions of $Z$ for every $R\ge R_{3/2}$. With the common axis
pressure $P_0=P_0^{\rm pre}$, Lemma~\ref{lem:moment-gluing} identifies
the pressure, radial velocity, inertial stress, and shear with the
fields computed from the right. In particular the terminal conditions
\eqref{eq:mc-targets} hold, and the exact heat exterior has
$\mathcal T=0$ for $R\ge R_h$.

The right-side cone estimates of Section~\ref{subsec:outer-waiting}
therefore apply to the completed profile on $[R_{3/2},R_h)$.

Together with the preceding intervals, this proves the relaxed cone on
$[R_{\rm useless},R_h)$ and the admissible cone on
$(R_{-\mu},R_h)$. The swirl is positive for $R>0$ and
$\mathcal S^\theta<0$ on $[R_{\rm ref},R_h]$. Up to $R_{\rm ref}$
the angular and axial velocities and all five moments retain their
reference values. The axis pressure is
$P_0^{\rm pre}$ by \eqref{eq:mc-pressure-preserved}, with the analytic
bounds of Lemma~\ref{lem:outer-axis-pressure}.

For later use in the inner construction of Section~\ref{sec:inner-construction},
we fix $P_*,\delta$ and the dimensionless
outer schedule. The selected $\tau$ and $P_0^{\rm pre}$ are then
independent of $R_{\rm ref}$. At each allowed $R_{\rm ref}$ we construct
the heat tail at that scale and solve the corresponding moment equations
with the uniform bounds \eqref{eq:mc-angular-coefficients},
\eqref{eq:mc-axial-end-coefficients}, and \eqref{eq:mc-pulse-amplitude}.
This completes the proof of Proposition~\ref{prop:road-outer-summary}.

\section{The inner construction: the stress-free core, connection, and moment matching}\label{sec:inner-construction}
\label{sec-inner-construction}

Let $U_o$ be the completed outer profile of Section~\ref{sec:outer-construction}.
We replace its reference continuation below $R_{\rm useless}=e^{-5}R_{\rm ref}$
by a regular stress-free core and a connection to the reference velocities
\begin{equation}\label{eq:inner-reference}
 U_{\rm ref}^\theta=\frac{P_*}{1+Z^2}
 \left(\frac R{R_{\rm ref}}\right)^{1/10},\qquad U_{\rm ref}^z=4Z
 \qquad(0<R\le R_{\rm ref}).
\end{equation}
The core uses the analytic axis pressure supplied by the outer construction,
\begin{equation}\label{eq:inner-axis-pressure}
 P_0(Z)=-\int_0^\infty\frac{(U_o^\theta)^2}{2R}\,dR=P_0^{\rm pre}(Z),
 \qquad \|P_0\|_{C^1_Z}\le C_{\rm pres}P_*^2,
\end{equation}
by \eqref{eq:mc-pressure-preserved} and Lemma~\ref{lem:outer-axis-pressure},
with $C_{\rm pres}\ge1$ independent of $P_*,\delta,R_{\rm ref}$ and the
waiting length $\tau$. Here $C^k_Z$ denotes the
$C^k$ norm on $[-1,1]$.

We integrate the inner moments from the axis, so the pressure is
\begin{equation}\label{eq:inner-local-pressure}
 P(R,Z)=P_0(Z)+M^p(R,Z),\qquad
 M^p(R,Z)=\int_0^R\frac{(U^\theta(\rho,Z))^2}{2\rho}\,d\rho.
\end{equation}

Write $\mathbf M_{\rm ref}$ for the five moments \eqref{fiveM} of the reference
continuation \eqref{eq:inner-reference}.
We require
\begin{equation}\label{eq:inner-required-matching}
\boxed{\quad M^j(R_{\rm useless},Z)=M^j_{\rm ref}(R_{\rm useless},Z),
\qquad j=\theta,z,\theta z,z\theta,p.\quad}
\end{equation}
Direct integration of the reference continuation gives
\begin{equation}\label{eq:inner-target-moments}
\begin{aligned}
M^\theta_{\rm ref}(R_{\rm useless})&=\tfrac58 R_{\rm useless}\sqrt{2R_{\rm useless}}\,
 \frac{e^{-1/2}P_*}{1+Z^2},
\\
M^z_{\rm ref}(R_{\rm useless})&=4ZR_{\rm useless},\\
M^{\theta z}_{\rm ref}(R_{\rm useless})&=4Z M^\theta_{\rm ref}(R_{\rm useless}),\\
M^{z\theta}_{\rm ref}(R_{\rm useless})&=16Z^2R_{\rm useless}
 -\frac{5R_{\rm useless}e^{-1}P_*^2}{12(1+Z^2)^2},
\\
M^p_{\rm ref}(R_{\rm useless})&=\frac{5e^{-1}P_*^2}{2(1+Z^2)^2}.
\end{aligned}
\end{equation}
We follow Intervals I.1--I.4 of Subsection~\ref{roadin:subsection}.
After constructing the core on I.1, we use a small shear opposing the
inertial stress to leave the core on I.2. Once the angular inertial
stress is strong enough, we remove the axial shear and join the swirl
to the reference power on I.3. On I.4 we correct all five accumulated
moment defects on $(R_m,2R_m)$, where $R_m=e^{-6}R_{\rm ref}$, while
preserving the relaxed cone.

First choose the shift $j$ in $U_0^z=4Z+j$ by \eqref{roadin:j-budget},
using the tolerances for restoring $4Z$ and correcting the five moments.
Fix $P_*\ge1$,
$0<\delta\le1/200$ and the dimensionless outer schedule as in
Section~\ref{sec:outer-construction}. This fixes the waiting length $\tau$
and the axis pressure $P_0$ independently of $R_{\rm ref}$.
Fix a complex neighborhood of $[-1,1]$ on which the axis data are
holomorphic, then take $\Lambda$ large enough for the exit estimate
\eqref{eq:inner-axial-Lambda} and the comparison
in Lemma~\ref{prop:inner-comparison}. Finally choose $C_*$
sufficiently large and use the radii
\begin{equation}\label{eq:inner-family-radii}
 R_a=4/\Lambda,\qquad R_{\rm ref}=110(C_*P_*)^{10}.
\end{equation}
Lemma~\ref{lem:imc-prepared-family} verifies the required smallness of
the five moment defects for these choices.

\subsection{Interval I.1: The stress-free core and its linear model}
\label{sec:analytic-core}
\label{sec-core-construction}

We seek a smooth profile with $\mathcal{T}=0$ on $[0,R_a]\times[-1,1]$, where
$R_a=4/\Lambda$. We construct it up to $R=4.1/\Lambda$ so that it
extends across $R_a$ for the connection in Interval~I.2.
Figure~\ref{fig:road-core-exit} shows the radial decay of $F/F_0$
and the location of the exit within this domain.
We use the stress-free equations \eqref{eq:stress-free-sources}, which we recall here:
\begin{equation}\label{eq:core-system}
\begin{aligned}
2{L_\delta}(R\partial_R^2F+2\partial_RF)
&=W(F+R\partial_RF)
+\frac{\delta}{2}(1-2ZU^z)F+H\partial_ZF,\\
2{L_\delta}(R\partial_R^2U^z+\partial_RU^z)
&=WR\partial_RU^z
+\frac{1+\delta}{2}(1-2ZU^z)U^z
+H\partial_ZU^z\\
&\quad+L\partial_ZP-2(1+\delta)ZP-2ZRF^2.
\end{aligned}
\end{equation}
Here the two transport coefficients are
\begin{equation}\label{eq:core-transport}
W=1-\frac{(1-\delta)ZM^z+L\partial_ZM^z}{R},
\qquad
H=\frac{1-\delta}{2}Z+LU^z.
\end{equation}
This form follows from the source formulas \eqref{inviscid-stress-sources-Ur}
using \eqref{Ur-Uz}, $\partial_R M^z=U^z$, and $\partial_R P=F^2$.

The axis pressure $P_0=P_0^{\rm pre}$ is fixed by \eqref{eq:mc-pressure-preserved}.
It is even and extends holomorphically to the fixed rectangle
$\Omega_0$ of Lemma~\ref{lem:outer-axis-pressure}. In particular,
\begin{equation}\label{eq:pressure}
P_0(Z)\le-\frac{5P_*^2}{2(1+Z^2)^2},\qquad
Z\partial_ZP_0(Z)\ge\frac{10P_*^2 Z^2}{(1+Z^2)^3}.
\end{equation}
It remains to prescribe the angular and axial velocity data. By axis
regularity \eqref{smooth-axis}, $U^\theta(0,Z)=0$ and
$F=U^\theta/\sqrt{2R}$ extends smoothly to the axis. We prescribe
\begin{equation}\label{eq:data-to-choose}
 F_0(Z)=F(0,Z)=\lim_{R\downarrow0}\frac{U^\theta(R,Z)}{\sqrt{2R}}>0,
 \qquad U^z_0(Z)=U^z(0,Z).
\end{equation}
Since $M^z/R=\int_0^1 U^z(sR,Z)\,ds$, the transport coefficients extend
smoothly to the axis, with $H_0(Z)=(1-\delta)Z/2+LU^z_0(Z)$.
Setting $R=0$ in \eqref{eq:core-system} gives the first radial slopes
\begin{equation}\label{eq:core-axis-slopes}
\begin{aligned}
\partial_RF(0,Z)
&=\frac1{4{L_\delta}}\left[
\left(1+\frac{\delta}{2}-ZU^z_0-L\partial_ZU^z_0\right)F_0
+H_0\partial_ZF_0\right],\\
\partial_RU^z(0,Z)
&=\frac1{2{L_\delta}}\left[
\frac{1+\delta}{2}(1-2ZU^z_0)U^z_0
+H_0\partial_ZU^z_0+L\partial_ZP_0
-2(1+\delta)ZP_0\right].
\end{aligned}
\end{equation}
Comparing successive powers of $R$ in \eqref{eq:core-system} then determines
every higher radial Taylor coefficient of $F$ and $U^z$. Hence the axis data
uniquely determine the analytic solution near $R=0$
(see Appendix~\ref{app:nonlinear-core-construction}).

\paragraph{The key exit condition.}
We require $F>0$ and $\partial_R F<0$ throughout the core, and
\begin{equation}\label{eq:exit}
 \kappa(R_a,Z)>2\qquad (|Z|\le1).
\end{equation}
The shear intensity is
\begin{equation}\label{eq:kappa}
\mathcal{S}^\theta=2R\partial_RF,\qquad
\kappa=-\frac{|\mathcal{S}|^2}{F\mathcal{S}^\theta}
=-2R\frac{\partial_RF}{F}
+\frac{(\partial_RU^z)^2}{-F\partial_RF}.
\end{equation}
Over most of $[-1,1]$ we seek $-\partial_R F/F>1/R_a$ at the exit,
so the angular term exceeds $2$. Where the angular decay weakens, we keep
$|\partial_R U^z|$ bounded below and reduce the angular amplitude.
The linear model will give a lower bound of $3$ for the corresponding
term in each region.

\paragraph{The linear model and its regular solution.}\label{subsec:linear-model}
To choose $F_0$ and $U_0^z$, we first solve a linear model with these
axis data left arbitrary. For fixed smooth data and a regular
solution, as $R\to0$ we have
\[
\begin{aligned}
U^z&=U^z_0+O(R),&
\partial_ZU^z&=\partial_ZU^z_0+O(R),\\
P&=P_0+O(R),&
\partial_ZP&=\partial_ZP_0+O(R),\\
\frac{\partial_ZF}{F}
&=\frac{\partial_ZF_0}{F_0}+O(R).
\end{aligned}
\]
In the angular equation write $\partial_ZF=F\,\partial_Z\log F$.
Freeze the coefficients and the axial forcing at their axis values and
drop terms with an explicit factor $R$ on the right of \eqref{eq:core-system}.
The constants in the Taylor expansions above depend on the axis data.
For the axis data chosen below, Appendix~\ref{app:analytic-core} shows that
the difference $\partial_Z\log F-\partial_Z\log F_0$ stays bounded throughout
the core and enters the rescaled angular equation with a factor $\Lambda^{-1}$.

The axial right-hand side of \eqref{eq:core-system} at $R=0$ is
$2L_\delta\partial_R U^z(0,Z)$. The model uses this value at every radius.
Keeping radial diffusion therefore gives
$\partial_R(R\partial_R U^{z,\rm m})=\partial_R U^z(0,Z)$. For fixed $Z$,
integration from $0$ to $R$ gives $\partial_R U^{z,\rm m}(R,Z)=\partial_R U^z(0,Z)$ for $R>0$, since
$\rho\,\partial_\rho U^{z,\rm m}(\rho,Z)\to0$ at the axis.
The resulting linear model is
\begin{equation}\label{eq:model}
\boxed{\begin{aligned}
2{L_\delta}(R\partial_R^2F^{\rm m}+2\partial_RF^{\rm m})
&=\left(1+\frac{\delta}{2}-ZU^z_0-L\partial_ZU^z_0
+H_0\frac{\partial_ZF_0}{F_0}\right)F^{\rm m},\\
2{L_\delta}\partial_RU^{z,\rm m}
&=\frac{1+\delta}{2}(1-2ZU^z_0)U^z_0
+H_0\partial_ZU^z_0+L\partial_ZP_0
-2(1+\delta)ZP_0,\\
F^{\rm m}(0,Z)&=F_0(Z),\qquad
U^{z,\rm m}(0,Z)=U^z_0(Z).
\end{aligned}}
\end{equation}
Its regular solution for each fixed $Z$ is
\begin{equation}\label{eq:core-general-model-solution}
\begin{aligned}
 F^{\rm m}(R,Z)&=F_0(Z)\mathcal B(q),&
 q&=-\frac{2R\,\partial_RF(0,Z)}{F_0(Z)},\\
 U^{z,\rm m}(R,Z)&=U_0^z(Z)+R\,\partial_R U^z(0,Z),&
 \mathcal B(q)&=\sum_{n=0}^\infty\frac{(-q)^n}{n!(n+1)!}.
\end{aligned}
\end{equation}
For $q>0$, $\mathcal B(q)=J_1(2\sqrt q)/\sqrt q$, where $J_1$ is the
Bessel function of the first kind of order one. The series gives the
normalization $\mathcal B(0)=1$.
Writing $F^{\rm m}=\sum_{n\ge0}a_n R^n$ in \eqref{eq:model} gives
\[
 (n+1)(n+2)a_{n+1}=\frac{2\partial_RF(0,Z)}{F_0(Z)}a_n,\qquad a_0=F_0(Z).
\]
The resulting series converges for every $R$ by the ratio test.
Multiplying the angular equation by $R$ and integrating twice from the
axis gives the Volterra equation
\[
 F^{\rm m}(R)=F_0(Z)+\frac{2\partial_RF(0,Z)}{F_0(Z)}
 \int_0^R\left(1-\frac{\rho}{R}\right)F^{\rm m}(\rho)\,d\rho .
\]
For the difference $D$ of two regular solutions with the same axis value,
this gives $|D(R)|\le \dfrac{2|\partial_RF(0,Z)|}{F_0(Z)}\int_0^R|D(\rho)|\,d\rho$,
so Gronwall's inequality gives $D=0$. Since $\mathcal B(0)=1$ and $\mathcal B'(0)=-1/2$, the model has
the same axis values and first radial derivatives as a regular nonlinear
solution with these data. Thus $F-F^{\rm m}$ and $U^z-U^{z,\rm m}$ are
both $O(R^2)$ as $R\to0$ for fixed smooth data. The Bessel comparison
function and the weighted analytic framework used in Appendix~\ref{app:analytic-core}
already appear in Appendix~B of \cite{1}.

\paragraph{Choice of axis data.}
The core construction applies to every fixed $0<j\le1/20$, with the
choice \eqref{roadin:j-budget} for the subsequent matching.
For positive parameters $\Lambda,C_*$, prescribe the axis values
\begin{equation}\label{eq:axial-data}
 U_0^z(Z)=4Z+j,\qquad F_0(Z)=C_*^{-1}e^{-\Lambda G(Z)}.
\end{equation}
Set $\sigma_0=j/500$ and let $Z_0$ be the unique zero of
\begin{equation}\label{eq:H0}
 H_0(Z)=\frac{1-\delta}{2}Z+(1-Z^2)(4Z+j).
\end{equation}
The pressure derivative vanishes at $Z=0$. The shift $j>0$ moves $Z_0$
away from this point so that the pressure supplies axial shear near $Z_0$.
Away from $Z_0$ we use the angular transport term $H_0\partial_Z\log F_0$
to increase the radial decay of $F$. For this purpose set
\[
 G(Z)=\int_{Z_0}^Z\frac{L_\delta(w)H_0(w)}{H_0(w)^2+\sigma_0^2}\,dw .
\]
With $\chi_c=H_0^2/(H_0^2+\sigma_0^2)$ we have
$H_0\partial_Z\log F_0=-\Lambda L_\delta\chi_c$.
Write
\begin{align}
 \beta(Z)&=\frac{3-\delta/2+jZ}{L_\delta(Z)},\label{eq:beta}\\
 g(Z)&=-\frac{1+\delta}{2}(1-2ZU^z_0)U^z_0-4H_0
       +2(1+\delta)ZP_0-L\partial_ZP_0 .\label{eq:g}
\end{align}
Since $\partial_ZU_0^z=4$, \eqref{eq:core-axis-slopes} gives
\begin{equation}\label{eq:chosen-model}
 \frac{\partial_RF(0,Z)}{F_0(Z)}=-\frac{\Lambda\chi_c+\beta}{4},\qquad
 \partial_RU^z(0,Z)=-\frac{g}{2L_\delta}.
\end{equation}
Substituting these slopes into \eqref{eq:core-general-model-solution} yields
\begin{equation}\label{eq:continuation-model}
\begin{aligned}
 F^{\rm m}(R,Z)&=F_0(Z)\mathcal B(q),&
 q&=\tfrac R2(\Lambda\chi_c+\beta),\\
 U^{z,\rm m}(R,Z)&=4Z+j-\frac{g}{2L_\delta}R .
\end{aligned}
\end{equation}
The first equation in \eqref{eq:model} governs the normalized swirl
$F^{\rm m}=U^{\theta,\rm m}/\sqrt{2R}$. With the chosen axis data its
right-hand side is $-L_\delta(\Lambda\chi_c+\beta)F^{\rm m}$.
At $Z_0$ the term $\Lambda\chi_c$ vanishes. The slope $\partial_ZU_0^z=4$
gives $\beta>0$, so $\partial_RF^{\rm m}(0,Z_0)=-\beta(Z_0)F_0(Z_0)/4<0$.
Recall from \eqref{eq:kappa} the two terms in the definition of $\kappa_{\rm m}$.
For fixed $\Lambda$, the angular term $-2R\,\partial_R F^{\rm m}/F^{\rm m}$
is independent of $C_*$, while the axial term
$(\partial_RU^{z,\rm m})^2/(-F^{\rm m}\partial_RF^{\rm m})$ is proportional to $C_*^2$.
After fixing $\Lambda$, we take $C_*$ sufficiently large to make the
axial term exceed $3$ where the angular decay is weak.

\begin{unnumberedlemma}[Exit bound for the linear model]
Suppose that $0<j\le1/20$, $P_*\ge1$, $0<\delta\le1/100$, and that
the axis pressure satisfies \eqref{eq:pressure}. For the model
\eqref{eq:continuation-model}, take
\begin{equation}\label{eq:continuation-thresholds}
 \Lambda\ge\max\{500,j^{-2}\},\qquad C_*\ge\Lambda^2.
\end{equation}
At $R_a=4/\Lambda$ we have
\[
 \kappa_{\rm m}(R_a,Z)>3\qquad (|Z|\le1).
\]
The angular term in \eqref{eq:kappa} exceeds $3$ where $\chi_c\ge99/100$,
and the axial term exceeds $3$ where $\chi_c\le99/100$.
\end{unnumberedlemma}

\begin{proof}
For $U_0^z=4Z+j$ the part of the coefficient of $F^{\rm m}$ in
\eqref{eq:model} independent of $\partial_ZF_0$ is
\begin{equation}\label{eq:negative-base}
1+\frac{\delta}{2}-ZU^z_0-L\partial_ZU^z_0
=-3+\frac{\delta}{2}-jZ
\le-\frac52<0
\qquad(-1\le Z\le1).
\end{equation}
Using \eqref{eq:beta}, ${L_\delta}(Z)=1-\delta Z^2$, $0<j\le1/20$, and
$0<\delta\le1/100$, we obtain
\[
\begin{aligned}
\beta(Z)&\ge3-\frac{\delta}{2}-j
\ge3-\frac1{200}-\frac1{20}>\frac{29}{10},\\
\beta(Z)&\le\frac{3+j}{1-\delta}
\le\frac{3+1/20}{1-1/100}<\frac{31}{10}.
\end{aligned}
\]

The quotient
$H_0/L=4Z+j+(1-\delta)Z/[2(1-Z^2)]$
is strictly increasing from $-\infty$ to $+\infty$ on
$(-1,1)$, since
\[
\partial_Z\left(\frac{H_0}{L}\right)
=4+\frac{(1-\delta)(1+Z^2)}{2(1-Z^2)^2}>0.
\]
Thus $H_0$ has one simple zero $Z_0$.
For $0<\delta\le1/100$ and $0<j\le1/20$, we have
\[
\begin{aligned}
H_0(-j/4)&=-\frac{(1-\delta)j}{8}<0,\\
H_0(-j/5)&=j\left(\frac{1+\delta}{10}-\frac{j^2}{125}\right)
>\frac{2j}{25}>0.
\end{aligned}
\]
Consequently,
\begin{equation}\label{eq:root}
-\frac j4<Z_0<-\frac j5<0,\qquad H_0(0)=j.
\end{equation}

Since ${L_\delta}>0$ and $H_0$ is negative to the left of $Z_0$
and positive to its right, $G$ attains its minimum at
$Z_0$. Thus $G(Z)\ge G(Z_0)=0$ on $[-1,1]$.

Away from the zero of $H_0$, the argument $q$ has size one
when $R$ has size $\Lambda^{-1}$. We therefore set
\begin{equation}\label{eq:continuation-scale}
\widetilde R=\Lambda R,
\qquad R_a=\frac4\Lambda,
\qquad q=\frac{\widetilde R}{2}\left(\chi_c+\frac{\beta}{\Lambda}\right).
\end{equation}
The upper endpoint $\widetilde R=4.1$ leaves a
short interval beyond the exit $\widetilde R=4$ on which $F^{\rm m}$
remains positive. At the exit the angular profile has already decreased
sufficiently where $\chi_c$ is close to one.

Since $G\ge0$ and $C_*\ge\Lambda^2$, the axis amplitude satisfies
\[
0<F_0\le\Lambda^{-2}.
\]

On the rescaled interval $0\le\widetilde R\le4.1$, we have
\[
0\le q
=\frac{\widetilde R}{2}
\left({\chi_c}+\frac{\beta}{\Lambda}\right)
\le\frac{41}{20}\left(1+\frac4{500}\right)
<\frac{21}{10}.
\]

In the Taylor expansion of $\mathcal B(q)$ in \eqref{eq:core-general-model-solution},
the ratio of successive term magnitudes is $q/((n+1)(n+2))$.
The tail after degree three therefore decreases in magnitude
on $0\le q\le21/10$. The term magnitudes of $-(\partial_q\mathcal B)$ have ratio
$q/((n+1)(n+3))<1$. For
$P_3(q)=1-q/2+q^2/12-q^3/144$, one also has
\[
 \partial_q P_3(q)=-\frac{(q-4)^2+8}{48}<0,\qquad
 P_3(21/10)=\frac{4051}{16000}>\frac14.
\]
Hence, for $0\le q\le21/10$,
\begin{equation}\label{eq:continuation-B}
\begin{aligned}
1\ge\mathcal B(q)
&\ge1-\frac q2+\frac{q^2}{12}-\frac{q^3}{144}
>\frac14,\\
\frac18
&<\frac12-\frac q6
\le-(\partial_q\mathcal B)(q)\le\frac12.
\end{aligned}
\end{equation}
The upper bound for $\mathcal B$ follows from
$\mathcal B(0)=1$ and $(\partial_q\mathcal B)<0$.

Using $F^{\rm m}=F_0\mathcal B(q)$ and
$\partial_Rq=(\Lambda{\chi_c}+\beta)/2$, we obtain
\begin{equation}\label{eq:continuation-signs}
\frac14F_0<F^{\rm m}\le F_0,\qquad
\frac1{16}F_0(\Lambda{\chi_c}+\beta)
\le-\partial_RF^{\rm m}
\le\frac14F_0(\Lambda{\chi_c}+\beta).
\end{equation}
In particular, $\partial_RF^{\rm m}<0$ everywhere,
including at $Z_0$, where ${\chi_c}=0$ and the positive
coefficient $\beta$ supplies the strict inequality.
Moreover,
\begin{equation}\label{eq:logarithmic-shear-size}
\frac1{16}(\Lambda{\chi_c}+\beta)
\le-\frac{\partial_RF^{\rm m}}{F^{\rm m}}
\le\Lambda{\chi_c}+\beta.
\end{equation}
Thus the logarithmic radial derivative has size comparable
to $\Lambda$ wherever ${\chi_c}$ is bounded below by a fixed
positive constant, and size comparable to one at $Z_0$.

At $R=R_a=4/\Lambda$ we distinguish two cases according to $\chi_c(Z)$.
First suppose that $\chi_c(Z)\ge99/100$. Then
\[
q=2\left({\chi_c}+\frac{\beta}{\Lambda}\right),
\qquad
\frac{19}{10}<q<\frac{21}{10}.
\]

The term magnitudes in the series for $3\mathcal B+2q\mathcal B'$ have consecutive ratio
\[
 \frac{2n+5}{2n+3}\frac{q}{(n+1)(n+2)}\le\frac{49}{100}
 \qquad(n\ge1,\ 0\le q\le21/10).
\]
Thus the degree-four truncation gives the upper bound
\[
3\mathcal B(q)+2q(\partial_q\mathcal B)(q)
\le p(q):=
3-\frac52q+\frac7{12}q^2-\frac1{16}q^3
+\frac{11}{2880}q^4.
\]

\noindent On $[19/10,21/10]$ we have $p'(q)\le-5/2+7q/6<0$,
so $p(q)\le p(19/10)<-1/50$. Since $\mathcal B(q)>0$, it follows that
\begin{equation}\label{eq:continuation-angular-exit}
-\frac{2R_a\partial_RF^{\rm m}}{F^{\rm m}}
=-\frac{2q(\partial_q\mathcal B)(q)}{\mathcal B(q)}>3.
\end{equation}

Next suppose that $\chi_c\le99/100$. To estimate the axial term in
\eqref{eq:kappa}, we first evaluate $g$ at $Z_0$. The relation $H_0(Z_0)=0$ gives
\[
U^z_0(Z_0)=-\frac{(1-\delta)Z_0}{2L(Z_0)},
\qquad
1-2Z_0U^z_0(Z_0)
=1+\frac{(1-\delta)Z_0^2}{L(Z_0)}
=\frac{{L_\delta}(Z_0)}{L(Z_0)}.
\]
Consequently,
\[
-\frac{1+\delta}{2}
\bigl(1-2Z_0U^z_0(Z_0)\bigr)U^z_0(Z_0)
=\frac{(1+\delta)(1-\delta)Z_0{L_\delta}(Z_0)}
       {4L(Z_0)^2}.
\]
Substituting into \eqref{eq:g} yields
\begin{equation}\label{eq:g-at-root}
g(Z_0)=(1+\delta)Z_0
\left[
2P_0(Z_0)+\frac{(1-\delta){L_\delta}(Z_0)}{4L(Z_0)^2}
\right]
-L(Z_0)\partial_ZP_0(Z_0)>\frac{4j}{5}.
\end{equation}
To verify the last inequality, use $|Z_0|<j/4\le1/80$.
Since $P_*\ge1$ and $(1-\delta){L_\delta}(Z_0)\le1$, the bracket
is at most
\[
-\frac5{(1+Z_0^2)^2}
+\frac1{4(1-Z_0^2)^2}
<-\frac92+\frac12=-4.
\]
Also $Z_0<0$ and $\partial_ZP_0(Z_0)<0$. Hence
\[
g(Z_0)>4(1+\delta)|Z_0|
>4|Z_0|>\frac{4j}{5},
\]
where the last inequality uses $|Z_0|>j/5$.
If the shift were zero, the zero of $H_0$ would be $Z=0$,
where evenness of $P_0$ would give $g(0)=0$;
the two shear mechanisms would then weaken at the same point.

For a bound on the entire region $\{\chi_c\le99/100\}$, set
\[
 b=\frac{2j}{5},\qquad I_0=\left[-\frac j4,-\frac j5\right].
\]
To see that $\chi_c\le99/100$ implies $Z\in I_0$, consider the complement of $I_0$.
For $Z\le-j/4$, we have $4Z+j\le0$, and hence
\[
H_0(Z)\le\frac{1-\delta}{2}Z
\le-\frac{(1-\delta)j}{8}<-\frac j{50}.
\]
On $[-j/5,1/2]$,
\[
\partial_ZH_0
=\frac{9-\delta}{2}-12Z^2-2jZ>0,
\qquad
H_0(-j/5)>\frac{2j}{25}.
\]
On $[1/2,1]$, both terms in \eqref{eq:H0} are
nonnegative, and
\[
H_0(Z)\ge\frac{1-\delta}{4}>\frac j{50}.
\]
Thus $|H_0|>j/50$ outside $I_0$. On the other hand,
\[
{\chi_c}(Z)\le\frac{99}{100}
\quad\Longrightarrow\quad
|H_0(Z)|\le\sqrt{99}\,\sigma_0<\frac j{50}.
\]
It remains to bound $g$ on $I_0$.
For $Z\in I_0$, we have
\[
\frac j5\le |Z|\le\frac j4\le\frac1{80},\qquad
0\le U^z_0(Z)\le\frac j5,\qquad
H_0(Z)\le U^z_0(Z)\le\frac j5.
\]
Consequently,
\[
\frac{1+\delta}{2}(1-2ZU^z_0)U^z_0
\le\frac{1+\delta}{2}\left(1+\frac{j^2}{10}\right)\frac j5
\le\frac j5.
\]
The derivative bound in \eqref{eq:pressure} gives
\[
-L(Z)\partial_ZP_0(Z)
\ge\frac{10P_*^2|Z|(1-Z^2)}{(1+Z^2)^3}
\ge8|Z|.
\]
Since $ZP_0(Z)>0$ on $I_0$, we obtain from \eqref{eq:g}
\[
g(Z)\ge-\frac j5-\frac{4j}{5}+8|Z|
\ge\frac{3j}{5}>b
\qquad(Z\in I_0).
\]
Therefore
\begin{equation}\label{eq:width}
{\chi_c}(Z)\le\frac{99}{100}
\quad\Longrightarrow\quad
Z\in I_0
\quad\Longrightarrow\quad
g(Z)\ge b.
\end{equation}
This bound holds for each fixed $j$, uniformly in $P_*\ge1$ and
$0<\delta\le1/100$. Since $0<L_\delta\le1$,
\[
|\partial_RU^{z,\rm m}|=\frac{g}{2{L_\delta}}\ge\frac b2.
\]
Also, \eqref{eq:continuation-signs} gives
\[
0<-F^{\rm m}\partial_RF^{\rm m}
\le\frac14F_0^2(\Lambda{\chi_c}+\beta).
\]
Combining these bounds with $F_0\le\Lambda^{-2}$,
${\chi_c}\le1$, and $\beta<4$, we obtain
\begin{equation}\label{eq:continuation-axial-exit}
\begin{aligned}
\frac{(\partial_RU^{z,\rm m})^2}
{-F^{\rm m}\partial_RF^{\rm m}}
&\ge
\frac{b^2}{F_0^2(\Lambda{\chi_c}+\beta)}\\
&\ge
\frac{b^2\Lambda^3}{1+4/\Lambda}
\ge\frac{2j^2\Lambda^3}{25}
\ge\frac{2}{25j^4}>3.
\end{aligned}
\end{equation}
Here we used $b^2=4j^2/25$, $1+4/\Lambda\le2$,
$\Lambda\ge j^{-2}$, and $0<j\le1/20$.

The two regions cover $[-1,1]$, and both contributions to $\kappa_{\rm m}$
are nonnegative by \eqref{eq:continuation-signs}. This proves the exit bound.
\end{proof}

\paragraph{The nonlinear core.}
With $j$ and $\sigma_0=j/500$ fixed, choose a bounded convex complex
neighborhood $\Omega$ of $[-1,1]$
whose closure lies in a common domain of analyticity of $P_0$,
$1/L_\delta$, and $1/(H_0^2+\sigma_0^2)$. This also gives a
holomorphic extension of $G$ to a neighborhood of $\overline\Omega$.
Put $r_\Omega=\operatorname{dist}([-1,1],\mathbb C\setminus\Omega)$,
fix $0<h_{\rm an}<r_\Omega/4$, and set
$A_\Omega=\max_{\overline\Omega}(-\operatorname{Re}G)$.
These choices may depend on $j$ and precede the choice of $\Lambda$
and $C_*$. For $F_0=C_*^{-1}e^{-\Lambda G}$ we require
\begin{equation}\label{eq:continuation-complex-bound}
 C_*\ge\Lambda^2e^{\Lambda A_\Omega}.
\end{equation}
In particular, $\sup_\Omega |F_0|\le\Lambda^{-2}$.

\begin{proposition}[The analytic core and its exit data]\label{prop:continuation-core}
Fix $0<j\le1/20$, the outer parameter $\delta$, the axis pressure
$P_0=P_0^{\rm pre}$, and the complex neighborhood $\Omega$ and analytic
width $h_{\rm an}$ chosen above. There exists
$\Lambda_0\ge\max\{500,j^{-2}\}$, depending only on these data,
such that for each $\Lambda\ge\Lambda_0$
and every $C_*$ satisfying \eqref{eq:continuation-complex-bound},
the stress-free system \eqref{eq:core-system} has a jointly analytic
solution $(F_c,U_c^z)$ on a neighborhood of $[0,4.1/\Lambda]\times[-1,1]$
with axis data
\[
 F_c(0,Z)=C_*^{-1}e^{-\Lambda G(Z)},\qquad
 U_c^z(0,Z)=4Z+j,\qquad P(0,Z)=P_0(Z).
\]
On this rectangle $\mathcal{T}_c=0$,
$F_c>0$ and $-\partial_R\log F_c\ge1/64$. At the exit $R_a=4/\Lambda$,
\begin{equation}\label{eq:continuation-exit}
 \kappa_c(R_a,Z):=
 \left[-2R\frac{\partial_RF_c}{F_c}
       +\frac{(\partial_R U_c^z)^2}{-F_c\partial_RF_c}\right]_{R=R_a}
 \ge\frac52\qquad (|Z|\le1).
\end{equation}
The angular term alone exceeds $5/2$ where $\chi_c\ge99/100$.
Where $\chi_c\le99/100$ the axial slope satisfies
$-2L_\delta\,\partial_R U_c^z(R_a,Z)\ge3j/10$.

The axial velocity and its radial average at the exit satisfy
\begin{equation}\label{eq:inner-core-axial-budget}
 \rho_{\rm core}:=\max\left\{\|U_c^z(R_a,\cdot)-U_0^z\|_{C_Z^2},
       \left\|\frac{M_c^z(R_a,\cdot)}{R_a}-U_0^z\right\|_{C_Z^2}\right\}
 \le \frac{B}{\Lambda}.
\end{equation}
At the same radius the axial derivative of $\log F_c$ satisfies
\begin{equation}\label{eq:inner-core-input}
 \left|\partial_Z\log F_c(R_a,Z)+\Lambda G'(Z)\right|
 \le B\bigl(\sqrt{\chi_c(Z)}+\Lambda^{-1}\bigr)
 \qquad (|Z|\le1).
\end{equation}
The constant $B$ depends only on $j,\delta,\Omega,h_{\rm an}$ and
the holomorphic pressure $P_0$ on $\Omega$.
For fixed $\Lambda$ and each integer $k\ge0$, uniformly in $C_*$,
\begin{equation}\label{eq:inner-prepared-normalized}
 \|C_*F_c\|_{C^k}+\|(C_*F_c)^{-1}\|_{C^k}+\|U_c^z\|_{C^k}\le B_{k,\Lambda}.
\end{equation}
The $C^k$ norm in \eqref{eq:inner-prepared-normalized} is taken in
$(R/R_a,Z)\in[0,e^{1/100}]\times[-1,1]$.
\end{proposition}
The contraction about the model \eqref{eq:continuation-model} in
Appendix~\ref{app:analytic-core} gives \eqref{eq:core-normalized-Ck}--\eqref{eq:continuation-comparison}.
On the fixed rectangle $(\widetilde R,Z)\in[0,4.1]\times[-1,1]$, the
$C^2$ norms of $(F_c-F^{\rm m})/F_0$ and $\Lambda(U_c^z-U^{z,\rm m})$
are $O(\Lambda^{-1})$, uniformly for $\Lambda\ge\Lambda_0$ and all
$C_*$ satisfying \eqref{eq:continuation-complex-bound}.
The stated core and exit bounds follow from this comparison, the weighted
slope estimates \eqref{eq:continuation-weighted}--\eqref{eq:continuation-relative},
and the argument in Subsection~\ref{app:core-exit}.

\subsection{Interval I.2: Leaving the core and making the shear small}
\label{sec:inner-leaving-core}

We first construct a smooth comparison profile from the exit data of
Proposition~\ref{prop:continuation-core}. Its normalized inertial stress
prescribes the direction of the shear on $R_a<R\le100$. Integrating
the radial derivatives then gives the actual velocity.
Figure~\ref{fig:inner-core-exit-stages} locates the admissible collar
within this connection.

\begin{figure}[H]
\centering
\begin{tikzpicture}[x=1cm,y=1cm,>=Stealth,
  every node/.style={font=\small}]
  \draw[->] (0,0)--(14.15,0) node[right] {$R$};
  \draw[line width=2pt,black!55] (.2,0)--(3.1,0);
  \draw[line width=2pt,blue!65!black] (3.1,0)--(6.4,0);
  \draw[line width=2pt,teal!70!black] (6.4,0)--(13.3,0);
  \foreach \x/\lab in {.2/0,3.1/R_a,4.25/R_{ac},6.4/R_ae^{h_b},13.3/100}
    {\draw (\x,.09)--(\x,-.09)
      node[below=3pt,font=\footnotesize] {$\lab$};}
  \node[align=center] at (1.65,.95)
    {analytic core\\$\mathcal{T}=0$};
  \draw[decorate,decoration={brace,amplitude=4pt},blue!65!black]
    (3.1,.38)--(6.4,.38);
  \node[align=center,blue!65!black] at (4.75,1.05)
    {reduce the shear\\$0<h_b\ll1$};
  \node[align=center,teal!70!black] at (9.85,1.05)
    {keep $|\mathcal{S}|/F$ small\\$\mathcal{S}$ nearly opposite to $\mathcal{I}$};
  \draw[<->,blue!65!black] (3.1,-.8)--(4.25,-.8);
  \node[align=center,blue!65!black,font=\footnotesize]
    at (3.675,-1.2) {admissible collar};
  \draw[<->] (3.1,-1.65)--(13.3,-1.65);
  \node[below] at (8.2,-1.65) {relaxed cone, $\mathcal{T}\ne0$ for $R>R_a$};
\end{tikzpicture}
\caption{Leaving the stress-free core; radial distances are schematic.
The admissible collar $(R_a,R_{ac}]$ lies inside the short interval
where the shear is reduced. Beyond $R_ae^{h_b}$ the shear is small,
while its direction remains nearly opposite to the inertial stress.
At $R=100$, the strict margin $\mathcal{I}^\theta/F>3>2$ prepares the removal
of the axial shear in the next interval.}
\label{fig:inner-core-exit-stages}
\end{figure}

\paragraph{The smooth comparison.}
Put $y=\log(R/R_a)$, where $R_a=4/\Lambda$.
At the exit write $f=F_c(R_a,\cdot)$ and $v=U_c^z(R_a,\cdot)$.
Use the data sizes
\begin{equation}\label{eq:inner-data-sizes}
\begin{aligned}
 A={}&10+\|\log(C_*F_c)\|_{C^3}+\|U_c^z\|_{C^3},\\
 K={}&10^6+100+R_a^{-1}+C_*+P_*+A
 +\|F_c\|_{C^3}+\|1/F_c\|_{C^3}
 +\|P_0\|_{C_Z^3}+\sum_i\|M_c^i\|_{C^3}.
\end{aligned}
\end{equation}
For each core field $w$ in \eqref{eq:inner-data-sizes}, $\|w\|_{C^3}$
is the joint $C^3$ norm of $(s,Z)\mapsto w(R_as,Z)$ on
$[0,e^{1/100}]\times[-1,1]$. The sum runs over
$i\in\{\theta,z,\theta z,z\theta,p\}$, and $\|P_0\|_{C_Z^3}$ is taken on $[-1,1]$.
For fixed $\Lambda$, \eqref{eq:inner-prepared-normalized} and the five
moment integrals give $A=O_\Lambda(1)$ and $C_*\le K=O_\Lambda(C_*)$.

\hypertarget{inner-bridge-width}{Choose} $K_1\ge1$ to bound the absolute
constants in the velocity, moment and stress estimates for the bridge and
the two short switches below. These constants are uniform for $0<c_*\le1$.
Fix $c_*$ independently of the core and set
\[
 0<c_*\le\min\left\{\frac{\epsilon_0}{10^{10}K_1},
                      \frac{\eta_{\rm tol}}{12K_1}\right\},
 \qquad h_b=c_*K^{-100}.
\]
For the actual profile, the shear multiplier $\chi_b$ defined in \eqref{eq:inner-chi}
decreases from $1$ to $h_b$ on $0\le y\le h_b$.
For the comparison profile, the cutoff $\alpha$ decreases from $1$ to $0$
on $h_b\le y\le 2h_b$.
The core continuation covers both intervals since
$R_a e^{2h_b}<R_a e^{1/100}<4.1/\Lambda$.
Use the smooth step $\sigma$ of
\eqref{eq:road-sigma}, with $0<\sigma'(s)\le8$ on $0<s<1$.
We denote the smooth comparison profile by $(\bar F,\bar U^z)$.
An overbar here refers to this profile or to a quantity computed from it.
In Section~\ref{sec:outer-construction}, the overbar on $\bar U$ denoted
the uncorrected outer candidate.
For $0\le y\le2h_b$, set $\alpha(y)=1-\sigma((y-h_b)/h_b)$ and prescribe
\begin{equation}\label{eq:inner-comparison}
 \partial_y\log\bar F=\alpha\,\partial_y\log F_c,\qquad
 \partial_y\bar U^z=\alpha\,\partial_y U_c^z,\qquad
 (\bar F,\bar U^z)|_{y=0}=(f,v).
\end{equation}
The comparison equals the core up to $y=h_b$. At
$r_*=R_a e^{2h_b}$ both radial derivatives turn off smoothly.
Keep $\bar F,\bar U^z$ constant in $R$ from $r_*$ to $110$, and put
$\bar U^\theta=\sqrt{2R}\bar F$. Continue the five moments from
their core values at $R_a$ and take $\bar P=P_0+\bar M^p$.
Their inertial stress defines
\[
 \mathbf q=(\bar D,\bar E)=\bar{\mathcal I}/\bar F.
\]
Write $\mathcal S/F=(-a,b)$, $a>0$, and $(D,E)=\mathcal I/F$.
For $0<\kappa\le2$, the relaxed cone test \eqref{eq:inner-small-shear-test}
reads $D-(b/a)E>2$. Later, in the actual construction
\eqref{eq:inner-actual-bridge}, we will take $\mathcal S/F=-\chi_b\mathbf q$
on $R_a<R\le100$, with $0<\chi_b<1$. Then $b/a=-\bar E/\bar D$
wherever $\bar D>0$. Suppose for the moment that the actual normalized
inertial stress $(D,E)$ equals $\mathbf q=(\bar D,\bar E)$. The test then
reduces to $\bar D+\bar E^2/\bar D>2$, independently of $\chi_b$.
The lower bound $9/4$ in the following lemma leaves room for
the difference $(D,E)-\mathbf q$, which we estimate in \eqref{eq:inner-bridge-error}.
The angular component alone supplies this margin where $\bar D$ is large,
and the axial component supplies it on the remaining part of the exit.

\begin{lemma}[Inertial stress of the smooth comparison]
\label{prop:inner-comparison}
Fix $j$ as in \eqref{roadin:j-budget} and the analytic data of
Proposition~\ref{prop:continuation-core}.
For all sufficiently large $\Lambda$ and, with $\Lambda$ fixed,
all sufficiently large $C_*$, the comparison on $R_a\le R\le110$,
$|Z|\le1$, satisfies
\begin{equation}\label{eq:inner-comparison-bounds}
 \bar D>0,\qquad
 \bar H:=\frac{\bar D^2+\bar E^2}{\bar D}\ge\frac94,\qquad
 \bar D\ge4\quad (100\le R\le110).
\end{equation}
There is a constant $B_\Lambda$, uniform in $C_*$, such that for every
$R_a\le R\le110$,
\[
 \|\bar D\|_{C_Z^2}+\|1/\bar D\|_{C_Z^2}
               +C_*^{-1}\|\bar E\|_{C_Z^2}\le B_\Lambda .
\]
\end{lemma}

\begin{proof}
Throughout this proof, $W,a,\zeta,Q$ and $J$ are evaluated on the
comparison profile, with $V=\bar U^z$ and $U=\bar U^\theta$.
The radial equations
\eqref{eq:normalized-inertial-stresses}--\eqref{eq:normalized-Q-J}
use centrifugal balance
$\partial_R\bar P=\bar F^2$ and the radial velocity determined by
incompressibility, so they apply with the inherited pressure
$\bar P=P_0+\bar M^p$. Set
\[
 H_V=\tfrac{1-\delta}{2}Z+LV,\quad
 \mathcal A w=(1-\delta)Zw+L\partial_Zw,\quad
 \mathcal P w=2(1+\delta)Zw-L\partial_Zw .
\]
Thus $W=1-\mathcal A(\bar M^z/R)$, $a=-2\partial_y\log\bar F$,
$\zeta=\partial_Z\log\bar F$. Write
\[
 \bar D=RQ/L_\delta,\qquad
 N=U^2J=L_\delta\sqrt{2/R}\,\bar{\mathcal I}^z,\qquad
 \bar E=\frac{\sqrt{R/2}}{L_\delta\bar F}\,N .
\]
We propagate $N$ so that the axial equation has no factor $1/\bar F$.
Put
\[
 g[V]=-\tfrac{1+\delta}{2}(1-2ZV)V-H_V\partial_ZV+\mathcal P P_0,
\]
so that $g[4Z+j]=g$ of \eqref{eq:g}.
Since $\partial_y U^2=(1-a)U^2$ and $\partial_yR=R$, the two
radial equations give
\begin{subequations}\label{eq:inner-comparison-ode}
\begin{align}
 \partial_yQ+(2-a/2)Q&=S_Q,\qquad
 S_Q=-W(1-a/2)-\tfrac{\delta}{2}(1-2ZV)-H_V\zeta,
 \label{eq:inner-Q-equation}\\
 \partial_y\bar D+(1-a/2)\bar D&=R S_Q/L_\delta,
 \label{eq:inner-D-equation}\\
 \partial_yN+N&=g[V]+\mathcal P(\bar P-P_0)+ZU^2-W\partial_yV.
 \label{eq:inner-N-equation}
\end{align}
\end{subequations}

\emph{Short transition.}
In this proof $B$ may depend on $j$ and the fixed analytic data,
and $B_\Lambda$ may also depend on $\Lambda$. Both are uniform in
$C_*$. The core estimates
\eqref{eq:inner-prepared-normalized} and integration over
$0\le y\le2h_b$ give
\[
 \|\log(\bar F/f)\|_{C_Z^2}+\|V-v\|_{C_Z^2}\le B_\Lambda h_b .
\]
Write $(w,\bar w)=(\log(C_*F_c),\log(C_*\bar F))$ or $(U_c^z,V)$
in the coordinates $(y,Z)$. For $0\le y\le2h_b$, integration by parts
in \eqref{eq:inner-comparison} gives
\[
 \bar w(y,Z)=\alpha(y)w(y,Z)+\int_0^y[-\alpha'(t)]w(t,Z)\,dt .
\]
Since $\alpha\ge0$, $-\alpha'\ge0$ and
$\alpha(y)+\int_0^y[-\alpha'(t)]\,dt=\alpha(0)=1$, this is a weighted
average of core values. The weights are independent of $Z$, so the same
average gives $\partial_Z^k\bar w$ in terms of $\partial_Z^k w$ for
$0\le k\le3$. Thus $C_*\bar F$, $(C_*\bar F)^{-1}$ and $V$ have
bounded $C_Z^3$ norms on the short layer and on the following constant
segment. Here and below the $Z$ derivatives are taken at fixed radius.
In particular $\bar M^z/R$ is bounded in $C_Z^3$, since it is the
radial average of $V$ and of the preceding core velocity. From
$\bar P-P_0=\int_0^R\bar F^2\,d\rho$, with the core used up to $R_a$,
we have on the whole interval
\[
 \|\bar P-P_0\|_{C_Z^3}+\|U^2\|_{C_Z^2}\le B_\Lambda C_*^{-2}.
\]
On the layer $a=-2\alpha\partial_y\log F_c$ and
$\partial_yV=\alpha\partial_yU_c^z$. The coefficients and right-hand
sides in \eqref{eq:inner-comparison-ode} therefore have $C_Z^2$ norms
bounded by $B_\Lambda$, using at most one radial derivative of the
core. The cutoff enters these equations through its value $\alpha$.

At $R_a$ the stress-free core gives
\[
 \bar D(R_a)=-2R_a\,\partial_R\log F_c(R_a),\qquad
 N(R_a)=-2L_\delta\,\partial_RU_c^z(R_a).
\]
The initial $C_Z^2$ norms are bounded by \eqref{eq:inner-prepared-normalized}.
Integration of the two linear equations, also after up to two $Z$
derivatives, gives
\[
 \|\bar D(R)-\bar D(R_a)\|_{C_Z^2}+
 \|N(R)-N(R_a)\|_{C_Z^2}\le B_\Lambda h_b
 \qquad(R_a\le R\le r_*=R_a e^{2h_b}).
\]
By Proposition~\ref{prop:continuation-core}, $\bar D(R_a)\ge R_a/32$,
$\bar D(R_a)>5/2$ where $\chi_c\ge99/100$, and $N(R_a)\ge3j/10$
where $\chi_c\le99/100$. Put $b_0=2j/5$.
For fixed $\Lambda$, $h_b\le c_*C_*^{-100}$, so we can take $C_*$
large enough that $B_\Lambda h_b\le\min\{R_a/64,1/8,j/30\}$.
Subtracting these errors from the three exit bounds gives, throughout the short transition,
\begin{equation}\label{eq:inner-comparison-transition}
 \bar D\ge R_a/64,\qquad
 \bar D\ge19/8\quad(\chi_c\ge99/100),\qquad
 N\ge 2b_0/3\quad(\chi_c\le99/100).
\end{equation}

\emph{Angular propagation on the constant segment.}
For $r_*\le R\le110$, $\bar F$ and $V$ retain their values at $r_*$, so
\[
 \bar F(R,Z)=f_*(Z):=\bar F(r_*,Z),\qquad
 V(R,Z)=v_*(Z):=V(r_*,Z).
\]
Thus $a=0$ and $\partial_yV=0$. With $m_z=M_c^z(R_a)$, the axial average is
\[
 \frac{\bar M^z(R)}R=\frac{R_a}{R}\,\frac{m_z}{R_a}+
                 \frac1R\int_{R_a}^R V(\rho)\,d\rho .
\]
By \eqref{eq:inner-core-axial-budget} and the short transition estimate,
both $\bar M^z/R$ and $v_*$ differ from $4Z+j$ by at most
$\rho_{\rm core}+B_\Lambda h_b$ in $C_Z^1$. Substitution in
$W=1-\mathcal A(\bar M^z/R)$ gives
\[
 -W-\tfrac{\delta}{2}(1-2Zv_*)
     =3-\tfrac{\delta}{2}+jZ+O(\rho_{\rm core})+O_\Lambda(h_b).
\]
At the exit, \eqref{eq:inner-core-axial-budget} gives
$H_v-H_0=L(v-U_0^z)=O(\Lambda^{-1})$. The definition of $G$ gives
$H_0G'=L_\delta\chi_c$ and $|G'|\le B\sqrt{\chi_c}$ for fixed $j$.
Hence \eqref{eq:inner-core-input} yields
\[
 \left|H_v\partial_Z\log f+\Lambda L_\delta\chi_c\right|
 \le |H_v|\,|\partial_Z\log f+\Lambda G'|
      +\Lambda|H_v-H_0|\,|G'|
 \le B\sqrt{\chi_c}+B/\Lambda .
\]
Since $L_\delta\ge1-\delta$, we have
$B\sqrt{\chi_c}\le\tfrac14\Lambda L_\delta\chi_c+B^2/[\Lambda(1-\delta)]$.
This gives the \hypertarget{inner-exit-transport}{angular transport bound}
\[
 -H_v\partial_Z\log f\ge\tfrac34\Lambda L_\delta\chi_c-B/\Lambda .
\]
The short transition bounds both $H_{v_*}-H_v=L(v_*-v)$ and
$\partial_Z\log(f_*/f)$ by $B_\Lambda h_b$ in absolute value.
Since $\|\partial_Z\log f\|_\infty\le B_\Lambda$, on $r_*\le R\le110$ we obtain
\[
 -H_V\zeta=-H_{v_*}\partial_Z\log f_*
 \ge\tfrac34\Lambda L_\delta\chi_c-B/\Lambda-B_\Lambda h_b .
\]
Choose $\Lambda$ sufficiently large that $B/\Lambda\le1/100$ and
the $O(\rho_{\rm core})$ error is at most $1/8$. For this $\Lambda$
and all sufficiently large $C_*$ the transition errors total at most $1/8$.
As $|j|\le1/20$, \eqref{eq:inner-Q-equation} now gives
\[
 S_Q\ge2+\tfrac34\Lambda L_\delta\chi_c.
\]
Since $a=0$ on this segment, \eqref{eq:inner-D-equation} gives
$\partial_R(R\bar D)=R S_Q/L_\delta\ge2R$.
The function $d(R)=R-r_*+R_a/64$ satisfies $d(r_*)\le\bar D(r_*)$
and $\partial_R(Rd)=2R-r_*+R_a/64\le2R$. Thus $\bar D\ge d\ge R_a/64$
up to $110$, with $\bar D>99$ on $100\le R\le110$.
On $\chi_c\ge99/100$, the same equation and $R\ge R_a=4/\Lambda$ give
\[
 \partial_R(R\bar D)=R S_Q/L_\delta
 \ge\tfrac34 R\Lambda\chi_c\ge3\chi_c>9/4 .
\]
Since $\bar D(r_*)\ge19/8$ by \eqref{eq:inner-comparison-transition},
the constant $9/4$ is a lower solution throughout the strong angular region.

\emph{Axial propagation and derivative bounds.}
On $\chi_c\le99/100$, \eqref{eq:width} gives $g[4Z+j]\ge b_0$.
The $C_Z^1$ distance of $v_*$ from $4Z+j$ is at most
$\rho_{\rm core}+B_\Lambda h_b$. The definition of $g[V]$
thus gives $g[v_*]\ge7b_0/8$ for
sufficiently large $\Lambda$ and then $C_*$. The pressure and
$ZU^2$ terms have total size at most $B_\Lambda C_*^{-2}\le b_0/8$.
Thus $\partial_yN+N\ge3b_0/4$ on this region from $r_*$ to $110$.
The constant $2b_0/3$ is a lower solution, and
\eqref{eq:inner-comparison-transition} gives this lower bound up to $r_*$.
Thus $N\ge2b_0/3$ on $R_a\le R\le110$. Since $\bar F\le B_\Lambda/C_*$
and $L_\delta\le1$, we obtain on this whole interval, for $\chi_c\le99/100$,
\[
 \bar E=\frac{\sqrt{R/2}}{L_\delta\bar F}\,N
 \ge\frac{2b_0}{3B_\Lambda}\sqrt{\frac{R_a}{2}}\,C_*
 \ge\frac98
\]
for all sufficiently large $C_*$ with $\Lambda$ fixed.
On this region $\bar H=\bar D+\bar E^2/\bar D\ge2|\bar E|\ge9/4$.
On the complementary angular region $\bar H\ge\bar D\ge9/4$.

Finally, on $r_*\le R\le110$ the coefficients and sources in
\eqref{eq:inner-comparison-ode} have bounded $C_Z^2$ norms by the
velocity, average, and pressure bounds above. The logarithmic
length is at most $\log(110/R_a)$, fixed with $\Lambda$.
Integrating the equations with the transition data bounds
$\|\bar D\|_{C_Z^2}+\|N\|_{C_Z^2}$ by $B_\Lambda$. The lower bound
$\bar D\ge R_a/64$ controls $1/\bar D$ and its first two $Z$
derivatives. The formula for $\bar E$ and the bound for
$(C_*\bar F)^{-1}$ give the stated estimate for $C_*^{-1}\bar E$.
\end{proof}

To join the axial velocity to $4Z$, use $j=\eta_{\rm tol}/8$ from
\eqref{roadin:j-budget}. With these data fixed, take $\Lambda$
sufficiently large for Lemma~\ref{prop:inner-comparison} and for
the bound in \eqref{eq:inner-core-axial-budget} to give
\begin{equation}\label{eq:inner-axial-Lambda}
 \rho_{\rm core}\le B/\Lambda\le\eta_{\rm tol}/4.
\end{equation}
With this $\Lambda$ fixed, every $C_*$ above a common lower bound
satisfies the requirements of the lemma and gives
\[
 \bar D\ge (2K)^{-1},\qquad |\mathbf q|\le2K^3,\qquad
 \|\mathbf q\|_{C_Z^2}\le K^7,\qquad \bar H\le K^{10}.
\]
Indeed, $K\ge C_*$, $\bar D^{-1}=O_\Lambda(1)$ and
$\|\mathbf q\|_{C_Z^2}=O_\Lambda(1+C_*)$. The integration in
\eqref{eq:inner-comparison} also gives
$\|\log(\bar F/f)\|_{C_Z^2}+\|\bar U^z-v\|_{C_Z^2}\le BKh_b$
directly from the core norms in \eqref{eq:inner-data-sizes}.
Consequently $\|\bar F\|_{C_Z^2}+\|1/\bar F\|_{C_Z^2}
+\|\bar U^z\|_{C_Z^2}\le BK$. In these elementary bridge estimates
$B$ denotes an absolute constant covered by $K_1$.

\paragraph{The actual velocity.}
Set
\begin{equation}\label{eq:inner-chi}
	\chi_b(y)=1-(1-h_b)\sigma(y/h_b),
\end{equation}
and, starting with the core values at $R_a$, solve the two explicit
integration equations
\begin{equation}\label{eq:inner-actual-bridge}
	\partial_y\log F=-\frac{\chi_b}2\bar D,
	\qquad
	\partial_y U^z=-\chi_b\sqrt{R/2}\,F\bar E,
	\qquad R_a<R\le100.
\end{equation}
First integrate the angular equation; then its positive solution
$F$ is a known coefficient in the axial equation. Writing
the functions in the coordinate $y$, this is the explicit ansatz
\[
\begin{aligned}
F(y,Z)&=f(Z)\exp\left[-\frac12\int_0^y
\chi_b(t)\bar D(t,Z)\,dt\right],\\
U^z(y,Z)&=v(Z)-\int_0^y\chi_b(t)\sqrt{R_a e^t/2}\,
F(t,Z)\bar E(t,Z)\,dt.
\end{aligned}
\]
These prescriptions give exactly $\mathcal{S}=-\chi_b F\mathbf q$.
They retain both components of the comparison stress near the
core, where the axial component may supply the exit margin.

For $y\le h_b$ the comparison is the stress-free core, and the
prescription with $\chi_b=1$ is exactly its velocity equation.
Since $1-\chi_b$ is flat at zero, all derivatives match the core
at $R_a$. For every $y>0$, $0<\chi_b<1$.

\paragraph{The error estimate and the collar.}
We compare the actual normalized inertial stress $\mathcal I/F$ with
$\mathbf q$ to verify the relaxed cone on $(R_a,100]$ and the admissible
cone on a collar immediately outside $R_a$. The actual and comparison
moments start from the same core values at $R_a$.
Write
\begin{equation}\label{eq:inner-ray-notation}
 \mathcal{S}=-\chi_b F\mathbf q,\qquad \mathcal{I}/F=\mathbf q+\mathbf e,\qquad
 \mathbf q=(\bar D,\bar E),\qquad \bar H=|\mathbf q|^2/\bar D .
\end{equation}
We prove
\begin{equation}\label{eq:inner-bridge-error}
	|\mathbf e|\le 2K_1K^{20}h_b(1-\chi_b).
\end{equation}
Since $\mathcal T/F=(1-\chi_b)\mathbf q+\mathbf e$, the leading term
vanishes as $R\downarrow R_a$. We retain the factor $1-\chi_b$ in the
velocity differences and carry it through the moment formulas to the stress error.

For $0\le y\le h_b$, put $\omega=1-\chi_b$ and
$\ell=\log(F/\bar F)$. Here the comparison equals the stress-free
core. Its equations are
$\partial_y\log F_c=-\bar D/2$ and
$\partial_y U_c^z=-\sqrt{R/2}F_c\bar E$. Thus subtraction gives the
exact identities
\[
 \partial_y\ell=\tfrac12\omega\bar D,\qquad
 \partial_y(U^z-\bar U^z)=(\chi_b e^\ell-1)\partial_y U_c^z,
 \qquad \ell(0)=0,\quad (U^z-\bar U^z)(0)=0.
\]
The weight $\omega$ is nondecreasing and independent of $Z$, so
\[
 \int_0^y\omega(t)\,dt\le h_b\omega(y),\qquad
 \|\ell(y)\|_{C_Z^2}\le BKh_b\omega(y).
\]
For the second bound we used
$\bar D=-2\partial_y\log F_c$ on this interval and the mixed
$C^3$ core bound. The elementary identity
\[
 \chi_b e^\ell-1=-\omega+\chi_b(e^\ell-1)
\]
and the product rule now give
$\|\chi_b e^\ell-1\|_{C_Z^2}\le B\omega(y)$, since $Kh_b$ is small.
Integrating the axial difference and using
$\|\partial_y U_c^z\|_{C_Z^2}\le BK$ therefore yields
\[
 \|U^z-\bar U^z\|_{C_Z^2}\le BKh_b\omega(y),\qquad
 \|F-\bar F\|_{C_Z^2}
 =\|\bar F(e^\ell-1)\|_{C_Z^2}\le BK^2h_b\omega(y).
\]
This proves the required weighted control on the flat part of the
cutoff, including the derivatives needed in the moment formulas.

For $h_b\le y\le\log(100/R_a)$ we have
$\chi_b=h_b$ and $\omega=1-h_b\ge1/2$.
The integral of the cutoff satisfies
\[
 \int_0^y\chi_b(t)\,dt
 \le h_b(1+\log(100/R_a))\le Kh_b.
\]
Combining this with the explicit angular formula and
$\|\mathbf q\|_{C_Z^2}\le BK^7$ gives
\[
 \|\log(F/f)\|_{C_Z^2}\le BK^8h_b.
\]
This is small, so $\|F\|_{C_Z^2}+\|1/F\|_{C_Z^2}\le BK$.
The explicit axial integral in turn gives
\[
 \|U^z-v\|_{C_Z^2}\le BK^9h_b.
\]
Together with the comparison estimates above, these bounds
also show, on the entire interval,
\[
 \|\ell\|_{C_Z^2}\le BK^8h_b.
\]
In particular, the angular estimate is a relative one; it does not
replace the small core amplitude by an absolute error independent
of that amplitude.

To estimate the moment differences at $R=R_a e^y$, let $X(y)$ measure the
$C_Z^1$ differences between the actual and comparison profiles up to that radius.
Here $\Delta$ denotes the actual quantity minus its comparison value
(in the viscous terms, $\Delta$ denotes the Laplacian).
With $\Delta F=F-\bar F$ and $\Delta V=U^z-\bar U^z$, set
\[
 X(y)=\sup_{0\le t\le y}
 \bigl(\|\Delta F(t,\cdot)\|_{C_Z^1}+\|\Delta V(t,\cdot)\|_{C_Z^1}\bigr).
\]
The early bound retains $\omega(y)$ by monotonicity. The later
bound retains it because $\omega(y)\ge1/2$. With a spare power
of $K$, both can therefore be written as
\[
 X(y)\le BK^{10}h_b\omega(y).
\]
The moments agree at $R_a$. Subtracting the integrands in \eqref{fiveM}--\eqref{moment-F}
gives the following identities for $\Delta\mathbf M=\mathbf M-\bar{\mathbf M}$.
\[
\begin{aligned}
 \Delta M^\theta&=2\int_{R_a}^R \rho\Delta F\,d\rho,&
 \Delta M^z&=\int_{R_a}^R\Delta V\,d\rho,\\
 \Delta M^{\theta z}&=2\int_{R_a}^R\rho
          (U^z\Delta F+\bar F\Delta V)\,d\rho,&
 \Delta M^p&=\int_{R_a}^R(F+\bar F)\Delta F\,d\rho,\\
 \Delta M^{z\theta}&=\int_{R_a}^R
   [(U^z+\bar U^z)\Delta V-\rho(F+\bar F)\Delta F]\,d\rho .
\end{aligned}
\]
On $R_a\le R\le100$, the $C_Z^1$ product rule and
$P-\bar P=\Delta M^p$ give
\[
 \sum_j\|M^j-\bar M^j\|_{C_Z^1}
       +\|P-\bar P\|_{C_Z^1}\le BKX(y).
\]
In the moment representations
\eqref{Itheta-moments}--\eqref{Iz-moments}, retain the pressure as
$P_0+M^p$ in the last term, as follows from
\eqref{Iz-moments-direct} and \eqref{Ur-Uz}.
Divide each stress by its own regular swirl. The angular fraction then has denominator
$2L_\delta RF$, numerator of size $BK^2$, and numerator difference
at most $BKX$. For the axial fraction the corresponding bounds are
$L_\delta\sqrt{2R}F$, $BK^3$, and $BK^2X$. These estimates follow
from the moment differences just displayed, the common bound $BK$
on both velocities, and the common $P_0$. Now
$|1/F-1/\bar F|\le BK^2X$, $R^{-1}\le K$ and $L_\delta^{-1}\le2$.
For either numerator $n$ and its comparison value $\bar n$, use
\[
 \frac{n}{F}-\frac{\bar n}{\bar F}
 =\frac{n-\bar n}{F}+\bar n\left(\frac1F-\frac1{\bar F}\right).
\]
Including the radial factors, the angular fraction changes by at most
$B R^{-1}(K^2+K^4)X\le BK^5X$, and the axial fraction by at most
$B R^{-1/2}(K^3+K^5)X\le BK^{11/2}X$.
The remaining angular term changes by $\mathcal A\Delta M^z/L_\delta$,
which is bounded by $BKX$. Hence
\[
 |\mathbf e|\le BK^6X(y)
       \le BK^{16}h_b\omega(y).
\]
This proves \eqref{eq:inner-bridge-error}: the prior choice of
$K_1$ covers the fixed constants, and $K^{16}\le K^{20}$.

For $R_a<R\le100$, we have $0<\chi_b<1$ and $\kappa=\chi_b\bar H$.
From $|\mathbf q|\ge\bar D\ge(2K)^{-1}$,
$|\mathbf q|\le2K^3$, $\bar H\le K^{10}$, and
$2K_1K^{20}h_b=2K_1c_*K^{-80}$, the
\hyperlink{inner-bridge-width}{choice of $c_*$} and
\eqref{eq:inner-bridge-error} give
\begin{equation}\label{eq:inner-ray-error}
 |\mathbf e|\le \min\!\left\{
 \frac{\bar D}{16(1+|\mathbf q|)},\,\frac{(1-\chi_b)|\mathbf q|}{10\sqrt{1+\bar H}}\right\},
 \qquad \mathcal{I}^\theta/F=\bar D+e_1>0.
\end{equation}
The factor $1-\chi_b$ in the second bound is essential near $R_a$, where
the leading term in $\mathcal T/F=(1-\chi_b)\mathbf q+\mathbf e$ vanishes.
It keeps $|\mathbf e|$ small relative to $(1-\chi_b)|\mathbf q|$ even as
$R\downarrow R_a$, which controls the stress direction in the angle estimate below.
In particular,
\[
 \frac{|\mathcal T|}{F}\ge(1-\chi_b)|\mathbf q|-|\mathbf e|
 \ge\tfrac9{10}(1-\chi_b)|\mathbf q|>0
 \qquad (R_a<R\le100).
\]
For $\kappa\le2$, apply the small-shear test \eqref{eq:inner-small-shear-test}
with $(a,b)=\chi_b(\bar D,-\bar E)$. After division by $a$, its left
side is $\bar H+\mathbf e\cdot\mathbf q/\bar D\ge\bar H-1/16\ge35/16>2$.

For $\kappa>2$, use $\mathcal{T}/F=(1-\chi_b)\mathbf q+\mathbf e$ and
$\mathcal{S}=-\chi_b F\mathbf q$. By \eqref{eq:inner-ray-error},
$-\mathcal{T}\cdot \mathcal{S}\ge\tfrac9{10}\chi_b F^2(1-\chi_b)|\mathbf q|^2>0$ and
\begin{equation}\label{eq:inner-ray-angle}
 \frac{|\mathcal{T}\cdot \mathcal{S}^\perp|}{-\mathcal{T}\cdot \mathcal{S}}
 \le \frac{|\mathbf e|}{(1-\chi_b)|\mathbf q|-|\mathbf e|}
 \le \frac{1}{9\sqrt{1+\bar H}}<\sqrt{\frac{2}{\kappa-2}}.
\end{equation}
Thus the relaxed cone holds on $(R_a,100]$. For the inner collar, set
\[
t_{ac}=\sigma^{-1}\!\left(\frac{1}{16K^{10}}\right),
\qquad R_{ac}=R_a\exp(h_b t_{ac}).
\]
For $R_a<R\le R_{ac}$, $(1-\chi_b)\bar H\le1/16$ by the choice of
$t_{ac}$ and $\bar H\le K^{10}$. Hence $\kappa\ge35/16$
and the cone is admissible there. All core derivatives match at $R_a$.
Since $t_{ac}<1$, the shear reduction continues past $R_{ac}$ to
$R_a e^{h_b}$. For $R_a e^{h_b}\le R\le100$ the multiplier equals $h_b$, so
\[
 \frac{|\mathcal S|}{F}=h_b|\mathbf q|\le2c_*K^{-97},\qquad
\kappa=h_b\bar H\le c_*K^{-90}<1.
\]

For the inner-edge derivative estimates in Section~\ref{sec:lower-order},
fix the constructed profile and take $0<y<h_b/2$. By \eqref{eq:road-sigma},
$\omega(y)=(1-h_b)\sigma(y/h_b)$ is $e^{-h_b^2/y^2}$ times a smooth
factor with a positive value at $y=0$.
The weight is independent of $Z$, and each $y$ derivative introduces at
most a finite power of $y^{-1}$. The difference equations for $\ell$ and
$\Delta V$ have zero initial data. Integrating the first and using
$e^\ell-1=\ell\int_0^1e^{s\ell}\,ds$ in the second gives the same
exponential factor for $\Delta F=\bar F(e^\ell-1)$ and $\Delta V$,
including their $Z$ derivatives. Their $y$ derivatives satisfy the
corresponding bounds by the same equations. The moment differences start
from zero at $R_a$, and $F,\bar F$ are positive there, so the moment
integrals and the quotient identity above give, for each $i,k\ge0$ and
some finite $m_{i,k}$,
\begin{equation}\label{eq:inner-flat-error}
 |\partial_y^i\partial_Z^k\mathbf e|
       \le C_{i,k}e^{-h_b^2/y^2}y^{-m_{i,k}} .
\end{equation}
The constants may depend on the constructed profile and on $h_b$.

\paragraph{Data for Interval I.3 at $R=100$.}
Continue $F,U^z$ and all five cumulative moments from their values at
$R=100$, with the same axis pressure $P_0$ and $P=P_0+M^p$.
Near this radius the incoming shear is
$\mathcal{S}/F=-h_b(\bar D,\bar E)$, where the comparison
direction $(\bar D,\bar E)$ extends smoothly across $100$.
Write $(D,E)=\mathcal{I}/F$. At $R=100$, uniformly in $Z$,
\[
 D>3,\qquad \kappa=h_b\bar H<1,\qquad
 D+\frac{\bar E}{\bar D}E\ge\bar H-1/16\ge35/16
\]
by \eqref{eq:inner-comparison-bounds}--\eqref{eq:inner-ray-error}.
Thus both endpoints of the affine test $D+t\bar E E/\bar D$, $0\le t\le1$,
have a strict margin when the axial shear is turned off.
The velocity estimates above give
\[
 \left\|\log\frac{F(100)}{f}\right\|_{C_Z^2}
 \le BK^8h_b,\qquad
 \max\left\{\|U^z(100)-4Z\|_{C_Z^2},
             \left\|\frac{M^z(100)}{100}-4Z\right\|_{C_Z^2}\right\}
 <\epsilon_0.
\]
Indeed, $M^z(100)=M_c^z(R_a)+\int_{R_a}^{100}U^z(R,Z)\,dR$.
The two axial errors are bounded by
$j+\rho_{\rm core}+2K_1K^{10}h_b$. The last term is at most
$2K_1c_*$, so the \hyperlink{inner-bridge-width}{choice of $c_*$},
\eqref{roadin:j-budget} and \eqref{eq:inner-axial-Lambda} give the claim.

\subsection{Interval I.3: Joining the reference velocity}
\label{subsec:inner-reference-velocity}
\label{subsec:inner-reference-connection}

We first remove the axial shear and restore the angular power $R^{1/10}$
on $[100,110]$. We then reshape the swirl, wait with the axial velocity
fixed, and restore the axial velocity to $4Z$. All these steps retain
the relaxed cone.
Figure~\ref{fig:inner-reference-stages} records their order and the
role of the waiting interval.

\begin{figure}[H]
\centering
\begin{tikzpicture}[x=1cm,y=1cm,>=Stealth,
  every node/.style={font=\footnotesize}]
  \draw[->] (0,0)--(14.15,0) node[right] {$R$};
  \draw[line width=2pt,violet!75!black] (.3,0)--(2.3,0);
  \draw[line width=2pt,orange!85!black] (2.3,0)--(4.3,0);
  \draw[line width=2pt,black!40] (4.3,0)--(6.1,0);
  \draw[line width=2pt,red!70!black] (6.1,0)--(8.6,0);
  \draw[line width=2pt,teal!70!black] (8.6,0)--(11,0);
  \draw[line width=2pt,violet!75!black] (11,0)--(13.6,0);
  \foreach \x/\lab in {.3/100,2.3/100e^{h_b},4.3/100e^{2h_b},6.1/110,8.6/R_{\theta\mathrm{ok}},11/R_{4Z},13.6/eR_{4Z}}
    {\draw (\x,.07)--(\x,-.07)
      node[below=4pt,font=\scriptsize] {$\lab$};}
  \node[align=center,violet!75!black] at (1.3,1.1)
    {flatten $U^z$\\$b\to0$};
  \draw[->,violet!75!black] (1.3,.68)--(1.3,.12);
  \node[align=center,black!65] at (5.2,1.1)
    {short buffer\\$U^\theta\propto R^{1/10}$};
  \draw[->,black!65] (5.2,.68)--(5.2,.12);
  \node[align=center,teal!70!black] at (9.8,1.1)
    {long wait\\$U^z=U_1^z(Z)$, $U^\theta$ grows};
  \draw[->,teal!70!black] (9.8,.68)--(9.8,.12);
  \node[align=center,orange!85!black] at (3.3,-1.15)
    {restore the power $R^{1/10}$\\$a\to4/5$};
  \draw[->,orange!85!black] (3.3,-.73)--(3.3,-.12);
  \node[align=center,red!70!black] at (7.35,-1.15)
    {restore $Z$ dependence\\and amplitude};
  \draw[->,red!70!black] (7.35,-.73)--(7.35,-.12);
  \node[align=center,violet!75!black] at (12.3,-1.15)
    {restore $U^z$\\to $4Z$};
  \draw[->,violet!75!black] (12.3,-.73)--(12.3,-.12);
  \node at (7,-2.05)
    {$R_{\theta\mathrm{ok}}=110e^T,\qquad
      R_{4Z}=e^{-8}R_{\rm ref},\qquad eR_{4Z}=e^{-7}R_{\rm ref}$.};
\end{tikzpicture}
\caption{Restoring the reference velocities; radial distances are
schematic and the short transitions are enlarged. Axial shear is
removed before the angular power is restored. During the wait,
the growth of the reference swirl reduces the normalized pressure
and angular moment defects, as well as the angular energy contribution
to the remaining defects. The wait also controls the inherited axial
inertial stress before $U^z$ is restored.
Both prescribed velocity components have their reference values
from $eR_{4Z}$ onwards.}
\label{fig:inner-reference-stages}
\end{figure}

\paragraph{Removing the axial shear and setting the angular slope.}
Use the shear coordinates \eqref{eq:outer-cone-variables} for the actual
velocity: $\mathcal{S}=F(-a,b)$, $D=\mathcal{I}^\theta/F$, and $w=\mathcal{I}^z/\mathcal{I}^\theta$.
Thus $a=1-2\partial_y\log U^\theta$ and $b=2\partial_y U^z/U^\theta$, where
$\partial_y=R\partial_R$.
On $100\le R\le100e^{h_b}$, keep $a=h_b\bar D$ and replace
$b=-h_b\bar E$ by
\[
b=-h_b\bar E\,[1-\sigma(\log(R/100)/h_b)].
\]
On the next interval, $100e^{h_b}\le R\le100e^{2h_b}$, keep $b=0$
and interpolate $a$ from $h_b\bar D$ to $4/5$ with
the same flat step. Integrate $\partial_y\log F=-a/2$ and
$\partial_y U^z=b\sqrt{2R}F/2$ as before. Thereafter keep
$(a,b)=(4/5,0)$ up to $R=110$.

The two switches have total logarithmic length $2h_b$. On both,
$\bar F,\bar U^z$ are constant in $R$, and the prescriptions have
$\|a\|_{C_Z^2}+\|b\|_{C_Z^2}\le B$ by
$h_b\|\mathbf q\|_{C_Z^2}\le c_*K^{-93}$.
Integration changes $\log F$ by at most $Bh_b$ in $C_Z^2$ and
$U^z$ by at most $BKh_b$. The moment and normalized-stress
difference estimates used for \eqref{eq:inner-bridge-error} hence
continue to hold with size $BK^{16}h_b$ on both
switches. In particular, writing $\mathbf e=(D,E)-(\bar D,\bar E)$
here, we have $|\mathbf e|\le1/16$ by the same smallness choice.
Since Lemma~\ref{prop:inner-comparison} gives $\bar D\ge4$,
$D>3$ on both switches.

During the first short interval, the exact relaxed test is
\[
D+\frac{t\bar E E}{\bar D}>2,
\qquad E= \mathcal{I}^z/F,\qquad 0\le t\le1.
\]
Here $\bar E E=\bar E^2+\bar E e_2\ge-e_2^2/4\ge-\bar D/4$,
so the left-hand side exceeds $3-1/4>2$. Also
$\kappa=h_b(\bar D+t^2\bar E^2/\bar D)
\le h_b\bar H<1$. On the second interval $b=0$
and $a$ lies between $h_b\bar D$ and $4/5$,
both in $(0,1)$. The test is then $D>2$.

Write $U_1^\theta=U^\theta(110,\cdot)$ and $U_1^z=U^z(110,\cdot)$.
At the core exit the sum of the two axial errors is at most
$2(j+\rho_{\rm core})<\epsilon_0$ by the choice $j=\eta_{\rm tol}/8$
in \eqref{roadin:j-budget} and the bound $\rho_{\rm core}\le\eta_{\rm tol}/4$
in \eqref{eq:inner-axial-Lambda}. The comparison and the two short
switches change these errors by $O(K^{20}h_b)$. On the subsequent
constant axial branch, radial averaging preserves the resulting bounds.
For the angular transport, the \hyperlink{inner-exit-transport}{angular transport bound}
in the proof of Lemma~\ref{prop:inner-comparison} gives
$H_v\partial_Z\log f\le B/\Lambda\le1/100$.
Up to the end of the two short switches, the bounds for
$\log(F/f)$ and $U^z-v$ above show that their first $Z$ derivatives
and their values tend to zero for fixed $\Lambda$ as $C_*$ increases.
Thereafter $a=4/5,\ b=0$ preserve both $U^z$ and $\partial_Z\log F$.
Since $\|\partial_Z\log f\|_\infty\le B_\Lambda$, the change from
$H_v\partial_Z\log f$ to $H(110,\cdot)\partial_Z\log U_1^\theta$
also tends to zero. The choice of $j$, followed by sufficiently
large $\Lambda$ and $C_*$, gives
$H\partial_Z\log U^\theta\le1/10$ throughout
$[100e^{2h_b},110]$.

Apply the radial equations \eqref{eq:inner-comparison-ode} to the
actual field by taking $U=U^\theta$, $V=U^z$, $H=H_V$,
$W=1-\mathcal A(M^z/R)$ and $P=P_0+M^p$. The two normalizations
are $Q=L_\delta D/R$ and $N=(U^\theta)^2J$ with the actual
inertial stress in \eqref{eq:normalized-inertial-stresses}.

On $[100e^{2h_b},110]$, we have $(a,b)=(4/5,0)$.
The axial velocity and average estimates give
$|W-(-3+4\delta Z^2)|\le10\epsilon_0$. Hence
\eqref{eq:inner-Q-equation} yields
\[
 S_Q=-\frac35W-\frac\delta2(1-2ZU^z)
     -H\partial_Z\log U^\theta\ge\frac75,
\]
using $\delta\le1/200$, $|U^z|\le5$, and the transport bound above.
Equation \eqref{eq:inner-D-equation} for the actual field becomes
$\partial_yD+\frac35D=RS_Q/L_\delta$.
Since $R\ge100$, the constant $3$ is a strict subsolution.
The bound $D>3$ at the end of the two switches therefore persists
up to $110$, and the relaxed cone follows from $b=0$ and $a=4/5$.
For the swirl amplitude at $110$, use $U_1^\theta=\sqrt{220}\,F(110)$ and
$f=F_c(R_a,\cdot)$ to write
\[
\log\!\bigl(C_*(1+Z^2)U_1^\theta\bigr)
=\log(C_*f)+\log\!\bigl(\sqrt{220}(1+Z^2)\bigr)+\log(F(110)/f).
\]
The first term has $C_Z^2$ norm at most $A-10$ by \eqref{eq:inner-data-sizes},
and the second has norm less than $7$. At $100e^{2h_b}$ the bound for
$\log(F/f)$ is $O(K^8h_b)$ in $C_Z^2$ by the incoming estimate at $100$
and the two short switches, with $K^8h_b=c_*K^{-92}$. From there to $110$, $\partial_y\log F=-2/5$
adds only the constant $-\tfrac25\log(110/(100e^{2h_b}))$.
Together these bounds give
\begin{equation}\label{eq:inner-110-data}
	\begin{gathered}
		D(110)>3,\qquad (a,b)=(4/5,0)\quad(100e^{2h_b}\le R\le110),\\
		\|U_1^z-4Z\|_{C^2_Z}
		+\|M^z(110)/110-4Z\|_{C^2_Z}\le2\epsilon_0,\\
		\|\log(C_*(1+Z^2)U_1^\theta)\|_{C^2_Z}\le2A,\qquad
		H(110,\cdot)\partial_Z\log U_1^\theta\le\tfrac1{10}.
	\end{gathered}
\end{equation}
For the axial velocity, the comparison on $R_a<R\le100$ and the two short
shear changes, together with the core bound
\eqref{eq:inner-core-axial-budget}, give
\begin{equation}\label{eq:inner-bridge-axial-budget}
	\begin{aligned}
		\rho_{\rm br}:=\|U_1^z-v\|_{C^2_Z}
		&\le 2K_1K^{10}h_b+K_1K^{20}h_b\\
		&\le3K_1c_*K^{-80}\le\frac{\eta_{\rm tol}}4,\\
		\eta:=\|U_1^z-4Z\|_1
		&\le j+\rho_{\rm core}+\rho_{\rm br}
		\le\frac{5\eta_{\rm tol}}8
		<\frac{\mathfrak e_*}{100}.
	\end{aligned}
\end{equation}
Here $\|q\|_1=\sup|q|+\sup|\partial_Z q|$ is bounded by the summed
$C^2_Z$ norm.

To estimate the five moments inherited at $R=110$ in Interval I.4, we
also need the velocity on the whole interval up to this radius. With
$j$, the analytic data and $\Lambda$ fixed, uniformly for all sufficiently
large $C_*$ we have
\begin{equation}\label{eq:imc-family-short-velocity}
 \sup_{0\le R\le110}
 \bigl(C_*\|F(R,\cdot)\|_1+\|U^z(R,\cdot)\|_1\bigr)\le B_\Lambda .
\end{equation}
On the core this follows from \eqref{eq:inner-prepared-normalized}.
For $R_a\le R\le100$, integrating \eqref{eq:inner-actual-bridge} gives
$\|\log(F/f)\|_1+\|U^z-v\|_1\le BK^9h_b$.
Since $h_b=c_*K^{-100}$, the identity $C_*F=(C_*f)\exp(\log(F/f))$
preserves the bound after one $Z$ derivative. On the two switches
$\|a\|_1+\|b\|_1\le B$, and afterwards $(a,b)=(4/5,0)$ up to $110$.
Integrating
$\partial_y\log F=-a/2$ and then $\partial_y U^z=b\sqrt{2R}F/2$
over this fixed logarithmic interval gives the same bounds up to $110$.

By \eqref{eq:imc-family-short-velocity} and the moment definitions
\eqref{fiveM}--\eqref{moment-F}, $M^p(110)=\int_0^{110}F^2\,dR=O_\Lambda(C_*^{-2})$
in $C_Z^1$, while $M^z(110)$ and $M^{z\theta}(110)$ are bounded in
the same norm independently of $C_*$. Thus $P(110)=P_0+M^p(110)$
is bounded in $C_Z^1$, and \eqref{Ur-Uz} bounds $U^r(110)$ uniformly
in $Z$. Substitution into \eqref{Iz-moments-direct} gives
$N(110)=L_\delta\mathcal I^z(110)/\sqrt{55}=O_\Lambda(1)$.
Since $K\ge C_*$, taking $C_*$ sufficiently large gives the
\hypertarget{inner-inherited-bounds}{bounds at $R=110$}
\[
 \|M^p(110)\|_{C_Z^1}\le1,\qquad \sup_Z|N(110,Z)|\le20K.
\]

\paragraph{Reshaping the swirl.}
Set $T=400A$ and $R_{\theta\mathrm{ok}}=110e^T$ for the shaping interval,
and $R_{4Z}=e^{-8}R_{\rm ref}$ for the subsequent axial restoration.
The required scale conditions are
\begin{equation}\label{eq:inner-compatibility}
	\begin{gathered}
		R_{\rm ref}=110(C_*P_*)^{10},\qquad C_*\ge e^{4A},\\
		R_{\rm ref}\ge110e^{T+10}(1+A)^{10},\qquad
		R_{4Z}\ge110(1+K)^2/P_*^2.
	\end{gathered}
\end{equation}
The last inequality controls the inherited axial inertial stress when
$U^z$ is restored.
For $y=\log(R/110)$, $0\le y\le T$, prescribe
\begin{equation}\label{eq:inner-swirl-shape}
	U^\theta(R,Z)=e^{y/10}
	[U_1^\theta(Z)]^{1-\sigma(y/T)}
	\left[\frac1{C_*(1+Z^2)}\right]^{\sigma(y/T)},
	\qquad U^z(R,Z)=U_1^z(Z).
\end{equation}
Since $\sigma$ is flat at zero, this prescription matches the incoming
branch $U^\theta=e^{y/10}U_1^\theta$, $U^z=U_1^z$ to every radial
order at $R=110$.
If $B(Z)=\log(C_*(1+Z^2)U_1^\theta(Z))$, then
\[
a=\frac45+\frac2T(\partial_s\sigma)(y/T)B(Z),\qquad
\zeta:=\partial_Z\log U^\theta
=(1-\sigma(y/T))\partial_Z\log U_1^\theta
-\sigma(y/T)\frac{2Z}{1+Z^2}.
\]
By \eqref{eq:inner-110-data}, $\|B\|_{C^2_Z}\le2A$, so $T=400A$ and
$|\sigma'|\le8$ give $|a-4/5|\le 2/25$. To bound $H\zeta$, observe that
$H$ is fixed in $R$ while $U^z=U_1^z$. At the incoming slope,
$H\partial_Z\log U_1^\theta\le1/10$ by \eqref{eq:inner-110-data}.
For the reference slope, \eqref{eq:core-transport} gives
\[
-\frac{2Z}{1+Z^2}H
=-\frac{[(1-\delta)+8L]Z^2}{1+Z^2}
  -\frac{2ZL}{1+Z^2}(U^z-4Z)
\le |U^z-4Z|.
\]
This is at most $2\epsilon_0$. Since $\zeta$ is a convex combination of
the two slopes, we obtain
\begin{equation}\label{eq:inner-shape-slopes}
	\frac7{10}\le a\le\frac9{10},\qquad
	H\zeta\le\frac1{10}+4\epsilon_0.
\end{equation}
Thus the angular shear stays strictly negative, and no large
radial derivative is needed to change the $Z$ dependence.

Because $U^z=U_1^z$ is constant in $R$,
\[
\frac{M^z(R)}R=\frac{110}R\frac{M^z(110)}{110}
+\left(1-\frac{110}R\right)U_1^z.
\]
Consequently $|W-(-3+4\delta Z^2)|\le10\epsilon_0$.
Equations \eqref{eq:inner-shape-slopes} and
\eqref{eq:inner-Q-equation} give $S_Q\ge7/5$ and
$31/20\le2-a/2\le33/20$. At $D=3$, the right side of
\eqref{eq:inner-D-equation} is strictly larger than $3(1-a/2)$.
Starting from $D(110)>3$ in \eqref{eq:inner-110-data}, we obtain
\begin{equation}\label{eq:inner-Q-lower}
	D>3\quad(R\ge110),\qquad Q\ge\frac12\quad(R\ge110e),
\end{equation}
as long as the displayed source bounds hold. For the second
assertion, integrate
$\partial_y Q+(33/20)Q\ge7/5$ for one logarithmic unit:
$\frac{28}{33}(1-e^{-33/20})>1/2$.
Here $b=0$ and $a\le1$, so $D>3$ proves the relaxed cone.

At $R_{\theta\mathrm{ok}}$, the compatibility relation gives
\[
U^\theta(R,Z)=\frac{(R/110)^{1/10}}{C_*(1+Z^2)}
=U_{\rm ref}^\theta(R,Z).
\]
Continue this swirl through $R_{\rm useless}$. It matches smoothly at
$R_{\theta\mathrm{ok}}$ because the step is flat. Keep $U^z=U_1^z$ until $R_{4Z}$.

\paragraph{The waiting interval.}
On $[R_{\theta\mathrm{ok}},R_{4Z}]$, the angular velocity has its
reference form and grows as $R^{1/10}$, while $U^z=U_1^z(Z)$ is fixed.
From $R_{\theta\mathrm{ok}}$ onwards the swirl agrees with the reference,
so the pressure and angular moment differences are constant in $R$, even
while $U^z$ is restored. Their normalized sizes decrease as the reference
swirl grows. In particular, at the correction radius
$R_m=e^{-6}R_{\rm ref}>eR_{4Z}$ we have
\begin{equation}\label{eq:inner-wait-pressure-defect}
\frac{|(M^p-M^p_{\rm ref})(R_m,Z)|}{[U^\theta(R_m,Z)]^2}
=\frac{|(M^p-M^p_{\rm ref})(R_{\theta\mathrm{ok}},Z)|}
        {[U^\theta(R_{\theta\mathrm{ok}},Z)]^2}
 \left(\frac{R_{\theta\mathrm{ok}}}{R_m}\right)^{1/5}.
\end{equation}
The contribution of $-\tfrac12\int (U^\theta)^2\,dR$ to the
$M^{z\theta}$ defect is also constant once the swirl agrees.
For this contribution the normalization by $R[U^\theta(R,Z)]^2$
gives the factor $(R_{\theta\mathrm{ok}}/R_m)^{6/5}$.
The same waiting interval also controls the inherited axial inertial
stress before $U^z$ is restored, when $\partial_RU^z$ becomes nonzero.
The bounds below retain this inherited stress explicitly.

\paragraph{Restoring the axial reference velocity.}
On $R_{4Z}\le R\le eR_{4Z}$, put $t=\log(R/R_{4Z})$ and prescribe
\begin{equation}\label{eq:inner-axial-restore}
	U^z(R,Z)=U_1^z(Z)+(4Z-U_1^z(Z))\sigma(t).
\end{equation}
For $eR_{4Z}\le R\le R_{\rm useless}$, set $U^z=4Z$.
The swirl remains $U_{\rm ref}^\theta$. Both endpoints of this
interpolation are smooth, and
\[
|\partial_y U^z|\le16\epsilon_0,\qquad
\|U^z-4Z\|_{C^1_Z}\le2\epsilon_0.
\]
The average $M^z/R$ stays within $2\epsilon_0$ of $4Z$ in
$C^1_Z$, since it is an average of its initial value and the
subsequent axial velocities. Thus $|W-(-3+4\delta Z^2)|\le10\epsilon_0$.
On $[R_{\theta\mathrm{ok}},R_{\rm useless}]$ the swirl has its reference
slope $\zeta=-2Z/(1+Z^2)$, so the same calculation as in
\eqref{eq:inner-shape-slopes} gives $H\zeta\le |U^z-4Z|\le2\epsilon_0$. With $a=4/5$,
\eqref{eq:inner-Q-equation} again gives $S_Q\ge7/5$, and the lower
bounds for $D$ and $Q$ in \eqref{eq:inner-Q-lower} persist up to
$R_{\rm useless}$.

\paragraph{Pressure and axial inertial stress.}
We start from the \hyperlink{inner-inherited-bounds}{bound
$\|M^p(110)\|_{C_Z^1}\le1$} and keep the $P_*$ dependence explicit.
On the shaping interval $\partial_y\log U^\theta\ge1/20$. Thus
\[
\int_0^T [U^\theta(110e^y,Z)]^2\,dy\le10[U^\theta(R_{\theta\mathrm{ok}},Z)]^2.
\]
Also $|\zeta|\le3A$. The compatibility inequality implies
\[
(1+A)[U^\theta(R_{\theta\mathrm{ok}},Z)]^2
\le (1+A)P_*^2(R_{\theta\mathrm{ok}}/R_{\rm ref})^{1/5}
\le e^{-2}P_*^2/(1+A).
\]
The remaining reference branch has $|\zeta|\le1$ and its
pressure integral is at most $\frac52P_*^2$.
Applying the same bounds after one $Z$ derivative gives, safely,
\begin{equation}\label{eq:inner-pressure-bound}
	\|P(R)\|_{C^1_Z}\le(C_{\rm pres}+100)P_*^2
	\qquad(110\le R\le R_{\rm useless}).
\end{equation}

On $110\le R\le R_{\rm useless}$, the prescribed axial velocity satisfies
$|U^z|,|\partial_Z U^z|\le5$ and $|\partial_y U^z|\le16\epsilon_0$. In \eqref{eq:inner-N-equation}
the pressure terms combine to $\mathcal P P$ for the actual field.
The preceding pressure estimate and the integrating factor $R$ give
\begin{equation}\label{eq:inner-N-bound}
	|N(R)|\le\frac{110}R|N(110)|
	+[500+4C_{\rm pres}]P_*^2.
\end{equation}
By the \hyperlink{inner-inherited-bounds}{bound at $110$},
$|N(110)|\le20K\le(1+K)^2$. Thus the last condition in
\eqref{eq:inner-compatibility} makes the inherited term at most
$P_*^2$ for $R\ge R_{4Z}$. In particular,
\begin{equation}\label{eq:inner-J-bound}
	|N|\le K_NP_*^2,
	\qquad |J|\le20K_N
	\qquad(R_{4Z}\le R\le R_{\rm useless}),
\end{equation}
because $U^\theta\ge e^{-4/5}P_*/2$ there. The inherited stress has
been estimated, not reset to zero.

\paragraph{The relaxed cone during restoration.}
Here $a=4/5$ and
\[
|b|=\frac{2|\partial_y U^z|}{U^\theta}\le\frac{144\epsilon_0}{P_*},
\qquad
bw=\frac{2J\partial_y U^z}{Q},\qquad
|bw|\le1280K_N\epsilon_0<\frac1{100}.
\]
This is the useful cancellation: although $w$ can have size
$P_*$, its product with $b$ is controlled by the small axial
change, with no hidden factor of $P_*$.
We have $\kappa=4/5+5b^2/4<1$, and
\[
D(a-bw)>3\left(\frac45-\frac1{100}\right)
>\frac85=2a.
\]
Thus \eqref{eq:inner-small-shear-test} verifies the relaxed
cone throughout the restoration. Afterwards $b=0$, so the
test returns to $D>2$.

\begin{figure}[htbp]
	\centering
	\begin{tikzpicture}[x=.95cm,y=.8cm,>=stealth,
		every node/.style={font=\small}]
		\fill[blue!5] (.2,0) rectangle (2.5,2.9);
		\fill[orange!8] (7.6,0) rectangle (9.7,2.9);
		\fill[green!7] (12.3,0) rectangle (14.2,2.9);
		\draw[->] (0,0)--(14.85,0) node[right] {$R$};
		\foreach \x/\lab in {.2/R_a,1.15/R_{ac},2.5/R_ae^{h_b},4/100,5.25/100e^{h_b},6.55/100e^{2h_b},7.6/110,9.7/R_{\theta\mathrm{ok}},11/R_{4Z},12.3/eR_{4Z},14.2/R_{\rm useless}}
		{\draw (\x,.06)--(\x,-.06) node[below,font=\scriptsize] {$\lab$};}
		\draw[gray!55,densely dotted] (5.25,-.58)--(5.25,-2.45);
		\draw[gray!55,densely dotted] (11,-.58)--(11,-2.45);
		\draw[gray!55,densely dotted] (6.55,.09)--(6.55,1.4);
		\draw[very thick,blue!70!black]
		(.2,2.1)..controls(.65,2.1)and(.96,2.1)..(1.15,2.04)
		..controls(1.66,1.92)and(1.94,.57)..(2.5,.57)
		..controls(3.45,.57)and(4.54,.65)..(5.25,.66)
		..controls(5.71,.67)and(6.14,1.4)..(6.55,1.4)
		--(7.6,1.4)..controls(8.04,1.4)and(8.13,1.13)..(8.65,1.13)
		..controls(9.17,1.13)and(9.26,1.4)..(9.7,1.4)--(14.2,1.4);
		\node[blue!70!black,anchor=west] at (0,3.35) {angular shear parameter $a$};
		\node[font=\footnotesize] at (1.35,2.58) {reduce the shear};
		\node[align=center,font=\footnotesize] at (8.65,2.3) {reshape\\the swirl};
		\node at (10.85,1.82) {$a=4/5$};
		\node[align=center,font=\footnotesize] at (13.25,2.26) {reference\\velocities};
		\draw[very thick,red!65!black]
		(.2,-1.82)..controls(.8,-1.92)and(1.4,-2.14)..(2.5,-2.22)
		..controls(3.05,-2.26)and(3.62,-2.32)..(4,-2.36)
		..controls(4.45,-2.41)and(4.77,-2.45)..(5.25,-2.45)
		--(11,-2.45)..controls(11.53,-2.45)and(11.77,-3.05)..(12.3,-3.05)
		--(14.2,-3.05);
		\node[red!65!black,anchor=west] at (0,-1.22)
		{axial velocity $U^z$ at a fixed $Z$};
		\node at (8.15,-2.83) {$U^z=U_1^z(Z)$};
		\node at (13.25,-2.67) {$U^z=4Z$};
	\end{tikzpicture}
	\caption{The inner connection at a fixed $Z$, with the short intervals
		enlarged and the amplitudes schematic. The axial shear is turned off on
		$[100,100e^{h_b}]$, so $U^z=U_1^z(Z)$ from $100e^{h_b}$ to $R_{4Z}$.
		The second switch brings $a$ to $4/5$ at $100e^{2h_b}$.
		Both velocities agree with the reference from $eR_{4Z}$ onwards,
		before the five-moment correction of Subsection~\ref{sec:inner-moment-corrections}.}
	\label{fig:inner-connection}
\end{figure}

\paragraph{Data for Interval I.4.}
The velocity joins the core smoothly at $R_a$. The prescribed components
$U^\theta,U^z$ agree with the reference on $[e^{-7}R_{\rm ref},R_{\rm useless}]$.
On $R_a<R\le R_{\rm useless}$,
\[
U^\theta>0,\qquad  \mathcal{S}^\theta<0,\qquad
 \mathcal{T}\ne0,
\]
and the relaxed cone holds for the actual moments and the pressure
\eqref{eq:inner-local-pressure}. The cone is admissible on
$R_a<R\le R_{ac}$, the collar already constructed in
Subsection~\ref{sec:inner-leaving-core}.
The moment correction in Subsection~\ref{sec:inner-moment-corrections}
is supported strictly beyond this collar.
In the subsequent shear modification, one may choose
$R_a<R_{-c}<R_{ac}$, so the velocity, moments, pressure, and stress
remain unchanged on the smaller collar $(R_a,R_{-c}]$.

The correction interval $(R_m,2R_m)$ lies inside this reference segment.
With $R_m=e^{-6}R_{\rm ref}$, the scale condition
\eqref{eq:inner-compatibility} gives $R_m\ge110e^{T+4}>16$.
The axial average under \eqref{eq:inner-axial-restore} and
\eqref{eq:inner-Q-lower}, \eqref{eq:inner-pressure-bound},
\eqref{eq:inner-J-bound} give the estimates used by the moment correction:
\begin{equation}\label{eq:imc-input-bounds}
 \begin{gathered}
 \sup_{R_m\le R\le R_{\rm useless}}\|M^z(R)/R-4Z\|_1\le2\epsilon_0,
 \qquad
 \sup_{R_m\le R\le R_{\rm useless}}\|P(R)\|_1\le(C_{\rm pres}+100)P_*^2,\\
 Q(R_m,Z)\ge\tfrac12,\qquad |N(R_m,Z)|\le K_NP_*^2 .
 \end{gathered}
\end{equation}
Since both prescribed velocities equal the reference on this segment,
the five moment discrepancies are constant in $R$. Write
\begin{equation}\label{eq:inner-defects}
 \Delta_j(Z)=M^j(R_m,Z)-M_{\rm ref}^j(R_m,Z)
            =M^j(R_{\rm useless},Z)-M_{\rm ref}^j(R_{\rm useless},Z).
\end{equation}
Append the given outer
\emph{velocity} and continue the moments from the axis. For
$R\ge R_{\rm useless}$ its moment differences from the supplied outer
construction remain the constant-in-$R$ functions $\Delta_j(Z)$.
In particular,
\[
P-P_o=\Delta_p,
\qquad
U^r-U_o^r=-\frac{\mathcal A\Delta_z}{L_\delta\sqrt{2R}}.
\]
The inertial stress differences follow from
\eqref{Itheta-moments-direct}--\eqref{Iz-moments-direct} and need not vanish.
Interval I.4 cancels the five discrepancies under the smallness
condition $\mathfrak e\le\mathfrak e_*$ on their centered, normalized
values \eqref{eq:imc-defect}. The axial error already satisfies
$\eta=\|U_1^z-4Z\|_1<\mathfrak e_*/100$ by
\eqref{eq:inner-bridge-axial-budget}. Interval I.4 uses
\eqref{eq:imc-family-short-velocity} for the moments inherited at
$R=110$ and the separation of $R_{\theta\mathrm{ok}}$ from $R_m$
to verify the full defect bound in Lemma~\ref{lem:imc-prepared-family}.

\subsection{Interval I.4: Correcting the five moments and preserving the relaxed cone}
\label{sec:inner-moment-corrections}\label{subsec:inner-corrected-cone}

We correct the five defects \eqref{eq:inner-defects} on $(R_m,2R_m)$,
using the entry estimates \eqref{eq:imc-input-bounds} from Interval I.3
and keeping the axis pressure $P_0$ fixed. On this interval the
prescribed velocities are $U^\theta=A_m^\theta x^{1/10}$ and $U^z=4Z$, where
\begin{equation}\label{eq:imc-scales}
 R_m=e^{-6}R_{\rm ref},\qquad R_{\rm useless}=eR_m,\qquad
 x=R/R_m,\qquad A_m^\theta(Z)=\frac{e^{-3/5}P_*}{1+Z^2}.
\end{equation}
Use the norm $\|q\|_1=\sup_{[-1,1]}|q|+\sup_{[-1,1]}|\partial_Z q|$,
summing over components and taking the supremum over the stated radial
interval when $R$ or $x$ is present. Define the centered, normalized defect by
\begin{equation}\label{eq:imc-defect}
 \mathbf d=\left(
 \frac{\Delta_z}{R_m},\quad
 \frac{\Delta_{\theta z}-4Z\Delta_\theta}{\sqrt2R_m^{3/2}A_m^\theta},\quad
 \frac{\Delta_\theta}{\sqrt2R_m^{3/2}A_m^\theta},\quad
 \frac{\Delta_{z\theta}-8Z\Delta_z}{R_m(A_m^\theta)^2},\quad
 \frac{\Delta_p}{(A_m^\theta)^2}\right),
 \qquad \mathfrak e=\|\mathbf d\|_1.
\end{equation}
The centering separates the two axial corrections from the three angular
corrections in the linearized moment map.

Retain the smooth nonnegative bump $\beta=\beta_{1/40}$, supported in
$[-1/40,1/40]$ with $\int\beta=1$. The constants $C_A,C_Q,C_S\ge1$ below
bound the inverse of the moment matrix, the quadratic remainder, and
the changes in velocity and stress sources, respectively. Fix them to
cover all the bounds in the proof below. They depend only on this bump
and the fixed reference exponent and radius ratios. Before choosing $j$
in \eqref{roadin:j-budget}, set
\begin{equation}\label{eq:imc-thresholds}
 t_*=\frac{1}{10^4K_NC_S},\qquad
 \mathfrak e_* =\min\left\{\frac{1}{8C_A^2C_Q},\frac{t_*}{2C_A}\right\}.
\end{equation}

\begin{proposition}[Inner moment correction]\label{prop:imc}
For a smooth connection as in Interval I.3 satisfying
\eqref{eq:imc-input-bounds}, assume $\mathfrak e\le\mathfrak e_*$.
Then $U^\theta$ and $U^z$ admit smooth corrections compactly supported
in $(R_m,2R_m)$ such that the corrected cumulative moments satisfy
\begin{equation}\label{eq:imc-output-moments}
 \widehat M^j(R,Z)=M^j_{\rm ref}(R,Z)
 \quad(2R_m\le R\le R_{\rm useless}),\qquad
 j=\theta,z,\theta z,z\theta,p.
\end{equation}
Recover $\widehat U^r$ from $\widehat U^z$ and $\widehat M^z$ by
\eqref{Ur-Uz}, and set $\widehat P=P_0+\widehat M^p$. The corrected inner
fields then match the supplied outer fields on a neighborhood of $R_{\rm useless}$.
On $R_a<R\le R_{\rm useless}$ they satisfy
\[
 \widehat U^\theta>0,\qquad \widehat{\mathcal{S}}^\theta<0,
 \qquad\widehat{\mathcal{T}}\ne0,
 \qquad\text{the relaxed cone condition}.
\]
All fields on $R\le R_m$ are unchanged, including the core and its
admissible collar. On the correction interval the construction gives
$\widehat\kappa<1$.
The size of the correction satisfies
\begin{equation}\label{eq:imc-output-size}
 \sup_{R_m\le R\le2R_m}
 \left\|
 \frac{\widehat U^\theta-U^\theta}{U^\theta},\ \widehat U^z-U^z,\
 R\partial_R\left(\frac{\widehat U^\theta-U^\theta}{U^\theta}\right),\
 R\partial_R(\widehat U^z-U^z)
 \right\|_1\le 2C_AC_S\mathfrak e.
\end{equation}
\end{proposition}

\begin{proof}
\emph{The moment correction.}
Set
\[
 s_1=\tfrac54,\quad s_2=\tfrac32,\quad s_3=\tfrac74,
 \qquad \gamma_j(x)=\beta(x-s_j),\qquad
 b_1=\gamma_1,\quad b_2=\gamma_3.
\]
For five coefficient functions $h=(c_1,c_2,\xi_1,\xi_2,\xi_3)$, write
\[
 f_h=\sum_{j=1}^3\xi_j\gamma_j,\qquad
 g_h=\sum_{i=1}^2c_i b_i,
\]
and define, on $R_m<R<2R_m$,
\begin{equation}\label{eq:imc-bumps}
 \widehat U^\theta=A_m^\theta(x^{1/10}+f_h),\qquad
 \widehat U^z=4Z+g_h.
\end{equation}
Extend the correction by zero. Figure~\ref{fig:imc-supports} shows
the supports of the two axial and three angular bumps.

\begin{figure}[H]
\centering
\begin{tikzpicture}[x=6.4cm,y=.80cm,>=stealth,
 every node/.style={font=\small}]
 \fill[green!7] (1,-.3) rectangle (2,2.7);
 \draw[->] (.84,0)--(2.76,0) node[right] {$x=R/R_m$};
 \foreach \x/\lab in {1/1,2/2,2.60/e}
  {\draw (\x,.08)--(\x,-.08) node[below] {$\lab$};}
 \foreach \x in {1.25,1.5,1.75}{
  \draw[very thick,blue!70!black]
    (\x-.085,1.1)..controls(\x-.040,1.1)and(\x-.045,1.75)..(\x,1.75)
    ..controls(\x+.045,1.75)and(\x+.040,1.1)..(\x+.085,1.1);}
 \foreach \x in {1.25,1.75}{
  \draw[very thick,red!65!black]
    (\x-.085,.30)..controls(\x-.040,.30)and(\x-.045,.85)..(\x,.85)
    ..controls(\x+.045,.85)and(\x+.040,.30)..(\x+.085,.30);}
 \node[blue!70!black,anchor=west] at (1.03,2.25)
  {three angular bumps};
 \node[red!65!black,anchor=west] at (1.02,-.83)
  {two axial bumps};
 \node[align=center] at (2.32,1.35)
  {all five moments\\already matched};
 \draw[thin,->] (2.32,.7)--(2.32,.1);
 \node[align=center] at (.93,3.32)
  {core and earlier\\connection unchanged};
 \node[align=center] at (2.60,3.32)
  {join the supplied\\outer fields};
\end{tikzpicture}
\caption{The correction uses a fixed portion of the terminal reference
interval. Bump widths and heights are enlarged; their coefficients
may have either sign. The tick labeled $e$ marks $R_{\rm useless}$ schematically.}
\label{fig:imc-supports}
\end{figure}

On this interval the exact identities
\[
 \widehat U^\theta(\widehat U^z-4Z)=\widehat U^\theta\,g_h,
 \qquad
 (\widehat U^z)^2-(4Z)^2-8Z(\widehat U^z-4Z)=g_h^2
\]
remove the unwanted linear coupling. Applying the five definitions
in \eqref{fiveM}, with exactly the normalization
\eqref{eq:imc-defect}, gives the terminal moment change
\begin{equation}\label{eq:imc-map}
 \mathbf F_Z(h)=
 \begin{pmatrix}
 \displaystyle\int g_h\\[1mm]
 \displaystyle\int x^{3/5}g_h+\int\sqrt{x}\,f_hg_h\\[1mm]
 \displaystyle\int\sqrt{x}\,f_h\\[1mm]
 \displaystyle-\int x^{1/10}f_h+
                 \int\left((A_m^\theta)^{-2}g_h^2-\tfrac12f_h^2\right)\\[1mm]
 \displaystyle\int x^{-9/10}f_h+\int\frac{f_h^2}{2x}
 \end{pmatrix}
 =\mathsf A h+\mathcal Q_Z(h,h).
\end{equation}
All integrals in this formula are over $1<x<2$. The only remaining
input-parameter factor is $(A_m^\theta)^{-2}$ in one quadratic term.

The axial block has weights $(1,x^{3/5})$ on the two separated
supports. Its determinant is positive, since it is the integral of
\[
 b_1(x_1)b_2(x_2)(x_2^{3/5}-x_1^{3/5}),\qquad x_1<x_2.
\]
The angular block has weights
$(x^{1/2},-x^{1/10},x^{-9/10})$ on three ordered supports.
Its evaluation determinant has a fixed nonzero sign. Indeed,
divide each column by $x_j^{-9/10}$: the three rows become
$(x_j^{7/5},-x_j,1)$. The graph of the strictly convex function
$x^{7/5}$ has strictly increasing secant slopes, so three ordered
points cannot be collinear. Integrating this determinant against
$\gamma_1(x_1)\gamma_2(x_2)\gamma_3(x_3)$ preserves its sign.
Thus $\mathsf A$ is a fixed invertible matrix, as in
Lemma~\ref{lem:bump-moment-matrix}. This gives the choice
\[
 C_A=\max\{1,\|\mathsf A^{-1}\|_{\ell^1\to\ell^1}\}.
\]
The explicit amplitude in \eqref{eq:imc-scales} gives
\[
 \|(A_m^\theta)^{-2}\|_1
 =\frac{12e^{6/5}}{P_*^2}<\frac{40}{P_*^2}\le40.
\]
The product inequality therefore gives a fixed constant $C_Q\ge1$
such that, for the symmetric bilinear form in \eqref{eq:imc-map},
\begin{equation}\label{eq:imc-quadratic-bound}
 \|\mathcal Q(h,k)\|_1\le C_Q\|h\|_1\|k\|_1.
\end{equation}

As in Lemma~\ref{lem:quadratic-moment-correction}, solve
$\mathbf F_Z(h)=-\mathbf d$ by iterating
\[
 h\longmapsto-\mathsf A^{-1}\bigl[\mathbf d+\mathcal Q(h,h)\bigr]
\]
on the closed $C^1_Z$ ball of radius $2C_A\mathfrak e$. If
$\mathfrak e\le(8C_A^2C_Q)^{-1}$, this map preserves the ball and has
Lipschitz constant at most $4C_A^2C_Q\mathfrak e\le1/2$. Hence there is
a solution satisfying
\begin{equation}\label{eq:imc-coefficients}
 t:=\|h\|_1\le2C_A\mathfrak e.
\end{equation}
The same estimate makes the derivative of the pointwise moment
map invertible. The implicit function theorem then shows that
$h$ is smooth in $Z$, including at $Z=\pm1$, for smooth input.

Returning to the original five moments, for $2R_m\le R\le R_{\rm useless}$
we have
\[
 \widehat M^j-M^j=-\Delta_j,\qquad
 \widehat M^j-M^j_{\rm ref}=0.
\]
The bumps cancel the incoming discrepancy $M^j-M^j_{\rm ref}=\Delta_j$
in each component. The coefficient bound \eqref{eq:imc-coefficients}
gives $t\le t_*$ for the cone estimate.

\noindent\emph{Preservation of the relaxed cone.}
To control $bw$ in the relaxed cone test \eqref{eq:inner-small-shear-test},
where $w=\mathcal{I}^z/\mathcal{I}^\theta$, we bound $Q$ from below and
$|N|$ from above. In the radial stress equations \eqref{eq:inner-Q-equation} and
\eqref{eq:inner-N-equation}, use the corrected velocities and their
actual cumulative moments, including partial bump integrals.
Write $P=P_0+M^p$ and omit hats throughout this calculation.
The fixed bumps and their derivatives give the following bounds
throughout $R_m\le R\le R_{\rm useless}$ whenever $C_St\le1/100$:
\begin{equation}\label{eq:imc-elementary}
\begin{gathered}
 P_*/8\le U^\theta\le P_*,\qquad 7/10\le a\le9/10,\\
 \|U^z-4Z\|_1+|\partial_y U^z|\le C_St,\qquad
 |\zeta+2Z/(1+Z^2)|\le C_St,\\
 |W-(-3+4\delta Z^2)|\le10\epsilon_0+C_St,\qquad
 \|P\|_1\le(C_{\rm pres}+100+C_St)P_*^2,\qquad
 S_Q\ge7/5.
\end{gathered}
\end{equation}
The change in $M^z/R$ is
\[
 \frac1x\int_1^{\min\{x,2\}}g_h(s,Z)\,ds,
\]
so its $C^1_Z$ norm is bounded by a fixed multiple of $t$.
The pressure change is
\[
 (A_m^\theta)^2\int_1^{\min\{x,2\}}
       \left(s^{-9/10}f_h+\frac{f_h^2}{2s}\right)ds,
\]
whose $C^1_Z$ norm is bounded by a fixed multiple of $P_*^2t$.
The formulas for $a$ and $\zeta$ follow by differentiating
$\log[A_m^\theta(x^{1/10}+f_h)]$. For $2\le x\le e$ the upper limit in both integrals is fixed at $2$,
so the same bounds hold there.

For the last inequality in \eqref{eq:imc-elementary}, the reference
values are
\[
 a=4/5,\quad U^z=4Z,\quad
 \zeta=-\frac{2Z}{1+Z^2},\quad W_0=-3+4\delta Z^2.
\]
Since $H_{4Z}\zeta\le0$, substitution gives
\[
 -W_0(1-a/2)-\tfrac\delta2(1-8Z^2)-H_{4Z}\zeta
 \ge\frac95-\frac\delta2.
\]
The error is bounded by $10\epsilon_0+C_St$ since $C_S$ includes
the constant in this source estimate. This proves
$S_Q\ge7/5$. With $C_S$ also bounding the fixed bumps and their
derivatives in \eqref{eq:imc-bumps}, the left side of
\eqref{eq:imc-output-size} is bounded by $C_St\le 2C_AC_S\mathfrak e$.

At $R_m$ no field or moment has changed. Thus $Q(R_m)\ge1/2$.
At a point where $Q=1/2$, \eqref{eq:inner-Q-equation} and \eqref{eq:imc-elementary} give
\[
 \partial_y Q\ge\frac75-\frac{33}{20}\frac12>0.
\]
The elementary scalar comparison therefore yields
\begin{equation}\label{eq:imc-Q-bound}
 Q\ge\tfrac12,\qquad D:=\mathcal{I}^\theta/F=RQ/L_\delta\ge R_m/2\ge8,
 \qquad F=U^\theta/\sqrt{2R}.
\end{equation}
For \eqref{eq:inner-N-equation}, the bounds $|U^z|,|\partial_Z U^z|\le5$,
$|\partial_y U^z|\le1$, $|W|\le4$, and \eqref{eq:imc-elementary} give
\[
 |Z(U^\theta)^2+\mathcal P P-W\partial_y U^z
   -\tfrac{1+\delta}{2}(1-2ZU^z)U^z-H\partial_Z U^z|
 \le(4C_{\rm pres}+500)P_*^2\le K_NP_*^2.
\]
Integrating $\partial_y N+N$ and using its unchanged entry value proves
\begin{equation}\label{eq:imc-N-bound}
 |N|\le K_NP_*^2\qquad(R_m\le R\le R_{\rm useless}).
\end{equation}

Since $t\le t_*$,
\[
 |b|\le\frac{16C_St}{P_*},\qquad
 \kappa=a+\frac{b^2}{a}<1.
\]
The product in the cone test satisfies the exact identity and bound
\begin{equation}\label{eq:imc-bw}
 bw=\frac{2N\partial_y U^z}{(U^\theta)^2Q},\qquad
 |bw|\le256K_NC_St<\frac1{10}.
\end{equation}
The factor $P_*^2$ in $N$ has canceled against $(U^\theta)^2$.
For $\kappa\le2$, the relaxed cone is exactly
$D(a-bw)>2a$, as in \eqref{eq:inner-small-shear-test}.
Here it follows at once from
\[
 D(a-bw)\ge8\left(\frac7{10}-\frac1{10}\right)
 =\frac{24}{5}>\frac95\ge2a.
\]
This also proves $\mathcal{T}\cdot \mathcal{S}<0$, hence nonzero
stress, and $a>0$ proves $\mathcal{S}^\theta<0$.

On $[2R_m,R_{\rm useless}]$ the corrections to $U^\theta$ and $U^z$
vanish, so these components agree with the reference. The five moments
also agree by \eqref{eq:imc-output-moments}. With the common axis pressure
$P_0$, Lemma~\ref{lem:moment-gluing} gives the asserted agreement of
the reconstructed fields there.
\end{proof}

\paragraph{Smallness of the five incoming defects.}
It remains to verify the hypothesis $\mathfrak e\le\mathfrak e_*$
of Proposition~\ref{prop:imc} for the connection just constructed.
Let $\mathbf d^0$ be given by \eqref{eq:imc-defect} with the moment
differences evaluated at $R=110$, keeping the same denominators
$R_m,A_m^\theta$. Set $\mathfrak e_0=\|\mathbf d^0\|_1$ and write
\begin{equation}\label{eq:imc-preparation-data}
 \alpha=R_{\theta\mathrm{ok}}/R_m,\qquad
 s=eR_{4Z}/R_m=e^{-1},\qquad \eta=\|U_1^z-4Z\|_1 .
\end{equation}
On the incoming connection the swirl agrees with the reference
from $R_{\theta\mathrm{ok}}$ to $R_m$.
For $110\le R\le R_{\theta\mathrm{ok}}$ put
$\ell=\log(R_{\theta\mathrm{ok}}/R)$ and $q=1-\sigma(1-\ell/T)$.
By \eqref{eq:inner-swirl-shape},
\[
 \frac{U^\theta}{A_m^\theta}
 =\alpha^{1/10}e^{-\ell/10+qB},\qquad
 B=\log\bigl(C_*(1+Z^2)U_1^\theta\bigr).
\]
Since $T=400A$, $|\sigma'|\le8$ and $|B|,|\partial_ZB|\le2A$,
we have $|qB|,|q\partial_ZB|\le\ell/25$. Thus
\begin{equation}\label{eq:imc-shape-decay}
 \left\|\left(\frac{U^\theta}{A_m^\theta}\right)^k\right\|_1
 \le\alpha^{k/10}e^{-3k\ell/50}(1+k\ell/25),\qquad k=1,2 .
\end{equation}

In $x=R/R_m$ the normalized angular, energy, and pressure moments
use the measures $\sqrt{x}\,dx$, $dx/2$, and $dx/(2x)$ for the
swirl terms. With $x=\alpha e^{-\ell}$,
\eqref{eq:imc-shape-decay} and the reference swirl $x^{1/10}$
give respective bounds $2\alpha^{8/5}$, $\alpha^{6/5}$, and
$10\alpha^{1/5}$ for their increments from $R=110$.
The constants follow by integrating the exponential bounds to
$\ell=\infty$, using $\int_0^\infty e^{-a\ell}(1+b\ell)\,d\ell
=a^{-1}+ba^{-2}$. In particular, the pressure moment has the
slowest decay because its measure is $dx/(2x)$.

On $[110,R_m]$ we have $\|U^z-4Z\|_1\le\eta$, with $U^z=4Z$
for $x\ge s$. The axial difference contributes at most $\eta s$ to the
first defect. For the centered mixed moment, integrate
$\sqrt{x}(U^\theta/A_m^\theta)(U^z-4Z)$, using
\eqref{eq:imc-shape-decay} up to $x=\alpha$ and $x^{1/10}$
for the normalized swirl on $[\alpha,s]$. Since $\alpha<s$ by
\eqref{eq:inner-compatibility}, the bound is $2\eta s^{8/5}$.
The axial part of the centered energy is
$(U^z-4Z)^2/(A_m^\theta)^2$, with integral at most
$40\eta^2s/P_*^2$ in $\|\cdot\|_1$. Summing the five estimates gives
\begin{equation}\label{eq:imc-defect-certificate}
 \boxed{\quad
 \mathfrak e\le\mathfrak e_0+16\alpha^{1/5}
                      +3\eta+40\eta^2/P_*^2 .
 \quad}
\end{equation}
The axial restoration ends at $sR_m$ with $s=e^{-1}$ fixed. The bound
$\eta<\mathfrak e_*/100$ from \eqref{eq:inner-bridge-axial-budget}
controls the two terms involving $\eta$. With the outer data and
$\Lambda$ fixed, the following lemma shows that both the inherited
defect $\mathfrak e_0$ and the swirl contribution $\alpha^{1/5}$
are $O(C_*^{-2})$.

\begin{lemma}[Choice of $C_*$ for the five-moment correction]
\label{lem:imc-prepared-family}
Fix $j$ as in \eqref{roadin:j-budget}, the outer data $P_*,\delta$
and dimensionless schedule as in Section~\ref{sec:outer-construction},
and $\Lambda$ and $c_*$ as in Interval I.2. For
$R_{\rm ref}=110(C_*P_*)^{10}$ and all sufficiently large $C_*$,
the connection satisfies \eqref{eq:inner-compatibility} and
\[
  \mathfrak e_0+\alpha^{1/5}=O(C_*^{-2}),\qquad
  \mathfrak e<\mathfrak e_* .
\]
Consequently Proposition~\ref{prop:imc} applies to this connection.
\end{lemma}

\begin{proof}
Subsection~\ref{subsec:outer-cone-closure} constructs the outer
profile at each $R_{\rm ref}$ with the same waiting length $\tau$
and analytic pressure $P_0=P_0^{\rm pre}$. Its complex neighborhood
and analytic bound are therefore fixed as $C_*$ increases. For fixed $\Lambda$, the
bounds following \eqref{eq:inner-data-sizes} give $A=O(1)$ and
$K=O(C_*)$, so $T=400A$ stays bounded. Since $R_{\rm ref}$,
$R_{4Z}$ and $R_m$ grow as $C_*^{10}$, all three inequalities
in \eqref{eq:inner-compatibility} and $R_m\ge16$ hold for large
$C_*$. We may increase $C_*$ also to meet the core and comparison
bounds in Proposition~\ref{prop:continuation-core} and
Lemma~\ref{prop:inner-comparison}.

On $[0,110]$ we have $\|F\|_1=O(C_*^{-1})$ and $\|U^z\|_1=O(1)$
by \eqref{eq:imc-family-short-velocity}. The reference profile satisfies
\[
 U_{\rm ref}^\theta=\frac{(R/110)^{1/10}}{C_*(1+Z^2)},\qquad
 \int_0^{110}\frac{(U_{\rm ref}^\theta)^2}{2R}\,dR
 =\frac{5}{2C_*^2(1+Z^2)^2}.
\]
The five centered numerators in \eqref{eq:imc-defect}, evaluated at
$R=110$, therefore have $C_Z^1$ sizes
$O(1),O(C_*^{-1}),O(C_*^{-1}),O(1),O(C_*^{-2})$, respectively.
Since $R_m=e^{-6}R_{\rm ref}$ and $A_m^\theta$ and its reciprocal
have fixed $C_Z^1$ bounds, their normalized sizes are bounded by a
constant times
$(C_*^{-10},C_*^{-16},C_*^{-16},C_*^{-10},C_*^{-2})$.
This proves $\mathfrak e_0=O(C_*^{-2})$. Also,
\[
 \alpha^{1/5}
 =e^{(T+6)/5}(C_*P_*)^{-2}=O(C_*^{-2}).
\]
By \eqref{eq:inner-bridge-axial-budget}, $\eta<\mathfrak e_*/100$.
Since $P_*\ge1$ and $\mathfrak e_*<1$, the last two terms in
\eqref{eq:imc-defect-certificate} sum to less than
$0.1\mathfrak e_*$. Taking $C_*$ sufficiently large makes the
first two terms smaller than $\mathfrak e_*/2$. Interval I.3 supplies
the other input bounds \eqref{eq:imc-input-bounds}, so the five
moments can be matched while preserving the relaxed cone.
\end{proof}

\subsection{Gluing the inner and outer profiles}
\label{subsec:inner-outer-gluing}

On $[2R_m,R_{\rm useless}]$ the corrected inner velocities
$(U_i^\theta,U_i^z)$ and their five cumulative moments agree with
the reference continuation of the outer profile $(U_o^\theta,U_o^z)$,
by Proposition~\ref{prop:imc} and Lemma~\ref{lem:imc-prepared-family}.
Both profiles use the same axis pressure $P_0$. We join them at
$R_{\rm useless}=eR_m$ by setting
\[
 (U^\theta,U^z)(R,Z)=
 \begin{cases}
 (U_i^\theta,U_i^z)(R,Z),&0\le R\le R_{\rm useless},\\
 (U_o^\theta,U_o^z)(R,Z),&R\ge R_{\rm useless}.
 \end{cases}
\]
Use the cumulative moments of this joined pair to set $P=P_0+M^p$,
and recover $U^r$ by \eqref{Ur-Uz}.

\begin{proposition}[The global profile with the relaxed cone]
\label{prop:global-relaxed-profile}
The profiles $U^\theta,U^z$ are bounded, and
$F=U^\theta/\sqrt{2R}$ and $U^z$ are smooth on
$[0,\infty)\times[-1,1]$. The pressure also satisfies
\eqref{P-Utheta}. The associated physical fields
are smooth away from $(t,x)=(1,0)$, and the velocity is exactly
divergence-free. All five terminal conditions
\eqref{eq:moment-conditions} hold for every $Z\in[-1,1]$. The complete
fields agree with those of the outer profile, including its reference
continuation, for $R\ge 2R_m$. In particular they agree with the exact
heat flow for $R\ge R_h$.

The stress $\mathcal{T}$ vanishes for $R\le R_a$ and $R\ge R_h$ and is
nonzero for $R_a<R<R_h$. One has $U^\theta>0$ for $R>0$ and
$\mathcal{S}^\theta<0$ on $[R_a,R_h]$. The relaxed cone holds on
$(R_a,R_h)$, and the cone is admissible on
$(R_a,R_{ac}]\cup[R_{+c},R_h)$.
\end{proposition}

\begin{proof}
On both sides of $R_{\rm useless}$ the prescribed velocities are given
locally by the same smooth outer profile, since the inner pair already
agrees with its reference continuation from $2R_m$ onwards. Thus the
join is smooth. At $2R_m$ the five moments agree as functions of $Z$.
With the common axis pressure, Lemma~\ref{lem:moment-gluing} gives
equality of the moments, pressure, radial velocity, inertial stress, and
shear for every $R\ge 2R_m$. The five terminal conditions, including the
renormalized angular moment limit, now follow from
Proposition~\ref{prop:road-outer-summary}(ii).
The total pressure moment equals $-P_0$, so
\begin{equation}\label{eq:global-pressure-normalization}
 P(R,Z)=P_0(Z)+M^p(R,Z)
 =-\int_R^\infty\frac{(U^\theta(\rho,Z))^2}{2\rho}\,d\rho
 \qquad(R\ge0).
\end{equation}
Near the axis, $F$ and $U^z$ are analytic by
Proposition~\ref{prop:continuation-core}. Since
$M^z/R=\int_0^1 U^z(sR,Z)\,ds$, formula \eqref{Ur-Uz} gives a smooth
extension of $U^r/\sqrt R$ to $R=0$. The pressure is smooth there by
$P=P_0+\int_0^R F^2\,d\rho$. The prescribed velocities are bounded on
the compact inner and connecting regions and, by \eqref{exterior-profile},
on the heat exterior. Formula \eqref{Ur-Uz} enforces incompressibility.
The similarity coordinates extend smoothly to $t=1$ at $z\ne0$, as shown
in Subsection~\ref{subsec:variables}. At the remaining terminal points
with $r>0$ and $z=0$, one has $R\to\infty$, so the velocity and pressure
extend smoothly by the exact heat formulas \eqref{exterior-ansatz}
and \eqref{eq:heat-physical-profile}.

The core and heat exterior have zero stress. On $(R_a,R_h)$ the signs
and cone conditions, including the two admissible collars, follow from
Proposition~\ref{prop:imc} and Proposition~\ref{prop:road-outer-summary}(ii).
At $R_a$, the core bounds $F_c>0$ and $-\partial_R\log F_c\ge1/64$
give $\mathcal S^\theta(R_a,Z)=2R_a\partial_R F_c(R_a,Z)<0$.
At $R_h$, the heat shear formula \eqref{eq:heat-shear} gives the
same strict sign.
\end{proof}

Section~\ref{sec:shear-modification} modifies the shear between the
two collars to obtain the admissible cone and restores the five moments
in the reserved outer interval $[R_c,2R_c]\subset I_1$.

\section{The admissible cone by modification of the radial shear}
\label{sec:shear-modification}

We modify the radial shear of the global profile in
Proposition~\ref{prop:global-relaxed-profile} to obtain the admissible cone.
We use rapid shear modulation as in \cite[Appendix C]{1}.
The Poisson kernel gives a periodic shear loop with the prescribed mean,
so the velocity stays close to the supplied profile. We then restore the
five moments to their values before modulation, preserving the core and
the exact heat exterior.

Let $0<\delta\le\frac12$. The input consists of bounded profiles
$(U^\theta,U^z)$, smooth for $R>0$, their moments $\mathbf M$ from \eqref{fiveM}, and their
pressure $P=P_0+M^p$, where $P_0(Z)=P(0,Z)$.
We impose the following requirements.
\begin{enumerate}[(i)]
\item
The profile has the axis regularity stated in \eqref{smooth-axis},
the prescribed exterior heat flow, and the five terminal moment
conditions \eqref{eq:moment-conditions}. Its stress is zero outside
$[R_a,R_h]$ and nonzero on $(R_a,R_h)\times[-1,1]$.
We have $U^\theta>0$ for every $R>0$ and
$\mathcal{S}^\theta<0$ on the closed stress annulus.
\item
The relaxed cone condition \eqref{cone-relaxed} holds wherever
the stress is nonzero. There are radii
\[
R_a<R_{-c}<R_{+c}<R_c<2R_c<R_h
\]
such that the input is already admissible on
$((R_a,R_{-c}]\cup[R_{+c},R_h))\times[-1,1]$, and on neighborhoods
of the two radial boundaries of $[R_{-c},R_{+c}]\times[-1,1]$.
\item
There is a reserved interval $R_c\le R\le2R_c$ on which
\begin{equation}\label{eq:sm-patch-input}
U^\theta(R,Z)=A_c^\theta(Z)x^{-1/2-\mu},\qquad
U^z(R,Z)=0,\qquad
x=R/R_c,\qquad 0<\mu\le\tfrac12,
\end{equation}
with $A_c^\theta>0$ smooth on $[-1,1]$. This interval will be used
only to correct moments. Its input shear has
$a=2+2\mu$, $b=0$, and hence $\kappa=2+2\mu>2$.
\end{enumerate}
The step $\sigma$ and the five bump functions used for moment
restoration are fixed independently of the input. Constants denoted
by $C$ depend only on these functions. For a given input profile and
requested accuracy $\varepsilon$, we compute the shear loop and its
cone margins before choosing the frequency $N$.

\begin{proposition}[Shear modification]\label{prop:sm}
Under (i)--(iii), for every $0<\varepsilon\le1$ and every integer
$N$ satisfying \eqref{eq:sm-frequency-choice}, there are profiles
$(\widehat U^\theta,\widehat U^z)$, smooth for $R>0$ and with the
same axis regularity, having the following properties.
\begin{enumerate}[(1)]
\item
The angular and axial velocities equal the input outside
$(R_{-c},R_{+c})\cup(R_c,2R_c)$. All five moments agree with the
input for $R\ge2R_c$. The same axis pressure is retained.
Consequently $\widehat P$, $\widehat U^r$, and
$\widehat{\mathcal{T}}$ agree with the input for $R\le R_{-c}$
and for $R\ge2R_c$.
\item
We have $\widehat U^\theta>0$ for $R>0$,
$\widehat{\mathcal{S}}^\theta<0$ on the stress annulus,
and the admissible cone condition wherever
$\widehat{\mathcal{T}}\ne0$. The stress is nonzero precisely
on $(R_a,R_h)\times[-1,1]$.
\item
Writing $\|\cdot\|_{C_Z^1}$ for the sum of the uniform
norms of a function and its first $Z$ derivative, we have
\begin{equation}\label{eq:sm-output-close}
\sup_{R\ge0}\left(
\|\widehat U^\theta-U^\theta\|_{C_Z^1}
+\|\widehat U^z-U^z\|_{C_Z^1}
+\|\widehat{\mathbf M}-\mathbf M\|_{C_Z^1}
+\|\widehat P-P\|_{C_Z^1}\right)\le\varepsilon.
\end{equation}
No smallness of radial derivatives is asserted. Their change
is exactly what raises $\kappa$ above two.
\end{enumerate}
\end{proposition}

\subsection{A smooth admissible shear loop with the prescribed mean}
\label{subsec:sm-loop}

At each $(R,Z)$ we hold $F$ and $\mathcal I$ fixed. We seek a smooth
periodic shear with the prescribed mean and an admissible stress at
every phase.

\begin{example}\label{ex:sm-shear-mean}
Freeze $F=1$, $\mathcal{I}=(10,0)$, and $\mathcal{S}=(-1,0)$. Then $\kappa=1$ and
the relaxed condition is satisfied:
\[
 \mathcal{T}=(9,0),\qquad \mathcal{T}\cdot(-\mathcal{S})=9>2-\kappa=1.
\]
The two shears
\[
 \mathcal{S}_\pm=(-1,\pm\sqrt{\sqrt5-1}),\qquad
 \tfrac12(\mathcal{S}_++\mathcal{S}_-)=\mathcal{S},\qquad \kappa_\pm=\sqrt5>2
\]
are admissible. Indeed, with $\mathcal{T}_\pm=\mathcal{I}+\mathcal{S}_\pm$,
\[
 \mathcal{T}_\pm\cdot(-\mathcal{S}_\pm)=10-\sqrt5>0,\qquad
 (\kappa_\pm-2)(\mathcal{T}_\pm\cdot \mathcal{S}_\pm^\perp)^2
 =100(7-3\sqrt5)<2(10-\sqrt5)^2.
\]
A smooth period-one loop with the same mean is
\[
 a_L(\phi)=\frac{\sqrt5}{1+4\cos^2(2\pi\phi)},\qquad
 b_L(\phi)=-2a_L(\phi)\cos(2\pi\phi).
\]
Direct integration gives $\int_0^1(a_L,b_L)\,d\phi=(1,0)$, and
$\kappa_L=\sqrt5$. For $\mathcal{S}_L=(-a_L,b_L)$ and $\mathcal{T}_L=\mathcal{I}+\mathcal{S}_L$,
\[
 -\mathcal{T}_L\cdot \mathcal{S}_L=a_L(10-\sqrt5)>0,\qquad
 \frac{(\kappa_L-2)(\mathcal{T}_L\cdot \mathcal{S}_L^\perp)^2}{a_L^2}
 \le400(\sqrt5-2)<2(10-\sqrt5)^2.
\]
Thus admissible shears can average to a relaxed shear.
\end{example}

Use the notation of Lemma~\ref{lem:outer-cone-test}:
\[
F=\frac{U^\theta}{\sqrt{2R}},\qquad
a=1-2R\partial_R\log U^\theta,\qquad
b=\frac{2R\partial_RU^z}{U^\theta}.
\]
For this section only, also put
\begin{equation}\label{eq:sm-normalized}
p=(p_1,p_2)=\frac{\mathcal{I}}{F},\qquad
t_0=-\frac ba,\qquad
\kappa=a(1+t_0^2),\qquad
H(t)=p_1+p_2t,\qquad J(t)=p_2-p_1t.
\end{equation}
Thus $-\mathcal{S}/F=a(1,t_0)$ and
$\mathcal{T}/F=p-a(1,t_0)$.
More generally, a proposed shear of the form
\[
\mathcal{S}'=-F\,\frac{v}{1+t^2}(1,t),\qquad
\mathcal{T}'=\mathcal{I}+\mathcal{S}'
\]
has $\kappa'=v$. With $a'=v/(1+t^2)$, direct scalar products give
\begin{equation}\label{eq:sm-scalar-cone}
\begin{aligned}
-\mathcal{T}'\cdot \mathcal{S}'&=F^2a'[H(t)-v],\\
2(\mathcal{T}'\cdot \mathcal{S}')^2
-(v-2)(\mathcal{T}'\cdot(\mathcal{S}')^\perp)^2
&=F^4(a')^2\{2[H(t)-v]^2-(v-2)J(t)^2\}.
\end{aligned}
\end{equation}
In particular, the relaxed cone implies
\begin{equation}\label{eq:sm-H-positive}
H(t_0)>2.
\end{equation}
Indeed, if $\kappa\le2$, its second inequality is precisely
$a[H(t_0)-\kappa]>a(2-\kappa)$.
If $\kappa>2$, its first inequality gives
$H(t_0)>\kappa>2$.
We use the full quantitative inequality in the $\kappa\le2$
branch. The sign $\mathcal{T}\cdot \mathcal{S}<0$ alone would
only give $H(t_0)>\kappa$, which is insufficient here.

The useful observation is now quite simple. If we keep
$H(t)$ a definite distance above two and take $v$ only slightly
above two, both expressions in \eqref{eq:sm-scalar-cone} are
positive. We may therefore vary the direction $t$ substantially.
This variation will allow the average shear to remain unchanged.
Merely increasing the length of the original shear would not
have that property.

\paragraph{Explicit scales for the loop.}
Let $\mathcal K=[R_{-c},R_{+c}]\times[-1,1]$; all extrema in this
paragraph are over $\mathcal K$. Define
\begin{equation}\label{eq:sm-loop-scales}
\begin{aligned}
a_*&=\min a,& m&=\min[H(t_0)-2],& d_*&=m/4,\\
q_*&=\sqrt{\frac3{2a_*}},&
B_*&=\|t_0\|_\infty+2q_*+
       \frac{4\|p_2\|_\infty q_*^2}{d_*},&
J_*&=\|p_1\|_\infty B_*+\|p_2\|_\infty.
\end{aligned}
\end{equation}
Both $a_*$ and $m$ are positive. To ensure that the loop stays
constant near the radial boundaries, include their margin
\[
m_\partial=
\min_{Z\in[-1,1],\ R\in\{R_{-c},R_{+c}\}}[\kappa(R,Z)-2]>0
\]
and set, with no further choice,
\begin{equation}\label{eq:sm-eta}
\eta=\min\left\{\frac12,\frac m8,
\frac{m^2}{8(1+J_*^2)},\frac{m_\partial}{2}\right\}.
\end{equation}
The scales $d_*$ and $\eta$ are local to this loop construction;
$\eta$ here is unrelated to the axial error in
Subsection~\ref{subsec:inner-reference-velocity}.

Using the smooth step $\sigma$ of \eqref{eq:road-sigma}, define a nonnegative amplitude
\begin{equation}\label{eq:sm-q-v}
\begin{aligned}
q&=\sigma\left(\frac{2+\eta-\kappa}{\eta}\right)
\sqrt{\frac{2+2\eta-\kappa}{2a}}\quad(\kappa<2+\eta),
&q&=0\quad(\kappa\ge2+\eta),\\
v&=\kappa+2aq^2.
\end{aligned}
\end{equation}
The square root is separated from zero wherever the cutoff
can be nonzero; flatness of $\sigma$ makes $q$ smooth at the
joining level. If $\kappa\le2$, then $v=2+2\eta$.
If $2<\kappa<2+\eta$, then $\kappa\le v\le2+2\eta$.
Otherwise $q=0$ and $v=\kappa$. In particular,
\begin{equation}\label{eq:sm-v-bounds}
v\ge2+\frac{3\eta}{8},\qquad
q\le q_*,\qquad v\le3\ \hbox{where }q\ne0.
\end{equation}
For the first estimate, when $2\le\kappa\le2+\eta/2$,
the cutoff is at least $\sigma(1/2)=1/2$; when
$\kappa\ge2+\eta/2$, use $v\ge\kappa$.
By \eqref{eq:sm-eta}, $q=0$ on neighborhoods of both radial
boundaries of $\mathcal K$.

Take a temporary angle $\psi\in\mathbb R/(2\pi\mathbb Z)$ and write
$\langle\cdot\rangle_\psi$ for the average over one period. We seek a
direction $t(\psi)$ with mean $\langle t\rangle_\psi=t_0$ and variance
$\langle(t-t_0)^2\rangle_\psi=(v-\kappa)/a=2q^2$. These conditions give
$a\langle1+t^2\rangle_\psi=v$, the normalization needed to recover the
mean shear by weighting each direction by $1+t^2$. To keep $H(t)$ above
two at every angle as well, we use the Poisson kernel. Put
\begin{equation}\label{eq:sm-poisson-direction}
u=\frac{p_2q}{d_*},\qquad h=\sqrt{1+u^2},\qquad r=\frac uh,
\qquad
t(\psi)=t_0+\frac{2q}{h}
\frac{\cos\psi-r}{1-2r\cos\psi+r^2}.
\end{equation}
Here $|r|<1$; $r$ is the kernel parameter, not the radial coordinate.
This formula is smooth even at $p_2=0$ or $q=0$.
It is useful to recognize the Poisson kernel
\[
\mathsf P_r(\psi)=\frac{1-r^2}{1-2r\cos\psi+r^2}
=1+2\sum_{n=1}^{\infty}r^n\cos(n\psi).
\]
When $p_2\ne0$, the last formula in
\eqref{eq:sm-poisson-direction} reads
$t-t_0=d_*(\mathsf P_r-1)/p_2$.
Orthogonality of the cosines gives
\[
\langle\mathsf P_r\rangle_\psi=1,\qquad
\langle(\mathsf P_r-1)^2\rangle_\psi
=\frac{2r^2}{1-r^2}=2u^2.
\]
Consequently, including $p_2=0$ by the explicit smooth formula,
\begin{equation}\label{eq:sm-direction-identities}
\boxed{\quad
\langle t\rangle_\psi=t_0,\qquad
\langle(t-t_0)^2\rangle_\psi=2q^2,\qquad
H(t)\ge H(t_0)-d_*\ge2+\frac{3m}{4}.
\quad}
\end{equation}
The last bound follows from $\mathsf P_r>0$ and
$p_2(t-t_0)=d_*(\mathsf P_r-1)$; it does not require a sign
assumption on $p_2$.
The pointwise size is also explicit:
\begin{equation}\label{eq:sm-t-bound}
|t-t_0|\le2q(h+|u|)
\le2q+\frac{4|p_2|q^2}{d_*},\qquad |t|\le B_*.
\end{equation}
For example, the first inequality follows by maximizing
$(c-r)/(1-2rc+r^2)$ over $-1\le c\le1$.

\paragraph{Phase reparametrization.}
The vectors $v(1,t)/(1+t^2)$ have the correct direction and
length parameter $v$, but their ordinary $\psi$ average need
not equal $a(1,t_0)$. We give each direction a weight proportional
to $1+t^2$, which cancels their denominator. Define
\begin{equation}\label{eq:sm-phase}
\phi(\psi)=\frac{a}{2\pi v}
\int_0^\psi[1+t(\vartheta)^2]\,d\vartheta.
\end{equation}
By \eqref{eq:sm-direction-identities} and \eqref{eq:sm-q-v},
\[
a\langle1+t^2\rangle_\psi
=a(1+t_0^2)+2aq^2=v.
\]
Thus $\phi(\psi+2\pi)=\phi(\psi)+1$, and $\partial_\psi\phi>0$.
It has a smooth inverse $\psi=\psi(R,Z,\phi)$.
Define the period-one loop by
\begin{equation}\label{eq:sm-loop}
(a_L,-b_L)(R,Z,\phi)
=\frac{v}{1+t(\psi(R,Z,\phi))^2}
       (1,t(\psi(R,Z,\phi))).
\end{equation}
Changing variables from $\phi$ to $\psi$ now gives exactly
\begin{equation}\label{eq:sm-loop-mean}
\int_0^1 a_L\,d\phi=a,\qquad
\int_0^1 b_L\,d\phi=b,\qquad
\frac{a_L^2+b_L^2}{a_L}=v.
\end{equation}
Where $q=0$, the loop is constant and equals $(a,b)$.
In particular it equals the input on the radial boundary
neighborhoods, with all derivatives.

\paragraph{Checking the cone.}
On $q\ne0$, equations \eqref{eq:sm-eta},
\eqref{eq:sm-direction-identities}, and \eqref{eq:sm-t-bound}
give
\[
H(t)-v\ge\frac m2,\qquad |J(t)|\le J_*,\qquad
2[H(t)-v]^2-(v-2)J(t)^2\ge\frac{m^2}{4}.
\]
Moreover,
\[
a_L\ge\frac2{1+B_*^2},\qquad |(a_L,b_L)|\le3.
\]
These estimates and \eqref{eq:sm-scalar-cone} verify the
admissible cone for the frozen pair
\[
\mathcal{S}_L=F(-a_L,b_L),\qquad
\mathcal{T}_L=\mathcal{I}+\mathcal{S}_L.
\]
On $q=0$, the loop equals the input and
$\kappa\ge2+\eta$, so the relaxed cone is already admissible.
This proves admissibility at every phase, including the
transition where the oscillation is turned off.

\begin{figure}[htbp]
\centering
\begin{tikzpicture}[scale=1.55,>=stealth]
  \draw[->] (-.12,0)--(2.65,0) node[right] {$a$};
  \draw[->] (0,-1.25)--(0,1.38) node[above] {$-b$};
  \fill[gray!12] (1,0) circle (1);
  \draw[gray,dashed] (1,0) circle (1);
  \draw[blue!65!black,thick,domain=-68:68,samples=100,variable=\t]
    plot ({1.1+1.1*cos(2*\t)},{1.1*sin(2*\t)});
  \draw[blue!65!black,->,thick] (2.17,.255)--(2.12,.410);
  \draw[blue!65!black,->,thick] (2.12,-.410)--(2.17,-.255);
  \fill (.80,0) circle (.026);
  \node[below] at (.80,-.04) {input mean};
  \node[gray!70!black] at (1.06,.56) {$\kappa=2$};
  \node[blue!65!black,right] at (2.22,.76) {$\kappa=v>2$};
  \node[align=center,font=\small] at (1.2,-1.65)
    {The loop runs out and back along the blue arc.\\
     Its weighted mean may lie inside the dashed circle.};
\end{tikzpicture}
\caption{Geometry in the $(a,-b)$ plane, drawn for $t_0=p_2=0$.
The circle $a^2+b^2=2a$ is the threshold $\kappa=2$.
Every point of the loop is outside it. Changing the time spent
at each point by \eqref{eq:sm-phase} gives the prescribed mean.
The drawing is schematic; admissibility also uses the
stress estimates proved above.}
\label{fig:sm-loop}
\end{figure}

\subsection{Putting the loop into the velocity profiles}
\label{subsec:sm-modulation}

We integrate the shear loop to obtain a small velocity perturbation.

Set $y=\log(R/R_{-c})$. For functions depending also on $\phi$,
$\partial_y=R\partial_R$ always differentiates while holding $\phi$
fixed. Let $\mathcal A_L,\mathcal B_L$ be the zero-mean periodic
primitives specified by
\begin{equation}\label{eq:sm-primitives}
\partial_\phi\mathcal A_L=-\frac12(a_L-a),\qquad
\partial_\phi\mathcal B_L=\frac12U^\theta(b_L-b).
\end{equation}
They exist by \eqref{eq:sm-loop-mean}. They vanish on the radial
boundary neighborhoods, so their extensions by zero are smooth.
Define
\begin{equation}\label{eq:sm-velocity-modulation}
U_N^\theta=U^\theta\exp\left(\frac{\mathcal A_L(R,Z,Ny)}N\right),
\qquad
U_N^z=U^z+\frac{\mathcal B_L(R,Z,Ny)}N.
\end{equation}
The phase is independent of $Z$. Thus taking a $Z$ derivative
does not introduce a factor $N$.
After evaluation at $\phi=Ny$, the operator $R\partial_R$
acts as $\partial_y+N\partial_\phi$.
This gives the exact formulas
\begin{equation}\label{eq:sm-exact-shears}
a_N=a_L-\frac{2\partial_y\mathcal A_L}{N},\qquad
b_N=e^{-\mathcal A_L/N}
\left(b_L+\frac{2\partial_y\mathcal B_L}{NU^\theta}\right),
\end{equation}
where all terms on the right are evaluated at $\phi=Ny$.
The exponential factor in the second formula is needed because
$b_N$ is normalized by the new angular velocity.

Compute $\mathbf M_N$ from the axis and, throughout the intermediate
construction, set
\begin{equation}\label{eq:sm-pressure-convention}
P_N=P_0+M_N^p.
\end{equation}
This keeps the core pressure exactly unchanged. Before moment
restoration, $P_N$ need not vanish at infinity. We use the
moment formulas for $\mathcal{I}_N$ with this $P_N$; there is
no need to impose the final normalization prematurely.

\paragraph{One explicit norm for the estimates.}
Put $\mathcal D=[R_{-c},2R_c]\times[-1,1]$. For a function on
$\mathcal D$, let $\|f\|_1=\sup_R\|f(R,\cdot)\|_{C_Z^1}$.
For a function also depending on $\phi$, take the supremum over
$(R,\phi)$ instead. For tuples use the sum of component norms.
These norms obey $\|fg\|_1\le\|f\|_1\|g\|_1$.
Define the following derived number, independent of $N$:
\begin{equation}\label{eq:sm-K}
\begin{aligned}
K={}&2+R_{-c}+R_{-c}^{-1}+2R_c+(2R_c)^{-1}\\
&+\|U^\theta,U^z,(U^\theta)^{-1},\mathbf M,P_0\|_{1,\mathcal D}
 +\|A_c^\theta,(A_c^\theta)^{-1}\|_{C_Z^1}\\
&+\|a,b,p_1,p_2\|_{L^\infty(\mathcal D)}\\
&+\|\mathcal A_L,\mathcal B_L,\partial_y\mathcal A_L,\partial_y\mathcal B_L,
          a_L,b_L\|_{1,\mathcal K\times(\mathbb R/\mathbb Z)}.
\end{aligned}
\end{equation}
All entries are known after constructing the loop. This norm
records, in particular, possible losses from a small angular
velocity, small radial scales, and small input cone margins.
For a family of inputs, $K$ must be evaluated or bounded for
that family; it is not an absolute constant.

For $N\ge2K$, the elementary inequality
$e^s-1\le2s$ for $0\le s\le1/2$, applied in the above product
norm, gives
\begin{equation}\label{eq:sm-first-estimates}
\begin{aligned}
\|U_N^\theta-U^\theta,U_N^z-U^z\|_1
&\le\frac{3K^2}{N},\\
\|a_N-a_L(\cdot,Ny),b_N-b_L(\cdot,Ny)\|_\infty
&\le\frac{8K^2}{N},\\
\|\mathbf M_N-\mathbf M,P_N-P\|_{1,\mathcal D}
&\le\frac{CK^5}{N}.
\end{aligned}
\end{equation}
Here and below $a_L,b_L$ are extended as the original $a,b$
off $\mathcal K$. To check the last bound, use the five
integrands in \eqref{fiveM}. The interval length is at most
$K$, both weights $\sqrt{2R}$ and $(2R)^{-1}$ are at most
$2K$, and a difference of quadratic products is at most
$CK$ times the velocity difference. These statements hold
also in $C_Z^1$. There are no large radial derivatives in
these integrands.

\subsection{Restoring all five moments}
\label{subsec:sm-moments}

We correct the five moment defects in the reserved interval, where
the velocity still equals \eqref{eq:sm-patch-input}. Before this correction,
$U_N=U$ for $R\ge R_{+c}$, so the moment integrands agree there and
$\mathbf M_N-\mathbf M$ is constant in $R$.
Write $D=\mathbf M_N(R_c,\cdot)-\mathbf M(R_c,\cdot)$ and normalize it as follows:
\begin{equation}\label{eq:sm-defect}
\mathbf D=\left(
\frac{D^z}{R_cA_c^\theta},\quad
\frac{D^{\theta z}}{\sqrt2R_c^{3/2}(A_c^\theta)^2},\quad
\frac{D^\theta}{\sqrt2R_c^{3/2}A_c^\theta},\quad
\frac{D^{z\theta}}{R_c(A_c^\theta)^2},\quad
\frac{D^p}{(A_c^\theta)^2}\right),\qquad
e=\|\mathbf D\|_{C_Z^1}.
\end{equation}
Equation \eqref{eq:sm-K} bounds all five normalization factors
in $C_Z^1$ by $CK^4$. Hence
\begin{equation}\label{eq:sm-defect-bound}
e\le CK^9/N.
\end{equation}

Use two nonnegative, nonzero smooth bumps $\beta_1,\beta_2$ and three
such bumps $\gamma_1,\gamma_2,\gamma_3$, compactly supported in $(1,2)$
with disjoint, ordered supports within each family. We take the
$\gamma_j$ from Interval I.4 and $(\beta_1,\beta_2)=(\gamma_1,\gamma_3)$,
as in Subsection~\ref{roadin:shear-modification}. Write
\[
f_h(x,Z)=\sum_{j=1}^3d_j(Z)\gamma_j(x),\qquad
g_h(x,Z)=\sum_{i=1}^2c_i(Z)\beta_i(x),\qquad
h=(c_1,c_2,d_1,d_2,d_3).
\]
On the reserved interval make the correction
\begin{equation}\label{eq:sm-moment-bumps}
\widehat U^\theta=U_N^\theta+A_c^\theta f_h,
\qquad \widehat U^z=U_N^z+A_c^\theta g_h.
\end{equation}
Use $\widehat U=U_N$ elsewhere. In the normalization
\eqref{eq:sm-defect}, the exact change in the terminal moments is
\begin{equation}\label{eq:sm-moment-map}
\mathbf F_\mu(h)=
\begin{pmatrix}
\displaystyle\int g_h\\[1mm]
\displaystyle\int x^{-\mu}g_h+\int\sqrt{x}\,f_hg_h\\[1mm]
\displaystyle\int\sqrt{x}\,f_h\\[1mm]
\displaystyle-\int x^{-1/2-\mu}f_h+\int(g_h^2-f_h^2/2)\\[1mm]
\displaystyle\int x^{-3/2-\mu}f_h+\int f_h^2/(2x)
\end{pmatrix}
=A_\mu h+Q_\mu(h,h),
\end{equation}
where every integral is over $1<x<2$.
In particular the factors $R_c$ and $A_c^\theta$ have disappeared
from this finite-dimensional map.

\begin{lemma}[Cost of the moment correction]\label{lem:sm-moments}
There are constants $C_A,C_Q\ge1$, depending only on the fixed
bumps, such that
\[
\|A_\mu^{-1}\|\le C_A/\mu,\qquad
\|Q_\mu(h,k)\|_{C_Z^1}\le
C_Q\|h\|_{C_Z^1}\|k\|_{C_Z^1}.
\]
If
\begin{equation}\label{eq:sm-moment-smallness}
e\le\frac{\mu^2}{8C_A^2C_Q},
\end{equation}
then $\mathbf F_\mu(h)=-\mathbf D$ has a smooth solution with
\begin{equation}\label{eq:sm-h-bound}
\|h\|_{C_Z^1}\le\frac{2C_Ae}{\mu}
\le\frac{CK^9}{\mu N}.
\end{equation}
\end{lemma}

\begin{proof}
The linear map has an axial block with weights $(1,x^{-\mu})$
and an angular block with weights
$(x^{1/2},-x^{-1/2-\mu},x^{-3/2-\mu})$.
For the axial block, its determinant is the integral over the
two ordered supports of
\[
\beta_1(x_1)\beta_2(x_2)
       (x_2^{-\mu}-x_1^{-\mu}).
\]
Its absolute value is bounded above and below by fixed positive
multiples of $\mu$, since $1<x_1<x_2<2$ and the supports
are separated. For example, divide by $\mu$ and use
\[
\frac{x_1^{-\mu}-x_2^{-\mu}}{\mu}
=\int_{\log x_1}^{\log x_2}e^{-\mu s}\,ds.
\]
For the angular block the three exponents are distinct for every
$0\le\mu\le1/2$. After reordering the rows and removing the minus sign
from the second weight, Lemma~\ref{lem:bump-moment-matrix} gives invertibility
for the three ordered bumps $\gamma_j$. By continuity in $\mu$ the inverse
is uniformly bounded on this closed interval.
This proves the asserted bound for $A_\mu^{-1}$.
The quadratic bound follows directly from
\eqref{eq:sm-moment-map} and the product inequality in $C_Z^1$.
At this point $C_A,C_Q$ are fixed independently of all input data.

The map
\[
h\longmapsto-A_\mu^{-1}
                  [\mathbf D+Q_\mu(h,h)]
\]
maps the closed ball of radius $2C_Ae/\mu$ in $C_Z^1$ to
itself and has Lipschitz constant at most
$4C_A^2C_Qe/\mu^2\le1/2$.
Its fixed point gives the moment identities.
The same inequality makes the derivative with respect to $h$
invertible. The pointwise implicit function theorem then gives
smoothness in $Z$, including at $Z=\pm1$.
\end{proof}

The factor $\mu^{-1}$ is real: as $\mu\downarrow0$, the two
axial weights both approach $1$. It cannot be absorbed into
an auxiliary constant. This is the reason for the $\mu^{-2}$
term in the frequency requirement below.

By construction,
\begin{equation}\label{eq:sm-moments-restored}
\widehat{\mathbf M}(R,Z)=\mathbf M(R,Z)\qquad(R\ge2R_c).
\end{equation}
Define $\widehat P=P_0+\widehat M^p$. Equality of the pressure moments
restores the pressure normalization at infinity. All five moments agree
at $2R_c$ as functions of $Z$, the angular and axial velocities agree
for $R\ge2R_c$, and the two profiles use the same axis pressure $P_0$.
Lemma~\ref{lem:moment-gluing} therefore gives equality of pressure, radial
velocity, inertial stress, and shear for $R\ge2R_c$.

\subsection{The actual stress and completion of the leading profile}
\label{subsec:sm-actual-cone}

On $\mathcal D=[R_{-c},2R_c]\times[-1,1]$ we compare the actual shear
$(\widehat a,\widehat b)$ with $(a_L,b_L)$ at phase $Ny$, and the actual
normalized inertial stress $\widehat p=\widehat{\mathcal I}/\widehat F$
with the input $p=\mathcal I/F$. Together these estimates transfer the
strict cone inequalities for $(a_L,b_L,p)$ to the corrected profile.
For $X=(a,b,p_1,p_2)$ set
\begin{equation}\label{eq:sm-polynomial-gaps}
\begin{aligned}
G_1(X)&=a,\\
G_2(X)&=a^2+b^2-2a,\\
G_3(X)&=p_1a-p_2b-a^2-b^2,\\
G_4(X)&=2aG_3(X)^2-G_2(X)(p_1b+p_2a)^2.
\end{aligned}
\end{equation}
The four inequalities $G_j>0$ are equivalent to $a>0$ and
the admissible cone. Indeed,
\[
\kappa-2=G_2/a,\qquad
-\mathcal{T}\cdot \mathcal{S}/F^2=G_3,\qquad
\mathcal{T}\cdot \mathcal{S}^\perp/F^2=-(p_1b+p_2a).
\]

Let $\mathcal X$ be the compact family consisting of
$(a_L,b_L,p_1,p_2)$ over $\mathcal K\times(\mathbb R/\mathbb Z)$
and the input $(a,b,p_1,p_2)$ over
$[R_{+c},2R_c]\times[-1,1]$.
Set
\begin{equation}\label{eq:sm-stability}
H_*=2+\max_{X\in\mathcal X}|X|_\infty,\qquad
g_* =\min_{X\in\mathcal X,\ 1\le j\le4}G_j(X)>0,
\qquad
\rho_* =\min\left\{1,\frac{g_*}{512H_*^5}\right\}.
\end{equation}
The margin $g_*$ may depend on all the input data.
For $|X|_\infty\le H_*$, differentiation of
\eqref{eq:sm-polynomial-gaps} gives
\[
\max_j\sum_{\nu\in\{a,b,p_1,p_2\}}
 |(\partial_\nu G_j)(X)|\le148H_*^5.
\]
Thus any perturbation of a state in $\mathcal X$ by at most
$\rho_*$ in the maximum norm retains $G_j\ge g_*/2$.
This also guarantees nonzero stress.

\paragraph{Bounds after the moment correction.}
For $N$ sufficiently large, the modulation and moment estimates give
\begin{equation}\label{eq:sm-final-errors}
\begin{aligned}
\|\widehat U-U\|_{1,\mathcal D}
&\le\frac{CK^{10}}{\mu N},\\
\|\widehat{\mathbf M}-\mathbf M,\widehat P-P\|_{1,\mathcal D}
&\le\frac{CK^{13}}{\mu N},\\
\|\widehat{\mathcal{I}}-\mathcal{I}\|_{L^\infty(\mathcal D)}
&\le\frac{CK^{15}}{\mu N},\\
\|\widehat a-a_L(\cdot,Ny),
\widehat b-b_L(\cdot,Ny),\widehat p-p\|_{L^\infty(\mathcal D)}
&\le\frac{CK^{18}}{\mu N}.
\end{aligned}
\end{equation}
Here $\widehat U-U$ denotes the pair of angular and axial
differences. We use the norm
$\|\cdot\|_1=\sup_R\|\cdot\|_{C_Z^1}$ throughout. Write $\Delta$ for corrected minus input
quantities, and set
\[
E_U=\|\Delta U\|_{1,\mathcal D},\qquad
E_M=\|\Delta\mathbf M,\Delta P\|_{1,\mathcal D}.
\]
\begin{itemize}
\item
The modulation estimate \eqref{eq:sm-first-estimates} and the
correction \eqref{eq:sm-moment-bumps} give
\[
E_U\le \frac{3K^2}{N}+CK\|h\|_{C_Z^1}
    \le\frac{CK^{10}}{\mu N}.
\]
For $E_U\le1$, a difference of quadratic velocity products is
bounded in $C_Z^1$ by $CKE_U$. Each of the five moment integrals
has length at most $K$ and weights bounded by $2K$, as in the
proof of \eqref{eq:sm-first-estimates}. Since
$\Delta P=\Delta M^p$, this yields
\[
E_M\le CK^3E_U\le\frac{CK^{13}}{\mu N}.
\]
\item
Here is the exact passage from moment differences to stress
differences. Recall $L=1-Z^2$ and $L_\delta=1-\delta Z^2$, and write
\[
B=-R+(1-\delta)ZM^z+L\partial_ZM^z,\qquad
\Delta B=(1-\delta)Z\Delta M^z+L\partial_Z\Delta M^z.
\]
Expanding the products in \eqref{Itheta-moments} and \eqref{Iz-moments} gives
\[
\begin{aligned}
\Delta \mathcal{I}^\theta
&=\frac{B\Delta U^\theta+\widehat U^\theta\Delta B}
       {L_\delta\sqrt{2R}}\\
&\quad+\frac1{2L_\delta R}\bigl[(1-\delta/2)\Delta M^\theta
 -\tfrac{1-\delta}{2}Z\partial_Z\Delta M^\theta\\
&\hspace{32mm}-L\partial_Z\Delta M^{\theta z}
 +(2\delta-1)Z\Delta M^{\theta z}\bigr],\\[1mm]
\Delta \mathcal{I}^z
&=\frac1{L_\delta\sqrt{2R}}\bigl[B\Delta U^z
 +\widehat U^z\Delta B
 +\tfrac{1-\delta}{2}(\Delta M^z-Z\partial_Z\Delta M^z)\\
&\hspace{23mm}+2\delta Z\Delta M^{z\theta}
 -L\partial_Z\Delta M^{z\theta}\\
&\hspace{23mm}+R\{2(1+\delta)Z\Delta P
                        -L\partial_Z\Delta P\}\bigr].
\end{aligned}
\]
By \eqref{eq:sm-K}, $|B|\le CK$,
$\|\Delta B\|_\infty\le CE_M$, and
$\|\widehat U\|_\infty\le K+1$. Together with $L_\delta^{-1}\le2$
and $R^{\pm1}\le K$, these give, term by term,
\[
\|\Delta \mathcal{I}\|_\infty
\le CK^2(E_U+E_M)\le\frac{CK^{15}}{\mu N}.
\]
The moment formulas involve only $Z$ derivatives, and the phase
is independent of $Z$.
\item
On the modulation interval the exact identities
\eqref{eq:sm-exact-shears} and their estimate in
\eqref{eq:sm-first-estimates} give a shear error at most $8K^2/N$.
On the correction interval,
\[
\widehat a=1-2x\partial_x\log(x^{-1/2-\mu}+f_h),\qquad
\widehat b=\frac{2x\partial_xg_h}{x^{-1/2-\mu}+f_h}.
\]
Here $1<x<2$, $0<\mu\le1/2$, and the bumps are fixed. For
$\|h\|_{C_Z^1}$ sufficiently small, the denominator stays
bounded below by a fixed positive constant. Differentiating
these finite sums of bumps therefore gives
\[
|\widehat a-(2+2\mu)|+|\widehat b|
\le C\|h\|_{C_Z^1}\le\frac{CK^9}{\mu N}.
\]
These are the errors from the extended loop on this interval,
where $(a_L,b_L)=(2+2\mu,0)$.

Finally, division by the swirl is quantified by
\[
\widehat p-p=\frac{\Delta \mathcal{I}-p\,\Delta F}{\widehat F},
\qquad
\Delta F=\frac{\Delta U^\theta}{\sqrt{2R}},\qquad
F^{-1}=\frac{\sqrt{2R}}{U^\theta}\le2K^2.
\]
In fact $|\Delta F|/F=|\Delta U^\theta|/U^\theta\le KE_U$.
Thus $CK^{11}/(\mu N)$ small ensures $\widehat F\ge F/2$.
Using $|p|\le K$, $|\Delta F|\le CKE_U$, and the preceding
stress estimate now gives
\[
\|\widehat p-p\|_\infty
\le CK^4(E_U+E_M)
\le\frac{CK^{17}}{\mu N}.
\]
The shear and normalized stress estimates prove the last line
of \eqref{eq:sm-final-errors}, with a spare power of $K$.
\end{itemize}
Choose $C_{\mathrm{sm}}$ sufficiently large, depending only on the fixed
step and bumps, to cover the preceding estimates. For any requested
accuracy $0<\varepsilon\le1$, take an integer $N$ satisfying
\begin{equation}\label{eq:sm-frequency-choice}
\boxed{\qquad
N\ge N_0:=
\left\lceil C_{\mathrm{sm}}K^{20}
\left(\mu^{-2}+(\mu\varepsilon)^{-1}
                         +(\mu\rho_*)^{-1}\right)\right\rceil.
\qquad}
\end{equation}
It implies $N\ge2K$, \eqref{eq:sm-moment-smallness},
$\widehat U^\theta>0$, the error bound
\eqref{eq:sm-output-close}, and a normalized state error at most
$\rho_*$. Equations \eqref{eq:sm-polynomial-gaps}--
\eqref{eq:sm-stability} therefore give the admissible cone
throughout $\mathcal D$.

The positive margin $g_*$ is used only on $\mathcal D$. For $R\le R_{-c}$
the velocities, the axis pressure, and all moments integrated from the
axis agree with the input. For $R\ge2R_c$ exact agreement follows from
\eqref{eq:sm-moments-restored} and Lemma~\ref{lem:moment-gluing}.
Thus the input cone inequalities hold throughout both open endpoint
collars. The stress still vanishes on the core and the heat exterior,
with any prescribed flat vanishing at $R_a$ and $R_h$ preserved.
This proves Proposition~\ref{prop:sm}.

\begin{proof}[Proof of Theorem~\ref{thm:leading}]
\leavevmode\phantomsection\label{subsec:sm-completion}%
Take the global profile of Proposition~\ref{prop:global-relaxed-profile}
for a sufficiently large member of the family in
Lemma~\ref{lem:imc-prepared-family}. Its axis pressure, velocities,
and moments are now fixed.
Choose $R_{-c}$ strictly inside the admissible collar $(R_a,R_{ac})$
furnished by the inner connection, and retain
\[
 R_{+c}=eR_{-\mu},\qquad R_c=e^{-24}R_{mz}.
\]
The input is admissible on $(R_a,R_{-c}]\cup[R_{+c},R_h)$.
The strict inequalities, smoothness, and compactness of $[-1,1]$
also give admissibility on neighborhoods of both endpoints of
$[R_{-c},R_{+c}]$.

The radius ordering is explicit. By \eqref{eq:road-mu-scales},
$R_{mz}=\mu^{-60}R_{+c}$, and hence
\[
 \frac{R_c}{R_{+c}}=e^{-24}\mu^{-60}>1,\qquad
 2R_c=2e^{-24}R_{mz}<R_{mz}<R_h.
\]
Here $\mu\le1/60$ and $R_{\rm useless}<R_{\rm ref}<R_{+c}$.
Thus $R_a<R_{-c}<R_{+c}<R_c<2R_c<R_h$.
Moreover, $[R_c,2R_c]$ lies strictly inside
$I_1=(e^{-25}R_{mz},e^{-20}R_{mz})$, where the earlier corrections
leave the pure-power swirl and zero axial velocity unchanged;
see Subsection~\ref{subsec:outer-buffer} and
\eqref{eq:road-reserved-patches}.
This is exactly the moment-restoration input
\eqref{eq:sm-patch-input}, with $A_c^\theta=U^\theta(R_c,\cdot)$.
Thus the same profile satisfies hypotheses (i)--(iii) of
Proposition~\ref{prop:sm}. Apply it with any $0<\varepsilon\le1$ and
a finite $N$ as in \eqref{eq:sm-frequency-choice}. The two
modifications are supported in
$(R_{-c},R_{+c})$ and $(R_c,2R_c)$.
They preserve the inner collar up to $R_{-c}$ and, after moment
restoration, all fields for $R\ge2R_c$.
Since $2R_c<e^{-14}R_{mz}=\inf I_2<\inf I_3$,
both later reserved intervals remain available.

Set $U^\theta_{(0)}=\widehat U^\theta$ and
$U^z_{(0)}=\widehat U^z$, and recover $P_{(0)}$ and $U^r_{(0)}$
by \eqref{P-Utheta} and \eqref{Ur-Uz}.
The modified shear has the admissible cone throughout the open
stress annulus. The five terminal moments, the prescribed axis
pressure, the regular stress-free core, and the exact heat exterior
are preserved. Smoothness of $F_{(0)},U^z_{(0)}$ and boundedness of
$U^\theta_{(0)},U^z_{(0)}$ follow from the input and the smooth
compactly supported modifications.
The identities for pressure, radial velocity, and physical residuals
then follow from Sections~\ref{sec:axi} and
\ref{sec:leading-order-system}. In particular, the two tangential
residuals have the divergence representation \eqref{eq:residual-form};
the remaining terms have the formal orders already computed in
Subsection~\ref{sec:orders}. This proves Theorem~\ref{thm:leading}.
\end{proof}

\section{The lower-order corrections and the complete background}
\label{sec:higher-order-correction}\label{sec:lower-order}

Fix the leading profile completed in Subsection~\ref{subsec:sm-completion}.
Here and below, \emph{tangential} refers to the angular and axial
components, in the directions $e_\theta$ and $e_z$ tangent to the
cylinders $r=\mathrm{constant}$.
The radii of the leading profile and $\delta$ are fixed throughout.
The leading index is $n=0$. For every coefficient index, use the regular
variables
\[
 F_{(n)}=\frac{U^\theta_{(n)}}{\sqrt{2R}},\qquad
 V_{(n)}=\sqrt{2R}\,U^r_{(n)}.
\]
Here $V_{(n)}$ is the radial primitive. Write $P_{(n)}$ for the pressure coefficient and
$\mathcal{T}_{(n)}=(\mathcal{T}^\theta_{(n)},\mathcal{T}^z_{(n)})$ for the normalized stress.
All these profiles depend on $(R,Z)$.
For $n\ge0$ put
\begin{equation}\label{lo:coefficient-weights}
 e_n=2n\delta,\qquad b_n=-2-\delta+e_n,\qquad
 c_n=-1-\delta+e_n,\qquad \pi_n=-2-2\delta+e_n.
\end{equation}
Since $r=\lambda\sqrt{2R}$, the physical coefficient fields are
\begin{equation}\label{lo:physical-coefficients}
\begin{aligned}
 u^\theta_{(n)}&=\lambda^{c_n}U^\theta_{(n)}
              =r\lambda^{b_n}F_{(n)},&
 u^z_{(n)}&=\lambda^{c_n}U^z_{(n)},\\
 r u^r_{(n)}&=\lambda^{e_n}V_{(n)},&
 p_{(n)}&=\lambda^{\pi_n}P_{(n)},&
 \mathsf T_{(n)}&=\lambda^{-2-\delta+e_n}\mathcal{T}_{(n)}.
\end{aligned}
\end{equation}
We seek the formal expansion
$(u,p,\mathsf T)\sim\sum_{n\ge0}(u_{(n)},p_{(n)},\mathsf T_{(n)})$.
Thus consecutive orders differ by a factor $\lambda^{2\delta}$.
Its smooth realization is constructed after the coefficient induction.

We first write the matching equations for these physical fields in
$(t,r,z)$ coordinates. Every derivative below acts on the full physical
field, including its $\lambda$ weight. Put
$\Delta_r=\partial_r^2+r^{-1}\partial_r$.
All indices in the sums are nonnegative, and negative-index fields are zero.
At each order $n\ge1$, incompressibility requires
\begin{equation}\label{lo:physical-divergence}
 \partial_r(ru^r_{(n)})+\partial_z(ru^z_{(n)})=0.
\end{equation}
The angular and axial matching equations are
\begin{equation}\label{lo:physical-tangential-matching}
\begin{aligned}
 \partial_tu^\theta_{(n)}
 &+\sum_{i+j=n}\left[
 u^r_{(i)}\left(\partial_r+\frac1r\right)u^\theta_{(j)}
 +u^z_{(i)}\partial_zu^\theta_{(j)}\right]\\
 &-\left(\Delta_r-\frac1{r^2}\right)u^\theta_{(n)}
 -\partial_z^2u^\theta_{(n-1)}
 =-\left(\partial_r+\frac2r\right)\mathsf T^\theta_{(n)},\\[1mm]
 \partial_tu^z_{(n)}
 &+\sum_{i+j=n}\left[
 u^r_{(i)}\partial_ru^z_{(j)}+u^z_{(i)}\partial_zu^z_{(j)}\right]
 +\partial_zp_{(n)}\\
 &-\Delta_ru^z_{(n)}-\partial_z^2u^z_{(n-1)}
 =-\left(\partial_r+\frac1r\right)\mathsf T^z_{(n)}.
\end{aligned}
\end{equation}
For the radial equation, define the source at index $k\ge0$ by
\begin{equation}\label{lo:physical-radial-source}
\begin{aligned}
 \mathcal H^r_{(k)}={}&\partial_tu^r_{(k)}
 +\sum_{i+j=k}\left[
 u^r_{(i)}\partial_ru^r_{(j)}+u^z_{(i)}\partial_zu^r_{(j)}\right]\\
 &-\left(\Delta_r-\frac1{r^2}\right)u^r_{(k)}
 -\partial_z^2u^r_{(k-1)}.
\end{aligned}
\end{equation}
The pressure and centrifugal terms must satisfy
\begin{equation}\label{lo:physical-radial-matching}
 \partial_rp_{(n)}-\frac1r\sum_{i+j=n}u^\theta_{(i)}u^\theta_{(j)}
 =-\mathcal H^r_{(n-1)}.
\end{equation}
At order $n$, $\mathcal H^r_{(n-1)}$ depends only on the completed lower orders.
The two index shifts reflect the scaling of the physical equations.
Axial viscosity has an extra factor $\lambda^{2\delta}$ relative to
the time derivative and radial viscosity, so the tangential equations
use the axial viscosity of order $n-1$.
In the radial equation, the order-$n$ pressure and centrifugal terms
balance the time, transport, and radial-viscosity terms of order $n-1$;
their axial-viscosity source comes from order $n-2$.
In particular, $\mathcal H^r_{(0)}$ excludes $-\partial_z^2u^r_{(0)}$,
which is matched at order two.
In each quadratic sum of total order $n$, only the pairs $(0,n)$ and
$(n,0)$ contain a new field, so the current matching equations are linear.

On the inner interval we impose $\mathsf T_{(n)}=0$ for $n\ge1$.
Subsection~\ref{lo:coefficient-equations} derives the corresponding
profile equations by substituting \eqref{lo:physical-coefficients}
into these physical matching equations and removing their common weights.
For the first correction we solve the linear equations for
$U^\theta_{(1)},U^z_{(1)},P_{(1)}$ on an inner interval extending beyond $R_a$.
As in Figure~\ref{fig:lo-cutoff-moments}, we cut off the two velocity
profiles before $R_{-c}$ and add three angular and two axial bumps in $I_2$
to cancel the five moments. The five equations are linear in the bump
coefficients and can be solved for defects of any size. From the corrected velocities we recover
the radial velocity and pressure, then integrate the two tangential residuals
to obtain the stress. We repeat this construction on the same intervals
at every positive order and sum the resulting profiles to prove
Theorem~\ref{thm:background}, following Section~5 of \cite{1}.

\begin{figure}[!ht]
\centering
\begin{tikzpicture}[x=1cm,y=.68cm,>=Latex,font=\small]
 \fill[roadblue!6] (0,-.55) rectangle (2.25,3.4);
 \fill[roadorange!9] (5.0,-.55) rectangle (8.05,3.4);
 \fill[roadpurple!7] (8.65,-.55) rectangle (9.65,3.4);
 \draw[->] (0,0)--(13.35,0) node[right] {$R$};
 \foreach \x/\lab in {0/0,1.45/R_a,3.35/R_{\rm in},4.05/R_{-c},10.55/R_{0z},12.9/R_h}
  {\draw (\x,.08)--(\x,-.08);\node[below] at (\x,-.08) {$\lab$};}
 \node at (6.53,-.29) {$I_2$};
 \node at (9.15,-.29) {$I_3$};
 \draw[dashed,gray] (2.25,0)--(2.25,2.8);
 \draw[dashed,gray] (3.35,0)--(3.35,2.45);
 \node[above,align=center] at (1.12,2.8) {inner solution\\and collar};
 \node[above,align=center] at (3.25,2.8) {smooth\\cutoff};
 \node[above,align=center] at (6.53,2.8) {three angular bumps\\and two axial bumps};
 \foreach \y in {1.35,.35} {\draw[gray,dashed] (0,\y)--(12.9,\y);}
 \draw[roadblue,thick] (0,1.35)
  .. controls (.35,1.35) and (.8,2.05) .. (1.45,2.05)
  --(2.25,2.05)
  .. controls (2.7,2.05) and (2.65,1.35) .. (3.25,1.35)
  --(5.1,1.35);
 \draw[roadblue,thick] (5.1,1.35)
  .. controls (5.22,1.35) and (5.2,2.15) .. (5.35,2.15)
  .. controls (5.5,2.15) and (5.48,1.35) .. (5.6,1.35)
  --(6.3,1.35)
  .. controls (6.42,1.35) and (6.4,2.15) .. (6.55,2.15)
  .. controls (6.7,2.15) and (6.68,1.35) .. (6.8,1.35)
  --(7.5,1.35)
  .. controls (7.62,1.35) and (7.6,2.15) .. (7.75,2.15)
  .. controls (7.9,2.15) and (7.88,1.35) .. (8,1.35)
  --(12.9,1.35);
 \draw[roadorange,thick] (0,.35)
  .. controls (.35,.35) and (.8,.85) .. (1.45,.85)
  --(2.25,.85)
  .. controls (2.7,.85) and (2.65,.35) .. (3.25,.35)
  --(5.7,.35)
  .. controls (5.82,.35) and (5.8,1.08) .. (5.95,1.08)
  .. controls (6.1,1.08) and (6.08,.35) .. (6.2,.35)
  --(6.9,.35)
  .. controls (7.02,.35) and (7,1.08) .. (7.15,1.08)
  .. controls (7.3,1.08) and (7.28,.35) .. (7.4,.35)
  --(12.9,.35);
 \node[roadblue,above] at (11.5,1.35) {$U^\theta_{(n)}=0$};
 \node[roadorange,above] at (11.5,.35) {$U^z_{(n)}=0$};
\end{tikzpicture}
\caption{Supports of the positive-order velocity corrections. The radial
cutoff is one on a fixed inner interval extending beyond $R_a$ and
vanishes before the inner-solution endpoint $R_{\rm in}<R_{-c}$.
Five fixed bumps in $I_2$ cancel the five moments. All positive-order
velocities vanish for $R\ge\sup I_2$. The recovered pressure correction
may persist on $I_3$ but vanishes for $R\ge R_{0z}$.
Only the first stress correction can reach $R_h$.}
\label{fig:lo-cutoff-moments}
\end{figure}

\subsection{Coefficient equations}
\label{lo:coefficient-equations}

We now pass from the physical matching equations to the profile coefficients.
The prescribed axis pressure is
$P_0(Z)=P_{(0)}(0,Z)$, distinct from the full profile $P_{(0)}(R,Z)$.
We keep this axis datum unchanged at every order by imposing
$P_{(n)}(0,Z)=0$ for $n\ge1$; only the pressure away from the axis
receives higher-order corrections.

Recall $L=1-Z^2$ and $L_\delta=1-\delta Z^2$. For a real weight $a$,
define
\begin{equation}\label{lo:weighted-operators}
\begin{aligned}
 \mathscr T_a G
 &=\frac1{L_\delta}\left(-\frac a2G+\frac{1-\delta}{2}Z\partial_ZG
                        +R\partial_RG\right),\\
 \mathscr Z_a G
 &=\frac1{L_\delta}\left(aZG+L\partial_ZG-2ZR\partial_RG\right),\\
 \mathscr Z_a^{[2]}G
 &=\mathscr Z_{a-1+\delta}\bigl(\mathscr Z_aG\bigr).
\end{aligned}
\end{equation}
The order of the two operators in the last line matters. The chain
rule from Section~\ref{sec:axi} gives the exact identities
\begin{equation}\label{lo:weighted-derivatives}
\begin{aligned}
 \partial_t(\lambda^aG)&=\lambda^{a-2}\mathscr T_aG,&
 \partial_z(\lambda^aG)&=\lambda^{a-1+\delta}\mathscr Z_aG,\\
 \partial_z^2(\lambda^aG)&=\lambda^{a-2+2\delta}
                                     \mathscr Z_a^{[2]}G.
\end{aligned}
\end{equation}
In particular, the second axial derivative increases the expansion
order by one relative to the radial and time derivatives.

Set
\[
 M^z_{(n)}(R,Z)=\int_0^R U^z_{(n)}(x,Z)\,dx,\qquad
 \mathcal A_RU^z_{(n)}=\frac{M^z_{(n)}}R
                    =\int_0^1U^z_{(n)}(tR,Z)\,dt.
\]
The last expression supplies the regular value at zero.
Incompressibility and the regular axis value $V_{(n)}(0,Z)=0$ give
\begin{equation}\label{lo:incompressibility}
\begin{aligned}
 \partial_RV_{(n)}&=-\mathscr Z_{c_n}U^z_{(n)},\\
 V_{(n)}
 &=\frac{2ZR U^z_{(n)}-(1-\delta+e_n)Z M^z_{(n)}
                       -L\partial_ZM^z_{(n)}}{L_\delta}.
\end{aligned}
\end{equation}
Indeed, integration by parts gives
$\int_0^R x\partial_xU^z_{(n)}\,dx=RU^z_{(n)}-M^z_{(n)}$,
and $c_n+2=1-\delta+e_n$.

At the inductive step $n\ge1$, the unknowns are $F_{(n)},U^z_{(n)},P_{(n)}$,
with $V_{(n)}$ given by \eqref{lo:incompressibility}.
We keep the completed leading profile fixed and assume that the profiles
of orders $1\le j<n$ have already been extended to $R\ge0$ and corrected
to satisfy their five moment conditions, as described in
Subsection~\ref{lo-ext:moments}. In a product with $i+j=n$, only $(0,n)$ and $(n,0)$
involve a new coefficient, and each pairs it with the fixed leading
profile. Thus the equations below are linear in the new unknowns.
Negative-index coefficients are taken to be zero.

On the inner interval, cancellation of the angular coefficient gives
\begin{equation}\label{lo:angular-recursion}
\begin{aligned}
 2\bigl(R\partial_R^2F_{(n)}+2\partial_RF_{(n)}\bigr)
 ={}&\mathscr T_{b_n}F_{(n)}
 +\sum_{i+j=n}\left[
 V_{(i)}\left(\partial_RF_{(j)}+\frac{F_{(j)}}R\right)
       +U^z_{(i)}\mathscr Z_{b_j}F_{(j)}\right]\\
 &-\mathscr Z_{b_{n-1}}^{[2]}F_{(n-1)}.
\end{aligned}
\end{equation}
The axial coefficient equation is
\begin{equation}\label{lo:axial-recursion}
\begin{aligned}
 2\bigl(R\partial_R^2U^z_{(n)}+\partial_RU^z_{(n)}\bigr)
 ={}&\mathscr T_{c_n}U^z_{(n)}
 +\sum_{i+j=n}\left[
 V_{(i)}\partial_RU^z_{(j)}
       +U^z_{(i)}\mathscr Z_{c_j}U^z_{(j)}\right]\\
 &+\mathscr Z_{\pi_n}P_{(n)}
       -\mathscr Z_{c_{n-1}}^{[2]}U^z_{(n-1)}.
\end{aligned}
\end{equation}
The radial pressure equation, imposed globally after extension, is
\begin{equation}\label{lo:pressure-recursion}
 \partial_RP_{(n)}=\sum_{i+j=n}F_{(i)}F_{(j)}
                         -\frac{\Omega_{(n-1)}}{2R},
\end{equation}
where the known radial source is
\begin{equation}\label{lo:radial-source}
\begin{aligned}
 \Omega_{(k)}={}&\mathscr T_{e_k}V_{(k)}
 +\sum_{i+j=k}\left[
 V_{(i)}\left(\partial_RV_{(j)}-\frac{V_{(j)}}{2R}\right)
       +U^z_{(i)}\mathscr Z_{e_j}V_{(j)}\right]\\
 &-2R\partial_R^2V_{(k)}
       -\mathscr Z_{e_{k-1}}^{[2]}V_{(k-1)}.
\end{aligned}
\end{equation}
For $k=0$ the final term is zero. We impose the positive-order axis data
\begin{equation}\label{lo:axis-values}
 F_{(n)}(0,Z)=U^z_{(n)}(0,Z)=P_{(n)}(0,Z)=0
 \qquad(n\ge1).
\end{equation}
In particular, no correction changes the prescribed axis pressure $P_0$.

We now derive \eqref{lo:angular-recursion}--\eqref{lo:pressure-recursion}
in more detail from the physical matching equations, checking the
radial differential operators and the powers of $\lambda$.
First, the radial parts of the angular and axial Laplacians applied
to \eqref{lo:physical-coefficients} are
\[
\begin{aligned}
 \left(\partial_r^2+\frac1r\partial_r-\frac1{r^2}\right)
       (r\lambda^{b_n}F_{(n)})
 &=2r\lambda^{b_n-2}
       (R\partial_R^2F_{(n)}+2\partial_RF_{(n)}),\\
 \left(\partial_r^2+\frac1r\partial_r\right)
       (\lambda^{c_n}U^z_{(n)})
 &=2\lambda^{c_n-2}
       (R\partial_R^2U^z_{(n)}+\partial_RU^z_{(n)}).
\end{aligned}
\]
For $i+j=n$, the angular transport terms have the common factor
$r\lambda^{b_n-2}$, and the axial transport and pressure terms
have the factor $\lambda^{c_n-2}$.
The preceding axial viscosity has the same respective factors because
$b_{n-1}+2\delta=b_n$ and $c_{n-1}+2\delta=c_n$.
In the radial equation,
\begin{equation}\label{lo:radial-balance-factors}
 r\partial_rp_{(n)}-\sum_{i+j=n}u^\theta_{(i)}u^\theta_{(j)}
 =2R\lambda^{\pi_n}
       \left(\partial_RP_{(n)}-\sum_{i+j=n}F_{(i)}F_{(j)}\right).
\end{equation}
The physical radial source \eqref{lo:physical-radial-source} satisfies
\[
 r\mathcal H^r_{(k)}=\lambda^{e_k-2}\Omega_{(k)}.
\]
Since $\pi_n=e_{n-1}-2$, the physical radial balance
\eqref{lo:physical-radial-matching} therefore gives \eqref{lo:pressure-recursion}.
Omitting this source would leave the preceding radial defect uncancelled.

The terms involving only lower orders are
\begin{equation}\label{lo:separated-sources}
\begin{aligned}
 N^\theta_{(n)}={}&\sum_{i=1}^{n-1}
 \left[V_{(i)}\left(\partial_RF_{(n-i)}+\frac{F_{(n-i)}}R\right)
       +U^z_{(i)}\mathscr Z_{b_{n-i}}F_{(n-i)}\right]
       -\mathscr Z_{b_{n-1}}^{[2]}F_{(n-1)},\\
 N^z_{(n)}={}&\sum_{i=1}^{n-1}
 \left[V_{(i)}\partial_RU^z_{(n-i)}
       +U^z_{(i)}\mathscr Z_{c_{n-i}}U^z_{(n-i)}\right]
       -\mathscr Z_{c_{n-1}}^{[2]}U^z_{(n-1)},\\
 N^p_{(n)}={}&\sum_{i=1}^{n-1}F_{(i)}F_{(n-i)}
                         -\frac{\Omega_{(n-1)}}{2R}.
\end{aligned}
\end{equation}
An empty sum is zero. To reconstruct the pressure after extending the
velocities, we denote its prescribed radial derivative by
\begin{equation}\label{lo:pressure-source}
 G_n:=2F_{(0)}F_{(n)}+N^p_{(n)}=\partial_RP_{(n)}.
\end{equation}
In Subsection~\ref{lo-ext:moments}, integration of $G_n$ from the axis
will recover $P_{(n)}$. Its integrals $\int_0^\infty G_n\,dR$ and
$\int_0^\infty R G_n\,dR$ enter the two pressure-dependent moment
conditions: the first removes the constant pressure tail, and the
negative of the second gives the pressure contribution to the integrated
axial flux once the first integral vanishes.
For completed lower indices write $V_{(j)}=Rv_{(j)}$.
The identities
\[
 \mathscr T_a(Rg)=R\mathscr T_{a-2}g,\qquad
 \mathscr Z_a(Rg)=R\mathscr Z_{a-2}g
\]
show that every term of $\Omega_{(n-1)}$ contains a factor $R$.
Consequently $\Omega_{(n-1)}/R$ and the apparent angular quotients
$V_{(i)}F_{(j)}/R$ are smooth at the axis, with all fixed axial
derivatives. These regular expressions are used in the inner solve.

\subsection{The inner coefficients on a common interval}
\label{sec:lower-order-inner}\label{lo-core:subsection}

We solve every positive order on the same interval $[0,R_{\rm in}]$,
which extends beyond $R_a$. We first record the regularity of the
leading coefficients and then prove that the radial endpoint need not
shrink with the order.
After fixing the leading profile, we may choose any
$R_{\rm in}\in(R_a,R_{-c})$, for example $(R_a+R_{-c})/2$, and keep it
fixed throughout the coefficient induction. The upper restriction
keeps the interval before the shear modification, where the required
$Z$-analyticity is available. The strict lower inequality leaves room
to retain the inner solution on a neighborhood of $[0,R_a]$ when it is
later cut off.

\begin{lemma}[Analytic inner data for the coefficient construction]
\label{prop:lo-prepared-input}\label{lo-core:analytic-input}
For the completed leading profile, fix a radius
\begin{equation}\label{lo-core:fixed-radii}
 R_a<R_{\rm in}<R_{-c}.
\end{equation}
The scalar profiles $F_{(0)},U^z_{(0)},P_{(0)}$ and $V_{(0)}/R$
are smooth in $R$ on $[0,R_{\rm in}]$. They and every fixed radial
derivative extend holomorphically in $Z$ to a common neighborhood
of $[-1,1]$, with bounds uniform on that radial interval.
The neighborhood, bounds, and radius may depend on the final leading
profile.
\end{lemma}
\begin{proof}
Proposition~\ref{prop:continuation-core} gives a jointly analytic core on
a neighborhood of $[0,4.1/\Lambda]\times[-1,1]$.
The inner collar satisfies $R_a=4/\Lambda<R_{-c}<R_{ac}<R_ae^{h_b}<4.1/\Lambda$.
The inner moment correction is supported beyond $R_m>R_{ac}$, and
the final shear modification starts at $R_{-c}$. On $[0,R_{\rm in}]$
the completed profile therefore consists of the core up to $R_a$ and
the inner connection beyond it.

Put $y=\log(R/R_a)$. For $0\le y\le\log(R_{\rm in}/R_a)<h_b$,
the comparison fields in \eqref{eq:inner-comparison} agree with the
analytic core. Their moments, pressure, and inertial stresses are
obtained by finite radial integrals and $Z$ derivatives.
Recall that the normalized inertial stresses of the comparison profile are
\[
 \bar D=\frac{\bar{\mathcal I}^\theta}{\bar F},\qquad
 \bar E=\frac{\bar{\mathcal I}^z}{\bar F},\qquad
 \bar F=\frac{\bar U^\theta}{\sqrt{2R}}.
\]
Choose a common smaller complex neighborhood on which these
expressions are holomorphic and the denominators $L_\delta$ and
$\bar F$ have no zeros, uniformly in the compact radial interval.
Consequently $\bar D$ and $\bar E$ are holomorphic there.

The actual bridge is given by the explicit integration formulas
\eqref{eq:inner-actual-bridge}. The cutoff $\chi_b$ depends only on
the real radial variable. Integration over that variable therefore
preserves holomorphy in $Z$, even though $\chi_b$ is only smooth.
Differentiating these formulas any fixed number of times in $R$
preserves the same complex domain. On a smaller neighborhood with
compact closure, each such derivative has a uniform bound.
At $R=R_a$, the infinite-order flatness of $1-\chi_b$ makes all
radial derivatives agree with the core. The agreement also holds
on the complex neighborhood, since the holomorphic jets agree for
real $Z$.

Pressure and the regular radial-velocity quotient are recovered by
\[
\begin{aligned}
 P_{(0)}(R,Z)&=P_0(Z)+\int_0^R F_{(0)}(s,Z)^2\,ds,\\
 \frac{V_{(0)}}R
 &=\frac{2ZU^z_{(0)}-(1-\delta)Z\mathcal A_RU^z_{(0)}
                    -L\partial_Z\mathcal A_RU^z_{(0)}}{L_\delta}.
\end{aligned}
\]
The averaging formula preceding \eqref{lo:incompressibility}
has no singular division. These formulas preserve the same
regularity, including at $R=0$.
\end{proof}

\begin{lemma}[A common radius at every positive order]
\label{lo-core:fixed-radius}
Fix $n\ge1$ and suppose all lower orders have already been completed.
Suppose their scalar profiles, including $V_{(j)}/R$, are smooth
on $[0,R_{\rm in}]\times[-1,1]$. Assume that each fixed radial
derivative is holomorphic near $[-1,1]$ in $Z$, with a bound
uniform on that radial interval.
Then \eqref{lo:incompressibility}--\eqref{lo:radial-source} have a
unique regular order-$n$ solution on $[0,R_{\rm in}]\times[-1,1]$
with
\begin{equation}\label{lo-core:zero-axis}
 F_{(n)}(0,Z)=U^z_{(n)}(0,Z)=P_{(n)}(0,Z)=0.
\end{equation}
The solution has the same smoothness and parameter regularity as
the lower orders. More precisely, for each fixed $n\ge1$ and $k\ge0$,
the $k$th radial derivatives of the regular scalar profiles extend
holomorphically in $Z$ to a neighborhood of $[-1,1]$, with bounds
uniform for $0\le R\le R_{\rm in}$. This neighborhood and these bounds
may depend on $n$ and $k$, whereas $R_{\rm in}$ is fixed independently
of both. Uniqueness is among regular solutions having
a bounded holomorphic parameter extension on some neighborhood of
$[-1,1]$, uniformly on this radial interval.
\end{lemma}
\begin{proof}
We write the equations as a first-order system in $\xi=\sqrt R$.
The $Z$ derivatives act only on the first four unknowns and enter only
the last two equations. The radial inverse is diagonal, so two consecutive
derivative factors in the integral iteration give zero. A $k$-fold
iteration thus costs at most $\lceil k/2\rceil$ derivatives in $Z$.
The factor $1/k!$ from the ordered radial integrals absorbs this loss
on $[0,\sqrt{R_{\rm in}}]$, even if the complex $Z$ neighborhood
shrinks with $n$.

\medskip
\noindent\textbf{Step 1: separate the current unknowns from the sources.}
Keep the leading coefficients fixed and put
\begin{equation}\label{lo-core:leading-coefficients}
\begin{gathered}
 f=F_{(0)},\qquad w=U^z_{(0)},\qquad v=V_{(0)}/R,\qquad
 v_{(j)}=V_{(j)}/R,\\
 H_c=\frac{1-\delta}{2}Z+Lw,\qquad
 B_c=v+\frac{1-2Zw}{L_\delta},\qquad E_c=Zw-\frac12,\\
 C_f=f+R\partial_Rf,\qquad q_n=c_n+2=1-\delta+e_n.
\end{gathered}
\end{equation}
Only the pairs $(0,n)$ and $(n,0)$ contain current unknowns.
The sources \eqref{lo:separated-sources} use completed lower orders.
In particular, the regular form of the radial source is
\begin{equation}\label{lo-core:regular-radial-source}
\begin{aligned}
 \frac{\Omega_{(k)}}R
 ={}&\mathscr T_{e_k-2}v_{(k)}
 +\sum_{i+j=k}\left[
 v_{(i)}\left(\frac12v_{(j)}+R\partial_Rv_{(j)}\right)
       +U^z_{(i)}\mathscr Z_{e_j-2}v_{(j)}\right]\\
 &-4\partial_Rv_{(k)}-2R\partial_R^2v_{(k)}
       -\mathscr Z_{e_{k-1}-2}^{[2]}v_{(k-1)}.
\end{aligned}
\end{equation}
The final term is absent at $k=0$.
Each fixed-order source uses finitely many derivatives of completed
lower profiles. Their parameter neighborhoods have a common smaller
neighborhood, with bounds uniform on the full radial interval.
The apparent angular division is also regular, since
$V_{(i)}(\partial_RF_{(j)}+F_{(j)}/R)
=v_{(i)}(R\partial_RF_{(j)}+F_{(j)})$.

\medskip
\noindent\textbf{Step 2: eliminate the average and display the system.}
Write $a=\sqrt{R_{\rm in}}$ and set
\[
 K_n=\mathcal A_R U^z_{(n)}-U^z_{(n)},\qquad
 \mathbf Y_n=(F_{(n)},U^z_{(n)},K_n,P_{(n)},
          \partial_\xi F_{(n)},\partial_\xi U^z_{(n)})^{\mathsf T}.
\]
Every entry is evaluated at $(R,Z)=(\xi^2,Z)$.
The last two entries are unknown radial derivatives, not extra data.
Incompressibility becomes
\begin{equation}\label{lo-core:radial-elimination}
 v_{(n)}=
 \frac{-c_nZU^z_{(n)}-q_nZK_n
             -L\partial_ZU^z_{(n)}-L\partial_ZK_n}{L_\delta},
\end{equation}
where $2-q_n=-c_n$. Subtracting the local velocity from its radial
average leaves $2K_n/\xi$ as the only singular term in its equation,
\begin{equation}\label{lo-core:average-equation}
 \partial_\xi K_n+\frac2\xi K_n=-\partial_\xi U^z_{(n)}.
\end{equation}
Conversely, for a regular solution of this equation, put
$A=U^z_{(n)}+K_n$. Then $\partial_R(RA)=U^z_{(n)}$.
Since $A$ is bounded, $RA$ vanishes at the axis, and integration
gives $A=\mathcal A_RU^z_{(n)}$. Thus no extra integration constant
has been introduced.

The radial second-order operators become
\[
 2(R\partial_R^2+2\partial_R)
 =\tfrac12(\partial_\xi^2+3\xi^{-1}\partial_\xi),\qquad
 2(R\partial_R^2+\partial_R)
 =\tfrac12(\partial_\xi^2+\xi^{-1}\partial_\xi).
\]
Substitute $\partial_RP_{(n)}=2fF_{(n)}+N^p_{(n)}$
in the axial equation. The three coefficient equations are
equivalent to
\begin{equation}\label{lo-core:first-order-system}
 \partial_\xi\mathbf Y_n+\xi^{-1}\mathsf D\mathbf Y_n
 =A_{0,n}\mathbf Y_n+A_1\partial_Z\mathbf Y_n+\mathbf f_n,
 \qquad
 \mathbf Y_n(0,Z)=0,
 \quad \mathsf D=\operatorname{diag}(0,0,2,0,3,1).
\end{equation}
The zero-order matrix and source are
\begin{equation}\label{lo-core:zero-order-matrix}
 A_{0,n}=
 \begin{pmatrix}
 0&0&0&0&1&0\\
 0&0&0&0&0&1\\
 0&0&0&0&0&-1\\
 4\xi f&0&0&0&0&0\\
 a_{51}&a_{52}&a_{53}&0&\xi B_c&0\\
 a_{61}&a_{62}&a_{63}&a_{64}&0&\xi B_c
 \end{pmatrix},\qquad
 \mathbf f_n=
 \begin{pmatrix}
 0\\0\\0\\2\xi N^p_{(n)}\\
 2N^\theta_{(n)}\\
 2N^z_{(n)}-4ZRN^p_{(n)}/L_\delta
 \end{pmatrix},
\end{equation}
where
\begin{equation}\label{lo-core:zero-order-entries}
\begin{aligned}
 a_{51}&=2\left(v+\frac{b_nE_c}{L_\delta}\right),&
 a_{52}&=2\mathscr Z_{b_0}f-\frac{2c_nZC_f}{L_\delta},\\
 a_{53}&=-\frac{2q_nZC_f}{L_\delta},&
 a_{61}&=-\frac{8ZRf}{L_\delta},\\
 a_{62}&=2\left(\frac{c_nE_c}{L_\delta}+\mathscr Z_{c_0}w
                         -\frac{c_nZR\partial_Rw}{L_\delta}\right),\\
 a_{63}&=-\frac{2q_nZR\partial_Rw}{L_\delta},&
 a_{64}&=\frac{2\pi_n Z}{L_\delta}.
\end{aligned}
\end{equation}
For example, the pressure substitution contributes
$-8ZRfF_{(n)}/L_\delta-4ZRN^p_{(n)}/L_\delta$
to the last row. The first three rows give the derivative
coordinates and \eqref{lo-core:average-equation}.
The only possibly nonzero entries of the derivative matrix are
\begin{equation}\label{lo-core:derivative-blocks}
\begin{aligned}
 (A_1)_{51}&=\frac{2H_c}{L_\delta},&
 (A_1)_{52}=(A_1)_{53}&=-\frac{2LC_f}{L_\delta},\\
 (A_1)_{62}&=\frac{2(H_c-LR\partial_Rw)}{L_\delta},&
 (A_1)_{63}&=-\frac{2LR\partial_Rw}{L_\delta},&
 (A_1)_{64}&=\frac{2L}{L_\delta}.
\end{aligned}
\end{equation}
This matrix is independent of $n$.

\medskip
\noindent\textbf{Step 3: invert the radial operator and count derivatives.}
Split the unknown into the two blocks
\[
 X_n=(F_{(n)},U^z_{(n)},K_n,P_{(n)})^{\mathsf T},\qquad
 W_n=(\partial_\xi F_{(n)},\partial_\xi U^z_{(n)})^{\mathsf T}.
\]
The derivative term $A_1\partial_Z\mathbf Y_n$ reads only $X_n$
and contributes only to the $W_n$ equations. The radial inverse below
preserves these blocks. Thus a second derivative factor cannot act
until a zero-order factor has returned some component to the $X_n$ block.
We make this observation precise in the integral formulation.
Both $A_1$ and $\partial_Z A_1$ have nonzero entries only in
the last two rows and first four columns by \eqref{lo-core:derivative-blocks}.
Thus for every diagonal $D_0$ and any radii $\xi,s$,
\begin{equation}\label{lo-core:nilpotence}
 A_1(\xi)D_0A_1(s)=0,\qquad
 A_1(\xi)D_0\partial_ZA_1(s)=0.
\end{equation}

For $d_i=\mathsf D_{ii}$, the regular zero-axis inverse is
\begin{equation}\label{lo-core:regular-inverse}
 (\mathscr G\mathbf g)_i(\xi,Z)
 =\int_0^\xi\left(\frac{s}{\xi}\right)^{d_i}g_i(s,Z)\,ds
 =\xi\int_0^1t^{d_i}g_i(t\xi,Z)\,dt.
\end{equation}
For $d_i>0$, regularity excludes the homogeneous solution
$\xi^{-d_i}$; for $d_i=0$, the zero axis value excludes the constant.
The integral kernel is diagonal, independent of $Z$, and bounded
by one for $0\le s\le\xi$.
Put
\[
 \mathscr K_{0,n}=\mathscr G A_{0,n},\qquad
 \mathscr K_1=\mathscr G A_1\partial_Z,\qquad
 \mathscr K^{(n)}=\mathscr K_{0,n}+\mathscr K_1.
\]
Then \eqref{lo-core:first-order-system} is equivalent to
\begin{equation}\label{lo-core:integral-equation}
 \mathbf Y_n=\mathscr K^{(n)}\mathbf Y_n+\mathscr G\mathbf f_n.
\end{equation}
Both terms in the double integral for $\mathscr K_1^2\mathbf Y$
vanish by \eqref{lo-core:nilpotence}, including the term in which
the parameter derivative acts on a coefficient.
Thus $\mathscr K_1^2=0$.
A word of length $k$ in $\mathscr K_{0,n}$ and $\mathscr K_1$
can be nonzero only if no two $\mathscr K_1$ factors are adjacent.
If it contains $p$ derivative factors, then
\begin{equation}\label{lo-core:word-count}
 p\le\nu_k:=\lceil k/2\rceil.
\end{equation}
There are at most $2^k$ words.
Separated derivative factors need not vanish, since $A_{0,n}$
may map the $W_n$ block back to the $X_n$ block. In particular,
$\mathscr K_1\mathscr K_{0,n}\mathscr K_1$ need not be zero,
so the iteration need not terminate after finitely many terms.

\medskip
\noindent\textbf{Step 4: distribute the Cauchy loss and sum on the full interval.}
Choose a bounded tube
$S_\rho=\{Z\in\mathbb C:\operatorname{dist}(Z,[-1,1])<\rho\}$
whose closure lies in a common coefficient and source neighborhood.
Let $C_{\mathrm{coef},n}\ge1$ bound $A_{0,n},A_1,\mathbf f_n$
on $[0,a]\times S_\rho$, using the maximum vector norm and its
induced matrix norm.
Fix $0<\rho'<\rho$ and put $\Delta=\rho-\rho'$.
For a word with $p\ge1$ derivative factors, apply Cauchy's estimate
with $p$ equal losses $\Delta/p$, starting at the rightmost factor.
Each derivative costs at most $p/\Delta$; multiplication costs
at most $C_{\mathrm{coef},n}$ and needs no loss of radius.
Cauchy's estimate applies to the entire expression to the right of
each derivative, so it includes all Leibniz terms.
For $p=0$ there is no Cauchy loss.

The $k$ matrix factors and the final source inverse give $k+1$
ordered radial integrations, whose simplex has volume
$\xi^{k+1}/(k+1)!$.
All diagonal kernels have norm at most one.
Summing the words and using \eqref{lo-core:word-count}, with
$C_n=2C_{\mathrm{coef},n}$, gives
\begin{equation}\label{lo-core:picard-bound}
 \left\|\bigl(\mathscr K^{(n)}\bigr)^k\mathscr G\mathbf f_n
       \right\|_{[0,a]\times S_{\rho'}}
 \le\frac{C_n^{k+1}a^{k+1}}{(k+1)!}
       \bigl(\max\{1,\nu_k/\Delta\}\bigr)^{\nu_k}.
\end{equation}
At $k=0$ the last factor is one.
The factorial wins over this derivative loss for every finite $aC_n$.
Indeed, for $k\ge1$, $\nu_k\le(k+1)/2$ and
$(k+1)!\ge((k+1)/e)^{k+1}$.
With $C_\Delta=\max\{1,(2\Delta)^{-1}\}$, the right side is at most
\[
 \left(\frac{eC_na\sqrt{C_\Delta}}{\sqrt{k+1}}\right)^{k+1},
\]
whose $k$th root tends to zero. Hence
\begin{equation}\label{lo-core:picard-series}
 \mathbf Y_n=\sum_{k=0}^\infty
                 \bigl(\mathscr K^{(n)}\bigr)^k\mathscr G\mathbf f_n
\end{equation}
converges uniformly on $[0,a]\times S_{\rho'}$.
The sum is continuous in $\xi$, holomorphic in $Z$, and zero at
the axis. To pass to \eqref{lo-core:integral-equation}, first use
the same convergence on an intermediate tube
$S_{\rho''}$, $\rho'<\rho''<\rho$; Cauchy's estimate then gives
uniform convergence of the parameter derivatives on $S_{\rho'}$.
Multiplication and radial integration preserve this convergence.

The same argument proves uniqueness. The difference $\mathbf Y$
of two regular solutions satisfies
$\mathbf Y=(\mathscr K^{(n)})^k\mathbf Y$ for every $k$.
On a common smaller parameter tube, the word estimate now has
radial volume $a^k/k!$ and gives
\[
 \|\mathbf Y\|_{[0,a]\times S_{\rho'}}
 \le\frac{(C_na)^k}{k!}
       \bigl(\max\{1,\nu_k/\Delta\}\bigr)^{\nu_k}
       \|\mathbf Y\|_{[0,a]\times S_\rho}.
\]
Letting $k\to\infty$ gives $\mathbf Y=0$; the identity theorem
extends equality to the common connected parameter neighborhood.
The tube and $C_n$ may depend on $n$. Convergence for every finite
$aC_n$ leaves the radial endpoint unchanged.

\medskip
\noindent\textbf{Step 5: recover radial derivatives and the axis parity.}
For continuous $\mathbf g$, the regular inverse satisfies
\begin{equation}\label{lo-core:inverse-derivative}
 \partial_\xi(\mathscr G\mathbf g)_i(\xi,Z)
 =g_i(\xi,Z)-d_i\int_0^1t^{d_i}g_i(t\xi,Z)\,dt.
\end{equation}
For $\xi>0$ this follows from the first expression in
\eqref{lo-core:regular-inverse}; at zero both sides equal
$g_i(0,Z)/(d_i+1)$.
Take $\mathbf g=A_{0,n}\mathbf Y_n+A_1\partial_Z\mathbf Y_n+\mathbf f_n$.
On a slightly smaller tube this is continuous, so the formula
gives one radial derivative.
Suppose derivatives through order $m$ have been obtained, with
holomorphic parameter extensions and uniform bounds.
Shrink the tube once to estimate their parameter derivatives.
The coefficients and sources have the required fixed radial
derivatives, so $\mathbf g$ is now $C^m$ in $\xi$.
Differentiating \eqref{lo-core:inverse-derivative} $m$ times gives
$\mathbf Y_n\in C^{m+1}$ up to the axis.
Only finitely many neighborhood reductions are needed for each
fixed derivative count.

Extend the coefficients to $-a\le\xi\le a$ by their displayed
formulas in $R=\xi^2$.
The nonsingular second formula in \eqref{lo-core:regular-inverse}
and the same estimates with $|\xi|$ give existence and uniqueness
on the two-sided interval.
Put $S=\operatorname{diag}(1,1,1,1,-1,-1)$.
The explicit matrices and source satisfy
\[
 A_{0,n}(-\xi)=-SA_{0,n}(\xi)S,\qquad
 A_1(-\xi)=-SA_1(\xi)S,\qquad
 \mathbf f_n(-\xi)=-S\mathbf f_n(\xi).
\]
Also $S$ commutes with $\mathsf D$.
Thus $S\mathbf Y_n(-\xi,Z)$ solves the same equation with the
same zero axis values. Uniqueness proves that the first four
coordinates are even and the last two are odd.

An even smooth function of $\xi$ is a smooth function of $R=\xi^2$
on $R\ge0$. Here this statement preserves the parameter estimates:
for an even function $g$,
\begin{equation}\label{lo-core:even-radial-derivative}
 \frac{\partial_\xi g(\xi,Z)}{2\xi}
 =\frac12\int_0^1\partial_\xi^2g(t\xi,Z)\,dt.
\end{equation}
The right side is smooth and even at zero.
Iteration gives all $R$ derivatives; each fixed derivative uses
only finitely many $\xi$ derivatives and preserves holomorphy on
the corresponding smaller tube.
Apply this to $F_{(n)},U^z_{(n)},K_n,P_{(n)}$.

\medskip
\noindent\textbf{Step 6: recover the radial field and check the axis.}
The third row recovers the average, as proved in Step 2.
Equation \eqref{lo-core:radial-elimination} therefore gives exactly
the field in \eqref{lo:incompressibility}.
It has the form $V_{(n)}=Rv_{(n)}$, with $v_{(n)}$ smooth and with
the same parameter regularity.
All eliminations were reversible, so the original coefficient
equations hold.
The zero axis values and their parameter derivatives give
$K_n(0,Z)=v_{(n)}(0,Z)=0$; hence $V_{(n)}=O(R^2)$.
Evaluation of the coefficient equations at the axis gives
\begin{equation}\label{lo-core:axis-jets}
 \partial_RF_{(n)}(0,Z)=\tfrac14N^\theta_{(n)}(0,Z),\qquad
 \partial_RU^z_{(n)}(0,Z)=\tfrac12N^z_{(n)}(0,Z),\qquad
 \partial_RP_{(n)}(0,Z)=N^p_{(n)}(0,Z).
\end{equation}
Finally, \eqref{lo-core:regular-radial-source} with $k=n$ shows
that $\Omega_{(n)}/R$ is regular for the next order.
\end{proof}

Each coefficient is used as a source only after its radial extension
and moment correction. Subsection~\ref{lo-ext:moments} carries out this
extension while preserving the inner parameter regularity required by the lemma.

\subsection{Five-moment correction and global extension at each order}
\label{lo-ext:moments}
\label{lo-ext:fixed-intervals}

With $R_{\rm in}$ fixed by \eqref{lo-core:fixed-radii}, choose any radius
\begin{equation}\label{lo-ext:retained-radius}
 R_a<R_{\rm keep}<R_{\rm in}.
\end{equation}
There is no further restriction on this choice; for example, one may take
$R_{\rm keep}=(R_a+R_{\rm in})/2$. Fix it once for all coefficient orders.
We retain the inner solution up to this radius and use the remaining
interval before $R_{\rm in}$ for the cutoff transition. Choose a smooth
radial cutoff such that
\[
 0\le\chi_{\rm rad}\le1,\qquad
 \chi_{\rm rad}=1\text{ on }[0,R_{\rm keep}],\qquad
 \chi_{\rm rad}=0\text{ for }R\ge\tfrac12(R_{\rm keep}+R_{\rm in}).
\]
In the reserved interval $I_2$ of \eqref{eq:road-reserved-patches},
choose five nonnegative, nonzero bumps
$b^z_1,b^z_2,b^\theta_1,b^\theta_2,b^\theta_3\in C_c^\infty(I_2)$
with disjoint supports, ordered increasingly within each family.
Use this cutoff and these bumps at every positive order, as in
Figure~\ref{fig:lo-cutoff-moments}.

Fix $n\ge1$ with all lower orders completed. In each product of total
order $n$, a new coefficient can occur only in the pairs $(0,n)$ and
$(n,0)$, where it multiplies the fixed leading profile. The five moment
conditions are therefore affine in the order-$n$ bump coefficients,
even when the incoming defects are large. From the inner solution define
\[
 \widetilde F_{(n)}=\chi_{\rm rad}F_{(n)}^{\rm in},\qquad
 \widetilde U^z_{(n)}=\chi_{\rm rad}U^{z,\rm in}_{(n)},\qquad
 \widetilde U^\theta_{(n)}=\sqrt{2R}\widetilde F_{(n)},
\]
extending the products by zero beyond $R_{\rm in}$.
On $[0,R_{\rm keep}]$ the cutoff is one and all moment bumps vanish,
so the completed profiles still solve the inner coefficient equations.
Every positive-order tangential residual therefore vanishes there,
and integration from the axis gives
\[
 \mathcal T_{(n)}=0\qquad(0\le R\le R_{\rm keep},\ n\ge1).
\]
The coefficient equations successively cancel the lower-order errors.
By taking sufficiently many orders, the remaining error
$\mathsf{NS}(u,p)+\mathcal D\mathsf T_{(0)}$ can thus be made $O(\lambda^N)$ for any
prescribed $N$ on this interval; see Lemma~\ref{lo-ext:finite-residual}.
The smooth summation in Subsection~\ref{lo-sum:subsection} realizes this
cancellation as a remainder flat in $\lambda$ with all fixed physical derivatives, while retaining
$\mathcal T_B=\mathcal T_{(0)}$ on $[0,R_{\rm keep}]$.
The strict inequality $R_{\rm keep}>R_a$ leaves a common inner collar
on which the stress is unchanged, as needed in Subsection~\ref{lo-cone:subsection}.

On $I_2$ the leading velocities still satisfy
\begin{equation}\label{lo-ext:reserved-power}
 U^z_{(0)}=0,\qquad U^\theta_{(0)}=A_*(Z)(2R)^{-1/2-\mu},\qquad
 A_*(Z)=\frac{a_*}{1+Z^2}>0.
\end{equation}
No vanishing of the leading radial velocity on $I_2$ is assumed or needed.
For coefficient functions $\alpha_{n,j}(Z)$ and $\beta_{n,j}(Z)$ set
\begin{equation}\label{lo-ext:bump-ansatz}
\begin{aligned}
 U^z_{(n)}&=\widetilde U^z_{(n)}+\sum_{j=1}^2\alpha_{n,j}b^z_j(R),\\
 U^\theta_{(n)}&=\widetilde U^\theta_{(n)}+\sum_{j=1}^3\beta_{n,j}b^\theta_j(R),
 \qquad F_{(n)}=U^\theta_{(n)}/\sqrt{2R}.
\end{aligned}
\end{equation}
The bumps vanish near the axis, so the last quotient retains its smooth
extension there.

For this tentative extension let
$G_n=2F_{(0)}F_{(n)}+N^p_{(n)}$ be the pressure derivative prescribed
by \eqref{lo:pressure-source}. The five moments at order $n$ are
\begin{equation}\label{lo-ext:five-moments}
\begin{aligned}
 m_{n,1}&=\int_0^\infty U^z_{(n)}\,dR,&
 m_{n,2}&=\int_0^\infty\sqrt{2R}\,U^\theta_{(n)}\,dR,\\
 m_{n,3}&=\int_0^\infty G_n\,dR,&
 m_{n,4}&=\int_0^\infty\sqrt{2R}
                   \sum_{i+j=n}U^z_{(i)}U^\theta_{(j)}\,dR,\\
 m_{n,5}&=\int_0^\infty
       \left(\sum_{i+j=n}U^z_{(i)}U^z_{(j)}-RG_n\right)dR.
\end{aligned}
\end{equation}
In the notation of \eqref{fiveM}, these moments correspond, in order, to
$M^z,M^\theta,M^p,M^{\theta z},M^{z\theta}$, with the known radial source
$\Omega_{(n-1)}$ contributing to the third and fifth entries through $G_n$.
All indices in these sums are nonnegative. We set $m_{n,1}=m_{n,3}=0$
to remove the exterior tails of the radial velocity and pressure,
respectively. In particular, once $m_{n,3}=0$, the pressure primitive
$P_{(n)}(R,Z)=\int_0^R G_n(x,Z)\,dx$ vanishes beyond the support of $G_n$,
so integration by parts gives $\int_0^\infty P_{(n)}\,dR
=-\int_0^\infty R G_n\,dR$. Thus $m_{n,5}$ is the integral of the
order-$n$ axial momentum flux $\sum_{i+j=n} U^z_{(i)}U^z_{(j)}+P_{(n)}$.
In \eqref{lo-ext:integrated-profile-residuals} below, the vanishing of
$m_{n,2},m_{n,4}$ removes the angular stress tail and the vanishing of
$m_{n,1},m_{n,5}$ removes the axial stress tail, using the preceding
order's moment identities in both cases.

For example, the first mixed moment condition reads
\[
 m_{1,4}=\int_0^\infty\sqrt{2R}\,
 \bigl(U^\theta_{(0)}U^z_{(1)}+U^z_{(0)}U^\theta_{(1)}\bigr)\,dR=0,
\]
with both first-order velocities given by \eqref{lo-ext:bump-ansatz}.
For the two moments involving pressure, we have
\begin{equation}\label{lo-ext:affine-pressure-moments}
\begin{aligned}
 m_{n,3}&=2\int_0^\infty F_{(0)}F_{(n)}\,dR
                              +\int_0^\infty N^p_{(n)}\,dR,\\
 m_{n,5}&=\int_0^\infty
       (2U^z_{(0)}U^z_{(n)}-U^\theta_{(0)}U^\theta_{(n)})\,dR\\
 &\quad+\int_0^\infty\left[
       \sum_{i=1}^{n-1}U^z_{(i)}U^z_{(n-i)}
       -\frac12\sum_{i=1}^{n-1}U^\theta_{(i)}U^\theta_{(n-i)}
       +\frac12\Omega_{(n-1)}\right]dR.
\end{aligned}
\end{equation}
Empty sums are zero. These formulas retain every lower-order cross term
and the known radial source. At $n=1$ the second integral in the last
formula is $\frac12\int\Omega_{(0)}\,dR$; it is not discarded.
The completed-field regularity gives $\Omega_{(n-1)}=R\omega_{n-1}$
near the axis with smooth $\omega_{n-1}$, and its exterior support follows
from the vanishing of the completed radial velocities.
Thus all five moments are well defined and smooth in $Z$.

On $I_2$, \eqref{lo-ext:reserved-power} makes the linear part block diagonal:
the two axial bumps change only $m_{n,1},m_{n,4}$, while the three angular
bumps change only $m_{n,2},m_{n,3},m_{n,5}$. Indeed, for variations
supported in $I_2$, with all lower profiles fixed,
\begin{equation}\label{lo-ext:moment-variations}
\begin{aligned}
 \dot m_{n,1}&=\int\dot U^z\,dR,&
 \dot m_{n,4}&=A_*\int(2R)^{-\mu}\dot U^z\,dR,\\
 \dot m_{n,2}&=\int(2R)^{1/2}\dot U^\theta\,dR,&
 \dot m_{n,3}&=2A_*\int(2R)^{-3/2-\mu}\dot U^\theta\,dR,\\
 \dot m_{n,5}&=-A_*\int(2R)^{-1/2-\mu}\dot U^\theta\,dR.
\end{aligned}
\end{equation}
The dot denotes variation of the current profiles with the sources fixed.

Define the fixed matrices
\begin{equation}\label{lo-ext:moment-matrices}
\begin{aligned}
 (\mathsf B^z)_{ij}&=\int_0^\infty(2R)^{a_i}b^z_j(R)\,dR,
       &(a_1,a_2)&=(0,-\mu),\\
 (\mathsf B^\theta)_{ij}&=\int_0^\infty(2R)^{s_i}b^\theta_j(R)\,dR,
       &(s_1,s_2,s_3)&=(1/2,-3/2-\mu,-1/2-\mu).
\end{aligned}
\end{equation}
The exponents in each block are distinct because $\mu>0$.
Lemma~\ref{lem:bump-moment-matrix}, after a row permutation and positive
row rescaling, shows that both matrices are invertible. They are
independent of $Z$ and of the coefficient order.

Let $m^0_{n,j}$ be the moments with all five bump coefficients zero,
including all known source terms, and put
\[
 d^z_n=\begin{pmatrix}m^0_{n,1}\\m^0_{n,4}/A_*\end{pmatrix},\qquad
 d^\theta_n=\begin{pmatrix}
 m^0_{n,2}\\m^0_{n,3}/(2A_*)\\-m^0_{n,5}/A_*
 \end{pmatrix}.
\]
The exact solution is
\begin{equation}\label{lo-ext:moment-solve}
 \alpha_n=-(\mathsf B^z)^{-1}d^z_n,\qquad
 \beta_n=-(\mathsf B^\theta)^{-1}d^\theta_n,
 \qquad m_{n,j}=0\quad(1\le j\le5).
\end{equation}
Since $A_*$ is bounded away from zero, this produces smooth signed
corrections for discrepancies of any finite size. No smallness of the
moment defects or of the bump coefficients is needed. The two solves
preserve each other and every earlier coefficient. Partial cumulative
moments within the bumps need not be zero; their actual values will be
used in the radial and stress primitives.

The supports and integration intervals are independent of $Z$.
Differentiating under the integrals and using the two constant inverses gives
\begin{equation}\label{lo-ext:repair-parameter-bound}
 \|\alpha_n\|_{C_Z^\ell}+\|\beta_n\|_{C_Z^\ell}
 \le C_\ell\sum_{j=1}^5\|m^0_{n,j}\|_{C_Z^\ell}\qquad(\ell\ge0).
\end{equation}
The constants $C_\ell$ depend only on the fixed matrices and
$\|A_*^{-1}\|_{C_Z^\ell}$, not on $n$. The discrepancies themselves
may grow with $n$.

\paragraph{Pressure and stress reconstruction.}\label{lo-ext:induction}
With the five bump coefficients fixed by \eqref{lo-ext:moment-solve},
we integrate the divergence and radial momentum equations from the axis to obtain
\begin{equation}\label{lo-ext:reconstruction}
\begin{aligned}
 M^z_{(n)}(R,Z)&=\int_0^R U^z_{(n)}(x,Z)\,dx,\\
 V_{(n)}(R,Z)&=-\int_0^R(\mathscr Z_{c_n}U^z_{(n)})(x,Z)\,dx,\\
 P_{(n)}(R,Z)&=\int_0^R
       \bigl(2F_{(0)}F_{(n)}+N^p_{(n)}\bigr)(x,Z)\,dx.
\end{aligned}
\end{equation}
In particular, the first pressure correction is
\begin{equation}\label{lo-ext:first-pressure}
 P_{(1)}(R,Z)=\int_0^R
 \bigl[2U^\theta_{(0)}U^\theta_{(1)}-\Omega_{(0)}\bigr](x,Z)
                      \,\frac{dx}{2x}.
\end{equation}
Here $U^\theta_{(1)}$ already includes its three moment corrections;
$\Omega_{(0)}$ is the known radial source of the completed leading field.
There is no further pressure constant: the positive-order axis value
is zero, as prescribed in \eqref{lo:axis-values}.
The regular quotient $\Omega_{(n-1)}/R$ makes these integrals smooth at
the axis. The divergence and pressure equations remain exact globally.
On $[0,R_{\rm keep}]$ the cutoff is one and the bumps vanish, so the
zero axis constants recover the solved inner $V_{(n)},P_{(n)}$ as well.

For $R\ge\sup I_2$, the corrected angular and axial velocities
vanish, and $M^z_{(n)}=m_{n,1}=0$ as a function of $Z$.
Equation \eqref{lo:incompressibility} then gives $V_{(n)}=0$ there too.
Since $\sup I_2<\inf I_3$, the velocity on $I_3$ is preserved.
The pressure can extend past $I_3$, since the first-order pressure
source still contains $\Omega_{(0)}$ beyond the bumps.

For $R\ge R_{0z}$, the leading axial velocity and its moment $M^z_{(0)}$
both vanish by the axial correction of Interval O.4. Hence $V_{(0)}=0$
by \eqref{lo:incompressibility}. All positive-order velocities already
vanish there. Thus every term in $\Omega_{(n-1)}$ vanishes, and each
product in $\sum_{i+j=n}F_{(i)}F_{(j)}$ has a vanishing positive-order
factor. This gives $G_n=\partial_RP_{(n)}=0$ there. The third moment
fixes the constant value of the pressure,
\begin{equation}\label{lo-ext:pressure-tail}
 P_{(n)}(R,Z)=\int_0^\infty G_n(x,Z)\,dx=m_{n,3}(Z)=0
 \qquad(R\ge R_{0z}).
\end{equation}
At the axis, the formulas
\[
 \frac{M^z_{(n)}}R=\int_0^1 U^z_{(n)}(tR,Z)\,dt,\qquad
 \frac{V_{(n)}}R=-\int_0^1(\mathscr Z_{c_n}U^z_{(n)})(tR,Z)\,dt
\]
prove the smoothness of both quotients. On the fixed inner interval
the cutoff has no $Z$ dependence and the bumps are absent. The product
rule and forward integration preserve the local holomorphic parameter
regularity of Lemma~\ref{lo-core:fixed-radius}; no value of a bump
coefficient from $I_2$ affects these inner integrals. This supplies
the regular source needed for the next coefficient.

We record the completed extension, including the stress estimates needed
for summation. Use the weight and edge distance from
\eqref{intro:edge-weight}, with the constants fixed by the endpoint collars:
\begin{equation}\label{lo-ext:flat-weight}
\begin{gathered}
 y_a=\log(R/R_a),\qquad y_h=\log(R_h/R),\qquad
 d_{\rm edge}=\min\{1,y_a,y_h\},\\
 c_a=h_b^2,\qquad c_h=4,\qquad
 \omega(R)=\exp\left(-\frac{c_a}{y_a^2}-\frac4{y_h^2}\right),
\end{gathered}
\end{equation}
extended by zero outside the annulus. Here $h_b$ is the logarithmic
width of the shear transition in \eqref{eq:inner-chi}.

\begin{proposition}[Extension with five vanishing moments]
\label{lo-ext:extension}
For every $n\ge1$, the construction above extends the order-$n$ inner
coefficient solution from Lemma~\ref{lo-core:fixed-radius} to $R\ge0$, $|Z|\le1$,
retaining the solved profiles on $[0,R_{\rm keep}]$ and satisfying the
divergence and pressure equations globally. The regular variables
$F_{(n)}$, $U^z_{(n)}$, $V_{(n)}/R$, $P_{(n)}$ and $M^z_{(n)}/R$ are smooth in $(R,Z)$
up to $R=0$. In particular, the transverse vector field
$U^\theta_{(n)}e_\theta+U^r_{(n)}e_r$ is smooth in Cartesian profile
coordinates at the axis, since $U^\theta_{(n)}=\sqrt{2R}\,F_{(n)}$ and
$U^r_{(n)}=\sqrt{2R}\,(V_{(n)}/(2R))$.
The radii and correction supports are independent of $n$.
For every $n\ge1$,
\begin{equation}\label{lo-ext:coefficient-support}
\begin{gathered}
 m_{n,j}=0\quad(1\le j\le5),\qquad
 U^\theta_{(n)}=U^z_{(n)}=M^z_{(n)}=V_{(n)}=0
                     \quad(R\ge\sup I_2),\\
 \operatorname{supp}_R P_{(n)}\subset[0,R_{0z}],
\end{gathered}
\end{equation}
and, for each $m\ge0$,
\begin{equation}\label{lo-ext:coefficient-bounds}
 \max_{k+\ell\le m}\sup_{R\ge0,\,|Z|\le1}
 \left|\partial_R^k\partial_Z^\ell
 (F_{(n)},U^z_{(n)},V_{(n)}/R,P_{(n)},M^z_{(n)}/R)\right|
 \le C_{n,m}<\infty.
\end{equation}
There are normalized stresses cancelling the order-$n$ tangential
residuals, with
\begin{equation}\label{lo-ext:stress-support}
 \operatorname{supp}_R \mathcal{T}_{(1)}\subset[R_{\rm keep},R_h],\qquad
 \operatorname{supp}_R \mathcal{T}_{(n)}\subset[R_{\rm keep},R_{0z}]
 \quad(n\ge2).
\end{equation}
For every multi-index $I$,
\begin{equation}\label{lo-ext:weighted-stress}
 |\partial_{R,Z}^I \mathcal{T}_{(n)}|
 \le C_{n,I}\omega d_{\rm edge}^{-N_{n,I}}.
\end{equation}
For $n\ge2$ one may take $N_{n,I}=0$.
The constants need not be uniform in the order.
\end{proposition}

\begin{proof}
The velocity, pressure, and parameter-regularity assertions have just
been proved. Smoothness and the common compact support give
\eqref{lo-ext:coefficient-bounds}. It remains to construct the stresses
and remove their exterior tails.

Substitute the corrected profiles into the tangential equations and
remove the common factor $\lambda^{-3-\delta+e_n}$. The signed residuals are
\begin{equation}\label{lo-ext:residual-profiles}
\begin{aligned}
 \mathsf r_{\theta,n}=\sqrt{2R}\bigg[&\mathscr T_{b_n}F_{(n)}
 +\sum_{i+j=n}\left\{V_{(i)}\left(\partial_RF_{(j)}+\frac{F_{(j)}}R\right)
                      +U^z_{(i)}\mathscr Z_{b_j}F_{(j)}\right\}\\
 &-2(R\partial_R^2F_{(n)}+2\partial_RF_{(n)})
                      -\mathscr Z_{b_{n-1}}^{[2]}F_{(n-1)}\bigg],\\
 \mathsf r_{z,n}={}&\mathscr T_{c_n}U^z_{(n)}
 +\sum_{i+j=n}\{V_{(i)}\partial_RU^z_{(j)}
                       +U^z_{(i)}\mathscr Z_{c_j}U^z_{(j)}\}
                       +\mathscr Z_{\pi_n}P_{(n)}\\
 &-2(R\partial_R^2U^z_{(n)}+\partial_RU^z_{(n)})
                       -\mathscr Z_{c_{n-1}}^{[2]}U^z_{(n-1)}.
\end{aligned}
\end{equation}
The inner equations make these residuals zero on $[0,R_{\rm keep}]$.
Set
\begin{equation}\label{lo-ext:stress-primitives}
 \mathcal{T}^\theta_{(n)}=-\frac1{2R}\int_0^R\sqrt{2x}\,\mathsf r_{\theta,n}(x,Z)\,dx,
 \qquad
 \mathcal{T}^z_{(n)}=-\frac1{\sqrt{2R}}\int_0^R\mathsf r_{z,n}(x,Z)\,dx.
\end{equation}
Differentiating gives
\begin{equation}\label{lo-ext:primitive-equations}
 \sqrt{2R}\left(\partial_R+\frac1R\right)\mathcal{T}^\theta_{(n)}=-\mathsf r_{\theta,n},
 \qquad
 \sqrt{2R}\left(\partial_R+\frac1{2R}\right)\mathcal{T}^z_{(n)}=-\mathsf r_{z,n}.
\end{equation}
Thus the physical stress $\mathsf T_{(n)}=\lambda^{-2-\delta+e_n}\mathcal{T}_{(n)}$
cancels the current tangential residual and is zero on the retained
inner interval. To identify the obstruction to compact exterior support,
define the two total residual integrals
\begin{equation}\label{lo-ext:total-residual-integrals}
 J_{\theta,n}(Z)=\int_0^\infty\sqrt{2R}\,\mathsf r_{\theta,n}(R,Z)\,dR,
 \qquad
 J_{z,n}(Z)=\int_0^\infty\mathsf r_{z,n}(R,Z)\,dR.
\end{equation}
Beyond the radial supports of the residuals, \eqref{lo-ext:stress-primitives} gives
\[
 \mathcal T^\theta_{(n)}(R,Z)=-\frac{J_{\theta,n}(Z)}{2R},
 \qquad
 \mathcal T^z_{(n)}(R,Z)=-\frac{J_{z,n}(Z)}{\sqrt{2R}}.
\]
Thus $J_{\theta,n}=J_{z,n}=0$ is precisely the condition that removes
both stress tails. The five incoming moment defects $m^0_{n,j}$ were
removed by the bump correction. We now use the resulting identities
$m_{n,j}=0$, together with the preceding order's moment identities,
to prove that these two residual integrals vanish as well.

For a divergence-free axisymmetric velocity, the conservative identities are
\begin{equation}\label{lo-ext:conservative-identities}
\begin{aligned}
 r^2\mathsf{NS}^\theta
 ={}&\partial_t(r^2u^\theta)+\partial_r(r^2u^ru^\theta)
    +\partial_z(r^2u^zu^\theta)
    -\partial_r(r^2\partial_ru^\theta-ru^\theta)
    -\partial_z^2(r^2u^\theta),\\
 r\mathsf{NS}^z
 ={}&\partial_t(ru^z)+\partial_r(ru^ru^z)
    +\partial_z\bigl(r((u^z)^2+p)\bigr)
    -\partial_r(r\partial_ru^z)-\partial_z^2(ru^z).
\end{aligned}
\end{equation}
At positive order, axis regularity and compact support remove the radial
boundary terms. Since $m_{n,3}=0$, the integrated axial flux is $m_{n,5}$
as observed after \eqref{lo-ext:five-moments}. At $n=1$ we use the
renormalized leading angular moment \eqref{lo-ext:leading-subtracted-moment}
for $m_{n-1,2}$ in the identity below.
The integrated equations express the total residual integrals as
differential combinations of the current and preceding moments:
\begin{equation}\label{lo-ext:integrated-profile-residuals}
\begin{aligned}
 J_{\theta,n}&=\mathscr T_{2-\delta+e_n}m_{n,2}
     +\mathscr Z_{1-2\delta+e_n}m_{n,4}
     -\mathscr Z_{2-\delta+e_{n-1}}^{[2]}m_{n-1,2}=0,\\
 J_{z,n}&=\mathscr T_{1-\delta+e_n}m_{n,1}
     +\mathscr Z_{-2\delta+e_n}m_{n,5}
     -\mathscr Z_{1-\delta+e_{n-1}}^{[2]}m_{n-1,1}=0.
\end{aligned}
\end{equation}
For $n\ge2$ every displayed moment converges and vanishes by the current
repair or the preceding induction. The weights follow, for example, from
$\int r^2u^\theta_{(n)}dr=\lambda^{2-\delta+e_n}m_{n,2}$ and
$\int r\sum_{i+j=n}u^z_{(i)}u^z_{(j)}dr+\int rp_{(n)}dr
=\lambda^{-2\delta+e_n}m_{n,5}$.
These identities hold before taking any physical derivative.

To justify the preceding-order term at $n=1$, recall that
\begin{equation}\label{lo-ext:leading-subtracted-moment}
 m_{0,2}^{\rm ren}=\int_0^\infty\sqrt{2R}
 \left(U^\theta_{(0)}-c_\infty R^{-(1+\delta)/2}\right)dR=0.
\end{equation}
This is the leading terminal moment, preserved in Section~\ref{sec:shear-modification}.
To justify its use, put $\alpha=(1+\delta)/2$. In physical variables it is
\[
 \int_0^\infty
 \left(r^2u^\theta_{(0)}-2^\alpha c_\infty r^{1-\delta}\right)dr=0.
\]
The difference converges: the heat correction to the power law is
$O(R^{-\alpha-1})$, and the subtracted power is locally integrable at
the axis. The subtraction is independent of physical $z$; past a bounded
physical radius the leading velocity equals the same $z$-independent
heat field on each compact neighborhood with $\lambda>0$.
Twice differentiating the convergent difference in $z$ therefore gives
the order-one angular-viscosity cancellation, without assigning a value
to the divergent unsubtracted moment. The leading axial moment is
convergent and equals zero, so it needs no subtraction.

Consequently both total residual integrals vanish at every positive order,
and the stresses can be integrated from infinity:
\begin{equation}\label{lo-ext:backward-stress}
 \mathcal{T}^\theta_{(n)}=\frac1{2R}\int_R^\infty\sqrt{2x}\,\mathsf r_{\theta,n}(x,Z)\,dx,
 \qquad
 \mathcal{T}^z_{(n)}=\frac1{\sqrt{2R}}\int_R^\infty\mathsf r_{z,n}(x,Z)\,dx.
\end{equation}
For $n\ge2$ every source, including the shifted viscosity of the completed
order $n-1$, is zero on $R\ge R_{0z}$. Both stresses vanish there.
For $n=1$ the only remaining source there is the axial viscosity of the
leading swirl, and
\begin{equation}\label{lo-ext:first-stress-exterior}
 \mathcal{T}^z_{(1)}=0,\qquad
 \mathcal{T}^\theta_{(1)}(R,Z)=-\frac1R\int_R^{R_h}
 x\bigl(\mathscr Z_{b_0}^{[2]}F_{(0)}\bigr)(x,Z)\,dx
 \quad(R\ge R_{0z}).
\end{equation}
The source is zero past $R_h$, where the heat swirl is independent of
physical $z$. This proves \eqref{lo-ext:stress-support}.

On the terminal collar the leading physical swirl is
$K(t,r)f_o(y)$, where $y=\log R$ and
$f_o=1-\varepsilon e^{-4/y_h^2}$; see
\eqref{eq:mc-terminal-positive-stress}. Since $K$ is independent of $z$,
\begin{equation}\label{lo-ext:collar-axial-viscosity}
 \partial_z^2(Kf_o)=K\lambda^{-2+2\delta}
 (a_z^2\partial_y^2f_o+b_z\partial_yf_o),\quad
 a_z=-\frac{2Z}{L_\delta},\quad
 b_z=\frac{-(1-\delta)Za_z+L\partial_Za_z}{L_\delta}.
\end{equation}
All fixed $Z$ derivatives of $a_z,b_z$ are bounded. The first two
derivatives of $f_o$ are bounded by $Ce^{-4/y_h^2}y_h^{-3}$ and
$Ce^{-4/y_h^2}y_h^{-6}$, respectively. Backward integration gains
three powers of the distance variable:
\begin{equation}\label{lo-ext:flat-integral}
 \int_0^v e^{-4/s^2}s^{-j}\,ds
 \le C_j e^{-4/v^2}v^{3-j}\qquad(0<v\le1).
\end{equation}
This follows by the substitution $t=4/s^2$ and the elementary
exponential-tail bound. Applying it to \eqref{lo-ext:first-stress-exterior}
and differentiating gives
\begin{equation}\label{lo-ext:first-order-edge}
 |\mathcal{T}_{(1)}|\le Ce^{-4/y_h^2}y_h^{-3},\qquad
 |\partial_{R,Z}^I \mathcal{T}_{(1)}|\le C_Ie^{-4/y_h^2}y_h^{-N_I}.
\end{equation}
On this collar $\omega/e^{-4/y_h^2}$ is bounded above and below by
positive constants. Near the inner edge every positive-order stress
is zero. On the remaining compact subannulus $\omega$ is positive.
These facts prove \eqref{lo-ext:weighted-stress}; for $n\ge2$ the entire
stress support is compactly contained in the annulus, so $N_{n,I}=0$ suffices.
\end{proof}

Starting with the leading field, solve the fixed-radius inner system,
apply the five-moment repair, and reconstruct $V_{(n)},P_{(n)},\mathcal{T}_{(n)}$
at each order before proceeding to the next. This preserves every earlier
coefficient and supplies the regular completed source $\Omega_{(n)}/R$.
It therefore produces the full formal coefficient sequence.

\paragraph{Residual of a finite truncation.}
For a physical tuple define
\begin{equation}\label{lo-ext:finite-tuple}
 \mathcal F(u,p,\mathsf T)=\mathsf{NS}(u,p)+\mathcal D\mathsf T,
 \qquad
 \mathbb U^{[N]}=\left(\sum_{n=0}^Nu_{(n)},\sum_{n=0}^Np_{(n)},
                      \sum_{n=0}^N\mathsf T_{(n)}\right),
\end{equation}
where $\mathcal D$ is \eqref{intro:stress-operator} and the physical
coefficients are \eqref{lo:physical-coefficients}.
No cutoff in $\lambda$ is used in this finite construction.
Write $|G|_m=\max_{|\beta|+b\le m}|\partial_x^\beta\partial_t^bG|$.

\begin{lemma}[Finite-order residual]
\label{lo-ext:finite-residual}
Every finite velocity sum is exactly divergence-free. For every fixed
$R_{\max}<\infty$ and $N,m\ge0$,
\begin{equation}\label{lo-ext:finite-residual-bound}
 |\mathcal F(\mathbb U^{[N]})|_m
 \le C_{N,m}\lambda^{2\delta(N+1)-3-2\delta-2m},
\end{equation}
on $0<\lambda\le1$, $0\le R\le R_{\max}$, $|Z|\le1$.
The constant may depend on $R_{\max}$ and the fixed construction.
The loss of at most two powers of $\lambda$ per derivative is independent of $N$.
\end{lemma}

\begin{proof}
Each coefficient satisfies the exact divergence equation.
At every retained order the pressure equation cancels the radial
residual, and the stress primitives cancel the tangential residuals.
At order zero these are the leading identities.
The remaining terms have total order at least $N+1$.
Here $0\le i,j\le N$ are coefficient indices in the finite truncation:
a linear term acts on index $j$, while a quadratic term combines
indices $i$ and $j$.
Indeed, relative to the tangential base power $\lambda^{-3-\delta}$,
time differentiation and radial viscosity have order $j$, transport
products have order $i+j$, and axial viscosity has order $j+1$.
Relative to the radial base power $\lambda^{-3-2\delta}$, the pressure
has order $j$, the centrifugal products have order $i+j$, radial time
and radial viscosity have order $j+1$, radial transport has order
$i+j+1$, and radial axial viscosity has order $j+2$.
Thus the finite remainder has the form
\begin{equation}\label{lo-ext:exact-finite-remainder}
\begin{aligned}
 \mathcal F(\mathbb U^{[N]})
 ={}&\sum_{k=N+1}^{\max\{2N,N+1\}}
 \lambda^{-3-\delta+e_k}(A^\theta_{k,N}e_\theta+A^z_{k,N}e_z)\\
 &+\sum_{k=N+1}^{\max\{2N+1,N+2\}}
 \lambda^{-3-2\delta+e_k}\frac{A^r_{k,N}}{\sqrt{2R}}e_r.
\end{aligned}
\end{equation}
These are finite sums of the completed coefficient products and
derivatives. Each $A^r_{k,N}$ is divisible by $R$: the centrifugal
products are $2RF_{(i)}F_{(j)}$, and each radial-source term has the
same factor. The apparent radial quotient is therefore smooth in
Cartesian profile coordinates. The angular terms are smooth for the
same reason: $U^\theta_{(j)}e_\theta=F_{(j)}(-X_2,X_1,0)$, where
$X_\perp=x_\perp/\lambda$; also
$(V_{(j)}/\sqrt{2R})e_r=(V_{(j)}/2R)(X_1,X_2,0)$.

The coordinate identities, with $L_\delta\ge1-\delta>0$, give for
any smooth Cartesian profile $g$ on a fixed compact range
\begin{equation}\label{lo-ext:physical-chain-bound}
 \left|\partial_{x_\perp}^\beta\partial_z^k\partial_t^\ell
       [\lambda^b g(X_\perp,Z)]\right|
 \le C_m(1+|b|)^m\|g\|_{C^m}
       \lambda^{b-|\beta|-(1-\delta)k-2\ell},\qquad
 m=|\beta|+k+\ell.
\end{equation}
To see the order-independent loss, one transverse derivative loses
one power of $\lambda$, one axial derivative loses $1-\delta$ powers,
and one time derivative loses two. Repeated differentiation of the
weight contributes at most $C_m(1+|b|)^m$; all coordinate coefficients
and their fixed derivatives are bounded. In particular, $L=0$ at
$Z=\pm1$ introduces no denominator.
Every summand of \eqref{lo-ext:exact-finite-remainder} has
$e_k\ge2\delta(N+1)$. Its profile derivatives are bounded by
\eqref{lo-ext:coefficient-bounds} and the fixed leading bounds.
At most $m$ physical derivatives lose at most $2m$ powers of $\lambda$.
Summing the finitely many terms gives \eqref{lo-ext:finite-residual-bound}.
\end{proof}

We impose no bound on the growth of the coefficient constants $C_{n,m}$
in \eqref{lo-ext:coefficient-bounds} as $n\to\infty$ and do not assert
convergence of the unmodified infinite series.
The smooth summation below absorbs these order-dependent constants
while preserving incompressibility and the exterior support.

\subsection{Smooth summation}
\label{lo-sum:subsection}

We sum the coefficient sequence to obtain smooth physical fields.
The cutoff construction preserves incompressibility. We use the
finite-order estimates and support properties proved above.

An estimate on a compact profile range means a uniform estimate for
$0\le R\le K$, $|Z|\le1$, and $0<\lambda\le1$.
Constants may depend on $K$, the derivative order, and the fixed
coefficient sequence. A residual is \emph{flat in $\lambda$} if its
$|\cdot|_m$ seminorm is bounded by every power of $\lambda$ on each
such range. All endpoint derivatives are understood by one-sided
extension.

The cumulative axial moment provides the physical streamfunction.
Set
\begin{equation}\label{lo-sum:potential}
 \Psi_{(n)}=\lambda^{1-\delta+2n\delta}M^z_{(n)},\qquad
 A_{(n)}=\frac{\Psi_{(n)}}r e_\theta
 =\frac12\lambda^{-1-\delta+2n\delta}
       \frac{M^z_{(n)}}R(-x_2,x_1,0).
\end{equation}
The regularity of $M^z_{(n)}/R$ makes this a smooth Cartesian vector
at the axis. The definition of $M^z_{(n)}$ and incompressibility
\eqref{lo:incompressibility} give
\begin{equation}\label{lo-sum:potential-identities}
 \frac1r\partial_r\Psi_{(n)}=u^z_{(n)},\qquad
 -\partial_z\Psi_{(n)}=r u^r_{(n)},\qquad
 \operatorname{curl}A_{(n)}=u^r_{(n)}e_r+u^z_{(n)}e_z.
\end{equation}
The first quotient extends to the axis by the smoothness of $U^z_{(n)}$.

Choose a fixed smooth function $0\le\chi\le1$, equal to one on
$[0,1/2]$ and zero on $[1,\infty)$. For scales $a_{n+1}\ge2a_n$,
to be chosen below, write $\chi_n=\chi(a_n\lambda)$ and define
\begin{equation}\label{lo-sum:realized-fields}
\begin{aligned}
 u_B&=u_{(0)}+\sum_{n\ge1}
     \left[\operatorname{curl}(\chi_n A_{(n)})
                         +\chi_n u^\theta_{(n)}e_\theta\right],\\
 p_B&=p_{(0)}+\sum_{n\ge1}\chi_n p_{(n)},\qquad
 \mathsf T_B=\mathsf T_{(0)}+\sum_{n\ge1}\chi_n\mathsf T_{(n)},\\
 \mathcal{T}_B&=\lambda^{2+\delta}\mathsf T_B
       =\mathcal{T}_{(0)}+\sum_{n\ge1}\chi_n\lambda^{2n\delta}\mathcal{T}_{(n)}.
\end{aligned}
\end{equation}
The five moment identities hold at each formal order.
In a quadratic term with $i,j\ge1$, the product of cutoffs
$\chi_i\chi_j$ can differ from $\chi_{i+j}$. We therefore estimate
the nonlinear residual by comparison with the uncut finite sums.
The curl in \eqref{lo-sum:realized-fields} preserves incompressibility exactly.
Since $\partial_r\lambda=0$, the angular and axial components have
the direct cutoff-sum form
\begin{equation}\label{lo-sum:tangential-series}
 (u_B^\theta,u_B^z)=\lambda^{-1-\delta}\left[
 (U^\theta_{(0)},U^z_{(0)})+
 \sum_{n\ge1}\chi(a_n\lambda)\lambda^{2n\delta}
                         (U^\theta_{(n)},U^z_{(n)})\right].
\end{equation}
The cutoff adds no axial component. Its radial contribution is explicit:
\begin{equation}\label{lo-sum:cutoff-radial}
 \lambda u^r_{(n),\mathrm{cut}}
 =\lambda^{2n\delta}\left[
    \chi(a_n\lambda)U^r_{(n)}
    -\frac{Z}{L_\delta}(a_n\lambda)\chi'(a_n\lambda)
                         \frac{M^z_{(n)}}{\sqrt{2R}}\right].
\end{equation}
Thus multiplying each velocity component by the same cutoff would
omit a term needed for exact incompressibility.

\begin{proposition}[Realization of the coefficient sequence]
\label{lo-sum:background}
The scales $a_n$ can be chosen so that \eqref{lo-sum:realized-fields}
defines smooth axisymmetric fields for $\lambda>0$, with
\begin{equation}\label{lo-sum:background-equation}
 \nabla\cdot u_B=0,\qquad
 \mathsf{NS}(u_B,p_B)=-\mathcal D\mathsf T_B+\mathsf E_B,
\end{equation}
where $\mathsf E_B$ is flat in $\lambda$ with all Cartesian space--time
derivatives. The fields have the prescribed coefficient expansion.
Recall that $\mathscr D^{\le m}$ denotes compositions of at most $m$
operators from
\[
 \{\lambda\partial_\lambda,\partial_\rho,\partial_Z\},\qquad
 \rho=\sqrt{2R},\qquad \partial_\rho=\sqrt{2R}\,\partial_R,
\]
where $(\lambda,\rho,Z)$ are treated as independent variables.
On every compact profile range one has
\begin{equation}\label{lo-sum:normalized-velocity}
 \mathscr D^{\le m}\left(
   \lambda^{1+\delta}u_B^\theta-U^\theta_{(0)},\,
   \lambda^{1+\delta}u_B^z-U^z_{(0)},\,
   \lambda u_B^r-U^r_{(0)}\right)
       =O_{m,K}(\lambda^{2\delta}).
\end{equation}
For the weight $\omega$ and edge distance $d_{\rm edge}$ in
\eqref{intro:edge-weight}, with $c_a=h_b^2$ and $c_h=4$,
\begin{equation}\label{lo-sum:normalized-stress}
 |\mathscr D^{\le m}(\mathcal{T}_B-\mathcal{T}_{(0)})|
   \le C_m\lambda^{2\delta}\omega d_{\rm edge}^{-q_m}.
\end{equation}
The stress has a smooth zero extension outside $[R_a,R_h]$.
The background equals the leading velocity on $I_3$ and retains
the exact heat velocity and pressure for $R\ge R_h$; in particular,
$\mathsf E_B=0$ there. Its axis pressure satisfies
\begin{equation}\label{lo-sum:axis-pressure}
 \lambda^{2+2\delta}p_B(t,0,z)=P_0(Z).
\end{equation}
All fixed physical derivatives extend to $Z=\pm1$ for $\lambda>0$.
\end{proposition}

\begin{proof}
We give the cutoff estimates used to choose the scales.
For functions of $(\lambda,Z)$, the coordinate identities read
\[
 \partial_z=\frac{\lambda^{-1+\delta}}{L_\delta}
       (Z\lambda\partial_\lambda+L\partial_Z),\qquad
 \partial_t=\frac{\lambda^{-2}}{2L_\delta}
       (-\lambda\partial_\lambda+(1-\delta)Z\partial_Z).
\]
Since $L_\delta\ge1-\delta>0$, induction gives
$|\partial_z^a\partial_t^b\lambda|
 \le C_{a,b}\lambda^{1-a(1-\delta)-2b}$.
The Cartesian profile estimate
\eqref{lo-ext:physical-chain-bound} applies at the axis as well:
each physical derivative loses at most two powers of $\lambda$,
independently of the coefficient index.

With $s=a_n\lambda$, the operator $\lambda\partial_\lambda$ acts on
the cutoff as $s\partial_s$. Thus
\begin{equation}\label{lo-sum:cutoff-derivative}
 |\partial_z^a\partial_t^b\chi_n|
       \le C_{a,b}\lambda^{-a(1-\delta)-2b}.
\end{equation}
For $a+b>0$ these derivatives are supported where $1/2\le s\le1$,
and every $(s\partial_s)^j\chi(s)$ is bounded there independently
of $a_n$. Transverse derivatives vanish.
In particular,
\[
 \operatorname{curl}(\chi_n A_{(n)})
 =\chi_n\operatorname{curl}A_{(n)}
           +a_n\chi'(a_n\lambda)\nabla\lambda\times A_{(n)}.
\]
The coefficient bounds and the product rule therefore bound the
$n$th summand $W_n$ in the physical tuple
$(u_B,p_B,\mathsf T_B)$ by
\begin{equation}\label{lo-sum:summand-bound}
 |W_n|_m\le B_{n,m,K}\lambda^{2n\delta-\ell_m},
\end{equation}
where $B_{n,m,K}$ is independent of $a_n$ and $\ell_m$ depends only
on $m$ and the fixed $\delta$. Increasing $\ell_m$
absorbs the single additional derivative in a curl. The same
argument bounds normalized derivatives of each velocity increment,
including the bracket in \eqref{lo-sum:cutoff-radial}, by
$B_{n,m,K}\lambda^{2n\delta}$.

Choose $a_1>2$ and successively $a_n\ge2a_{n-1}$ for $n\ge2$,
in each case large enough that
\begin{equation}\label{lo-sum:diagonal-choice}
 B_{n,m,K}\lambda^{n\delta}\le2^{-n}
       \quad(0<\lambda\le a_n^{-1})
\end{equation}
for all physical and normalized seminorms of order $m\le n$
on the first $n$ members of a compact exhaustion. Include also,
for $n\ge2$, the corresponding bounds for $\mathcal{T}_{(n)}/\omega$
on the closed stress annulus. These quotients have smooth zero
extensions by their common interior support
\eqref{lo-ext:stress-support}. At each stage there are only
finitely many conditions, each with a positive power $n\delta$.
The constants do not depend on $a_n$, so a sufficiently large
$a_n$ satisfies them all. A single sequence of scales consequently
works for all the seminorms and all compact profile ranges.

A compact subset of the region $\lambda>0$ has $\lambda\ge c>0$,
so all terms with $a_n>c^{-1}$ vanish there. This proves local
finiteness and smoothness. Every curl has zero divergence, and
$\nabla\cdot(\chi_nu^\theta_{(n)}e_\theta)=0$ by axisymmetry.
Thus incompressibility holds exactly, including the cutoff regions.

Let $\mathbb U^{[J]}$ be the uncut physical partial sum through
order $J$, and write $\mathbb U_B=(u_B,p_B,\mathsf T_B)$.
If $\lambda<(2a_J)^{-1}$, the first $J$ cutoffs are one. For
$J$ at least the derivative and exhaustion indices,
\begin{equation}\label{lo-sum:tail}
 |\mathbb U_B-\mathbb U^{[J]}|_m
 \le\sum_{n>J}2^{-n}\lambda^{n\delta-\ell_m}
 \le2^{-J}\lambda^{(J+1)\delta-\ell_m}.
\end{equation}
We compare the nonlinear residual \eqref{lo-ext:finite-tuple} with
that of the finite sum.
For $(e,\pi,\tau)=\mathbb U_B-\mathbb U^{[J]}$ and
$v=u^{[J]}$, the exact difference is
\begin{equation}\label{lo-sum:nse-difference}
\begin{split}
 \mathcal F(\mathbb U_B)-\mathcal F(\mathbb U^{[J]})
 ={}&\partial_te-\Delta e+(v\cdot\nabla)e+(e\cdot\nabla)v
       +(e\cdot\nabla)e+\nabla\pi+\mathcal D\tau.
\end{split}
\end{equation}
Every stress is supported where $R\ge R_a$, so the coefficients
$1/r$ in $\mathcal D$ and their derivatives cost only fixed powers
of $\lambda^{-1}$. The product rule yields, for fixed $m$,
\[
 |\mathcal F(\mathbb U_B)-\mathcal F(\mathbb U^{[J]})|_m
 \le C_m\lambda^{-H_m}|(e,\pi,\tau)|_{m+2}
             (1+|v|_{m+2}+|(e,\pi,\tau)|_{m+2}).
\]
Here $H_m$ is independent of $J$. Likewise
$|v|_{m+2}\le C_{J,m,K}\lambda^{-K'_m}$ with $K'_m$
independent of $J$, since every positive-order coefficient has
an additional factor $\lambda^{2n\delta}$. Given $m,N$, choose
$J$ so large that \eqref{lo-sum:tail}, inserted in this last
estimate, is $O(\lambda^N)$ and the finite-order residual
\eqref{lo-ext:finite-residual-bound} is also $O(\lambda^N)$.
Keep $J$ fixed and then take $\lambda$ sufficiently small.
This proves flatness of $\mathsf E_B=\mathcal F(\mathbb U_B)$.
No bound uniform in $J$ for the finite constants is required.
Increasing constants on a compact strip bounded away from
$\lambda=0$ extends each estimate to the fixed range
$0<\lambda\le\lambda_*$.

The normalized diagonal bounds give
\eqref{lo-sum:normalized-velocity}: at a fixed derivative order,
the finite initial block is $O(\lambda^{2\delta})$, and the
remaining terms, starting at an index at least two, are bounded by
$\sum 2^{-n}\lambda^{n\delta}=O(\lambda^{2\delta})$.
For stress, put
\[
 H_B=\sum_{n\ge2}\chi_n\lambda^{2n\delta}\mathcal{T}_{(n)}/\omega.
\]
The same argument gives $|\mathscr D^{\le m}H_B|
\le C_m\lambda^{2\delta}$ and the exact decomposition
\begin{equation}\label{lo-sum:weighted-decomposition}
 \mathcal{T}_B-\mathcal{T}_{(0)}=\chi_1\lambda^{2\delta}\mathcal{T}_{(1)}+\omega H_B.
\end{equation}
The first-order bound \eqref{lo-ext:weighted-stress}, together
with $|\partial_R^j\omega|\le C_j\omega d_{\rm edge}^{-q_j}$,
proves \eqref{lo-sum:normalized-stress} by the product rule.
In particular, all fixed derivatives still vanish at the radial edges.
For any formal order $N$, take $J\ge2N+2$ and above the relevant
derivative and exhaustion indices, and make the first
$J$ cutoffs equal to one. The finite block $N<n\le J$ starts
at order $2\delta(N+1)$, and the remaining tail has at least
that order by \eqref{lo-sum:tail}. Thus
\begin{equation}\label{lo-sum:full-expansion}
 |\mathbb U_B-\mathbb U^{[N]}|_m
       \le C_{N,m,K}\lambda^{2\delta(N+1)-\widetilde\ell_m},
\end{equation}
with a fixed loss $\widetilde\ell_m$ independent of $N$;
normalized derivatives retain the original formal orders.

Finally, all positive-order velocity profiles and $M^z_{(n)}$
vanish for $R\ge\sup I_2$, before $I_3$.
Equations \eqref{lo-sum:realized-fields} and
\eqref{lo-sum:cutoff-radial} therefore preserve the velocity on
$I_3$ exactly, including the radial cutoff term containing $M^z_{(n)}$.
By \eqref{lo-ext:pressure-tail} the pressure corrections vanish
for $R\ge R_{0z}$. The velocity and pressure hence
retain the exact heat exterior, where the stress and residual
are zero. At the axis, $P_{(n)}(0,Z)=0$ for $n\ge1$, proving
\eqref{lo-sum:axis-pressure}.
On a compact set with $\lambda>0$, the locally finite sum has
smooth coefficient extensions to $Z=\pm1$. The Jacobian of
$(\lambda,Z)\mapsto(t,z)$ is
$-2\lambda^{2-\delta}L_\delta\ne0$ there. Consequently all
fixed Cartesian derivatives have continuous one-sided extensions.
The remaining physical points away from $(1,0)$ with $\lambda=0$
have $t=1$, $z=0$ and $r>0$. A neighborhood of each such point
lies in the exact heat exterior whenever $\lambda>0$ is small,
since $R=r^2/(2\lambda^2)$. The heat formulas extend smoothly
to $t=1$ for $r>0$, completing the claimed physical extension.
\end{proof}

\subsection{Stress bounds, preservation of the cone, and proof of Theorem~\ref*{thm:background}}
\label{lo-cone:subsection}

We prove that the realized stress stays in the admissible cone on
the whole open annulus. A weighted absolute bound alone does not
give this conclusion where the leading stress tends to zero. We
first compare the stresses relative to their leading magnitude.

\begin{lemma}[Relative stress control at both edges]
\label{lo-cone:relative-stress}
For the completed coefficient construction, there is a constant $C$
such that, for all sufficiently small $\lambda>0$,
\begin{equation}\label{lo-cone:global-relative}
 |\mathcal{T}_B-\mathcal{T}_{(0)}|\le C\lambda^{2\delta}|\mathcal{T}_{(0)}|
       \qquad(R_a<R<R_h,\ |Z|\le1).
\end{equation}
Moreover, with $c_a=h_b^2$ and $c_h=4$ in
\eqref{intro:edge-weight},
\begin{equation}\label{lo-cone:leading-weighted}
 |\mathcal{T}_{(0)}|\ge c\omega,\qquad
 |\mathscr D^{\le m}\mathcal{T}_{(0)}|
                 \le C_m\omega d_{\rm edge}^{-q_m}.
\end{equation}
\end{lemma}

\begin{proof}
At the inner edge every positive-order stress vanishes for
$R\le R_{\rm keep}$, so the difference in
\eqref{lo-cone:global-relative} is exactly zero there. The velocity
corrections can still change the shear on this collar.
Fix a terminal collar $R_{\rm out}\le R<R_h$ with
$R_{\rm out}>R_{0z}$, sufficiently short that
\eqref{eq:outer-stress-factorization} holds with $b_\theta\ge c_0>0$.
Writing $y_h=\log(R_h/R)$ as in \eqref{lo-ext:first-order-edge}, we have
\begin{equation}\label{lo-cone:leading-edge}
 \mathcal{T}^\theta_{(0)}=e^{-4/y_h^2}y_h^{-3}b_\theta(y_h,Z),\qquad
 \mathcal{T}^z_{(0)}=e^{-4/y_h^2}y_h^3b_z(y_h,Z),\qquad b_\theta\ge c_0.
\end{equation}
The coefficients $b_\theta,b_z$ are smooth and bounded up to $y_h=0$.
All positive-order velocities and pressures vanish for
$R\ge R_{0z}$. At order one the only tangential residual in this
collar is the axial viscosity of the leading swirl. The two axial
derivatives in \eqref{lo-ext:collar-axial-viscosity} give a normalized
source of size $O(e^{-4/y_h^2}y_h^{-6})$. Integration inward from $R_h$
gains three powers of $y_h$ by \eqref{lo-ext:flat-integral}. The resulting
bound has exactly the weight of the leading angular stress,
\[
 \mathcal{T}^z_{(1)}=0,\qquad
 |\mathcal{T}^\theta_{(1)}|\le C e^{-4/y_h^2}y_h^{-3}
                         \le (C/c_0)\mathcal{T}^\theta_{(0)}.
\]
At order $n\ge2$, every transport contribution contains a
positive-order factor, and the axial-viscosity source uses order
$n-1\ge1$. Both tangential residuals therefore vanish in this
collar. The backward stress formulas give $\mathcal{T}_{(n)}=0$ there.
Consequently
\[
 \mathcal{T}_B-\mathcal{T}_{(0)}=\lambda^{2\delta}\chi_1\mathcal{T}_{(1)}
                       \qquad(R\ge R_{\rm out}),
\]
which proves the relative estimate at the outer edge. On the
remaining compact interval $[R_{\rm keep},R_{\rm out}]$, the
continuous leading stress has a positive minimum in magnitude,
uniformly for $|Z|\le1$. The normalized stress estimate
\eqref{lo-sum:normalized-stress} gives the relative bound there.

For the inner weight, put $y=\log(R/R_a)$ and recall
$\mathbf q=(\bar D,\bar E)$. On the collar retained in Section~\ref{sec:shear-modification},
\eqref{eq:inner-ray-notation} and \eqref{eq:inner-ray-error} give
\begin{equation}\label{lo-cone:inner-factorization}
 \mathcal{T}_{(0)}=F_{(0)}[(1-\chi_b)\mathbf q+\mathbf e],\qquad
 |\mathbf e|\le\tfrac1{10}(1-\chi_b)|\mathbf q|.
\end{equation}
Both $F_{(0)}$ and $|\mathbf q|$ have positive minima on the closed collar.
Thus $|\mathcal T_{(0)}|\asymp1-\chi_b$, and
\eqref{eq:inner-chi} and \eqref{eq:road-sigma} give
\[
 1-\chi_b=(1-h_b)\sigma(y/h_b)
                 \asymp e^{-h_b^2/y^2}\quad(y\downarrow0).
\]
For every fixed derivative order, \eqref{eq:inner-flat-error} bounds
the mixed derivatives of $\mathbf e$ by $C e^{-h_b^2/y^2}y^{-q}$
for a finite $q$. The derivatives of $1-\chi_b$ satisfy the same
bound by \eqref{eq:road-sigma}. Since $F_{(0)}$ and $\mathbf q$ are
smooth up to $R_a$, the product rule in \eqref{lo-cone:inner-factorization}
gives the same bound for $\mathcal{T}_{(0)}$.
The outer derivative bounds follow directly from \eqref{lo-cone:leading-edge}.
The factor in $\omega$ coming from the opposite edge is smooth
and bounded above and below on either collar. In between, both
$\omega$ and $|\mathcal{T}_{(0)}|$ have positive minima. This proves
\eqref{lo-cone:leading-weighted}; $\partial_\rho$ and
$\partial_R$ have equivalent fixed-order bounds on the annulus.
\end{proof}

\begin{proposition}[Cone preservation on the full open annulus]
\label{lo-cone:preservation}
For a sufficiently small $\lambda_*>0$, the complete background
satisfies the weighted stress bounds
\eqref{intro:weighted-stress} and the uniform admissible cone
\eqref{intro:uniform-cone}. In particular, $\mathcal{T}_B\ne0$ throughout
$R_a<R<R_h$, uniformly for $|Z|\le1$.
\end{proposition}

\begin{proof}
Choose $\lambda_*$ so that $C\lambda_*^{2\delta}\le1/2$ in
\eqref{lo-cone:global-relative}. Then $\mathcal T_B\ne0$ on the open
annulus, and its unit direction $\widehat t_B$ and the leading unit
direction $\widehat t_0$ satisfy
\begin{equation}\label{lo-cone:magnitude-and-direction}
 \tfrac12|\mathcal{T}_{(0)}|\le|\mathcal{T}_B|\le\tfrac32|\mathcal{T}_{(0)}|,
 \qquad |\widehat t_B-\widehat t_0|
                  \le4C\lambda^{2\delta}.
\end{equation}
The lower weighted bound follows from
\eqref{lo-cone:leading-weighted}; the upper derivative bound
follows by combining it with \eqref{lo-sum:normalized-stress}.

For a unit stress direction $t$, the three quantities in
\eqref{intro:uniform-cone} are
\begin{equation}\label{lo-cone:three-tests}
 \kappa-2,\qquad \frac{-t\cdot\mathcal S}{|\mathcal S|},\qquad
 \sqrt{\frac{\kappa-2}{2}}\,
 \frac{|t\cdot\mathcal S^\perp|}{-t\cdot\mathcal S}.
\end{equation}
For the leading profile we bound the first two away from zero and the
third strictly below one. On the inner collar, write
$\bar H=|\mathbf q|^2/\bar D$. By Lemma~\ref{prop:inner-comparison}
and the choice of $R_{ac}$, $\bar H\ge9/4$ and $(1-\chi_b)\bar H\le1/16$.
Hence $\kappa_0=\chi_b\bar H\ge35/16$. The direction estimate
\eqref{eq:inner-ray-angle} gives
\[
 \frac{|\widehat t_0\cdot\mathcal S_{(0)}^\perp|}
      {-\widehat t_0\cdot\mathcal S_{(0)}}
 \le\frac1{9\sqrt{1+\bar H}}.
\]
Since $\widehat t_0$ is a unit vector and $-\widehat t_0\cdot\mathcal S_{(0)}>0$,
the second quantity in \eqref{lo-cone:three-tests} is at least
$(1+1/81)^{-1/2}$. The third is at most $1/(9\sqrt2)$ since
$\kappa_0\le\bar H$.

On the outer collar the shear is purely angular with
$\mathcal S_{(0)}^\theta<0$, and the heat-collar bounds give
$\kappa_0\ge2+\delta/2$ and a uniform upper bound for $\kappa_0$.
By \eqref{lo-cone:leading-edge}, $\mathcal T_{(0)}^\theta>0$ and
$|\mathcal T_{(0)}^z/\mathcal T_{(0)}^\theta|\le C y_h^6$.
Thus the second quantity tends to one and the third tends to zero
uniformly in $Z$ as $R\uparrow R_h$. On the closed subannulus between
these collars, strict admissibility and compactness give positive margins
there as well. Consequently, on the whole open annulus the first two
leading quantities are at least $\eta_*>0$ and the third is at most
$1-\eta_*$, for a common $\eta_*$.

Radial differentiation commutes with every cutoff. Thus the normalized
velocity estimate \eqref{lo-sum:normalized-velocity}, through one radial
derivative, gives
\begin{equation}\label{lo-cone:shear-close}
 |F_B-F_{(0)}|+|\mathcal{S}_B-\mathcal{S}_{(0)}|\le C\lambda^{2\delta}
                 \quad\hbox{on }[R_a,R_h]\times[-1,1].
\end{equation}
On the closed annulus $F_{(0)}$ and $-\mathcal S_{(0)}^\theta$ are bounded
above and away from zero. By \eqref{lo-cone:shear-close} the same holds
for $F_B$ and $-\mathcal S_B^\theta$ for sufficiently small $\lambda_*$.
The leading triples $(F_{(0)},\mathcal S_{(0)},\widehat t_0)$ have compact
closure. On this closure $\kappa_0-2\ge\eta_*$ and
$-\widehat t_0\cdot\mathcal S_{(0)}\ge\eta_*|\mathcal S_{(0)}|$.
Hence the three functions in \eqref{lo-cone:three-tests} are uniformly
Lipschitz near the leading triples. By \eqref{lo-cone:magnitude-and-direction}
and \eqref{lo-cone:shear-close}, each changes by at most $C\lambda^{2\delta}$.
Taking $C\lambda_*^{2\delta}\le\eta_*/2$ gives
\eqref{intro:uniform-cone} with $c_{\rm cone}=\eta_*/2$.
The physical regularized swirl, shear and stress pair have the
common factor $\lambda^{-2-\delta}$ relative to $F_B,\mathcal{S}_B,\mathcal{T}_B$.
Thus these are also the admissible inequalities for the physical fields.
\end{proof}

\begin{proof}[Proof of Theorem~\ref{thm:background}]
\leavevmode\phantomsection\label{lo-cone:conclusion}\label{lo:conclusion}%
Apply Proposition~\ref{lo-sum:background} to the completed coefficient
sequence and choose $\lambda_*$ as in
Proposition~\ref{lo-cone:preservation}.

The equation, exact incompressibility and endpoint smoothness
are proved in Proposition~\ref{lo-sum:background}. Its residual
estimate on $R\le R_h$, together with $\mathsf E_B=0$ outside, proves
\eqref{intro:flatness} on the whole stated domain. The normalized
velocity bounds follow from \eqref{lo-sum:normalized-velocity}
on $R\le R_h$ and the explicit heat formulas for $R\ge R_h$;
their fixed normalized derivatives are bounded there as well.
The weighted stress bounds and the cone follow from
Proposition~\ref{lo-cone:preservation}. The exact velocity on
$I_3$ and the unchanged heat exterior give conclusions (v)
and (vi). Thus all six assertions hold for the same background.

Annular support and the edge estimates make the tensor in
\eqref{intro:full-tensor} smooth. Its radial divergence is
$\partial_z\mathsf T_B^z-(r\partial_z\mathsf T_B^z)/r=0$,
and its tangential divergence is $\mathcal D\mathsf T_B$.
Finally, multiplying the tangential equations by $r^2$ and $r$
and integrating gives the exact identities \eqref{intro:force-moments},
since the stress primitives vanish at both ends.
\end{proof}

\appendix
\section{The analytic construction of the stress-free core}
\label{app:analytic-core}
We prove the analytic existence, uniqueness and exit estimates used
in Proposition~\ref{prop:continuation-core} by comparing the core with the
explicit linear model \eqref{eq:continuation-model} on the fixed
interval $0\le\widetilde R\le4.1$.
Subsection~\ref{subsec:common-analytic-space} places the model in a common
analytic space and constructs the inverses needed for its nonlinear
perturbation. Subsection~\ref{app:nonlinear-core-construction} constructs
the nonlinear core and transfers the model's slope and exit bounds to it.
Subsection~\ref{app:core-exit} then verifies the quantitative outputs
\eqref{eq:inner-core-axial-budget}--\eqref{eq:inner-prepared-normalized}
used in the main construction. The endpoint $4.1$ leaves room beyond
the core exit $\widetilde R=4$; in particular, it covers the collar
$R\le R_a e^{1/100}$ since $4e^{1/100}<4.1$.

\subsection{The linear model in the common analytic space}
\label{subsec:common-analytic-space}
Fix $0<j\le1/20$, the outer parameter $\delta$, and the even axis
pressure $P_0$ from Subsection~\ref{sec:analytic-core}. Put
$\sigma_0=j/500$ and $U_0^z=4Z+j$. The pressure satisfies the real
bounds in Lemma~\ref{lem:outer-axis-pressure}. Keep $\Omega$ and $h_{\rm an}$
fixed as in that subsection, independently of $\Lambda$ and $C_*$.
On $\Omega$ the functions $P_0$, $1/L_\delta$ and
$1/(H_0^2+\sigma_0^2)$ are holomorphic and bounded, as are the
fixed coefficients $\chi_c$, $\beta$ and $g/L_\delta$.
Constants in this subsection may depend on these fixed analytic data,
including $j,\delta,P_0,\Omega$ and $h_{\rm an}$, but not on $\Lambda$ or $C_*$.

We use the rescaled radius $\widetilde R=\Lambda R$ from
\eqref{eq:continuation-scale} and write $\varepsilon=\Lambda^{-1}$ in this appendix.
Normalize the linear profiles by
\[
F^{\rm m}=F_0\Phi^{\rm m},\qquad
U^{z,\rm m}=U_0^z+\varepsilon\Psi^{\rm m}.
\]
Their regular equations and axis values are
\begin{equation}\label{eq:linear-analytic-system}
\begin{gathered}
2\bigl(\widetilde R\partial_{\widetilde R}^2
       +2\partial_{\widetilde R}\bigr)\Phi^{\rm m}
 +({\chi_c}+\varepsilon\beta)\Phi^{\rm m}=0,\\
2\bigl(\widetilde R\partial_{\widetilde R}^2
       +\partial_{\widetilde R}\bigr)\Psi^{\rm m}
 =-g/{L_\delta},\\
\Phi^{\rm m}(0,Z)=1,\qquad \Psi^{\rm m}(0,Z)=0.
\end{gathered}
\end{equation}
The second equation is equivalent, for a regular solution,
to the first-order axial equation of the model: integration
of $\partial_{\widetilde R}(\widetilde R\partial_{\widetilde R}\Psi^{\rm m})
=-g/(2{L_\delta})$ removes the singular integration constant and gives
$\partial_{\widetilde R}\Psi^{\rm m}=-g/(2{L_\delta})$.
These normalized equations do not contain $C_*$. We construct their
solution in the same analytic space used for the nonlinear perturbation,
on the fixed radial interval and with bounds uniform as $\varepsilon\to0$.
This requires control of products and differentiated nonlinear sources
after radial inversion.

For functions of $(\widetilde R,Z)$, define the radial average and primitive
\[
\mathcal Mf(\widetilde R,Z)
=\int_0^1f(s\widetilde R,Z)\,ds,\qquad
\mathcal Vf(\widetilde R,Z)
=\int_0^{\widetilde R}f(v,Z)\,dv.
\]
These formulas include the regular values at the axis.
They will recover the axial moment and pressure from the normalized
velocities in Subsection~\ref{app:nonlinear-core-construction}.

\noindent\textbf{The analytic space and radial inversion.}
For a series $f=\sum_{n\ge0}f_n(Z)\widetilde R^n$ with smooth
coefficients on $[-1,1]$, define
\begin{equation}\label{eq:continuation-norm}
\|f\|_{h_{\rm an}}=
\sup_{\substack{n,m\ge0\\Z\in[-1,1]}}
\frac{20^n h_{\rm an}^m(n+1)^2(m+1)^2|\partial_Z^mf_n(Z)|}
{m!\binom{n+m}{m}}.
\end{equation}
Let $\mathcal X_{h_{\rm an}}$ be the space of such series with finite norm.
The radial weight $20^n$ leaves a fixed margin beyond the
interval $0\le\widetilde R\le4.1$.
The factors $(n+1)^2(m+1)^2$ make the coefficient convolutions
summable in both indices.
The binomial weight allows an axial derivative to be
controlled after radial inversion, as verified below.

The space $\mathcal X_{h_{\rm an}}$ is complete.
Indeed, a Cauchy sequence converges uniformly coefficientwise
with every axial derivative. The limits of successive
derivatives are compatible by the fundamental theorem of
calculus, and the weighted bounds pass to the limit.
More explicitly, if $\|f^{(a)}-f^{(b)}\|_{h_{\rm an}}\le\eta$ for
$a,b\ge N$, then passing to the coefficientwise limit as
$b\to\infty$ gives $\|f^{(a)}-f\|_{h_{\rm an}}\le\eta$.
Thus convergence also holds in the norm.

For each fixed radial degree, the derivative bounds imply
analyticity of the coefficient function. The resulting joint
Taylor series converges normally when
\[
|\widetilde R|<a,\qquad |Z-Z_*|<\eta,\qquad
\frac a{20}+\frac\eta {h_{\rm an}}<1,
\quad Z_*\in[-1,1].
\]
Indeed, its absolute value is bounded by $\|f\|_{h_{\rm an}}$ times
the convergent majorant
\[
\sum_{n,m\ge0}\binom{n+m}{m}
\left(\frac a{20}\right)^n
\left(\frac\eta {h_{\rm an}}\right)^m
=\frac1{1-a/20-\eta/{h_{\rm an}}}.
\]
The local extensions agree on overlaps because they agree
with the same coefficient functions on the real interval.
Taking $a=5$ and $\eta={h_{\rm an}}/2$ gives a common complex
neighborhood of the real rectangle, with a bound
$4\|f\|_{h_{\rm an}}$ for the above majorant.
Cauchy's estimates on smaller neighborhoods therefore give,
for each fixed integer $k\ge0$,
\begin{equation}\label{eq:nonlinear-embedding}
\|f\|_{C^k([0,4.1]\times[-1,1])}
\le K_{k,{h_{\rm an}}}\|f\|_{h_{\rm an}}.
\end{equation}
All derivatives in this norm use $(\widetilde R,Z)$.

For $\nu=1,2$, define the regular zero-axis inverse by
\[
(\mathcal J_\nu f)_0=0,\qquad
(\mathcal J_\nu f)_{n+1}
=\frac{f_n}{(n+1)(n+\nu)}.
\]
Equivalently, on the real radial interval,
\[
(\mathcal J_\nu f)(\widetilde R,Z)
=\int_0^{\widetilde R}s^{-\nu}
  \int_0^s t^{\nu-1}f(t,Z)\,dt\,ds.
\]
The inner integral cancels the apparent singularity at
$s=0$. Differentiation gives
\[
\partial_{\widetilde R}
\left(\widetilde R^\nu
\partial_{\widetilde R}\mathcal J_\nu f\right)
=\widetilde R^{\nu-1}f,
\]
and hence
\[
(\widetilde R\partial_{\widetilde R}^2
+\nu\partial_{\widetilde R})\mathcal J_\nu f=f,
\qquad
(\mathcal J_\nu f)(0,Z)=0.
\]
Thus the angular and axial radial diffusion terms in
\eqref{eq:linear-analytic-system} use $\tfrac12\mathcal J_2$ and
$\tfrac12\mathcal J_1$, respectively. The zeroth-order angular
coefficient $\chi_c+\varepsilon\beta$ will be treated separately below.

For the nonlinear sources, we shall use the following product and
integration estimates. The differentiated products are estimated after
radial inversion; no boundedness of axial differentiation alone in
$\mathcal X_{h_{\rm an}}$ is required.
\begin{equation}\label{eq:continuation-calculus}
\begin{aligned}
\|fk\|_{h_{\rm an}}&\le K\|f\|_{h_{\rm an}}\|k\|_{h_{\rm an}},\\
\|\mathcal Mf\|_{h_{\rm an}}&\le\|f\|_{h_{\rm an}},\qquad
\|\mathcal Vf\|_{h_{\rm an}}\le80\|f\|_{h_{\rm an}},\\
\|\mathcal J_\nu((\partial_Zf)k)\|_{h_{\rm an}}
&+\|\mathcal J_\nu(f\widetilde R
                    \partial_{\widetilde R}k)\|_{h_{\rm an}}\\
&+\|\mathcal J_\nu((\partial_Zf)\widetilde R
                    \partial_{\widetilde R}k)\|_{h_{\rm an}}
\le K_{h_{\rm an}}\|f\|_{h_{\rm an}}\|k\|_{h_{\rm an}}.
\end{aligned}
\end{equation}

Here is the coefficient verification. Write
$\widehat f_{n,m}=\partial_Z^mf_n/m!$ and put
\[
w_{n,m}
=\frac{\binom{n+m}{m}}{(n+1)^2(m+1)^2}.
\]
Then
\[
|\widehat f_{n,m}|
\le\|f\|_{h_{\rm an}}20^{-n}{h_{\rm an}}^{-m}w_{n,m}.
\]
Splitting the sum at $n/2$ gives
\[
\sum_{i=0}^n\frac1{(i+1)^2(n-i+1)^2}
\le\frac8{(n+1)^2}\sum_{i=0}^\infty\frac1{(i+1)^2}
\le\frac{16}{(n+1)^2}.
\]
For $i+j=n$ and $p+q=m$, we also have
\[
\binom{i+p}{p}\binom{j+q}{q}
\le\binom{n+m}{m}.
\]
Indeed, the product on the left is one term in the
Vandermonde sum for the binomial coefficient on the right.
Leibniz's formula for normalized derivatives gives
\[
\widehat{(fk)}_{n,m}
=\sum_{i+j=n}\sum_{p+q=m}
\widehat f_{i,p}\widehat k_{j,q}.
\]
Applying the convolution estimate in both indices yields
\[
\sum_{i+j=n}\sum_{p+q=m}w_{i,p}w_{j,q}
\le256w_{n,m},
\]
which proves the algebra bound.

The most demanding derivative term is
$\mathcal J_\nu((\partial_Zf)
\widetilde R\partial_{\widetilde R}k)$.
Its normalized coefficient of radial degree $n+1$ and
axial order $m$ is
\[
\frac1{(n+1)(n+\nu)}
\sum_{i+j=n}\sum_{p+q=m}
(p+1)\widehat f_{i,p+1}\,j\widehat k_{j,q}.
\]
Use
\[
(p+1)\binom{i+p+1}{p+1}
=(i+1)\binom{i+p+1}{p},
\qquad
\binom{i+p+1}{p}\binom{j+q}{q}
\le\binom{n+m+1}{m}.
\]
The binomial coefficient on the right is precisely the one
in the output weight at radial degree $n+1$.
Moreover,
\[
\frac{(i+1)j}{(n+1)(n+\nu)}\le1.
\]
After multiplying by the norm weight at degree $(n+1,m)$,
the coefficient is therefore bounded by
\[
\begin{aligned}
&\frac{20}{h_{\rm an}}\|f\|_{h_{\rm an}}\|k\|_{h_{\rm an}}
(n+2)^2(m+1)^2
\sum_{i+j=n}\sum_{p+q=m}
\frac1{(i+1)^2(j+1)^2(p+2)^2(q+1)^2}\\
&\qquad\le K_{h_{\rm an}}\|f\|_{h_{\rm an}}\|k\|_{h_{\rm an}}.
\end{aligned}
\]
The last step uses the preceding convolution estimates,
$(p+2)^{-2}\le(p+1)^{-2}$, and
$(n+2)^2/(n+1)^2\le4$.
The factor $20/{h_{\rm an}}$ records the radial degree shift and the
single axial derivative.

For $\mathcal J_\nu((\partial_Zf)k)$, replace $j$ by $1$;
the factor $(i+1)/((n+1)(n+\nu))$ is still at most one.
For $\mathcal J_\nu(f\widetilde R\partial_{\widetilde R}k)$,
omit the axial shift and use
$j/((n+1)(n+\nu))\le1$.
This proves all derivative estimates in
\eqref{eq:continuation-calculus}.

The average divides radial coefficients by $n+1$.
For a shift from degree $n$ to $n+1$, the ratio of the
radial and binomial weights is
\[
20\frac{(n+2)^2}{(n+1)^2}
\frac{\binom{n+m}{m}}{\binom{n+m+1}{m}}
=
20\frac{(n+2)^2}{(n+1)^2}
\frac{n+1}{n+m+1}
\le80.
\]
This proves boundedness of multiplication by
$\widetilde R$. The primitive has the additional divisor
$n+1$, so $\|\mathcal Vf\|_{h_{\rm an}}\le80\|f\|_{h_{\rm an}}$ as well.
The same computation gives, on series starting at degree $k$,
\[
\|\mathcal J_\nu f\|_{h_{\rm an}}
\le\frac{80}{(k+1)(k+\nu)}\|f\|_{h_{\rm an}}.
\]

Finally, if $a(Z)$ is bounded and holomorphic on $\Omega$,
Cauchy's estimate on disks of radius $r_\Omega/2$ gives
\[
\frac{|\partial_Z^ma(Z)|}{m!}
\le\left(\frac2{r_\Omega}\right)^m\sup_\Omega|a|.
\]
Regarding $a$ as a series of radial degree zero, we obtain
\begin{equation}\label{eq:nonlinear-multiplier}
\|a\|_{h_{\rm an}}
\le
\sup_{m\ge0}(m+1)^2(2{h_{\rm an}}/r_\Omega)^m
\sup_\Omega|a|
\le K_{h_{\rm an}}\sup_\Omega|a|.
\end{equation}
The supremum is finite because $2{h_{\rm an}}/r_\Omega<1/2$.
Together with the algebra estimate, this proves the
multiplier property for the fixed analytic coefficients.

\medskip
\noindent\textbf{A factorial inverse for the linear model.}
The coefficient $\chi_c$ remains of order one as $\varepsilon\to0$.
We invert the full angular model using the gain of one radial degree
at each integration, which gives a factorial bound on powers of its
integral operator.
For $0\le\varepsilon\le1/500$, set
\begin{equation}\label{eq:linear-analytic-operator}
\mathcal A_\varepsilon f
=\frac12\mathcal J_2\bigl(({\chi_c}+\varepsilon\beta)f\bigr).
\end{equation}
The multiplier estimate gives a constant $M_{{\chi_c},\beta}$,
depending only on the fixed analytic data and on ${h_{\rm an}}$, such that
\[
\|({\chi_c}+\varepsilon\beta)f\|_{h_{\rm an}}
\le M_{{\chi_c},\beta}\|f\|_{h_{\rm an}}
\quad\text{for }0\le\varepsilon\le1/500.
\]
For example, the algebra constant times
$\|{\chi_c}\|_{h_{\rm an}}+\|\beta\|_{h_{\rm an}}/500$ is an admissible choice.
On the closed subspace of series starting at radial degree
$k$, the radial inverse estimate gives
\[
\|\mathcal A_\varepsilon f\|_{h_{\rm an}}
\le\frac{40M_{{\chi_c},\beta}}{(k+1)(k+2)}\|f\|_{h_{\rm an}}.
\]
Multiplication by a coefficient depending only on $Z$ preserves
this subspace, whereas $\mathcal J_2$ raises the radial degree
by one. Applying this bound successively, starting at degree
zero, yields
\[
\|\mathcal A_\varepsilon^n\|
\le\prod_{k=0}^{n-1}
\frac{40M_{{\chi_c},\beta}}{(k+1)(k+2)}
=\frac{(40M_{{\chi_c},\beta})^n}{n!(n+1)!}.
\]
Consequently,
\begin{equation}\label{eq:continuation-resolvent}
\begin{aligned}
\|\mathcal A_\varepsilon^n\|
&\le\frac{K_{h_{\rm an}}^n}{n!(n+1)!},\\
\mathcal R_\varepsilon
=(I+\mathcal A_\varepsilon)^{-1}
&=\sum_{n=0}^\infty(-\mathcal A_\varepsilon)^n.
\end{aligned}
\end{equation}
The constant in the first line can be chosen uniformly for
$0\le\varepsilon\le1/500$.
The series converges in operator norm, uniformly in this
parameter range. Indeed, its partial sums satisfy
\[
(I+\mathcal A_\varepsilon)
\sum_{n=0}^N(-\mathcal A_\varepsilon)^n
=I-(-\mathcal A_\varepsilon)^{N+1},
\]
and the remainder tends to zero in operator norm.
The same identity holds with the factors reversed.
Thus this is a two-sided inverse and
\[
\sup_{0\le\varepsilon\le1/500}
\|\mathcal R_\varepsilon\|
\le\sum_{n=0}^\infty\frac{K_{h_{\rm an}}^n}{n!(n+1)!}<\infty.
\]
In particular, no smallness of ${\chi_c}$ or of its multiplier
norm is required.
For a zero-axis angular correction, the full operator is inverted by
$\tfrac12\mathcal R_\varepsilon\mathcal J_2$: for $f\in\mathcal X_{h_{\rm an}}$,
\[
2\bigl(\widetilde R\partial_{\widetilde R}^2
       +2\partial_{\widetilde R}\bigr)y
 +(\chi_c+\varepsilon\beta)y=f,\qquad y(0,Z)=0
\]
has the unique regular solution $y=\tfrac12\mathcal R_\varepsilon\mathcal J_2f$
in this space. Here $\mathcal J_2$ inverts radial diffusion, and
$\mathcal R_\varepsilon$ accounts for the zeroth-order angular term.

\medskip
\noindent\textbf{Existence, uniqueness, and identification with the explicit solution.}
For either radial diffusion operator, its regular homogeneous solutions
are constant; the prescribed axis value therefore removes
all freedom when applying $\mathcal J_\nu$.
It follows that \eqref{eq:linear-analytic-system} is equivalent
in $\mathcal X_{h_{\rm an}}\times\mathcal X_{h_{\rm an}}$ to
\[
(I+\mathcal A_\varepsilon)\Phi^{\rm m}=1,
\qquad
\Psi^{\rm m}=-\frac12\mathcal J_1(g/{L_\delta}).
\]
The right sides belong to $\mathcal X_{h_{\rm an}}$, and the inverse
just constructed gives the unique solution in this space:
\begin{equation}\label{eq:linear-analytic-solution}
\Phi^{\rm m}=\mathcal R_\varepsilon1,
\qquad
\Psi^{\rm m}=-\frac{g}{2{L_\delta}}\widetilde R.
\end{equation}
Conversely, the embedding \eqref{eq:nonlinear-embedding}
allows termwise differentiation of these series on a
neighborhood of the real rectangle, so this pair solves the
differential equations and the axis conditions.
This proves existence and uniqueness by radial inversion
without using the explicit formula for $\mathcal B$.

To identify the solution, write
$\Phi^{\rm m}=\sum_{n\ge0}\phi_n(Z)\widetilde R^n$.
The first equation in \eqref{eq:linear-analytic-system} gives
\[
\phi_0=1,\qquad
2(n+1)(n+2)\phi_{n+1}
=-({\chi_c}+\varepsilon\beta)\phi_n.
\]
Induction yields
\[
\phi_n(Z)
=\frac{(-1)^n({\chi_c}(Z)+\varepsilon\beta(Z))^n}
       {2^n n!(n+1)!}.
\]
Thus the solution obtained in the function space is the model
\eqref{eq:continuation-model} in the variables $(\widetilde R,Z)$:
\[
\Phi^{\rm m}
=\mathcal B\bigl(({\chi_c}+\varepsilon\beta)\widetilde R/2\bigr),
\qquad
\Psi^{\rm m}=-\frac{g}{2{L_\delta}}\widetilde R.
\]

The uniform inverse bound \eqref{eq:continuation-resolvent} and the
explicit axial solution give a uniform bound for
$X^{\rm m}=(\Phi^{\rm m},\Psi^{\rm m})$ in
$\mathcal X_{h_{\rm an}}\times\mathcal X_{h_{\rm an}}$ for
$0\le\varepsilon\le1/500$.

\subsection{The nonlinear construction and its comparison estimates}
\label{app:nonlinear-core-construction}
\paragraph{The exact rescaled equations.}
We construct the nonlinear core as a perturbation of the linear model
on the fixed rectangle $0\le\widetilde R\le4.1$, $|Z|\le1$.
We first separate the model terms from the nonlinear corrections;
their factor $\varepsilon$ will yield a contraction after radial inversion.
Throughout this proof,
\[
L=1-Z^2,\qquad {L_\delta}=1-\delta Z^2,\qquad
\widetilde R=\Lambda R,\qquad \varepsilon=\Lambda^{-1}.
\]

The moment and pressure relations used together with
\eqref{eq:core-system} are
\begin{equation}\label{eq:nonlinear-recovery}
M^z(R,Z)=\int_0^R U^z(s,Z)\,ds,\qquad
P(R,Z)=P_0(Z)+\int_0^R F(s,Z)^2\,ds.
\end{equation}
Thus $M^z$ and $P$ are determined by the velocity profiles;
they are not additional unknowns in the analytic construction.
On the scale $R=O(\varepsilon)$, the model varies by $O(F_0)$ in
$F$ and by $O(\varepsilon)$ in $U^z$. Accordingly, write
\begin{equation}\label{eq:continuation-normalization}
F(R,Z)=F_0(Z)\Phi(\widetilde R,Z),\qquad
U^z(R,Z)=U^z_0(Z)+\varepsilon\Psi(\widetilde R,Z),
\end{equation}
with
\[
\Phi(0,Z)=1,\qquad \Psi(0,Z)=0.
\]
Here $\Psi$ is the normalized axial-velocity increment. The model
pair $(\Phi^{\rm m},\Psi^{\rm m})$ is uniformly bounded on the fixed
rectangle by the preceding subsection.

Using the radial average $\mathcal M$ and primitive $\mathcal V$
from the preceding subsection gives
\[
\frac{M^z(\varepsilon\widetilde R,Z)}
     {\varepsilon\widetilde R}
=U^z_0(Z)+\varepsilon\mathcal M\Psi(\widetilde R,Z),
\]
with the quotient interpreted by continuity at $\widetilde R=0$.
Recall that $W,H$ are the transport coefficients in
\eqref{eq:core-transport}, with axis values $W_0,H_0$.
The moment, pressure, and transport identities become
\begin{equation}\label{eq:continuation-moments}
\begin{aligned}
M^z(\varepsilon\widetilde R,Z)
&=\varepsilon\widetilde R U^z_0
  +\varepsilon^2\mathcal V\Psi,\\
P(\varepsilon\widetilde R,Z)
&=P_0+\varepsilon\mathcal P,
\qquad
\mathcal P=F_0^2\mathcal V(\Phi^2),\\
W
&=W_0-\varepsilon\bigl(
 (1-\delta)Z\mathcal M\Psi
 +L\partial_Z\mathcal M\Psi\bigr),\\
H&=H_0+\varepsilon L\Psi,
\end{aligned}
\end{equation}
where
\[
W_0=1-(1-\delta)ZU^z_0-4L.
\]
Using $U^z_0=4Z+j$ and $L=1-Z^2$, we obtain
\[
W_0+\frac{\delta}{2}(1-2ZU^z_0)
=1+\frac{\delta}{2}-ZU^z_0-4L
=-3+\frac{\delta}{2}-jZ=-{L_\delta}\beta.
\]
The coefficient $\chi_c=H_0^2/(H_0^2+\sigma_0^2)$ records the
strength of the axis angular transport and vanishes at $Z_0$.
The positive coefficient $\beta$ keeps the model angular slope negative
there, while $-g/L_\delta$ is the fixed axial forcing of the model.

Substitution into \eqref{eq:core-system} gives
\begin{equation}\label{eq:continuation-rescaled}
\begin{aligned}
2\bigl(\widetilde R\partial_{\widetilde R}^2\Phi
       +2\partial_{\widetilde R}\Phi\bigr)+{\chi_c}\Phi
&=\varepsilon E_\theta,\\
2\bigl(\widetilde R\partial_{\widetilde R}^2\Psi
       +\partial_{\widetilde R}\Psi\bigr)
&=-\frac g{L_\delta}+\varepsilon E_z,
\end{aligned}
\end{equation}
where
\begin{equation}\label{eq:continuation-remainders}
\begin{aligned}
E_\theta={}&
\frac1{L_\delta}\left[
 W\bigl(\Phi+\widetilde R\partial_{\widetilde R}\Phi\bigr)
 +\frac{\delta}{2}
   \bigl(1-2Z(U^z_0+\varepsilon\Psi)\bigr)\Phi
 +H\partial_Z\Phi
\right]
-\frac{LH_0}{H_0^2+\sigma_0^2}\Psi\Phi,\\
E_z={}&
\frac1{L_\delta}\left[
 W\widetilde R\partial_{\widetilde R}\Psi
 +\left(\frac{1+\delta}{2}(1-4ZU^z_0)+4L\right)\Psi
 -\varepsilon(1+\delta)Z\Psi^2
 +H\partial_Z\Psi \right.\\
&\hspace{3em}\left.
 +L\partial_Z\mathcal P
 -2(1+\delta)Z\mathcal P
 -2Z\widetilde R F_0^2\Phi^2
\right].
\end{aligned}
\end{equation}
To verify the powers of $\varepsilon$ and identify the model forcing
in these equations, first note that the radial derivatives transform as
\[
\begin{aligned}
\partial_RF
&=\varepsilon^{-1}F_0\partial_{\widetilde R}\Phi,
&
\partial_R^2F
&=\varepsilon^{-2}F_0\partial_{\widetilde R}^2\Phi,\\
\partial_RU^z
&=\partial_{\widetilde R}\Psi,
&
\partial_R^2U^z
&=\varepsilon^{-1}\partial_{\widetilde R}^2\Psi.
\end{aligned}
\]
The angular equation is divided by ${L_\delta}F_0/\varepsilon$,
whereas the axial equation is divided by ${L_\delta}$.
The quadratic axial source expands exactly as
\[
\begin{aligned}
\frac{1+\delta}{2}(1-2ZU^z)U^z
={}&\frac{1+\delta}{2}(1-2ZU^z_0)U^z_0
+\varepsilon\frac{1+\delta}{2}(1-4ZU^z_0)\Psi
-\varepsilon^2(1+\delta)Z\Psi^2,\\
H\partial_ZU^z
={}&4H_0+\varepsilon(4L\Psi+H\partial_Z\Psi).
\end{aligned}
\]
The terms independent of $\varepsilon$ in the axial
right-hand side sum to $-g$.
For the angular transport term, use
\[
\frac{\partial_ZF_0}{F_0}
=-\varepsilon^{-1}\frac{{L_\delta}H_0}{H_0^2+\sigma_0^2}.
\]
It follows that
\[
\begin{aligned}
\frac{\varepsilon H\partial_ZF}{{L_\delta}F_0}
&=\frac{\varepsilon H}{L_\delta}\partial_Z\Phi
-\frac{HH_0}{H_0^2+\sigma_0^2}\Phi\\
&=-{\chi_c}\Phi
+\varepsilon\left[
\frac H{L_\delta}\partial_Z\Phi
-\frac{LH_0}{H_0^2+\sigma_0^2}\Psi\Phi
\right].
\end{aligned}
\]
Thus the product with $H_0$ produces the term ${\chi_c}\Phi$
on the left, while the product with
$H-H_0=\varepsilon L\Psi$ produces the last term of
$E_\theta$. Neither contribution has been discarded.
These identities verify the exact rescaled equations.

\paragraph{The difference from the linear model.}
Write $X=(\Phi,\Psi)$ and $Y=X-X^{\rm m}$. Adding
$\varepsilon\beta\Phi$ to both sides of the angular equation in
\eqref{eq:continuation-rescaled} gives the same radial operator
as in \eqref{eq:linear-analytic-system}. Put
\begin{equation}\label{eq:core-compensated-source}
\begin{aligned}
G_\theta(X)
&=E_\theta(X)+\beta\Phi\\
&=\frac{1}{L_\delta}
  \left[(W-W_0-\varepsilon\delta Z\Psi)\Phi
     +W\widetilde R\partial_{\widetilde R}\Phi+H\partial_Z\Phi\right]
   -\frac{LH_0}{H_0^2+\sigma_0^2}\Psi\Phi .
\end{aligned}
\end{equation}
The second identity uses
$W_0+\delta(1-2ZU_0^z)/2=-L_\delta\beta$.
Thus $G_\theta$ cancels the fixed zero-order term $-\beta\Phi$
in $E_\theta$. Retaining $\varepsilon\beta$ in the comparison operator
will also allow a relative slope estimate where $\chi_c$ vanishes,
as shown below.
Substitution of the linear model into \eqref{eq:continuation-rescaled}
therefore gives the residual
$-\varepsilon(G_\theta(X^{\rm m}),E_z(X^{\rm m}))$
when the right-hand side is subtracted from the left.
Subtracting the model equations and applying the zero-axis inverses
gives the exact fixed-point equation
\begin{equation}\label{eq:continuation-fixed-point}
\begin{aligned}
Y&=\varepsilon\mathcal N_{\varepsilon,C_*}(X^{\rm m}+Y),\\
\mathcal N_{\varepsilon,C_*}(X)
&=\tfrac12\bigl(\mathcal R_\varepsilon\mathcal J_2G_\theta(X),
                    \mathcal J_1E_z(X)\bigr).
\end{aligned}
\end{equation}
Indeed, the angular difference satisfies
$(I+\mathcal A_\varepsilon)Y_\theta
=\varepsilon\mathcal J_2G_\theta(X)/2$, and the axial model
already accounts for the term $-g/L_\delta$.

We will bound this radially inverted map and its Lipschitz constant
uniformly for $0<\varepsilon\le1/500$ and all $C_*$ satisfying
\eqref{eq:continuation-complex-bound}. The analytic data
$j,\delta,P_0,\Omega,h_{\rm an}$ are already fixed. Such bounds permit
us to choose $\Lambda$ first and then any sufficiently large $C_*$.
Use the sum norm on the product space and let $\mathbb B$ be
the closed unit ball in its subspace of pairs with zero axis value.
Evaluation at $\widetilde R=0$ is continuous. Both radial inverses
have zero axis value, and $\mathcal R_\varepsilon$ preserves this
subspace by its series in \eqref{eq:continuation-resolvent}.
For $Y\in\mathbb B$, the uniform model bound controls
$\|X^{\rm m}+Y\|$ independently of $\varepsilon$ and $C_*$.

The only derivative in $W$ is
$-\varepsilon L\partial_Z\mathcal M\Psi$.
After radial inversion, its products with $\Phi$,
$\widetilde R\partial_{\widetilde R}\Phi$, and
$\widetilde R\partial_{\widetilde R}\Psi$ are controlled
by \eqref{eq:continuation-calculus}, with
$f=\mathcal M\Psi$.
Factors depending only on $Z$, such as $L/{L_\delta}$, commute
with $\mathcal J_\nu$ and are bounded multipliers.
The derivative-free part of $W$ uses the algebra estimate
and the estimate for a single radial derivative.
Similarly, the terms
$(H_0+\varepsilon L\Psi)\partial_Z\Phi$ and
$(H_0+\varepsilon L\Psi)\partial_Z\Psi$ use the integrated
estimate for one axial derivative.

The condition $C_*\ge\Lambda^2e^{\Lambda A_\Omega}$ in
\eqref{eq:continuation-complex-bound} compensates for the possible
growth of $e^{-\Lambda G}$ on the complex neighborhood $\Omega$.
Since $F_0=C_*^{-1}e^{-\Lambda G}$, it gives
\[
 \sup_{\Omega}|F_0|
 \le C_*^{-1}e^{\Lambda A_\Omega}\le\Lambda^{-2}=\varepsilon^2 .
\]
In particular, Cauchy's estimate on disks of radius $r_\Omega/2$
gives $\sup_{[-1,1]}|\partial_Z^k(F_0^2)|
\le k!(2/r_\Omega)^k\varepsilon^4$.
The multiplier bound \eqref{eq:nonlinear-multiplier} therefore gives

\[
\|F_0^2\|_{h_{\rm an}}\le K_{h_{\rm an}}\varepsilon^4,\qquad
\|\mathcal P\|_{h_{\rm an}}
\le K_{h_{\rm an}}\varepsilon^4\|\Phi\|_{h_{\rm an}}^2.
\]
Thus $\mathcal J_1(L\partial_Z\mathcal P/{L_\delta})$ is bounded
by the integrated axial derivative estimate with
$f=\mathcal P$ and $k=1$.
The other pressure and swirl terms contain only products,
a radial primitive, or multiplication by $\widetilde R$.
The term $LH_0\Psi\Phi/(H_0^2+\sigma_0^2)$ likewise uses
only the algebra and multiplier estimates.

We estimate $G_\theta$ using the second expression in
\eqref{eq:core-compensated-source}, where $\beta\Phi$ has already
canceled the term $-\beta\Phi$ in $E_\theta$. The uniform bound for $\mathcal R_\varepsilon$
then applies to the inverted angular source.

Each inverted remainder is a finite sum of bounded
multilinear expressions in $X$.
Their differences are estimated by replacing one factor
at a time. For example,
\[
\mathcal P(\Phi_1)-\mathcal P(\Phi_2)
=F_0^2\mathcal V\bigl(
(\Phi_1-\Phi_2)(\Phi_1+\Phi_2)\bigr),
\]
so
\[
\|\mathcal P(\Phi_1)-\mathcal P(\Phi_2)\|_{h_{\rm an}}
\le K_{h_{\rm an}}\varepsilon^4
(\|\Phi_1\|_{h_{\rm an}}+\|\Phi_2\|_{h_{\rm an}})\|\Phi_1-\Phi_2\|_{h_{\rm an}}.
\]

Thus for $X_i=X^{\rm m}+Y_i$, $Y_i\in\mathbb B$,
\begin{equation}\label{eq:nonlinear-map-bounds}
\|\mathcal N_{\varepsilon,C_*}(X_1)\|\le K_*,
\qquad
\|\mathcal N_{\varepsilon,C_*}(X_1)
-\mathcal N_{\varepsilon,C_*}(X_2)\|
\le K_*\|Y_1-Y_2\|.
\end{equation}
The constant depends only on the fixed analytic data. These estimates
apply after radial inversion, so they require no bound for the
uninverted axial derivatives in $\mathcal X_{h_{\rm an}}$.

For $\Lambda\ge\max\{500,j^{-2},2K_*\}$ and $C_*$ satisfying
\eqref{eq:continuation-complex-bound}, the map in
\eqref{eq:continuation-fixed-point} takes $\mathbb B$ into its
concentric ball of radius $1/2$ and has Lipschitz constant at most
$1/2$. Its unique fixed point satisfies
\begin{equation}\label{eq:continuation-first-error}
\|\Phi-\Phi^{\rm m}\|_{h_{\rm an}}
+\|\Psi-\Psi^{\rm m}\|_{h_{\rm an}}\le K_*\varepsilon .
\end{equation}

The map preserves real coefficients, so the solution is real on the
real rectangle. The joint analytic extension used in
\eqref{eq:nonlinear-embedding} gives, after rescaling, a solution on
a neighborhood of $[0,4.1/\Lambda]\times[-1,1]$.
Termwise differentiation on smaller complex neighborhoods verifies
the exact rescaled differential equations.

Uniqueness is not restricted to the contraction ball.
Indeed, write the radial series of any jointly analytic
solution in \eqref{eq:continuation-rescaled}.
At degree $n$, all terms on the right involve coefficients
of $\Phi,\Psi$ of degrees at most $n$ and their axial
derivatives. A radial primitive raises degree, a radial
average preserves it, and
$\widetilde R\partial_{\widetilde R}$ preserves it.
The term ${\chi_c}\Phi$ also involves only degree $n$.
The coefficients of degree $n+1$ are therefore determined
with nonzero divisors $2(n+1)(n+2)$ and $2(n+1)^2$.
Starting from the prescribed degree-zero data, induction
gives identical Taylor coefficients for any two solutions.
They agree near the axis, and the identity principle gives
agreement on the connected component of their common
analytic domain containing it.

\paragraph{Comparison and the weighted angular slope.}
By \eqref{eq:continuation-first-error} and
\eqref{eq:nonlinear-embedding}, for each fixed $k\ge0$,
\begin{equation}\label{eq:core-normalized-Ck}
\|\Phi-\Phi^{\rm m}\|_{C^k}
+\|\Psi-\Psi^{\rm m}\|_{C^k}\le K_k\varepsilon .
\end{equation}
The $C^k$ norms use $(\widetilde R,Z)$ on
$[0,4.1]\times[-1,1]$. Equivalently, for the original profiles
evaluated at $R=\widetilde R/\Lambda$, the left-hand side is
$\|(F-F^{\rm m})/F_0\|_{C^k}
+\Lambda\|U^z-U^{z,\rm m}\|_{C^k}$. In particular,
\begin{equation}\label{eq:continuation-comparison}
\left\|\frac{F-F^{\rm m}}{F_0}\right\|_{C^2}
+\Lambda\|U^z-U^{z,\rm m}\|_{C^2}\le\frac{K_2}{\Lambda}.
\end{equation}
Axial derivatives act on the entire normalized difference. Every
$K_k$ depends on the fixed analytic data, including $j$, and on $k$,
independently of
sufficiently large $\Lambda$ and every subsequent $C_*$
satisfying \eqref{eq:continuation-complex-bound}.

The $O(\varepsilon)$ comparison alone does not preserve the angular
slope near $Z_0$, where the model slope is itself of order
$\varepsilon$. We therefore estimate the slope error relative to
the scale $\chi_c+\varepsilon$.
Write $c=\chi_c+\varepsilon\beta$ and $D=\Phi-\Phi^{\rm m}$.
The model is $\Phi^{\rm m}=\mathcal B(c\widetilde R/2)$.
By the definition of $\chi_c$,
\[
 |H_0|=\sqrt{H_0^2+\sigma_0^2}\sqrt{\chi_c}\le K\sqrt{\chi_c},
 \qquad
 \frac{|H_0|}{H_0^2+\sigma_0^2}
   \le \sigma_0^{-1}\sqrt{\chi_c}.
\]
Since $\mathcal B'$ is bounded on the model interval, the identity
$\partial_{\widetilde R}\Phi^{\rm m}=c\mathcal B'(c\widetilde R/2)/2$
and the $C^1$ comparison
\eqref{eq:core-normalized-Ck} give
$|\partial_{\widetilde R}\Phi|\le K(\chi_c+\varepsilon)$.
The same comparison bounds $\Phi,\Psi,\partial_Z\Phi,\partial_Z\Psi$.
Hence \eqref{eq:continuation-moments} gives
$|W-W_0|\le K\varepsilon$ and $|H|\le K(\sqrt{\chi_c}+\varepsilon)$.
Substitution in \eqref{eq:core-compensated-source}
yields $|G_\theta(X)|\le K(\sqrt{\chi_c}+\varepsilon)$ throughout
the real rectangle. The angular difference equation now reads
\[
 2\bigl(\widetilde R\partial_{\widetilde R}^2D
       +2\partial_{\widetilde R}D\bigr)
 =-cD+\varepsilon G_\theta(X).
\]
Since $|D|\le K\varepsilon$ and $0\le\chi_c\le1$, the right-hand
side is bounded by $K\varepsilon(\sqrt{\chi_c}+\varepsilon)$.
Integrating
\[
 \partial_{\widetilde R}
   \bigl(\widetilde R^2\partial_{\widetilde R}D\bigr)
 =\frac{\widetilde R}{2}\bigl[-cD+\varepsilon G_\theta(X)\bigr]
\]
from the axis, using regularity there and continuity at
$\widetilde R=0$, gives
\begin{equation}\label{eq:continuation-weighted}
 \left|\partial_{\widetilde R}\Phi
       -\partial_{\widetilde R}\Phi^{\rm m}\right|
 \le K\bigl(\varepsilon\sqrt{\chi_c}+\varepsilon^2\bigr).
\end{equation}

By \eqref{eq:continuation-signs} and $\beta>1$,
\[
-\partial_{\widetilde R}\Phi^{\rm m}
\ge\frac1{16}({\chi_c}+\varepsilon\beta)
\ge\frac1{16}({\chi_c}+\varepsilon).
\]
Moreover,
$\partial_RF/\partial_RF^{\rm m}
=\partial_{\widetilde R}\Phi/
 \partial_{\widetilde R}\Phi^{\rm m}$.
Dividing \eqref{eq:continuation-weighted} by the preceding
lower bound yields
\begin{equation}\label{eq:continuation-relative}
\left|
\frac{\partial_RF}{\partial_RF^{\rm m}}-1
\right|
\le K\frac{\varepsilon\sqrt{\chi_c}+\varepsilon^2}
              {{\chi_c}+\varepsilon}
\le K\sqrt\varepsilon.
\end{equation}
Indeed,
\[
\frac{\varepsilon\sqrt{\chi_c}}{{\chi_c}+\varepsilon}
\le\frac{\sqrt\varepsilon}{2},
\qquad
\frac{\varepsilon^2}{{\chi_c}+\varepsilon}\le\varepsilon.
\]
After increasing $\Lambda$ so that
$K\sqrt\varepsilon\le1/2$, this proves
$\partial_RF<0$ everywhere, including at $Z_0$.
It also gives
\begin{equation}\label{eq:nonlinear-weighted-sign}
\frac1{32}({\chi_c}+\varepsilon)
\le-\partial_{\widetilde R}\Phi
\le2({\chi_c}+\varepsilon).
\end{equation}
For the upper bound, use
$-\partial_{\widetilde R}\Phi^{\rm m}
\le({\chi_c}+\varepsilon\beta)/4$, $\beta<4$, and the
relative bound by $3/2$.

\paragraph{The nonlinear exit and the extension.}

Since $1/4<\Phi^{\rm m}\le1$, increasing $\Lambda$ gives
\[
\frac18\le\Phi\le2,\qquad
\frac{F}{F^{\rm m}}
=\frac{\Phi}{\Phi^{\rm m}}=1+O(\varepsilon).
\]
Thus $F>0$. Together with
\eqref{eq:continuation-relative}, this proves
$F>0$ and $\partial_RF<0$ throughout the rectangle.
For $R>0$, it follows that
$U^\theta=\sqrt{2R}F>0$ and
$\mathcal{S}^\theta=2R\partial_RF<0$.

Combining the bounds for $\Phi$ with
\eqref{eq:nonlinear-weighted-sign} also gives, for a constant $c_0>0$
depending only on the fixed data,
\[
c_0(\Lambda{\chi_c}+1)
\le-\frac{\partial_RF}{F}
=-\Lambda\frac{\partial_{\widetilde R}\Phi}{\Phi}
\le K(\Lambda{\chi_c}+1).
\]
Thus the nonlinear logarithmic radial derivative has the
same scales as the model: order $\Lambda$ where ${\chi_c}$
has a fixed positive lower bound and order one at $Z_0$.

Recall from \eqref{eq:width} that
$I_0=[-j/4,-j/5]$, $b=2j/5$, and $g\ge b$ on $I_0$;
this interval contains $\{\chi_c\le99/100\}$.
The axial comparison gives
\begin{equation}\label{eq:axial-regions}
\partial_RU^z(R,Z)
=-\frac{g(Z)}{2{L_\delta}(Z)}+O(\varepsilon),\qquad
\partial_RU^z(R,Z)\le-\frac b4
\quad(Z\in I_0),
\end{equation}
after increasing $\Lambda$.
Indeed, on $I_0$ the model derivative is at most $-b/2$,
and the error $K/\Lambda$ can be made at most $b/4$ by
taking $\Lambda\ge4K/b$, with $K$ depending on the fixed core data.
This negativity is only asserted on $I_0$.
At $Z=0$, evenness of $P_0$ and $U^z_0(0)=H_0(0)=j$ give
$g(0)=-(9+\delta)j/2$. Thus $\partial_RU^z(R,0)>0$ after
increasing $\Lambda$.

It remains to transfer the model exit bound on its two regions.
We use the angular contribution where $\chi_c\ge99/100$ and the
axial contribution where $\chi_c\le99/100$; the comparisons above
control the relative change of the corresponding positive term.
Recall
\[
\kappa
=-\frac{2R\partial_RF}{F}
+\frac{(\partial_RU^z)^2}{-F\partial_RF}.
\]
Both contributions are nonnegative.
At $R=R_a$, on $\{{\chi_c}\ge99/100\}$, the ratio of the
nonlinear angular contribution to the model contribution is
\[
\frac{\partial_RF}{\partial_RF^{\rm m}}
\frac{F^{\rm m}}F=1+O(\sqrt\varepsilon).
\]
On $\{{\chi_c}\le99/100\}$, \eqref{eq:width} gives $g\ge b$,
and hence $|\partial_RU^{z,\rm m}|\ge b/2$.
The axial comparison therefore implies
\[
\left|
\frac{\partial_RU^z}{\partial_RU^{z,\rm m}}-1
\right|
\le\frac{2K}{b}\varepsilon.
\]
Since $b=2j/5>0$ is fixed before $\Lambda$ is chosen, this is an
$O(\varepsilon)$ estimate for the fixed core data; its constant
may depend on $j$. In particular, $\Lambda\ge40K/b$ makes the
displayed relative error at most $1/20$. The ratio of the nonlinear
axial contribution to the model contribution is
\[
\left(
\frac{\partial_RU^z}{\partial_RU^{z,\rm m}}
\right)^2
\frac{F^{\rm m}}F
\frac{\partial_RF^{\rm m}}{\partial_RF}
=1+O(\sqrt\varepsilon).
\]

Increase $\Lambda$, using only the fixed data, so that the angular
and axial ratios just obtained both exceed $5/6$ on their respective
regions. Each corresponding model contribution exceeds $3$ by
\eqref{eq:continuation-angular-exit} and \eqref{eq:continuation-axial-exit}.
Thus at $R=R_a=4/\Lambda$,
\begin{equation}\label{eq:nonlinear-angular-exit}
\begin{aligned}
 -2R_a\frac{\partial_RF}{F}&>\frac52
       &&(\chi_c\ge99/100),\\
 \frac{(\partial_R U^z)^2}{-F\partial_R F}&>\frac52
       &&(\chi_c\le99/100).
\end{aligned}
\end{equation}
Both terms are nonnegative everywhere, so
$\kappa(R_a,Z)\ge5/2$ for every $|Z|\le1$.
All lower bounds on $\Lambda$ in this construction depend only on
the fixed analytic data, including $j$ and $h_{\rm an}$.
These bounds and the exit estimates in Subsection~\ref{app:core-exit}
determine a finite $\Lambda_0$ independently of every subsequent $C_*$
satisfying \eqref{eq:continuation-complex-bound}.

The pressure and moment are recovered by
\eqref{eq:nonlinear-recovery}. They are jointly analytic,
and
\[
\frac{M^z(R,Z)}R=\int_0^1U^z(sR,Z)\,ds
\]
extends analytically through $R=0$.
Thus the transport coefficients are regular there, and
the analytic solution solves \eqref{eq:core-system} with the
prescribed axis data. By \eqref{inviscid-stress-radial} and
\eqref{eq:stress-free-sources}, $R\mathcal{T}^\theta$ and $\sqrt R\,\mathcal{T}^z$ are
constant in $R$. Both vanish at the axis by regularity, so $\mathcal{T}=0$.
This uses the local pressure with
axis value $P_0$; matching it to a global pressure moment
is a separate step.

Finally, the solution exists up to $4.1/\Lambda>R_a$.
On a closed interval to the right of $R_a$, the denominator
$-F\partial_RF$ is strictly positive, so $\kappa$ is
continuous. Uniform continuity in $Z$ and the exit margin
give, for each chosen solution, a number $\eta>0$ such that
\[
\kappa>2
\quad\hbox{on }[R_a,R_a+\eta]\times[-1,1],
\qquad
R_a+\eta\le4.1/\Lambda.
\]

\subsection{Estimates at the core exit}
\label{app:core-exit}

We complete the proof of Proposition~\ref{prop:continuation-core}.
Put $R_a=4/\Lambda$, $\varepsilon=\Lambda^{-1}$, $U_0^z=4Z+j$
and $\sigma_0=j/500$. At the exit write
\[
 f=F_c(R_a,\cdot),\qquad v=U_c^z(R_a,\cdot),\qquad
 m_z=M_c^z(R_a,\cdot).
\]
We use the stress-free continuation to $R_a e^{1/100}$, which is available
since $4e^{1/100}<4.1$.
Constants denoted by $B$ below may depend on the fixed analytic data,
including $j$, and are uniform for sufficiently large $\Lambda$ and
all $C_*$ satisfying \eqref{eq:continuation-complex-bound}.
A subscript $\Lambda$ allows additional dependence on the fixed $\Lambda$.
The lower bounds on $\Lambda$ below are included in the threshold
$\Lambda_0$ of Proposition~\ref{prop:continuation-core} for the given $j$.

\paragraph{Normalized core estimates.}
Recall the normalized fields $\Phi=F_c/F_0$ and
$\Psi=\Lambda(U_c^z-U_0^z)$, expressed in
$(\widetilde R,Z)$ with $\widetilde R=\Lambda R$.
By \eqref{eq:core-normalized-Ck},
for every fixed $k$,
\[
 \|\Phi-\Phi^{\rm m}\|_{C^k}+
 \|\Psi-\Psi^{\rm m}\|_{C^k}\le B_k\varepsilon
\]
on $[0,4.1]\times[-1,1]$ in $(\widetilde R,Z)$, with $B_k$ depending only on
$k$ and the fixed analytic data. Here
\[
 \Phi^{\rm m}(\widetilde R,Z)
 =\mathcal B((\chi_c+\varepsilon\beta)\widetilde R/2),
 \qquad \Psi^{\rm m}(\widetilde R,Z)
 =-\frac{g}{2L_\delta}\widetilde R .
\]
The fixed smooth coefficients
and the explicit model give uniform $C^k$ bounds for both components.
By \eqref{eq:continuation-signs}, $1/4<\Phi^{\rm m}\le1$.
Increasing the lower bound for $\Lambda$ gives $1/8\le\Phi\le2$.
For fixed $\Lambda$, $e^{-\Lambda G}$ and its reciprocal have bounded
derivatives of every fixed order on $[-1,1]$.
Since $C_*F_c=e^{-\Lambda G}\Phi$, $U_c^z=U_0^z+\varepsilon\Psi$,
and $\widetilde R=4R/R_a$, this proves \eqref{eq:inner-prepared-normalized}
on the stated core rectangle, including the reciprocal bound.

\paragraph{Exit slopes and the axial defect.}
By \eqref{eq:nonlinear-weighted-sign} and $\Phi\le2$,
$-\partial_R\log F_c
=-\Lambda\,\partial_{\widetilde R}\Phi/\Phi\ge1/64$
throughout the real core rectangle.
On $\{\chi_c\ge99/100\}$, \eqref{eq:nonlinear-angular-exit} gives
the stronger bound
\[
 -2R_a\,\partial_R\log F_c(R_a,Z)>\frac52.
\]
On $\{\chi_c\le99/100\}$ we have $g\ge b=2j/5$ by \eqref{eq:width}.
Equation \eqref{eq:axial-regions} gives
\[
 -2L_\delta\partial_R U_c^z(R_a,Z)
 =g+O(\Lambda^{-1})\ge3b/4
\]
for sufficiently large $\Lambda$.

At the exit the explicit axial model and \eqref{eq:continuation-comparison} give
\[
\begin{aligned}
	v&=4Z+j-\frac{2g}{\Lambda L_\delta}+O_{C^2_Z}(\Lambda^{-2}),\\
	\frac{m_z}{R_a}&=4Z+j-\frac{g}{\Lambda L_\delta}
	+O_{C^2_Z}(\Lambda^{-2}).
\end{aligned}
\]
Indeed, throughout $0\le R\le R_a$,
\[
 U_c^z(R,Z)=U_0^z(Z)-\frac{Rg(Z)}{2L_\delta(Z)}
                 +O_{C_Z^2}(\Lambda^{-2}),
\]
and $m_z/R_a$ is the radial average of this expression.
Recalling the axial defect from \eqref{eq:inner-core-axial-budget},
\[
 \rho_{\rm core}
 =\max\left\{\|v-U_0^z\|_{C_Z^2},
       \left\|\frac{m_z}{R_a}-U_0^z\right\|_{C_Z^2}\right\},
\]
we obtain
\[
 \rho_{\rm core}\le\frac{2}{\Lambda}
      \left\|\frac{g}{L_\delta}\right\|_{C_Z^2}
       +\frac{B}{\Lambda^2},
\]
which proves \eqref{eq:inner-core-axial-budget} for the fixed $j$.

\paragraph{Axial variation of the swirl at the exit.}
The axis factor $F_0=C_*^{-1}e^{-\Lambda G}$ contributes the term
$-\Lambda\partial_ZG$ to the axial logarithmic derivative. At the exit,
\[
 \partial_Z\log f+\Lambda\partial_ZG
 =\partial_Z\log\Phi(4,Z).
\]
We estimate this remaining variation with the weight $\sqrt{\chi_c}$,
which records its smaller size near the zero of $H_0$.
The identity $\chi_c=H_0^2/(H_0^2+\sigma_0^2)$ gives
$|\partial_Z\chi_c|\le B\sqrt{\chi_c}$ for the fixed $\sigma_0>0$.
Writing $q=(\chi_c+\varepsilon\beta)\widetilde R/2$
in the explicit model, we have
\[
 \partial_Z\log\Phi^{\rm m}
 =\frac{\mathcal B'(q)}{\mathcal B(q)}\,\frac{\widetilde R}{2}
   (\partial_Z\chi_c+\varepsilon\partial_Z\beta).
\]
The Bessel argument ranges over a fixed compact interval on which
$\mathcal B>1/4$, so $|\partial_Z\log\Phi^{\rm m}|
\le B(\sqrt{\chi_c}+\varepsilon)$.
The $C^1$ comparison and the lower bounds for $\Phi$ and $\Phi^{\rm m}$
give $|\partial_Z\log\Phi-\partial_Z\log\Phi^{\rm m}|\le B\varepsilon$.
Evaluating at $\widetilde R=4$ proves \eqref{eq:inner-core-input}.

\section*{Acknowledgement}
ZL was supported in part by NSFC Excellent Research Group Project
(No.62588101), NSFC (No.12494544), NSFC (No.12431007) and
New Cornerstone Science Foundation through the XPLORER PRIZE. XR was supported in part by NSFC (No.12601444) and the National Key R\&D Program of China (No. 2023YFA1010700).
We used AI tools for drafting and polishing the paper.

\end{document}